\documentclass[11pt]{article}
\usepackage[utf8]{inputenc}
\usepackage[margin=1.0in]{geometry}
\usepackage{mathtools} 
\usepackage{amsthm}
\usepackage{amssymb} 
\usepackage{amsfonts, dsfont} 
\usepackage{graphicx} 
\usepackage{color} 
\usepackage{subcaption}
\usepackage{enumerate} 
\usepackage{enumitem}
\usepackage{fancybox} 
\usepackage{float}  
\usepackage{verbatim}
\usepackage{stmaryrd}
\usepackage{empheq}
\usepackage{booktabs}
\usepackage{mathrsfs}
\usepackage{tikz}
\usetikzlibrary{cd, patterns,arrows,decorations.pathreplacing, decorations.markings,arrows.meta, shapes.multipart}
\usetikzlibrary{shapes.misc}
\usepackage{appendix}
\usepackage{bbold}
\usepackage{cases}
\usepackage{slashed}
\usepackage{authblk}
\usepackage{hyperref}

\newtheorem{theorem}{Theorem}[section]
\newtheorem{prop}[theorem]{Proposition}
\newtheorem{lemma}[theorem]{Lemma}
\newtheorem{remark}[theorem]{Remark}
\newtheorem{coro}[theorem]{Corollary}
\newtheorem{definition}[theorem]{Definition}

\newcommand{\R}{\mathbb{R}}
\newcommand{\Z}{\mathbb{Z}}
\newcommand{\N}{\mathbb{N}}
\newcommand{\C}{\mathbb{C}}
\newcommand{\dd}{\; \mathrm{d}}

\newcommand{\cotan}{\text{cotan}}
\newcommand{\supp}{\text{supp}}
\renewcommand{\div}{\text{\normalfont div}}
\newcommand{\grad}{\text{\normalfont grad}}
\newcommand{\vol}{\text{\normalfont vol}}
\newcommand{\eff}{\text{\normalfont eff}}
\renewcommand{\tt}{{\tilde{\mathfrak{t}}}}

\newcommand{\ul}{u_{l}}
\newcommand{\uh}{u_{h}}
\title{Optimal decay for solutions of the Teukolsky equation on the Kerr metric for the full subextremal range $\lvert a\rvert<M$}
\author{Pascal Millet}
\affil{Université Grenoble Alpes, Institut Fourier,\\
100 rue des Maths, 38610 Gières, France\\
pascal.millet@univ-grenoble-alpes.fr}
\date{}
\begin{document}
\maketitle
\sloppy
\begin{abstract}
We derive the large time asymptotics of initially regular and localized solutions of the Teukolsky equation on the exterior of a subextremal Kerr black hole for any half integer spin. More precisely, we obtain the leading order term (predicted by Price's law) in the large time regime assuming that the initial data have compact support and have enough (but finite) Sobolev regularity. For initial data with less spatial decay (typically decaying like $r^{-1-\alpha}$ with $\alpha\in (0,1)$), we prove that the solution has a pointwise decay of order $t^{-1-\alpha-s-\left|s\right|+\epsilon}$ on spatially compact regions.

In the proof, we adopt the spectral point of view and make use of recent advances in microlocal analysis and non elliptic Fredholm theory which provide a robust framework to study linear operators on black hole type spacetimes. 
\end{abstract}

\tableofcontents

\section{Introduction}
The study of wave propagation on black hole spacetimes has been the subject of intense research in the last decades. An important motivation is to understand the stability properties of explicit solutions of the Einstein equations. First works on this subject are \cite{friedrich1986existence} for the De-Sitter solution and \cite{christodoulou1993global} (see also more recently \cite{lindblad2010global}) for the Minkowski solution. More recently, stability results have been obtained for black hole solutions: in \cite{hintz2018global}, \cite{hintz2018non} for Kerr-De-Sitter and Kerr-Newman-De-Sitter, in \cite{klainerman2017global} and \cite{dafermos2021non} for Schwarzschild and in \cite{klainerman2022construction, klainerman2019effective, giorgi2020general, klainerman2021kerr, shen2022construction} for the slowly rotating Kerr solution.

All these results are based on a precise description of the propagation of perturbations on the underlying spacetime at the linear level. Linear stability results include \cite{dafermos2019linear},\cite{hung2020linear}, \cite{zbMATH07167670}, \cite{johnson2019linear} for Schwarzschild, \cite{andersson2019stability} and \cite{hafner2021linear} for slowly rotating Kerr and \cite{giorgi2020linear} for subextremal Reissner-Nordstr{\"o}m spacetimes. As shown in some of these approaches, the Teukolsky equation introduced in \cite{teukolsky1973perturbations} can be used to reduce the tensorial equations of linearized gravity to a scalar wave type equation involving gauge independent quantities. Moreover, the Teukolsky equation also encompass Maxwell equation and the scalar wave equation which are often considered as simplified model for linearized gravity. The equation depends on a spin parameter $s\in \frac{1}{2}\Z$ whose value changes according to the type of field being studied: $s=0$ for scalar waves, $s=\pm\frac{1}{2}$ for sourceless neutrino fields, $s = \pm 1$ for Maxwell fields and $s=\pm 2$ for linearized gravity.

There exists a large literature about the Teukolsky equation and the related linear wave equations. The optimal decay rate for solutions of the Teukolsky equation on a Schwarzschild spacetime was first conjectured by Price \cite{price1972nonsphericalI, price1972nonsphericalII} and later clarified by Price-Burko \cite{price2004late}. The conjecture was generalized to the Kerr case in \cite{hod1999mode, gleiser2008late}.
The most studied case has been the wave equation $s=0$ starting with pioneer works by Wald and Kay-Wald \cite{wald1979note, kay1987linear} followed by many authors. The question of optimal decay on subextremal Kerr black hole is now well understood. Tataru \cite{tataru2013local} (see also \cite{metcalfe2012price} for a generalization) obtained the optimal local decay for a family of stationary spacetimes including the subextremal Kerr family assuming a local integrated energy estimate (this estimate holds for the \emph{full subextremal range} of paramters as proved by Dafermos–Rodnianski–Shlapentokh-Rothman \cite{dafermos2016decay}). Using a different approach Donninger-Schlag-Soffer \cite{donninger2011proof} obtained the optimal local decay on Schwarzschild black holes. The global optimal decay (by global we mean uniform up to null infinity) was obtained in the spherically symmetric case by Angelopoulos-Aretakis-Gajic \cite{angelopoulos2018vector, angelopoulos2018late} where they also compute the precise asymptotic profile. The global optimal decay and leading order term are obtained by Hintz \cite{hintz2022sharp} in the subextremal Kerr case. On a Reissner–Nordström spacetime, we also mention the result of Luk-Oh \cite{luk2017proof} where the authors obtained a sharp decay for the wave equation and deduce the instability of the Cauchy horizon under scalar perturbations, highlighting the link between sharp decay of waves and the strong cosmic censorship.

For higher spin, Donninger-Schlag-Soffer \cite{donninger2012pointwise} (later refined by Metcalfe-Tataru-Tohaneanu \cite{metcalfe2017pointwise} in the case $s=\pm 1$ under an integrated energy decay assumption)  obtained an explicit but sub optimal polynomial decay rate for Teukolsky solutions with spin $\pm1, \pm2$ on a Scwharzschild background. On a slowly rotating Kerr black hole, integrated energy decay was proved for the Teukolsky equation for spin $s=\pm 1$ by Ma \cite{MaMorawetzI} and for spin $s =\pm 2$ by Ma \cite{MaMorawetz} and independently by Dafermos-Holzegel-Rodnianski \cite{dafermos2019boundedness}. On the Schwarzschild background, Price's law was obtained by Ma-Zhang for spin $s=\pm \frac{1}{2}, \pm 1, \pm 2$ in \cite{ma2022price, ma2022sharp}. Ma-Zhang \cite{ma2021sharp} further generalized their result to the Price's law (with computation of the leading term) for spin $s = \pm 1, \pm 2$ on the exterior region of a slowly rotating Kerr black hole $\left|a\right|\ll M$ (and conditionally in the case $\left|a\right|<M$). More recently, the first part of a boundedness and decay result in the full subextremal range has been released by Shlapentokh-Rothman-Teixeira da Costa \cite{shlapentokh2020boundedness}. In the current work, we prove the unconditional Price's law for the whole subextremal range of parameters ($\left|a\right|<M$) and for arbitrary spins $s\in \frac{1}{2}\Z$. Our approach relies on microlocal and spectral methods (and in particular on works by Vasy \cite{vasy2013microlocal}, \cite{vasy2020limiting}, \cite{vasy2020resolvent}, Wunsch-Zworski \cite{wunschZworski} and Dyatlov \cite{dyatlov2016spectral}). We make use of Melrose's $b$ and scattering pseudodifferential algebras \cite{melrose2020spectral, melrose1993atiyah}. A crucial point in the proof of our result is the analysis of the low energy limit of the resolvent which has been initiated in the euclidean context by the work of Jensen-Kato \cite{JensenKato}. Although we adopt the Vasy's point of view \cite{vasy2020resolvent} for the low energy analysis, we mention \cite{Bony2010semilinear, Bony2010low, Bony2013local, GuillarmouHasselI, GuillarmouHasselII, zbMATH06208251} for a different perspective on this problem. These methods have recently led to the previously mentioned linear stability result for Einstein's equations on Kerr black holes by Häfner-Hintz-Vasy in \cite{hafner2021linear} and to the sharp asymptotic description of scalar waves by Hintz \cite{hintz2022sharp}. A key ingredient of our proof is the mode stability result, originally obtained by Whiting in \cite{whiting1989mode} and further improved by Andersson-Ma-Paganini-Whiting in \cite{andersson2017mode} and Andersson-Häfner-Whiting in \cite{andersson2022mode}(see also \cite{casals2021hidden} for a partial result for Kerr-De-Sitter).

\subsection{Main results}
Let $s\in \frac{1}{2}\Z$. We fix a subextremal Kerr black hole spacetime with parameters $\left|a\right|<M$ (see section \ref{KerrMetric}).
Our main results concern the solution of the Cauchy problem with initial data on a spacelike hypersurface $\Sigma_0$ transversal to the future horizon and equal to the Boyer-Lindquist initial hypersurface near Boyer-Lindquist radius $r=+\infty$.
\begin{align}\label{cauchyPbmIntro}
\begin{cases}
T_su= 0\\
u_{|_{\Sigma_0}} = u_0\\
n_{\Sigma_0} u_{|_{\Sigma_0}} = u_1
\end{cases}
\end{align} where $T_s$ is the Teukolsky operator (see Section \ref{analyticFramework}) and $n_{\Sigma_0}$ is the unit normal vector field to $\Sigma_0$. We prove two different results depending on the assumptions on the initial data. In both cases, we do not intend to provide optimal regularity assumptions on the initial data.
We denote by $\Delta^{[s]}$ the spin weighted laplacian. To state the results, we introduce the function $\mathfrak{t}$ whose level sets are transversal to null infinity and to the future horizon (see Subsection \ref{KerrMetric} for the exact definition).
\begin{theorem}\label{mainTheoWeak}
Let $\alpha \in (0,1)$. Let $k\in \N$. 
There exists $N\in \N$ such that the following holds.
For all initial data $u_0,u_1$ such that for all $j\leq N$, $r^{-\frac{1}{2}+\alpha}(r\partial_r )^{2j}u\in L^2(\Sigma_0)$ and $r^{-\frac{1}{2}+\alpha}(\Delta^{[s]})^j u \in L^2(\Sigma_0)$, we have for all $p\leq k$ and $\mathfrak{t}\geq 1$:
\begin{align*}
\left|(\mathfrak{t}\partial_{\mathfrak{t}})^{p}u(\mathfrak{t},r)\right|\leq C_{k}\mathfrak{t}^{-\alpha+}\frac{r^{s+\left|s\right|+}}{(r+\mathfrak{t})^{1+s+\left|s\right|}}
\end{align*}

\end{theorem}

\begin{theorem}\label{mainTheoStrong}
Let $k\in \N$. There exists $N\in \N$ and $\epsilon>0$ such that the following holds:
For all compactly supported initial data $u_0,u_1 \in H^N(\Sigma_0)$, for all $p\leq k$, all $\delta>0$ and all $\mathfrak{t}\geq 0$ we have:
\begin{align*}
\left|(\mathfrak{t}\partial_{\mathfrak{t}})^p \left(u(\mathfrak{t},r)- \mathfrak{p}(\mathfrak{t},r^{-1},\omega)\right)\right|\leq C_{\delta,k} r^{-1+\left|s\right|+s+\delta}\mathfrak{t}^{-3-\left|s\right|+s-\epsilon}(t+r)^{-\left|s\right|-s}
\end{align*}
where $\mathfrak{p}$ is an explicit function depending on $s$ and on the initial data.
\end{theorem}
\begin{remark}
The $\delta>0$ in the previous theorem comes from the fact that we stated Theorem \ref{PreciseThmHighDecay} in a Sobolev space with respect to $r$ (and we lose a small power of $r$ using the Sobolev embedding). However we could get rid of it with some extra work starting from proposition \ref{radiationField}.
\end{remark}
\begin{remark}
For a more precise version of these theorems and the explicit definition of $\mathfrak{p}$, see Corollary \ref{preciseCauchyLowDecay} and Corollary \ref{PreciseCauchyHighDecay}.
\end{remark}

\subsection{Method of proof}
Theorems \ref{mainTheoWeak} and \ref{mainTheoStrong} are stated in terms of a Cauchy problem \eqref{cauchyPbmIntro}, but the forcing problem 
\begin{align*}
T_s v = f
\end{align*}
(where $v,f$ have bounded support in the past) is more convenient in the spectral analysis perspective. To get a forcing problem from \eqref{cauchyPbmIntro}, a natural idea is to take $v= \chi(\mathfrak{t})u$ for $\chi \in C^{\infty}(\R_{\mathfrak{t}},[0,1])$ equal to $1$ in a neighborhood of $+\infty$ and to $0$ in a neighborhood of $-\infty$. With this definition, $u$ and $v$ share the same asymptotic behavior near $\mathfrak{t} = +\infty$ and $f$ has compact support in $\mathfrak{t}$.  In the case of initial data with compact support, it is possible to choose $\mathfrak{t}$ such that $f$ also has a spatial compact support. Without the compact support hypothesis, we use energy estimates similar to what is done in \cite{hintz2020stability} to compute the behavior of $f$ near null infinity.

Then we take the Fourier-Laplace transform with respect to $\mathfrak{t}$ and obtain the following equation for $\Im(\sigma)$ large enough: $\hat{T}_s(\sigma) \hat{v} = \hat{f}$. If we are able to prove that $\hat{T}_s$ is invertible between suitable spaces, we obtain the following expression for $v$: 
\begin{align*}
v(\mathfrak{t})= \frac{1}{2\pi}\int_{\Im(\sigma)=C}e^{-i\mathfrak{t}\sigma}R(\sigma)\hat{f}(\sigma)\dd \sigma
\end{align*}
where $R(\sigma) = \hat{T}_s(\sigma)^{-1}$. A formal estimate of the right-hand side using integration by part, assuming $\partial_\sigma^k R(\sigma)\hat{f}(\sigma)$ is integrable provides:
\begin{align*}
\left|v(\mathfrak{t})\right|\leq e^{C\mathfrak{t}}\mathfrak{t}^{-k}\left\|\partial_\sigma^k R(\sigma)\hat{f}(\sigma)\right\|_{L^1(\R_{\sigma})}
\end{align*}
This formal computation leads to the intuition that better estimates are obtained when $C$ is small and inverse polynomial estimates correspond to $C=0$.
The previous observations suggest the following key points to address:
\begin{enumerate}
\item Prove that $\hat{T}_s(\sigma)$ is invertible on the upper half plane between suitable spaces. \label{inversibility}
\item Prove that $R(\sigma)$ admits a polynomial bound when $\left|\sigma\right|\to +\infty$ and $\Im(\sigma)$ remains in a compact set (note that since we impose high regularity on the initial data, $\hat{f}(\sigma)$ has a high polynomial decay with respect to $\sigma$). \label{highFreqIntro}
\item Prove that $R(\sigma)$ is holomorphic on the strictly upper half plane and continuous up to the real axis. \label{holomorphy}
\item Analyze precisely the regularity of $R(\sigma)$ on the real axis. Note that as shown in \cite{hintz2022sharp} for the wave equation, the leading order term can be obtained by computing the higher order singularity of $R(\sigma)f(\sigma)$ on the real axis (in our case, it is localized at $\sigma = 0$). \label{regularityReal}
\end{enumerate}

Points \ref{inversibility} can be subdivided into three steps: Proving that $\hat{T}_s(\sigma)$ is Fredholm, proving that its index is zero and proving that its kernel is trivial. The Fredholm property is obtained by gluing Fredholm estimates on different regions of phase space: Near radial points on the horizon, (using \cite{vasy2013microlocal}), in a small region inside the black hole (using a hyperbolic estimate as is done in \cite{zworski2015resonances}) and near radial points at spatial infinity using \cite{vasy2020limiting, vasy2020resolvent}. Note that in order to apply the results of \cite{vasy2020resolvent} we have to check the triviality of the kernel of the effective normal operator of $\hat{T}_s(\sigma)$ (see definition \ref{defNormalOp}) which involves theory of the confluent hypergeometric equation. The gluing process relies on elliptic estimates and on propagation of singularities. Therefore, we need a global analysis of the classical Hamiltonian flow of the operator which can be computed from its principal symbol. Triviality of the kernel follows from the mode stability result for the Teukolsky equation on a subextremal Kerr spacetime originally obtained by Whiting in \cite{whiting1989mode} and further improved by Andersson-Ma-Paganini-Whiting in \cite{andersson2017mode} (see also \cite{andersson2022mode} for $\sigma =0$). The index zero property will follows from the continuity of the index and the invertibility of $\hat{T}_s(\sigma)$ for large $\Re(\sigma)$.

We prove this invertibility together with the polynomial bound of point \ref{highFreqIntro} by introducing the semiclassical parameter $h = \frac{1}{\left|\sigma\right|}$, the operator $\hat{T}_{s,h}(z) := h^2\hat{T}_s(h^{-1}z)$ and proving a bound of the form $\left\|u\right\|\leq h^{-2}\left\|\hat{T}_{s,h}(z)u\right\|$ (see Proposition \ref{globalSemiclassicalEstimate} for the precise statement). As before, we glue semiclassical estimates obtained on different regions of phase space. This time, the analysis is driven by the semiclassical Hamiltonian flow whose global structure has to be computed. The semiclassical flow has a more complicated structure than the classical flow and in particular it contains a set of trapped trajectories. It appears that a global estimates can be obtained by gluing an estimate  near radial points on the horizon (using \cite{vasy2013microlocal}), an estimate in a small region inside the black hole (using a semiclassical version of the hyperbolic estimate), an estimate near radial points at spatial infinity (using \cite{vasy2020limiting}) and an estimate near the normally hyperbolic trapped set (based on \cite{wunschZworski} and \cite{dyatlov2016spectral}). 

Point \ref{holomorphy} follows from the resolvent identity $R(\sigma)-R(\sigma') = R(\sigma)\left(\hat{T}_s(\sigma')-\hat{T}_s(\sigma)\right)R(\sigma')$ once the mapping properties of $R(\sigma)$ (resulting from the global Fredholm estimate) have been clarified.

Concerning point \ref{regularityReal}, we obtain the high regularity of $R(\sigma)$ on the real axis outside of $\sigma= 0$ by using repeatedly the identity of the resolvent. The number of iterations is only limited by the regularity of the initial data (which is assumed to be high in this work). This contrast with the situation at zero where the number of iterations is limited by the spatial decay of $\hat{f}(\sigma)$. For Theorem \ref{mainTheoWeak}, the regularity we get this way and the resolvent's bounds we have obtained at zero and at infinity are sufficient to conclude after taking the inverse Fourier transform. Under the hypotheses of theorem \ref{mainTheoStrong}, we are able to go further (up to $2\left|s\right|+2$ iterations) using that $\hat{f}(\sigma)$ has more spatial decay. We obtain the expression of the higher order singularity at zero as $\sigma^{2\left|s\right|+2}R(\sigma)w$ with $w$ explicit independent of $\sigma$. Note that this step requires a precise knowledge of the kernel and cokernel of $\hat{T}_s(0)$ (in weaker spaces than the ones on which we have invertibility) which we compute using theory of hypergeometric equations. Adapting an idea of \cite{hintz2022sharp} which roughly consists in performing a last iteration with $\hat{T}_s(0)$ replaced by the effective normal operator $N_{\eff}(\hat{T}_s(\sigma))$ (see Definition \ref{defNormalOp}), which governs the transition between $\hat{T}_s(\sigma)$ and $\hat{T}_s(0)$ near $x =0$, we are then able to compute explicitly $\sigma^{2\left|s\right|+2}R(\sigma)w$ modulo terms, which are irrelevant for the late time expansion.

After all the points \ref{inversibility}-\ref{regularityReal} are achieve, we are ready to perform a contour deformation argument to write:
\begin{align*}  
v(\mathfrak{t}) = \frac{1}{2\pi}\int_{\R} e^{-i\sigma\mathfrak{t}} R(\sigma)\hat{f}(\sigma)\dd \sigma.
\end{align*}
We then use standard properties of the Fourier transform (together with careful Fourier transform computations) to conclude the proof. 
\subsection{Outline of the paper}
\begin{itemize}
\item In Sections \ref{geometricFramework} and \ref{analyticFramework}, we introduce the notations and classical properties of pseudodifferential operators needed for the proof. 
\item In Section \ref{secCauchy} we state the propositions needed to rephrase the Cauchy problem as a forcing problem while keeping track of the spatial decay of the forcing term.
\item In Section \ref{secFlow}, we compute precisely the classical and semiclassical Hamiltonian flows of the Teukolsky operator. As mentioned earlier, the structure of these flow is paramount to obtain global Fredholm and semiclassical estimates. 
\item In Section \ref{secFred} we first get the Fredholm and semiclassical estimates on the different problematic regions of phase space. Subsection \ref{subSecEstHorizon} handles the region near the horizon (with both Fredholm and semiclassical estimates), Subsection \ref{subSecEstSpatInf} handles the region near $r = +\infty$ (with two Fredholm estimates: one uniform for $\sigma$ in a compact not containing zero and one uniform for $\sigma$ in a small neighborhood of zero, and a semiclassical estimate) and Subsection \ref{subSecTrapp} handles the region near the trapped set in the semiclassical regime. In the beginning of Subsection \ref{subSecGlob}, we explain how a global Fredholm estimate implies the Fredholm property. We then obtain global Fredholm estimates and global semiclassical estimates. Subsection \ref{subSecIndexZero} is dedicated to the zero index property. The heart of the proof is the continuity of the index, but some care is required since spaces on which the operator is Fredholm depends on $\sigma$. The case of $\hat{T}_s(0)$ is treated separetly in Subsection \ref{subSecHatTZero} since we also need to compute its kernel and cokernel in weaker spaces to later compute the higher order singularity at $\sigma = 0$.
\item Section \ref{secAnalRes} contains the analysis of the resolvent which includes its mapping properties, uniform bounds when $\sigma\to 0$ and $\left|\sigma\right|\to +\infty$ and the precise computation of the higher order singularity at $\sigma= 0$.
\item In section \ref{secContour}, we use the previous result about the resolvent to perform the contour deformation argument and we finish the proof of the main theorems. 
\item Appendix \ref{semiclassicalHyperboEstSection} presents a quite general version of the semiclassical hyperbolic estimate that we use in the proof. In Appendix \ref{absenceOfKernel}, we obtain the absence of kernel for a class of operator including the effective normal operator associated to the Teukolsky equation. Finally, Appendix \ref{ProofCauchyConormal} presents in detail the energy estimate used to translate the Cauchy problem into a forcing problem.
\end{itemize}

\subsection*{Acknowledgments.} This work was carried out during my Ph.D thesis at the Institut Fourier. I am very grateful to Dietrich Häfner for numerous fruitful discussions and for his proofreading and corrections. I would also like to thank Peter Hintz for his valuable comments and suggestions.

\section{Geometric framework}\label{geometricFramework}
\subsection{Kerr metric}\label{KerrMetric}
The exterior of a rotating black hole of mass $M$ and of angular momentum $a$ with $\lvert a\rvert<M$ is described by the Kerr solution to the Einstein equation \cite{kerr1963gravitational}. We are interested in the exterior region, modeled by the smooth manifold $\mathcal{M} := \R_t \times (r_+, +\infty)\times \mathbb{S}^2$ where $r_+ := M+\sqrt{M^2-a^2}$.
The manifold $\mathcal{M}$ is endowed with the Kerr metric given in Boyer-Lindquist coordinates by:
\begin{align*}
g=& \frac{\Delta_r - a^2\sin^2\theta}{\rho^2}\dd t^2 + \frac{4Mar\sin^2\theta}{\rho^2}\dd t \dd\phi - \frac{\rho^2}{\Delta_r}\dd r^2\\
& - \rho^2\dd \theta^2 - \frac{\sin^2\theta}{\rho^2}((a^2+r^2)^2-a^2\Delta_r\sin^2\theta)\dd \phi^2\\
\end{align*}
where
\begin{align*}
\rho^2 :=& r^2+a^2\cos^2\theta\\
\Delta_r :=& r^2-2Mr+a^2
\end{align*}

Note that $\Delta_r$ vanishes when $r=r_+$ (the other root is $r_- := M-\sqrt{M^2-a^2}$). However, the singularity is merely a coordinate singularity. We introduce $\text{Kerr}_*$ coordinates: 
\begin{align*}(t_*, r, \theta, \phi_*) = (t + T(r), r, \theta, \phi + A(r))
\end{align*} where $T(r):= \int_{3M}^r \frac{a^2+r^2}{\Delta_r} \dd r$ and $A(r)=\int_{3M}^r \frac{a}{\Delta_r}\dd r$. In these coordinates, the metric $g$ writes:
\begin{align*}
g =& \frac{\Delta_r - a^2\sin^2\theta}{\rho^2}\dd t_*^2 -2\dd t_* \dd r + \frac{4Mar \sin^2\theta}{\rho^2}\dd t_* \dd \phi_*  + 2a\sin^2\theta\dd r \dd \phi_* - \rho^2\dd \theta^2\\
& - \frac{\sin^2\theta}{\rho^2}((a^2+r^2)^2-a^2\Delta_r\sin^2\theta)\dd \phi_*^2.
\end{align*}
Therefore, it can be extended analytically to a larger manifold.
Let $\epsilon>0$ such that $r_-<r_+ - \epsilon$. We consider here an extension of $\mathcal{M}\cap t^{-1}((0,+\infty))$ defined as $\mathcal{M}_\epsilon := \R_{t_*}\times (r_+-\epsilon, +\infty)\times \mathbb{S}^2_{\theta,\phi_*}$ endowed with the analytic extension of $g$. The future horizon is by definition $\mathfrak{H} := \R_{t_*}\times \left\{r_+\right\}\times \mathbb{S}^2_{\theta,\phi_*}$.

Let $\psi_1$, $\psi_2$ be smooth, non negative and monotonic functions with $\psi_1 =1$ on $\left\{r<3M\right\}$, $\psi_1 = 0$ on $\left\{r>4M\right\}$, $\psi_2=0$ on $\left\{r<5M\right\}$, $\psi_2 =1$ on $\left\{6M<r\right\}$.
In our analysis, we mainly use the smooth function $\mathfrak{t} := t_*-L(r)$ where the smooth function $L$ is defined by
\begin{align}\label{defLr}
L(r) := \begin{cases}1+T(r)+\int^{+\infty}_r \psi_1(r)\frac{a^2+r^2}{\Delta_r} \dd r+ \int_{-\infty}^r \psi_2(r)\frac{a^2+r^2}{\Delta_r}\dd r \text{ for $r>r_+$}\\
1 + \int_{3M}^{+\infty}\psi_1(r)\frac{a^2+r^2}{\Delta_r} \dd r \text{ for $r\leq r_+$}
\end{cases}.
\end{align} Note that $\mathfrak{t} = t_* + c_1$ on $\left\{r<3M\right\}$, $\mathfrak{t} = t-1$ on $\left\{4M<r<5M\right\}$ and $\mathfrak{t} = t-T(r)+c_2$ on $\left\{6M<r<+\infty\right\}$ where $c_1$ and $c_2$ are real constants which will play no role in the analysis. As a consequence of the definition, level sets of $\mathfrak{t}$ are smooth hypersurfaces transverse to the future horizon and transverse to null infinity. We also define $\varphi = \phi_* - A(r)(1-\psi_1(r))$ (so that $\varphi$ is equal to $\phi_*$ near the horizon and equal to $\phi$ far from the horizon).

We define $\Sigma^{\mathfrak{t}}_0:= \mathfrak{t}^{-1}(\left\{0\right\})$ and we have that $\mathcal{M}_{\epsilon}$ is diffeomorphic to $\R_{\mathfrak{t}}\times \Sigma^{\mathfrak{t}}_0$.
We introduce the coordinate $x := \frac{1}{r}$ and we define $\overline{\mathcal{M}_\epsilon}$ the smooth manifold with boundary $\R_{\mathfrak{t}}\times[0,\frac{1}{r_+-\epsilon})_x\times \mathbb{S}^2_{\theta,\varphi}$ and we call the boundary $\mathscr{I}^+:=\R_{\mathfrak{t}}\times \left\{0\right\}\times \mathbb{S}^2_{\theta,\varphi}$ future null infinity. Note that the metric does not extend smoothly to the boundary but blows up like $x^{-2}$. For this reason, we will often consider the conformally rescaled metric $\tilde{g}:= \rho^{-2}g$. We denote by $G$ the metric induced by $g$ on the cotangent bundle and by $\tilde{G}$ the metric induced by $\tilde{g}$.

\subsection{Spin weighted functions}
Let $s\in \frac{1}{2}\Z$. 
Let $U_N := \mathbb{S}^2 \setminus \left\{N\right\}$ and $U_S := \mathbb{S}^2 \setminus \left\{ S\right\}$ where $N$ is the north pole and $S$ the south pole. We define $\mathcal{B}_s$, the complex line bundle over $\mathbb{S}^2$ with projection $\pi$, with local trivializations $(U_S,\mathcal{T}_S)$ and $(U_N, \mathcal{T}_N)$ and transition function from $\mathcal{T}_S$ to $\mathcal{T}_N$ given by
\[
f^{S,N}_s: \begin{cases}
U_N \cap U_S\rightarrow GL_1(\C)\\
(\theta, \phi) \mapsto e^{-2is\phi}
\end{cases}
\]
In other words, for every $e \in \mathcal{B}_s$ with $\pi(e) = x\in U_N \cap U_S$, if we use the notation $(x, v_N) = \mathcal{T}_N(e)$ and $(x,v_S) = \mathcal{T}_S(e)$, we have $v_N = f_s^{S,N}(x)v_S$.
For conventional reason, we also introduce the (redundant) local trivialization on $U_N\cap U_S \setminus\left\{\phi = 0\right\}$, $\mathcal{T}_m$ which is such that the transition function from $\mathcal{T}_S$ to $\mathcal{T}_m$ is $f_s^{S,m}(\theta, \phi) = e^{-is\phi}$ (which has no continuous extension to $U_N\cap U_S$ if $s$ is not an integer).
 
We introduce the complex line bundle $\mathcal{B}(s,s):= \pi_3^*(\mathcal{B}_s)$ over $\overline{\mathcal{M}_\epsilon}$ where $\pi_3$ is the third projection ($\pi_3: \overline{\mathcal{M}_\epsilon}\rightarrow \mathbb{S}^2_{\theta, \varphi}$). We denote by $\tilde{\mathcal{T}}_{S}$ (respectively $\tilde{\mathcal{T}}_{N}$ and $\tilde{\mathcal{T}}_m$) the local trivialization of $\mathcal{B}(s,s)$ on $(\R_\mathfrak{t}\times [0, \frac{1}{r_+-\epsilon})_x\times U_S)$ (respectively on $(\R_\mathfrak{t}\times [0, \frac{1}{r_+-\epsilon})_x\times U_N)$ and $(\R_\mathfrak{t}\times [0, \frac{1}{r_+-\epsilon})_x\times U_m)$) induced by $\mathcal{T}_S$ (respectively $\mathcal{T}_N$ and $\mathcal{T}_m$). The sections of $\mathcal{B}(s,s)$ are spin-weighted functions. This definition may seem a bit artificial but will be enough for the purpose of this paper. For a more intrinsic and meaningful definition (but equivalent up to explicit isomorphism of vector bundle) and some explanation of how this bundle naturally appears in the study of tensorial equations, see \cite{millet2021geometric}.

If we have a vector bundle $E$ over a manifold with boundary (or with corners) $X$ (and $E'$ its dual), we denote by:
\begin{itemize}
\item $\Omega_X$ the bundle of densities over $X$.
\item $\Gamma(X,E)$ (or $\Gamma(E)$ for short) the space of smooth sections (smooth up to the boundary). 
\item $\Gamma^k(X,E)$ (or $\Gamma^k(E)$) the space of sections of regularity $C^k$ (up to the boundary).
\item $\Gamma(\overset{\circ}{X},E)$ (or $\Gamma^{\circ}(E)$ for short) the space of smooth sections on the boundle $\overset{\circ}{E}:= E_{|_{\overset{\circ}{X}}}$
\item $\dot{\Gamma}(X, E)$ (or $\dot{\Gamma}(E)$ for short) the space of smooth sections vanishing at infinite order at the boundary.
\item $\mathcal{D}'(\overset{\circ}{X},E)$ (or $\mathcal{D}'^{\circ}(E)$ for short) the space of distributional sections of $E_{|_{\overset{\circ}{X}}}$ (the dual space of $\Gamma_c^{\circ}(E'\otimes \Omega_X)$ where the index $c$ indicates that the sections have compact support in the base manifold, in this case $\overset{\circ}{X}$)
\item $\dot{\mathcal{D}}'(X,E)$ (or $\dot{\mathcal{D}}'(E)$ for short) the space of supported distribution (the dual space of $\Gamma_c(E^*\otimes\Omega_X)$)
\item $\mathcal{D}'(X,E)$ (or $\mathcal{D}'(E)$ for short) the space of extendible distributions (the dual space of $\dot{\Gamma}_c(E^*\otimes\Omega_X)$ the space of smooth compactly supported sections on $X$ vanishing up to infinite order at $\partial X$). The space of extendible distribution can be viewed as a quotient space of distributions on a larger manifold. It is naturally included in the space of distributions on $\overset{\circ}{X}$. Note that we have a natural application $\dot{\mathcal{D}}'(X)\rightarrow \mathcal{D}'(X)$ whose kernel is the space of distributions supported at the boundary. 
\end{itemize}
For more details about supported and extendible distributions, see \cite{hormander2007analysis} (appendix B).
If the boundary has several connected components, we can prescribe the extendible/ supported character of the distribution at each boundary component independently.

Note that since $\mathcal{B}(s,s)$ is a pullback of $\mathcal{B}_s$, we have a natural identification between elements of $\Gamma^{\circ}(\mathcal{B}(s,s))$ and $\C^{\infty}(\R_{\mathfrak{t}}\times (r_+-\epsilon, +\infty), \Gamma(\mathcal{B}_s))$ and between $\mathcal{D}'(\mathcal{B}(s,s))$ and $\mathcal{D}'(\R_{\mathfrak{t}}\times (r_+-\epsilon, +\infty), \mathcal{D}'(\mathcal{B}_s))$.

If $Y,X$ are smooth manifolds (with or without boundary) and $E$ is a vector bundle over $X$, we denote by $Y\times E$ the pullback bundle $\pi_2^*(E)$ over $Y\times X$ (where $\pi_2$ is the second projection on the product). With this notation, we have for example $\mathcal{B}(s,s) = \R_{\mathfrak{t}}\times [0, \frac{1}{r_+-\epsilon})_x \times \mathcal{B}_s$. When the context is clear, we sometimes omit the pullback in the notation, for example $\Gamma(Y\times X, E)$ will be a shortcut for $\Gamma(Y\times X, \pi_2^*(E))$.

Let $\Theta$ be the unique connection on $\mathcal{B}_s$ such that in trivialization $\mathcal{T}_m$:
\begin{align*}
\Theta_{\partial_\theta} =& \partial_\theta \\
\Theta_{\partial_\phi} =& \partial_\phi + is\cos\theta
\end{align*}
We also denote by $\Theta$ the pullback on $\mathcal{B}(s,s)$, concretely:
\begin{align*}
\Theta_{\partial_x} = \partial_x \\
\Theta_{\partial_\mathfrak{t}} = \partial_{\mathfrak{t}}
\end{align*}
This redundant notation is sometimes convenient to use Einstein's summation convention. To alleviate notations, we also sometimes omit the $\partial$. For example, $\Theta_{\mathfrak{t}}$ stands for $\Theta_{\partial_{\mathfrak{t}}}$.

If we have a Frechet space $F$, we define the Fourier transform $\mathcal{F}_{\mathfrak{t}}u$ of a section $u \in \mathscr{S}'\left(\R, F\right)$ by (for all $\phi \in \mathscr{S}(\R_{\mathfrak{t}})$):
\[
\mathcal{F}_{\mathfrak{t}} u (\phi) = u(\hat{\phi})
\]
where $\hat{\phi}(\sigma):= \int_{\R_{\mathfrak{t}}}e^{i\mathfrak{t}\sigma}\phi(\mathfrak{t})\dd \mathfrak{t}$.

\section{Analytic framework}
\label{analyticFramework}

The Teukolsky operator (originally introduced in \cite{teukolsky1973perturbations}) is a smooth second order differential operator on $\overset{\circ}{\mathcal{B}}(s,s)$. Its original expression was given in a local trivialization $\mathcal{T}_{kinner}$ defined on a dense open subset of $\mathcal{M}$ (this local trivialization is associated to the Kinnersley tetrad, in a way explained in \cite{millet2021geometric}). This local trivialization does not extend to a local trivialization on a dense open subset of $\mathcal{M}_{\epsilon}$, this is linked to the fact that the Kinnersley tetrad degenerate at the horizon. However, the Teukolsky operator itself extends analytically to $\mathcal{M}_{\epsilon}$ (it can also be defined more intrinsically directly on $\mathcal{M}_\epsilon$, see \cite{millet2021geometric} for details). The local trivialization $\mathcal{T}_m$ that we have introduced earlier (based on a renormalized tetrad which is smooth across the horizon) gives the following expression for $T_s$:
\begin{align}
(T_s)_{m} =& -a^2\sin^2\theta \partial_{t_*}^2 - 2a \partial_{t_*}\partial_{\phi_*} - \frac{1}{\sin^2\theta}\partial_{\phi_*}^2 - \frac{1}{\sin\theta}\partial_{\theta_*}\left(\sin\theta\partial_{\theta_*}\right) - \Delta_r^{-s}\partial_{r_*}\left(\Delta_r^{s+1}\partial_{r_*}\right) \notag \\
&- 2(a^2+r^2)\partial_{t_*}\partial_{r_*} - 2a\partial_{\phi_*}\partial_{r_*} + 4s(r-M)\partial_{r_*} - \frac{2is\cos\theta}{\sin^2\theta}\partial_{\phi_*}\notag\\
 &- 2\left((1-2s)r-ias\cos\theta\right)\partial_{t_*}+(s^2\cot^2\theta + s)\label{TeukolskyOperateurBasique}
\end{align}
We can check (by changing trivialization with $\mathcal{T}_S$ and $\mathcal{T}_N$) that this expression defines a unique second order smooth differential operator $T_s$ on $\overset{\circ}{\mathcal{B}}(s,s)$.
We can also write the operator in coordinates $(\mathfrak{t}, r, \theta, \varphi)$:
\begin{align}
(T_s)_{m} =& a_{\mathfrak{t},\mathfrak{t}}\partial_{\mathfrak{t}}^2 + a_{\mathfrak{t},\varphi}\partial_{\mathfrak{t}}\partial_{\varphi} + a_{\varphi,\varphi}\partial_{\varphi}^2 -\frac{1}{\sin\theta}\partial_{\theta}\sin\theta\partial_{\theta} - \Delta_r^{-s}\partial_r \Delta_r^{s+1}\partial_r + a_{\mathfrak{t},r}\partial_{\mathfrak{t}}\partial_r  \notag \\
 &+ a_{r,\varphi}\partial_r\partial_{\varphi} + 4s(r-M)\partial_r + a_{\varphi} \partial_{\varphi} + a_{\mathfrak{t}} \partial_{\mathfrak{t}} +s^2\cot^2\theta +s.
\end{align}
where the coefficients $a_{\mathfrak{t},\mathfrak{t}}$, $a_{\mathfrak{t}, \varphi}$, $a_{\varphi,\varphi}$, $a_{\mathfrak{t},r}$, $a_{r,\varphi}$, $a_\varphi$, $a_{\mathfrak{t}}$ are smooth, independent of $\mathfrak{t}$ and $\varphi$, with:
\[
\begin{array}{|c|c|c|c|}
\hline
 & r<3M & 4M<r<5M & 6M<r\\
\hline
a_{\mathfrak{t},\mathfrak{t}}& -a^2\sin^2\theta & \frac{(r^2+a^2)^2}{\Delta_r}-a^2\sin^2\theta & -a^2\sin^2\theta\\
\hline
a_{\mathfrak{t}, \varphi}& -2a & \frac{4Mar}{\Delta_r} & \frac{4Mar}{\Delta_r}\\
\hline
a_{\varphi,\varphi}& -\frac{1}{\sin^2\theta} &\frac{a^2}{\Delta_r}-\frac{1}{\sin^2\theta}&\frac{a^2}{\Delta_r}-\frac{1}{\sin^2\theta}\\
\hline
a_{\mathfrak{t},r}& -2(a^2+r^2) & 0& 2(a^2+r^2)\\
\hline
a_{r,\varphi} & -2a & 0 & 0\\
\hline
a_{\varphi} & -\frac{2is\cos\theta}{\sin^2\theta} & -2s\left(\frac{a(r-M)}{\Delta_r}+i\frac{\cos\theta}{\sin^2\theta}\right) & -2s\left(\frac{a(r-M)}{\Delta_r}+i\frac{\cos\theta}{\sin^2\theta}\right)\\
\hline
a_{\mathfrak{t}} & -2\left((1-2s)r-ias\cos\theta\right) & -2s\left(\frac{M(r^2-a^2)}{\Delta_r}-r-ia\cos\theta\right) & \frac{4sMa^2-2sr(a^2+r^2)}{\Delta_r}+2(s+1)r+2ias\cos\theta\\
\hline
\end{array}
\]
With the previous expression, we see that 
\begin{align*}T_s\left(\mathscr{S}'\left(\R_{\mathfrak{t}}, \mathcal{D}'((r_+-\epsilon, +\infty)\times \mathcal{B}_s)\right)\right)\subset \mathscr{S}'\left(\R_{\mathfrak{t}}, \mathcal{D}'((r_+-\epsilon, +\infty)\times \mathcal{B}_s)\right).
\end{align*}

We are interested in the following equation
\begin{align}
T_s u = f
\end{align}

We define the differential operator $\hat{T}_s$ on $\R_\sigma\times (r_+-\epsilon, +\infty)_r\times \mathcal{B}_s$ as the operator such that for all $u\in\mathscr{S}'(\R_{\mathfrak{t}}, \mathcal{D}'_{(r_+-\epsilon, +\infty)_r\times \mathcal{B}_s})$:
\begin{align*}
\hat{T}_s \mathcal{F}_{\mathfrak{t}}u = \mathcal{F}_{\mathfrak{t}}(T_s u).
\end{align*}
If we call $\sigma$ the conjugate variable of $\mathfrak{t}$, we get:
\begin{align*}
\label{That}
\left(\hat{T}_s(\sigma)\right)_{m} =& -a_{\mathfrak{t},\mathfrak{t}}\sigma^2 - ia_{\mathfrak{t},\varphi}\sigma\partial_{\varphi} + a_{\varphi,\varphi}\partial_{\varphi}^2 -\frac{1}{\sin\theta}\partial_{\theta}\sin\theta\partial_{\theta} - \Delta_r^{-s}\partial_r \Delta_r^{s+1}\partial_r - ia_{\mathfrak{t},r}\sigma\partial_r\\
& + a_{r,\varphi}\partial_r\partial_{\varphi} + 4s(r-M)\partial_r + a_\varphi \partial_\varphi - ia_{\mathfrak{t}} \sigma +s^2\cot^2\theta +s.
\end{align*}
Therefore, $(\hat{T}_s(\sigma))_{\sigma}$ is a family of differential operators on $(r_+-\epsilon, +\infty)\times \mathcal{B}_s$. 
Note that if we call $T_s(\sigma)'$ the operator obtained with the Fourier transform with respect to $\mathfrak{t}$ replaced by the Fourier transform with respect to $t$, we have the relation: $\hat{T}_s(\sigma) = e^{-i\sigma T(r)}T_s(\sigma)'e^{i\sigma T(r)}$ on $r>6M$.

Finally, we introduce the semiclassical rescaling: $\hat{T}_{s,h}(z) := h^{2}\hat{T}_s(h^{-1}z)$ where $z = \frac{\sigma}{\lvert\sigma\rvert}$ and the small semiclassical parameter is $h= \lvert \sigma\rvert$. It will be used to study the regime $\Re(\sigma)\to +\infty$ with $\Im(\sigma)$ bounded.

\subsection{Sobolev spaces} \label{SobolevSpaces}

We define the spatial manifolds $X := (r_+-\epsilon, +\infty)_r \times \mathbb{S}^2$ and $\overline{X} := [0,\frac{1}{r_+-\epsilon})_x \times \mathbb{S}^2$.
We define the volume form $\dd vol := r^2|\sin\theta|\dd r \dd \theta \dd \varphi$ on $X$. The choice of this exact volume form is not crucial in this paper for the radial point estimates since the subprincipal symbol at a radial point does not depend on the volume form. However, we have less freedom for the computation of the subprincipal symbol at the trapped set (see Remark \ref{remarkVolumeForm}).
To define the $L^2$ norm on $(r_+-\epsilon, +\infty)\times\mathcal{B}_s$, we also need a hermitian metric $\mathbb{m}$ on $\mathcal{B}_s$. We define it using the local trivialization $\mathcal{T}_m$. On this trivialization, for any $x\in U_m$ and $z_1,z_2\in \C$, $\mathbb{m}_x(z_1,z_2) := z_1\overline{z_2}$. This metric on $(\mathcal{B}_s)_{|_{U_m}}$ extends uniquely to a smooth metric on $\mathcal{B}_s$. We often write $\mathbb{m}(z)$ for $\mathbb{m}(z,z)$.

We can now define the following scalar product on $\Gamma_c((r_+-\epsilon, +\infty)\times \mathcal{B}_s)$. For $u,v \in \Gamma_c((r_+-\epsilon, +\infty)_r\times\mathcal{B}_s)$:
\begin{align*}
\left< u, v\right> = \int_{X} \mathbb{m}_{r(x)}(u(x),v(x)) \dd vol(x).
\end{align*}
We define $L^2_{(r_+-\epsilon, +\infty)_r\times\mathcal{B}_s}$ as the completion of $\Gamma_c((r_+-\epsilon, +\infty)_r\times\mathcal{B}_s)$ for the associated norm. We remark that we have the natural inclusion $L^2((r_+-\epsilon, +\infty)_r\times\mathcal{B}_s)\subset \mathcal{D}'_{(r_+-\epsilon, +\infty)_r\times\mathcal{B}_s}$ (we also have the equality $L^2_{(r_+-\epsilon, +\infty)_r\times\mathcal{B}_s} = L^2(X, \dd vol)\otimes_{C^{\infty}(X)}\Gamma((r_+-\epsilon, +\infty)_r\times\mathcal{B}_s)$, where $\otimes_{C^{\infty}(X)}$ denotes the tensor product of $C^{\infty}(X)$-modules).   

We now define the $b$-tangent bundle and the associate notations.
\begin{definition}
Let $\mathcal{N}$ be a manifold with boundary of dimension $n$.
We define ${}^b T\mathcal{N}$, the bundle of $b$-vectors on $\mathcal{N}$ as the bundle whose sections are smooth vector fields tangent to the boundary. In local coordinates $(y_i)_{i=0}^{n-1}$ near a point of the boundary, if $y_0$ is a defining function of the boundary, such vector fields write $a_0(y)y_0\partial_{y_0}+\sum_{i=1}^{n-1} a_i(y)\partial_{y_i}$ with $(a_i)_{i=0}^{n-1}$ a family of smooth functions (smooth up to the boundary). The dual bundle is denoted by ${}^bT^*\mathcal{N}$. These definitions extend to the case where $\mathcal{N}$ is a manifold with corners by defining $b$-vector fields as smooth vector fields tangent to each boundary face. The bundle ${}^{b}T\mathcal{N}$ is then used to construct other $b$-objects such as $b$-metrics which are metrics on ${}^b T\mathcal{N}$ or $b$-volume forms which are non vanishing sections of $\Lambda^n\left({}^bT^*\mathcal{N}\right)$. For $k\in \N$, we denote by $\text{Diff}^{k}_b(E)$ the algebra of differential operators generated (as a $C^{\infty}(\mathcal{N})$-module) by the set $\left\{\text{Id}\right\}\cup \left\{\Theta_{X_1}...\Theta_{X_j}, j\leq k, X_i \in \Gamma({}^{b}T\mathcal{N})\right\}$.
\end{definition}

We also need to define the scattering tangent bundle.
\begin{definition}
Let $\mathcal{N}$ be a manifold with boundary of dimension $n$. Let $y_0$ be a boundary defining function.
We define ${}^{sc} T\mathcal{N}$, the bundle of $sc$-vectors on $\mathcal{N}$ as the bundle whose sections are of the form $y_0Z$ for $Z\in \Gamma({}^bT\mathcal{N}$. In local coordinates $(y_i)_{i=0}^{n-1}$ near a point of the boundary, such vector fields write $a_0(y)y_0^2\partial_{y_0}+\sum_{i=1}^{n-1} a_i(y)y_0\partial_{y_i}$ with $(a_i)_{i=0}^{n-1}$ a family of smooth functions (smooth up to the boundary). The dual bundle is denoted by ${}^{sc}T^*\mathcal{N}$. The bundle ${}^{sc}T\mathcal{N}$ is then used to construct other $sc$-objects such as $sc$-metrics which are metrics on ${}^{sc} T\mathcal{N}$ or $sc$-volume forms which are non vanishing sections of $\Lambda^n\left({}^{sc}T^*\mathcal{N}\right)$. For $k\in \N$, we denote by $\text{Diff}^{k}_{sc}(E)$ the algebra of differential operators generated (as a $C^{\infty}(\mathcal{N})$-module) by the set $\left\{\text{Id}\right\}\cup \left\{\Theta_{X_1}...\Theta_{X_j}, j\leq k, X_i \in \Gamma({}^{sc}T\mathcal{N})\right\}$ where $\Theta$ is any connection on $E$.
\end{definition}

\begin{definition}
Let $\mathcal{N}$ be a compact manifold with boundary of dimension $n$ with a boundary defining function $y_0$. Let $E$ be a vector bundle over $\mathcal{N}$ with connection $\Theta$ and a metric $\mathbb{m}$. We fix a finite family of smooth vector fields $(Z_i)_{i=1}^N$ which generate ${}^bT\mathcal{N}$ as a $C^{\infty}(\mathcal{N})$ module and we fix a $b$ volume form.
We define $H^0_b(E)$ as the space $L_b^2(E)$ where the index $b$ on $L^2$ indicates that integration is performed against the $b$-volume form. For $\tilde{r}\in \N$, $\tilde{r}\geq 1$, we define recursively the $b$-Sobolev space $H_b^{\tilde{r}}(E)$ and its semiclassical version $H_{b,h}^{\tilde{r}}(E)$ by completion of $\Gamma_c^{\circ}(E)$ for the norms:
\begin{align*}
\left\|u\right\|_{H^{\tilde{r}+1}_b} :=& \left\|u\right\|_{H^{\tilde{r}}_b}+\sum_{i=1}^N \left\|\Theta_{Z_i} u\right\|_{H^{\tilde{r}}_{b}(E)}\\
\left\|u\right\|_{H^{\tilde{r}+1}_{b,h}} :=& \left\|u\right\|_{H^{\tilde{r}}_{b,h}}+\sum_{i=1}^N \left\|h\Theta_{Z_i} u\right\|_{H^{\tilde{r}}_{b,h}(E)}\\
\end{align*}
We then define $H_b^{\tilde{r}}(E)$ and $H_{b,h}^{\tilde{r}}(E)$ for $\tilde{r}\in \R$ by interpolation and duality. For $l\in \R$, we then define $H_b^{\tilde{r},l}:= y_0^{l}H_b^{\tilde{r}}$ and $H_{b,h}^{\tilde{r},l}:= y_0^{l}H_{b,h}^{\tilde{r}}$.

Similarly, we can define the spaces $H_{sc}^{\tilde{r}, l}$ and $H_{sc,h}^{\tilde{r},l}$ for $\tilde{r},l\in \R$ by replacing all the $b$ indices by $sc$ indices in the previous definition.
\end{definition}

Note that the previous definition makes sense even if the manifold has no boundary (then ${}^{b}T\mathcal{N} = {}^{sc}T\mathcal{N} = T\mathcal{N}$). We use it to define the Sobolev space $H^k(\mathcal{B}_s)$.

When the boundary has several connected component, we can specify the behavior of the section near each component of the boundary independently. In particular, in the case of sections of $[0, \frac{1}{r_+-\epsilon}]\times \mathcal{B}_s$, we define the following hybrid spaces:
\begin{definition}
We can view distributions in $\dot{\mathcal{D}}'((0,\frac{1}{r_+-\epsilon}]_x\times \mathcal{B}_s)$ as elements of $\mathcal{D}'((0,\frac{1}{r_+-\epsilon}+1)_x\times \mathcal{B}_s)$ (therefore, it makes sense to ask whether or not some $u \in \dot{\mathcal{D}}'((0,\frac{1}{r_+-\epsilon}]\times \mathcal{B}_s)$ is in $H^{k,l}_{b}([0,\frac{1}{r_+-\epsilon}+1]_x\times \mathcal{B}_s$). We denote by $\dot{H}_{b}^{k,l} := H^{k,l}_{b}\cap \dot{\mathcal{D}}'((0,\frac{1}{r_+-\epsilon}]_x\times \mathcal{B}_s)$ endowed with the induced norm.
We perform the same construction (with indices $b$ replaced by $sc$) to get $\dot{H}^{\tilde{r},l}_{sc}$. We also get the semiclassical version of these spaces by adding an index $h$ next to the index $b$ (or $sc$) in the definition.
\end{definition}

\begin{definition}
For distributions in $u\in\mathcal{D}'(\left(0,+\frac{1}{r_+-\epsilon}\right]_x\times \mathcal{B}_s)$, we denote by $\text{Ext}(u)$ the set of $\tilde{u} \in \mathcal{D}'((0,+\infty)_x\times \mathcal{B}_s)$ with support contained in some $x^{-1}((0,C))$ such that $\tilde{u}_{|_{(0,\frac{1}{r_+-\epsilon})}} = u$. 
Therefore, it makes sense to define the (possibly infinite) norm:
\begin{align*}
\left\|u\right\|_{\overline{H}^{k,l}_{b}} := \inf_{\tilde{u}\in \text{Ext}(u)} \left\|\tilde{u}\right\|_{H^{k,l}_{b}([0,C]_x\times\mathcal{B}_s)}
\end{align*}
We perform the same construction (with indices $b$ replaced by $sc$) to get $\overline{H}^{\tilde{r},l}_{sc}$. We also get the semiclassical version of these spaces by adding an index $h$ next to the index $b$ in the definition.
\end{definition}

To be consistent with the convention used in \cite{vasy2020limiting} and \cite{vasy2020resolvent} (where a scattering volume form is used in the definition of $b$-Sobolev spaces), it is useful to introduce $\overline{H}^{\tilde{r},l}_{(b)}:= \overline{H}_b^{\tilde{r},l+\frac{3}{2}}$ and $\dot{H}^{\tilde{r},l}_{(b)} = \dot{H}^{\tilde{r},l+\frac{3}{2}}_{b}$. If there is no parenthesis, it means that the definition are the same without this multiplication by $x^{\frac{3}{2}}$.

Note that the spaces $\overline{H}^{k,l}_{(b)}$ and $\dot{H}^{-k,-l}_{(b)}$ are dual to each other (using the volume form $\dd vol$ and the metric $\mathbb{m}$ for the identifications).  

For $k\in \R$, we define the Mellin transform of a function $f \in \dot{C}^{\infty}([0,+\infty), H^k(\mathcal{B}_s))$ which is Schwartz at infinity by 
\begin{align*}
Mf(\lambda) := \int_{0}^{+\infty} x^{-i\lambda}f(x)\frac{\dd x}{x}
\end{align*}
By Plancherel formula, $M$ can be extended to an isomorphism from $L^2_b([0,+\infty); H^k(\mathcal{B}_s))$ to\\ $L^2(\R_{\lambda}, H^k(\mathbb{S}^2; \mathcal{B}_s))$. Also note that $\lambda Mf(\lambda) = M(x\partial_x f)(\lambda)$ and therefore, $M$ extends to an isomorphism between $H^k_b([0,+\infty]_x\times\mathcal{B}_s)$ (where the boundary defining function of $x=+\infty$ is $x^{-1}$) and $(1+|\lambda|^2)^{-\frac{k}{2}}L^2(\R_{\lambda}, L^2(\mathbb{S}^2; \mathcal{B}_s))\cap L^2(\R_{\lambda}, H^k(\mathbb{S}^2; \mathcal{B}_s))$.

\subsection{Pseudodifferential algebras}
In this section, we introduce the definitions and properties of the various Pseudodifferential algebras involved. Since we do not provide proofs for standard properties of pseudodifferential operators and we refer to \cite[Chapter~4]{melrose1993atiyah}, \cite[Appendix~A]{hintz2018global}, \cite[Chapter~XVIII]{hormander2007analysis}, \cite{vasy2018minicourse}, \cite{zworski2022semiclassical} for details and proofs.

\subsubsection{Definitions}
Let $E$ be a vector bundle over a manifold $\mathcal{N}$ (of dimension $n$) and $F$ be a vector bundle over $\mathcal{N}'$. We use the notation $E\boxtimes F$ to denote the vector bundle over $\mathcal{N}\times \mathcal{N}'$ with fiber $E_x\otimes F_y$ over $(x,y)\in \mathcal{N}\times \mathcal{N}'$.
We use the notation $\Omega^s(E)$ to denote the real line bundle over $\mathcal{N}$ whose fibers are $\Omega^s(E)_x = \left\{ u:(\Lambda^n E_x)\setminus\left\{0\right\}\rightarrow \R, \forall t\in \R\setminus \left\{0\right\}, u(t\alpha) = |t|^s u(\alpha)\right\}$. Note that $\Omega^1(T\mathcal{N}) = \Omega(T\mathcal{N})$ is the usual density bundle over $\mathcal{N}$.

We define the fiber radial compactification of $E$ as the bundle obtained by adding a boundary at fiber infinity. More precisely, for $\mathfrak{m}$ a smooth positive definite metric on $E$, a boundary defining function of fiber infinity is $(x,\xi)\mapsto \sqrt{\mathfrak{m}_x(\xi, \xi)}^{-1}$. This definition does not depend on the choice of $\mathfrak{m}$ as two different metrics gives locally equivalent boundary defining functions. If $\mathcal{N}$ has a boundary, the definition is the same (but we require that $\mathfrak{m}$ is a smooth positive metric up to the boundary). If the boundary of $\mathcal{N}$ is not empty, the resulting bundle is a manifold with corners. We denote by $\overline{T}^*X$ the fiber radial compactification of $T^*X$ and by $S^*X$ fiber infinity.

Let $G$ be an other vector bundle over $\mathcal{N}$.
\begin{definition}
We define the space $S^m(E, G)$ of $G$-valued symbols of order $m$ on $E$ as the set of functions $u\in C^{\infty}(E,G)$ such that the following estimates hold: For all open subset $U$ on which $E$ is trivial, if $(y, \xi)$ are local coordinates and $K$ is a compact subset of $U$, 
for all $\alpha,\beta\in \N^{n}$, there exists $C>0$ such that (uniformly on $K$)
\begin{equation}
\label{estimateSymbol}
|\partial_y^{\alpha}\partial_{\xi}^{\beta} u(x,\xi)|_G\leq C\left<\xi\right>^{m-|\beta|}
\end{equation}
where $\left|.\right|_G$ is computed using a fixed (but arbitrary) metric on $G$. We similarly define the space of semiclassical symbols of order $m$ $S_h^{m}(E,G)$ as the set of smooth $h$-indexed families $(u_h)_{h\in[0,1)} \in C^{\infty}([0,1)_h, S^{m}(E,G))$ such that for all $h \in [0,1)$, $u_h \in C^{\infty}(E,G)$ and estimate \eqref{estimateSymbol} holds uniformly with respect to $h \in [0,1)$.
\end{definition}

In the case where $\mathcal{N}$ has a boundary with boundary defining function $y_0$, we also define the following space of symbols:
\begin{definition}
For $l\in \R$, we denote by $S^{m,l}(E,G)$ the space of $G$-valued symbols of order $(m,l)$ on $E$ as the set of functions $u\in C^{\infty}(\overset{\circ}{E},G)$ such that the following estimates hold: For all open subset $U$ of $\mathcal{N}$ on which $E$ is trivial, if $(y, \xi)$ are local coordinates (with $y_0$ the boundary defining function if $U$ intersect the boundary) and $K$ is a compact subset of $U$, 
for all $\alpha,\beta\in \N^{n}$, there exists $C>0$ such that (uniformly on $K$)
\begin{equation}
\label{estimateSymbol2}
|\partial_y^{\alpha}\partial_{\xi}^{\beta} u(x,\xi)|_G\leq C\left<\xi\right>^{m-|\beta|}y_0^{l-\alpha_0}
\end{equation}
where $\left|.\right|_G$ is computed using a fixed (but arbitrary) metric on $G$. We similarly define the space of semiclassical symbols of order $(m,l)$ $S_h^{m,l}(E,G)$ as the set of $h$-indexed families $(u_h)_{h\in (0,1)}$ such that for all $h \in (0,1)$, $u_h \in C^{\infty}(\overset{\circ}{E},G)$ and estimate \eqref{estimateSymbol2} holds uniformly with respect to $h \in (0,1)$.
\end{definition}

We recall the Schwartz kernel theorem applied to our setting (see theorem 4.14 in \cite{melrose1993atiyah} and theorem 5.2.1 in \cite{hormanderI}). In the statement of the theorem, $L$ denotes the set of continuous linear operators from $C^{\infty}_c(X, \mathcal{B}_s)$ to $\mathcal{D}'(X, \mathcal{B}_s)$.
\begin{theorem}
There is a one to one correspondence between $L$ and the set of distributions 
$\mathcal{D}'\left(X\times X, \mathcal{B}_s\boxtimes \left(\mathcal{B}'_s\otimes \Omega(TX)\right)\right)$. The correspondence is given by:
\[
S:\begin{cases}
\mathcal{D}'\left(X\times X, \mathcal{B}_s\boxtimes \left(\mathcal{B}'_s\otimes \Omega(TX)\right)\right) \rightarrow L\\
A \mapsto \left(u \in C^{\infty}_c(X, \mathcal{B}_s) \mapsto \left(v \mapsto <A, v\otimes u>\right)\right)
\end{cases}
\]
The distribution $A$ is the Schwartz kernel of the operator $S(A)$.
\end{theorem}

\begin{definition}
We denote by ${}^bT\overline{X}$ the vector bundle of $b$-vector fields on $\overline{X}$ (we denote by ${}^bT^*\overline{X}$ it dual). We denote by ${}^b\overline{T}^* \overline{X}$ the fiber radial compactification of ${}^b T^* \overline{X}$ and by ${}^b S^* \overline{X}$ fiber infinity.
We define the set $\overline{L}:=S\left(\mathcal{D}'\left(\overline{X}\times \overline{X}, \mathcal{B}_s\boxtimes \left(\mathcal{B}'_s\otimes \Omega({}^bT\overline{X})\right)\right)\right)$. We have the natural inclusion $\overline{L}\subset L$.
\end{definition}
We define $\overline{X}^2_b$ as the blow-up of the manifold with corner $\overline{X}^2$ with respect to the corner $\partial \overline{X}\times \partial \overline{X}$ with the blow-down map $\beta$. Concretely if we use primes to denote coordinates on the second factor and $\omega$ to denote coordinates on $\mathbb{S}^2$, $\overline{X}^2_b$ is obtained by the introduction of coordinates $(\rho := x+x',\tau := \frac{x-x'}{x+x'}, \omega, \omega')$ on $\overline{X}^2\setminus \left\{x=x'=0\right\}$ and the addition of the front face which corresponds to $\left\{0\right\}_{\rho}\times[-1,1]_{\tau}\times \mathbb{S}^2_{\omega}\times \mathbb{S}^2_{\omega'}$ in these new coordinates. The blow-down map is then the identity map on $\overline{X}^2\setminus \left\{x=x'=0\right\}$ and the map $(0,\tau,\omega,\omega')\mapsto ((0,\omega), (0, \omega'))\in (\partial\overline{X})^2$ on the front face (see \cite[Chapter~4]{melrose1993atiyah} for details about this construction). We denote by $\Delta$ the diagonal in $\overline{X}^2$. We define $\Delta_b:= \overline{\beta^{-1}\left(\Delta\setminus \partial {\overline{X}}\times\partial {\overline{X}}\right)}$ and the front face $f= \beta^{-1}(\partial {\overline{X}}\times\partial {\overline{X}})$.

\begin{definition}
If $\mathcal{E}$ is a vector bundle of order $m$ over $\overline{X}^2_b$, we say that a distribution $A\in \mathcal{D'}(\overline{X}^2_b, \mathcal{E})$  is conormal of order $m$ to $\Delta_b$ it is smooth on $\overline{X}^2_b\setminus \Delta_b$ and if for every $y \in \Delta_b$, there exist:
\begin{itemize}
\item Local coordinates $(\alpha^{i})_{i=0}^{5}$ and a local trivialization of $\mathcal{E}$ on a neighborhood $U$ of $y$ such that $\Delta_b\cap U = \left\{\alpha^{0}=\alpha^{1} = \alpha^{2} = 0\right\}$. We use the notation $\alpha' := (\alpha^0, \alpha^1, \alpha^2)$ and $\alpha'' := (\alpha^3, \alpha^4, \alpha^5)$.
\item A symbol $a\in S^{m}(\R^3\times \R^3, \R^m)$ such that $A (\alpha',\alpha'')= (2\pi)^{-3}\int_{\R^3} e^{i\xi\cdot \alpha'} a(\xi,\alpha'')\dd \xi$ on $U$ (where the local trivialization is used implicitly to identify $A$ with a $\R^m$-valued function.
\end{itemize}
\end{definition}
It can be checked (see for example \cite[Theorem~18.2.9]{hormander2007analysis}) that the element $a$ in the definition is invariantly defined in $S^{m}(N^*(\Delta_b), \mathcal{E}_{|_{\Delta_b}}\otimes \Omega(N^*\Delta_b))$ modulo $S^{m-1}(N^*(\Delta_b)\otimes\mathcal{E}_{|_{\Delta_b}}\otimes \Omega(N^*\Delta_b))$.
\begin{definition}
Following \cite{melrose1993atiyah}, we define the (small) algebra of b-operators $\Psi_b$ as the subset of $\overline{L}$ of operators whose Schwartz kernels (lifted to $\overline{X}^2_b$) are conormal to $\Delta_b$, vanish up to infinite order at $\partial \overline{X}^2_b \setminus f$ and are properly supported in $\overline{X}^2_b$.
\end{definition}
For a conormal kernel $A$ of order $m$ (we write $A\in \Psi^{m,0}_b$), the principal symbol is therefore given by an element of 
\[S^m\left(N^*\Delta_b, \left(\mathcal{B}_s\boxtimes \left(\mathcal{B}'_s\otimes \Omega({}^bT\overline{X})\right)\right)_{|\Delta_b}\otimes \Omega(N^*\Delta_b)\right)/S^{m-1}.\] We can use the map 
\[I:\begin{cases}N^*\Delta_b\rightarrow {}^bT^*\overline{X}\\
\alpha \mapsto \left(v \mapsto \alpha(v,0)\right)
\end{cases}\]
with inverse
\[
I^{-1}: \begin{cases}{}^bT^*\overline{X} \rightarrow N^*\Delta_b\\
\beta \mapsto \left((v,v')\mapsto \beta(v)-\beta(v')\right) \end{cases}
\]
to identify $N^*\Delta_b$ and ${}^bT^*\overline{X}$. Moreover, we have canonical isomorphisms 
$\Omega({}^bT\overline{X})\otimes \Omega({}^bT^*\overline{X})= \R \times \overline{X}$ and $\mathcal{B}_s\otimes \mathcal{B}'_s = \C\times \overline{X}$.
Using the identification $I$, the symbol of $A$ can be seen as an element of $S^m({}^bT^*\overline{X})/S^{m-1}$ which is called the principal symbol of the operator. 

Following \cite[Appendix~A3]{hintz2018global}, we introduce the space $\Psi^{m,0}_{b,h}$ of $b$ semiclassical operators Note that in this case, the principal symbol is invariantly defined in $S^m_h({}^bT^*\overline{X})/hS^{m-1}_h$.

Eventually, we define $\Psi^{m,l}_b$ as the set $x^{-l}\Psi^{m,0}_b$ (and similarly $\Psi^{m,l}_{b,h} := x^{-l}\Psi^{m,0}_{b,h}$).

We now define the local model for scattering pseudodifferential operators.
We consider the manifold $\mathcal{N} := \R_y^3$ and the trivial bundle $E =\C\times \mathcal{N}$. Let $\overline{\R^3_y}$ be the radial compactification of $\R^3_y$. For $p\in S^{m,l}(\R^3_{\xi}\times\overline{\R^3_y},\mathcal{L}(E,E))$, we define the operator $\text{Op}(p)$ by its action on $\phi \in C^{\infty}_c(\R^3, \C)$:
\begin{align*}
\text{Op}(p)u(y) = (2\pi)^{-3}\int e^{i\xi\cdot (y-y')}p(y, \xi)u(y') \dd y' \dd \xi
\end{align*}
We call $\Psi_{sc}^{m,l}(\C\times\mathcal{N})$ the operators obtained by this procedure.
We similarly define $Op_h(p_h)$ for a symbol in $p_h\in S^{m,l}(\R^3_{\xi}\times\overline{\R^3_y})$ by:
\begin{align*}
\text{Op}_h(p_h)u(y) = (2\pi h)^{-3}\int e^{ih^{-1}\xi\cdot (y-y')}p_h(y, \xi)u(y') \dd y' \dd \xi.
\end{align*} 
The corresponding operator space will be denoted by $\Psi_{sc,h}^{m,l}(\C\times \mathcal{N})$.

We can then use this local model to define the algebra of scattering pseudodifferential operators on $[0,\frac{1}{r_+-\epsilon})\times\mathcal{B}_s$ (see \cite[Section~5.3.2]{vasy2018minicourse} for more details about this construction in a more general context):
\begin{definition}
We define $\Psi^{m,l}_{sc}([0,\frac{1}{r_+-\epsilon})\times\mathcal{B}_s)$ ($b$-pseudodifferential operators) as the set of continuous linear operators $A:C^{\infty}_c(X, \mathcal{B}_s)\rightarrow \mathcal{D}'(X, \mathcal{B}_s)$ such that: 
\begin{itemize}
\item The Schwartz kernel of $A$ is properly supported.
\item For all $\chi_1, \chi_2$ smooth and compactly supported on some open sets of trivialization for $[0,\frac{1}{r_+-\epsilon})\times\mathcal{B}_s$, $\chi_1 A\chi_2 = \Psi_{sc}^{m,l}(\C\times \mathcal{N})$.
\end{itemize}
The definition of $\Psi^{m,l}_{sc,h}([0,\frac{1}{r_+-\epsilon})\times\mathcal{B}_s)$ is obtained by replacing $A$ by a family $(A_h)_{h\in(0,1)}$ of continuous linear operators from $C^{\infty}_c(X, \mathcal{B}_s)$ to $\mathcal{D}'(X, \mathcal{B}_s)$ and $\Psi_{sc,h}^{m,l}(\C\times \mathcal{N})$ by $\Psi_{sc, h}^{m,l}(\C\times \mathcal{N})$ in the previous definition.
\end{definition}

We denote by $\Psi^{m,l}_{b,c}$ the set of operator with Schwartz kernel supported on $\mathcal{M}_{\epsilon-\eta}^2$ for some $0<\eta<\epsilon$. Standard theory of b-pseudodifferential operators (see for example \cite[Paragraph 5.9]{melrose1993atiyah}) provides the composition rule: 
Let $A\in \Psi^{m,l}_{b,c}$ with principal symbol $a$ and $B\in \Psi^{m',l'}_{b,c}$ with principal symbol $b$, then $AB \in \Psi^{m+m',l+l'}_{b,c}$ and has principal symbol $ab$. The fact that we restrict to operators with Schwartz kernel supported in $\mathcal{M}_{\epsilon-\eta}^2$ enables to reduce to the case of b-pseudodifferential operators on a compact manifold with boundary. Sometimes, we also want to take the composition of an operator $A\in \Psi^{m,l}_{b,c}$ and a differential operator $B\in x^{l'}\text{Diff}_b^{m'}$. In this case we can use the locality of $B$ to make sense of the compositions $AB$ and $BA$ as elements of $\Psi^{m+m',l+l'}_{b,c}$. Indeed, for any cutoff $\chi$ equal to $1$ on $(r_+-\epsilon+\eta, +\infty)_r$ and with compact support in $(r_+-\epsilon, +\infty)_r$, $A\chi B\chi = AB$ and $\chi B \chi A = BA$ by locality of $B$. These compositions properties will be enough for our purpose. We have the same composition rules for semiclassical $b$-operators, scattering operators and semiclassical scattering operators (see for example \cite{vasy2018minicourse}).

We say that an operator in $\Psi^{m,l}_b$ is elliptic if its principal symbol $p\in x^l S^m({}^bT^*\overline{X})/S^{m-1}$ has an inverse in other words if there exists $q\in x^{-l} S^{-m}({}^bT^*\overline{X})/S^{-m-1}$ such that $pq = qp = [1]$ where $[1]$ is the class of $1$ in $S^0({}^bT^*\overline{X})/S^{-1}$. We can localize this definition near a point $(x,\xi) \in {}^bS^* \overline{X}$: we say that $P\in \Psi^{m,l}$ with principal symbol $p\in x^{l} S^{m}({}^bT^*\overline{X})/S^{m-1}$ is elliptic at $(x,\xi)$ if there exists $q \in x^{-l}S^{-m}({}^bT^*\overline{X})/S^{m-1}$ such that $pq = qp = [g]$ with $g = 1$ on a neighborhood of $(x,\xi)$. The subset of ${}^bS^*\overline{X}$ at which $P$ is elliptic is denoted $\text{Ell}(P)$. The complementary subset (in ${}^bS^*\overline{X}$) is the characteristic set denoted by $\text{Char}(P)$. The wavefront set of $P$ denoted $WF(P)$ is defined negatively:  $(x,\xi)\in {}^bS^*\overline{X}$ is not in the wavefront set of $P$ if there exists $A \in \Psi^{0,0}_{b,c}$ with $(x,\xi)\in \text{Ell}(A)$ such that $AP\in \Psi^{-\infty, l}_{b,c}$. All these definitions have an analog in the semiclassical and scattering cases (see \cite{vasy2018minicourse} and \cite{zworski2022semiclassical}).

In this paper, we denote the principal symbol map by $\mathfrak{s}$. In the semiclassical setting, we denote it by $\mathfrak{s}_h$.

\subsubsection{Mapping properties}
Let $A\in \Psi^{m_0,l_0}_b$ or $A\in \Psi^{m_0,l_0}_{sc}$. Then $A$ is bounded as an operator between the following spaces:
\begin{align*}
A: H^{m,l}_{b,c} \rightarrow H^{m-m_0, l+l_0}_{b,c}\\
A: H^{m,l}_{b,loc}\rightarrow H^{m,l+l_0}_{b,loc}
\end{align*}
Where $H^{m,l}_{b,c}$ is the subspace of $\overline{H}^{m,l}_b$ whose elements have compact support in $\overline{\mathcal{M}}_\epsilon$.
Note that differential operators in $\text{Diff}^{m_0,l_0}_b$, $\text{Diff}^{m_0,l_0}_{sc}$, $\Psi^{m_0,l_0}_{b,c}$ and $\Psi^{m_0,l_0}_{sc, c}$ are bounded from $\overline{H}^{m,l}_{b}$ to $\overline{H}^{m-m_0, l+l_0}_b$ and from $\dot{H}^{m,l}_b$ to $\dot{H}^{m-m_0, l+l_0}_b$. The same mapping properties holds for semiclassical pseudodifferential operators if the Sobolev spaces are replaced by their semiclassical version.

\subsubsection{Elliptic estimates} 
We state the standard estimate which is a consequence of Proposition 18.1.23 in \cite{hormander2007analysis}:
\begin{prop}
\label{ellipticEstimate}
Let $A,B \in \Psi^{0,l}_{b}$ be pseudodifferential operators with compactly supported Schwarz kernels in $X^2$ (in particular it vanishes near the boundary of $\overline{X}^2$). Let $P \in \Psi^{m,l}_b$ be a differential operator with characteristic set $\Sigma$ (subset of fiber infinity). We assume that $WF(A)\cap \Sigma = \emptyset$ and $WF(A)\subset Ell(B)$. Then we have the following estimates:
For every integers $N,M>0$ and every $s\in \R$, there exists $C>0$ such that:
\[
\left\|Au\right\|_{\overline{H}^{s+m,l}_{(b)}}\leq \left(\left\|B Pu\right\|_{\overline{H}^{s,l}_{(b)}} + \left\|u\right\|_{\overline{H}^{-N,-M}_{(b)}}\right)
\]
\end{prop}
\begin{remark}
Since the support of the Schwartz kernels of $A$ and $B$ are compactly supported in $X\times X$, the index $l$ and the extendible character at $r= r_+-\epsilon$ are irrelevant. For the same reason, the following estimates is also true:
\[
\left\|Au\right\|_{\dot{H}^{s+m,l}_{(b)}}\leq C\left(\left\|B Pu\right\|_{\dot{H}^{s,l}_{(b)}} + \left\|u\right\|_{\dot{H}^{-N,-M}_{(b)}}\right)
\]
\end{remark}

We also need the semiclassical version of the $b$-elliptic estimate:
\begin{prop}
Let $A,B \in \Psi^{0,l}_{b,h}$ with Schwartz kernels supported inside a compact subset of $\overline{X}^2$. Let $P \in \Psi^{m,l}_b$ be a differential operator with characteristic set $\Sigma$. We assume that $WF_h(A)\cap \Sigma = \emptyset$ and $WF_h(A)\subset Ell(B)$. Then we have the following estimate:
For every integers $N,M>0$ and every $s\in \R$, there exists $C>0$ such that:
\[
\left\|Au\right\|_{\overline{H}^{s+m,l}_{(b)}}\leq C\left(\left\|B Pu\right\|_{\overline{H}^{s,l}_{(b)}} + h\left\|u\right\|_{\overline{H}^{-N,l}_{(b)}}\right)
\]
\end{prop}
\begin{remark}
As before, by the support condition on the Schwartz kernel, the character at $r = r_+-\epsilon$ is irrelevant and we can replace $\overline{H}$ by $\dot{H}$.  
\end{remark}
\begin{remark}
The same estimate holds for semiclassical scattering operators. In the scattering case, we could even get an error term $h\left\|u\right\|_{\overline{H}^{-N,-N}_{(b)}}$ but this will not be needed.
\end{remark}

\subsubsection{Propagation of singularities}
We state the result in the case of differential operators. We first need to introduce the notion of Hamiltonian flow for a differential operator.
 
\begin{definition}
Let $p \in C^{\infty}(T^*X, \C)$.
The Hamiltonian vector field $H_p$ associated to $p$ is the unique vector field on ${}^bT^*\overline{X}$ such that, for every local coordinates $y$ on an open subset $U$ of $X$:
\begin{align*}
H_p = \sum_{i=1}^n \partial_{\xi_i}p\partial_{y_i}-\partial_{y_i}p\partial_{\xi_i}.
\end{align*}
where $(y,\xi)$ are the induced local coordinates on $T^*U$.
\end{definition}

Let $P\in \text{Diff}_b^{k,l}$. In particular $P\in \Psi^{k,l}_b$ and we take $x^l p \in C^{\infty}({}^bT^*\overline{X},\C)$ a representative of its principal symbol and we assume that $p$ has real values. Then, we can check that the Hamiltonian vector field $H_{x^lp}$ rescaled by the factor $\mu^{-k+1}$ extends to a vector field on ${}^b\overline{T}^*\overline{X}$ tangent to the boundary (where $\mu$ is a boundary defining function of fiber infinity). In particular, it defines a flow on fiber infinity which we call the Hamiltonian flow. It does not depend on the choice of the representative $p$.

In the case $P\in \text{Diff}_{sc}^{k,l}$, the rescaled Hamiltonian vector field extends to a smooth vector field on ${}^{sc}\overline{T}^*\overline{X}$ tangent to the boundaries and it defines a flow on fiber infinity but also on the face $\left\{x=0\right\}$.

Finally, we define the semiclassical Hamiltonian flow for operator $P_h \in \text{Diff}_{b,h}^{k,l}$ or $P_h\in\text{Diff}_{sc,h}^{k,l}$. We call $p_h$ the semiclassical principal symbol restricted to the face $h=0$ (by a slight abuse we will call this restriction the semiclassical principal symbol in the rest of this paper). Then we have $x^{l}p_h \in C^{\infty}({}^{b/sc}\overline{T}^*X, \C)$. We assume that $p_h$ is real valued. Then, the Hamiltonian vector field $H_{p_h}$ (after rescaling) extends to a vector field on ${}^{b/sc}\overline{T}^*X$ tangent to the boundaries. We call the flow of this vector field the semiclassical flow of the operator $P_h$.

The following proposition is a standard propagation of singularity estimate. For a proof in the scattering case (which can be adapted to treat the $b$ case as well), see \cite[Section~5.4]{vasy2018minicourse}.
\begin{prop}
Let $B_0,B_1,G \in \Psi^{0,0}_{b,c}$. Let $P\in \text{Diff}^{m,k}_b$ with real principal symbol. Assume that for every $x\in WF_b(B_1)$, there exists $t>0$ (resp. $t<0$) such that $e^{-tH_p}x \in \text{Ell}(B_0)$ and $(e^{-sH_p}x)_{s\in [0,t]}$ (resp. $(e^{-sH_p}x)_{s\in [t, 0]}$) remains in the elliptic set of $G$. For every $N>0$, there exists a constant $C>0$ such that we have the following estimate (which holds for every $u$ such that the right hand side is finite):
\begin{align*}
\left\|B_1 u\right\|_{\overline{H}^{\tilde{r},l}_b}\leq C\left(\left\|GPu\right\|_{\overline{H}_b^{\tilde{r}-m+1,l-k}}+\left\|B_0 u\right|_{\overline{H}^{\tilde{r}-N,l}_b}\right).
\end{align*}
\end{prop}
\begin{remark}
By the assumption on the support of the Schwartz kernels of $B_0,B_1$ and $G$, we see that the behavior at the end $\left\{r = r_+-\epsilon\right\}$ is irrelevant and we have the same estimate with $\dot{H}$ instead of $\overline{H}$.
\end{remark}

We also have a semiclassical version of this estimate:
\begin{prop}
Let $B_0,B_1,G \in \Psi^{0,0}_{b,h,c}$. Let $P\in \text{Diff}^{m,k}_{b,h}$ with real principal symbol. Assume that for every $x\in WF_{b,h}(B_1)$, there exists $t>0$ such that $e^{-tH_{p_h}}x \in \text{Ell}(B_0)$ and $(e^{-sH_{p_h}}x)_{s\in [0,t]}$ remains in the elliptic set of $G$. For every $N>0$, there exists a constant $C>0$ such that we have the following estimate:
\begin{align*}
\left\|B_1 u\right\|_{\overline{H}^{\tilde{r},l}_{b,h}}\leq C\left(h^{-1}\left\|GPu\right\|_{\overline{H}_{b,h}^{\tilde{r}-m+1,l-k}}+h\left\|B_0 u\right|_{\overline{H}^{\tilde{r}-N,l}_{b,h}}\right).
\end{align*}
\end{prop}
\begin{remark}
There is also a scattering version for the previous estimates (see \cite[Section~5.4]{vasy2018minicourse}).
Moreover, we stress the fact that there is also a second microlocal version of propagation of singularity estimates and elliptic estimates (including in the semiclassical regime) for an algebra $\Psi_{b,sc}^{\tilde{r},m,l}$ of operators refining the $b$ algebra. We refer to \cite{vasy2020limiting} for a detailed presentation of this algebra and its properties. We will make very little use of this second microlocal algebra here but we build on results of \cite{vasy2020limiting, vasy2020resolvent} in which it plays a central role.
\end{remark}

\section{Cauchy problem}\label{secCauchy}

To state the Cauchy problem, we consider a time coordinates $t_0 = t_*+h(r)$ whose level sets are transverse to the future event horizon (with $\dd t_0$ remaining timelike up to $r_+-2\epsilon$) and which is equal to the usual Boyer-Lindquist coordinate when $r$ is large (see \eqref{deft_0} for a possible concrete definition of $t_0$). 

 We consider the Cauchy problem with initial data on $\Sigma_0 := t_0^{-1}(\left\{ 0\right\})$. Note that the level sets of $t_0$ are naturally identified with $X$.
\subsection{Cauchy Problem for smooth compactly supported initial data}

\begin{prop}\label{CauchyFirst}
Let $\tilde{r}\geq 0$.
Let $u_0 \in H^{\tilde{r}+1}(\Sigma_0, \mathcal{B}_s)$, $u_1\in H^{\tilde{r}}(\Sigma_0, \mathcal{B}_s)$ be compactly supported. The Cauchy problem:
\[\begin{cases}
T_s u = 0\\
u_{|_{\Sigma_0}} = u_0\\
\nabla^\mu t_0 \partial_\mu u_{|_{\Sigma_0}} = u_1
\end{cases}\]
has a unique solution (in the sense of distributions) $u$ in $H^{\tilde{r}+1}_{loc}(\mathcal{M}_\epsilon, \mathcal{B}_s)$. Moreover, for all $k\leq \tilde{r}$ (where $k\in \N$) we have 
\begin{align*}x^{-1}u \in C^{0}([0,+\infty)_{\mathfrak{t}}, \overline{H}^{\tilde{r}+\frac{1}{2}}(X,\mathcal{B}_s))\cap C^{k}([0,+\infty)_{\mathfrak{t}}, \overline{H}^{\tilde{r}-k+\frac{1}{2}}(X,\mathcal{B}_s))
\end{align*} with the exponential bound ($C_{\tilde{r}}>0$ is independent of $T$):
\begin{align*}
\sup_{\mathfrak{t}\in [0,T]}\left\| x^{-1}u(\mathfrak{t}) \right\|_{\overline{H}^{\tilde{r}+\frac{1}{2}}(X,\mathcal{B}_s)} \leq C_{\tilde{r}}e^{C_{\tilde{r}} T}\left( \left\|u_0\right\|_{\overline{H}^{\tilde{r}+1}(\Sigma_0, \mathcal{B}_s)} + \left\|u_1\right\|_{\overline{H}^{\tilde{r}}(\Sigma_0, \mathcal{B}_s)}\right)
\end{align*}
\end{prop}

\begin{proof}
Recall that we defined $\tilde{g} := \rho^{-2}g$ and $\tilde{G}$ the associated metric on the cotangent bundle.
Since the principal symbol of $T_s$ is $\tilde{G}$ (see for example the expression in \eqref{TeukolskyOperateurBasique}), the metric to consider for hyperbolic theory is $\tilde{g}$ (which induces the same causal structure as $g$).
The existence and uniqueness of $u$ in $H^{\tilde{r}+1}_{loc}(\mathcal{M}_\epsilon, \mathcal{B}_s)$ follows from the classical hyperbolic theory on the globally hyperbolic Lorentzian manifold $\mathcal{M}_{\epsilon}$ (see \cite[Chapter~XXIII]{hormander2007analysis} and \cite[Section~12]{ringstrom2009cauchy} for an introduction to this theory). For $A\subset \mathcal{M}_{\epsilon}$, we denote by $J^+(A)$ the causal future of $A$ (for more detail about this notion see \cite[Section~10.2.4]{ringstrom2009cauchy}). By finite speed of propagation, the support of $u$ is contained in $J^+(K)$ where $K$ is a compact subset of $\Sigma_0$ whose interior contains the union of the supports of $u_0$ and $u_1$. As a consequence, we can find a strictly spacelike hypersurface $\Sigma_0'$ which is transverse to $\mathscr{I}^+$, which contains $K$ and such that $J^+(K)\cap \Sigma_0' = K$. The solution $u$ is then solution to the Cauchy problem with initial data on $\tilde{\Sigma_0}$. Note that $x^{-1}T_s x$ is smooth up to $x=0$, and we check that we can extend it analytically on a slightly larger manifold $\tilde{\mathcal{M}}_\epsilon := \R_{\mathfrak{t}}\times[-\eta, \frac{1}{r_+-\epsilon-\eta})_x\times\mathbb{S}^2$ for a small $\eta>0$. For $\eta$ sufficiently small, $\tilde{g}$ extends analytically as a Lorentzian metric on $\tilde{\mathcal{M}}_\epsilon$ and $\Sigma_0'$ extends as a spacelike hypersurface $\tilde{\Sigma}_0'$ such that $\Sigma_0'$ is relatively compact in $\tilde{\Sigma}_0'$. We denote by $\overline{\mathcal{M}}_\epsilon^{[0,T]} := \left(\overline{\mathcal{M}}_\epsilon \cup \left\{x = \frac{1}{r_+-\epsilon}\right\}\right)\cap \mathfrak{t}^{-1}([0,T])$. The Cauchy problem for the operator $x^{-1}Tx$ on the domain of dependence of $\tilde{\Sigma}_0'$ and, by classical hyperbolic theory, the solution $x^{-1}u$ belongs to $\overline{H}^{\tilde{r}+1}(\overline{\mathcal{M}}_\epsilon^{[0,T]})$. Moreover, we can find $u^{n} \in xC^\infty(\overline{\mathcal{M}}_\epsilon^{[0,T]})$ such that $\lim\limits_{n\to +\infty}x^{-1}u^{n} = x^{-1}u$ in $\overline{H}^{\tilde{r}+1}(\overline{\mathcal{M}}_\epsilon^{[0,T]})$ and $\lim\limits_{n\to +\infty}(u^{n}_0,u^{n}_1)= (u_0,u_1)$ in $H^{\tilde{r}+1}\times H^{\tilde{r}}$. Therefore, we can assume that $u$ is smooth when performing energy estimates on $\overline{\mathcal{M}}_\epsilon^{[0,T]}$. The standard energy estimate for hyperbolic partial differential equations gives (for some constant $C>0$ independent of $T$ because coefficients of $T_s$ do not depend on $\mathfrak{t}$):
\begin{align*}
\left\| x^{-1}u\right\|_{\overline{H}^{\tilde{r}+1}\left(\overline{\mathcal{M}}_\epsilon^{[0,T]}\right)}\leq Ce^{CT}\left(\left\|u_0\right\|_{H^{\tilde{r}+1}(\Sigma_0)}+\left\|u_1\right\|_{H^{\tilde{r}}(\Sigma_0)}\right)
\end{align*}
We use the trace theorem (see for example Theorem B.2.7 in \cite{hormander2007analysis}) to get $x^{-1}u \in C^{0}([0,+\infty)_{\mathfrak{t}}, \overline{H}^{\tilde{r}+\frac{1}{2}}(X,\mathcal{B}_s))$ and:
\begin{align*}
\sup_{\mathfrak{t}\in [0,T]}\left\| x^{-1}u(\mathfrak{t}) \right\|_{\overline{H}^{\tilde{r}+\frac{1}{2}}(X,\mathcal{B}_s)} \leq C_{\tilde{r}} e^{C_{\tilde{r}} T}\left( \left\|u_0\right\|_{\overline{H}^{\tilde{r}+1}(\Sigma_0, \mathcal{B}_s)} + \left\|u_1\right\|_{\overline{H}^{\tilde{r}}(\Sigma_0, \mathcal{B}_s)}\right)
\end{align*}
Eventually, still by the trace theorem, we obtain $x^{-1}u \in C^{k}([0,+\infty)_{\mathfrak{t}}, \overline{H}^{\tilde{r}-k+\frac{1}{2}}(X,\mathcal{B}_s))$.
\end{proof}

The following proposition translates the Cauchy problem into a forcing problem and records regularity and decay properties of the forcing term.
\begin{prop}\label{cutoffCauchy}
Let $k\in \N$ and let $\tilde{r}>\frac{1}{2}+k$.
Let $u$ be as in Proposition \ref{CauchyFirst}, then there exists $v \in xC^{1+k}(\R_{\mathfrak{t}}, \overline{H}^{\tilde{r}-\frac{1}{2}-k}(X,\mathcal{B}_s))$ and reals $\mathfrak{t}_0<\mathfrak{t}_1$ such that:
\begin{itemize}
\item $v(\mathfrak{t})=0$  if $\mathfrak{t}<\mathfrak{t}_0$
\item $v(\mathfrak{t})=u(\mathfrak{t})$  if $\mathfrak{t}_1<\mathfrak{t}$
\item $T_s v= f$ with $f \in C^k(\R_{\mathfrak{t}}, \overline{H}_b^{\tilde{r}-\frac{3}{2}-k, \infty})$, $\supp(f) \subset [\mathfrak{t}_0, \mathfrak{t}_1]$.
\end{itemize}
\end{prop}
\begin{proof}
By finite speed of propagation, there exists a compact subset $K\subset (r_+-\epsilon, +\infty)\times \mathbb{S}^2$ and a time interval $[\mathfrak{t}_0,\mathfrak{t}_1]$ such that $\supp(u)\cap \mathfrak{t}^{-1}([\mathfrak{t}_0,\mathfrak{t}_1])\subset [\mathfrak{t}_0,\mathfrak{t}_1]\times K$. We define $\chi \in C^{\infty}(\R,[0,1])$ equal to $0$ on $(-\infty, \mathfrak{t}_0]$ and equal to $1$ on $[\mathfrak{t}_1, +\infty)$ and $v = \chi u$. The fact that $v \in C^{1+k}(\R_{\mathfrak{t}}, \overline{H}^{\tilde{r}-\frac{1}{2}-k}(X,\mathcal{B}_s))$ follows from the same property on $u$. Since $T_s u = 0$, we get $T_s v = [T_s, \chi] u$ which has time support in $[\mathfrak{t}_0, \mathfrak{t}_1]$ and spatial support in $K$.
\end{proof}

\begin{remark}\label{ChangeOfCutoff}
The function $\chi$ introduced in the proof of Proposition \ref{cutoffCauchy} is not the only choice. In particular, it is not necessary to choose a function of $\mathfrak{t}$. The relevant properties for $\chi \in C^{\infty}(\mathcal{M}_\epsilon,[0,1])$ to get $v:= \chi u$ as in Proposition \ref{cutoffCauchy} are:
\begin{itemize}
\item $\chi =1$ on $\left\{\mathfrak{t} \geq \mathfrak{t}_1\right\}$ for some $\mathfrak{t}_1\in \R$.
\item $\chi =0$ on $\left\{\mathfrak{t}\leq \mathfrak{t}_0\right\}\cap \supp(u)$ for some $\mathfrak{t}_0<\mathfrak{t}_1$.
\item $\supp (\dd \chi )\cap \supp (u)$ is compact in $\mathcal{M}_\epsilon$.
\end{itemize}
In particular, it is possible to choose a function of the form $\chi(t_0)$ with $\chi \in C^\infty(\R, [0,1])$ equal to $1$ on $[\eta,+\infty)$ and equal to $0$ on  $(-\infty, 0]$ for $\eta>0$ small enough. This freedom will be useful to highlight the dependence of the forcing term on the initial data (see Remark \ref{constantWrtInitialData}).
\end{remark}

We are now able to take the Fourier-Laplace transform.
\begin{coro}\label{FourierTrans1}
We use the notation of proposition \ref{CauchyFirst} and assume $\tilde{r}>k+\frac{1}{2}$ for some $k\in \N$. We also fix $l\in \R$ (we can think of it as a large decay rate).
For all $\sigma\in \C$ such that $\Im(\sigma)>C_{\tilde{r}}$, we have the following equality between the Fourier-Laplace transforms (with respect to $\mathfrak{t}$): 
\begin{align*}
\hat{T}_s(\sigma)\hat{v}(\sigma) = \hat{f}(\sigma)
\end{align*}
Moreover, we have that $\hat{f}$ is holomorphic on $\C$ with value in $\overline{H}_b^{\tilde{r}-\frac{1}{2}-k,l}$ and there exists $D>0$ such that for all $j\in \N$ there exists $D_{j}>0$ such that:
\begin{align}
\label{boundOnf}
\left\|\partial_{\sigma_x}^{j}\hat{f}(\sigma_x + i\sigma_y)\right\|_{\overline{H}_b^{\tilde{r},l}}\leq D_{j}\left<\sigma_x\right>^{-k} e^{D\lvert \sigma_y\rvert}
\end{align}
\end{coro}

\begin{proof}
The fact that $\hat{u}(\sigma)$ is well defined for $\Im(\sigma)>C_{\tilde{r}}$ follows from the exponential estimate in \ref{CauchyFirst} and the equality follows from the definition of $\hat{T}$. The estimate on $\hat{f}$ follows from the Paley-Wiener-Schwartz theorem (see for example \cite{hormanderI}, Theorem 7.3.1) with the observation that $\partial_{\mathfrak{t}}^{k}(\mathfrak{t}^{j}f)$ is a compactly supported distribution of order zero.
\end{proof}

\subsection{More general initial data} \label{CauchyMoreGeneral}
The goal of this section is to translate the Cauchy problem on $\Sigma_0$ with data $u_0 \in \overline{H}_b^{\tilde{r}+1,1+\alpha}$ and $u_1 \in \overline{H}_{b}^{\tilde{r}, 1+\alpha}$ into a forcing problem and to specify the properties of the Fourier transform of the forcing term.  

Before stating the propositions, we need some geometric preparation. Following \cite{hintz2020stability}, we add another boundary on $\overline{\mathcal{M}}_{\epsilon}$ as follows:
On $\mathcal{U} := \left\{ \mathfrak{t}<0 \right\}$, we can define $\rho_0:= -\mathfrak{t}^{-1}$ and $\rho_I := -x\mathfrak{t}$ and we add the boundary: $I_0 := \left\{\rho_0 = 0\right\}$ to $\mathcal{U}$. Note that the set $\left\{\rho_I=0, \rho_0>0\right\}$ corresponds to $\mathscr{I}^+\cap \mathcal{U}$ and $\left\{\rho_0=0\right\}$ has been glued at the end $\mathfrak{t} = -\infty$, therefore the manifold has a corner at $\rho_0 = \rho_I = 0$. The closure of the hypersurface $\Sigma_0$ on this new manifold intersects $I_0$ transversally. Moreover, we consider $\mathcal{M}_{\epsilon}$ as included in $\mathcal{M}_{2\epsilon}$. It enables us to define extendible Sobolev estimates at $r=r_+-\epsilon$. We call $\textbf{M}_\epsilon := \overline{\mathcal{M}}_\epsilon \cup I_0 \cup \left\{r=r_+-\epsilon\right\}$. Note that in this section, we do not need the boundary $I_+$ and its boundary defining function $\rho_+$ which are introduced in \cite{hintz2020stability}. Indeed, it is enough for our purpose to prove a crude exponential bound with respect to $\mathfrak{t}$ which do not require a precise analysis near $I_+$ (more precise estimate will then follow from the analysis of the resolvent).
\begin{figure}
\begin{center}
\begin{tikzpicture}
\draw (-2.286, 3)--(-5.5, 0.5).. controls +(1, -0.5) and +(-1,0) .. (-2,0) node[yshift = -0.4cm] {$\Sigma_0$}--(1,0)--(2,1) node [midway, xshift = 0.2cm, yshift = -0.2cm] {$I_0$}--(0,3) node [midway, xshift = 0.2cm, yshift = 0.2cm] {$\mathscr{I}^+$};
\draw [dashed](-2,3)--(-5,0) node[yshift = -0.2cm] {$\mathfrak{H}$};
\end{tikzpicture}
\end{center}
\caption{Representation of $\textbf{M}_\epsilon\cap \left\{t_0\geq 0\right\}$}
\end{figure}
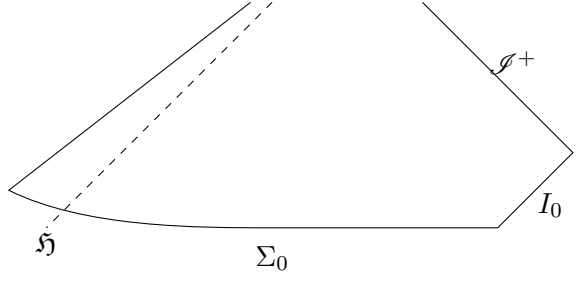
\begin{definition}
For $\tilde{r}\in \N$, we define the $b$ Sobolev space $\mathcal{E}^{\tilde{r}}$ of distributions which are (locally in $\textbf{M}_\epsilon\cap \left\{t_0\geq 0\right\}$): extendible (as a distribution of Sobolev order $\tilde{r}$) at $r=r_+-\epsilon$ and at $t_0 = 0$ and of b regularity $\tilde{r}$ near $I_0 \cup \mathscr{I}^+$. Concretely, $u\in \mathcal{E}^{\tilde{r}}$ means that for any family of smooth\footnote{up to the boundary of $\textbf{M}_{\epsilon}$} vector fields $(L_i)_{i=1}^{N}$ on $\textbf{M}_\epsilon\cap \left\{t_0\geq 0\right\}$ tangent to $\mathscr{I}^+\cup I_0$ with $N\leq \tilde{r}$, $L_1...L_Nu \in L^2_{b,loc}(\textbf{M}_{\epsilon}, \mathcal{B}_s)$. 
\end{definition}
We use the time coordinate $\tt$ on $\textbf{M}_\epsilon \cap \left\{t_0\geq 0\right\}$, which satisfies $c\leq \tilde{G}(\dd \tt, \dd \tt)\leq C$ for positive constants $c$ and $C$ and is smooth up to $\mathscr{I}^+$.

The main proposition of this Section is obtained by adapting energy estimates from \cite[Section~4.1]{hintz2020stability}:
\begin{prop}\label{basicConormalEstimate}
Let $\tilde{r}\in \N$. Let $a_I<0$ and $a_0>a_I$. Let $u_0 \in \overline{H}_b^{\tilde{r}+1, a_0}$ and $u_1 \in \overline{H}_b^{\tilde{r}, a_0}$.
The unique solution $u$ to the Cauchy problem:
\begin{align*}
\begin{cases}
T_s u = 0\\
u_{|_{\Sigma_0}} = u_0\\
\rho_0\nabla^\mu t_0\partial_{\mu}u_{|_{\Sigma_0}} = u_1
\end{cases}
\end{align*}
belongs to $\rho_0^{a_0}\rho_I^{a_I+1}\mathcal{E}^{\tilde{r}}$. Moreover if $k\leq \tilde{r}$ for some $k\in \N$, it also belongs to $C^k([0,+\infty)_{\mathfrak{t}}, \overline{H}^{\tilde{r}-k,a_I+1}_b)\cap C^0([0,+\infty)_{\mathfrak{t}}, \overline{H}^{\tilde{r},a_I+1}_b)$ and there exists $C_{\tilde{r}}>0$ such that for all $T>0$:
\begin{align*}
\sup_{\mathfrak{t}\in [0, T]}\left\|u(\mathfrak{t})\right\|_{\overline{H}^{\tilde{r},a_I+1}_b}\leq C_{\tilde{r}}e^{C_{\tilde{r}}T}\left(\left\|u_0\right\|_{\overline{H}^{\tilde{r}+1, a_0}_b}+\left\|u_1\right\|_{\overline{H}^{\tilde{r},a_0}}\right)
\end{align*}
\end{prop}

\begin{remark}
In the proof, we obtain a more precise result (namely $u\in\rho_0^{a_0}\rho_I^{a_I}\mathcal{E}_{\mathscr{I}}^{\tilde{r}+1}$, see the proof for the definition of $\mathcal{E}_{\mathscr{I}}^{\tilde{r}}$). However, since we do not aim to optimize the regularity assumptions here, we allow this loss in order to stay with simpler spaces $\mathcal{E}^{\tilde{r}}$.
\end{remark}

We then use an idea mentioned in \cite[Section~5.1]{hintz2020stability} to get:
\begin{prop}\label{refinementAtScri}
With the notations of Proposition \ref{basicConormalEstimate}, if we assume in addition that $a_0=1+\alpha$ with $\alpha \in (0,1)$ and $\tilde{r}> k+2$ with $k\in \N$, then there exists $\mathfrak{t}_0$,$\mathfrak{t}_1$ and $v\in C^{k}\left(\R_{\mathfrak{t}}, H_{b}^{\tilde{r}-2-k, 1-}\right)$ such that:
\begin{itemize}
\item $v(\mathfrak{t})=0$ if $\mathfrak{t}<\mathfrak{t}_0$
\item $v(\mathfrak{t})=u(\mathfrak{t})$ if $\mathfrak{t}>\mathfrak{t}_1$
\item For all $\mathfrak{t}\in \R$, $v(\mathfrak{t}) = xv_0(\mathfrak{t}) + \overline{H}_{b}^{\tilde{r}-2-k, 1+\alpha-}$ where $v_0 \in C^k(\R_{\mathfrak{t}}, H^{\tilde{r}-2-k}(\mathcal{B}_s))$
\item $T_s v = f$ with $f \in C^{k-1}\left(\R_{\mathfrak{t}}, \overline{H}_{(b)}^{\tilde{r}-3-k, -\frac{3}{2}+\alpha-}\right)$, $\mathfrak{t}(\supp(f))\subset [\mathfrak{t}_0,\mathfrak{t}_1]$
\end{itemize}
\end{prop}

Proposition \ref{basicConormalEstimate} and Proposition \ref{refinementAtScri} allow us to reduce the study of the asymptotic behavior of the solution to the Cauchy problem to the study of the solution $v$ of the forcing problem $T_s v= f$ provided by Proposition \ref{refinementAtScri}. Moreover, the properties of $u$, $v$ and $f$ stated in the propositions allow us to consider the Fourier transform with respect to $\mathfrak{t}$:
\begin{coro}\label{FourierTrans2}
We use the notation of Propositions \ref{basicConormalEstimate} and \ref{refinementAtScri} and assume $\tilde{r}>k+3$ for some $k\in \N$.
For all $\sigma\in \C$ such that $\Im(\sigma)>C_{\tilde{r}}$, we have the following equality between the Fourier-Laplace transforms: 
\begin{align*}
\hat{T}_s(\sigma)\hat{v}(\sigma) = \hat{f}(\sigma)
\end{align*}
Moreover, we have that $\hat{f}$ is holomorphic on $\C$ with values in $\overline{H}_{(b)}^{\tilde{r}-3-k, -\frac{3}{2}+\alpha-}$ and there exists $D>0$ such that for all $j\in \N$ there exists $D_{j}>0$ such that:
\begin{align}
\label{boundOnf2}
\left\|\partial_{\sigma_x}^{j}\hat{f}(\sigma_x + i\sigma_y)\right\|_{\overline{H}_{(b)}^{\tilde{r},-\frac{3}{2}+\alpha-}}\leq D_{j}\left<\sigma_x\right>^{-(k-1)} e^{D\lvert \sigma_y\rvert}
\end{align}
\end{coro}

Since the detailed proofs of Proposition \ref{basicConormalEstimate} and Proposition \ref{refinementAtScri} are quite long and not essential to understand the rest of the paper, they have been added in Appendix \ref{ProofCauchyConormal}.

\section{Analysis of the classical and semiclassical Hamiltonian flow of \texorpdfstring{$\hat{T}_s(\sigma)$}{ hat\{T\}\_ s}} \label{secFlow}
Before stating the Fredholm estimate, we analyse the structure of the classical and semiclassical Hamiltonian flow of $\hat{T}_s(\sigma)$. In particular, we identify the elliptic regions, the transport regions (i.e. the subset of the characteristic set where the Hamilton vector field does not vanish) and the radial points. In the semiclassical regime, we will see that some trajectories under the Hamiltonian flow remains in a compact region of $X$, a phenomenon known as trapping.
\subsection{Analysis of the classical flow}
\begin{prop}
Let $\eta>0$. The operator $\hat{T}_s(\sigma)$ is (classically) elliptic on every region of the form $\left\{r>2M+\eta\right\}\cap \overset{\circ}{K}$ where $K$ is a compact subset of $X$. 
\end{prop}
\begin{proof}
If we denote by $\hat{T}_s(\sigma)'$ the operator obtained by Fourier transform with respect to $t_*$ instead of $\mathfrak{t}$, we check that $\hat{T}_s(\sigma) = e^{-i\sigma L(r)}\hat{T}_s'(\sigma)e^{i\sigma L(r)}$ where $L$ is a smooth function on $X$. In particular, the classical principal symbols of $\hat{T}_s(\sigma)$ and $\hat{T}_s'(\sigma)$ are the same. In coordinates $(r,\theta,\phi_*)$, writing cotangent vectors as $\xi \dd r + \zeta \dd \phi_* + \eta \dd \theta$, we have:
\begin{align}
\label{symbPrinc}
\mathfrak{s}(\hat{T}_s(\sigma))=& \Delta_r \xi^2 + 2a\xi\zeta + \frac{\zeta^2}{\sin^2\theta} + \eta^2\\
=& \left(a^2\cos^2\theta + r(r-2M)\right)\xi^2 + \left(a\sin\theta \xi + \frac{\zeta}{\sin(\theta)}\right)^2  + \eta^2 \notag
\end{align}
where $\mathfrak{s}$ is the classical principal symbol map. Strictly speaking, this computation is only valid outside the rotation axis, however with a change of coordinates we can show that the operator is elliptic also in a neighborhood of the rotation axis (see \eqref{symbolAxis} below).
\end{proof}
We now investigate the characteristic set and the Hamiltonian flow on it (in the region $\left\{r\leq 2M\right\}$). The Hamiltonian vector field is
\begin{align}
\label{H}
H = \left(2\Delta_r \xi + 2a\zeta\right)\partial_r + \left(2a\xi + \frac{2\zeta}{\sin^2\theta}\right)\partial_\phi + 2\eta\partial_\theta - 2(r-M)\xi^2\partial_\xi + \frac{2\zeta^2\cos\theta}{\sin^3\theta}\partial_{\eta}
\end{align}
We now work in the radial compactification $\overline{T}^*X$ of $T^*X$ and we set $\tilde{\rho} := \frac{1}{\sqrt{\xi^2+\frac{\zeta^2}{\sin^2\theta}+\eta^2}}$ (note that $\tilde{\rho}$ extends smoothly to the rotation axis). We set $\tilde{\xi}:=\tilde{\rho} \xi$, $\tilde{\zeta}:=\tilde{\rho} \frac{\zeta}{\sin\theta}$, $\tilde{\eta}:= \tilde{\rho} \eta$. We use functions $(\tilde{\rho}, (\tilde{\xi}, \tilde{\zeta}, \tilde{\eta}))\in [0,+\infty)\times\mathbb{S}^2$ to parametrize the fibers of $\overline{T}^*\mathcal{M}_\epsilon\setminus\left\{0\right\}$ away from the rotation axis. The rescaled Hamiltonian vector field (in terms of $(r,\theta,\phi,\tilde{\rho},(\tilde{\xi},\tilde{\eta},\tilde{\zeta})$) is
\begin{align}
\tilde{H} = \tilde{\rho} H =& \left(2\Delta_r \tilde{\xi} + 2a\tilde{\zeta}\sin\theta\right)\partial_r +\left(2a\tilde{\xi} + \frac{2\tilde{\zeta}}{\sin\theta}\right)\partial_\phi + 2\tilde{\eta}\partial_\theta + 2(r-M)\tilde{\xi}^3\tilde{\rho}\partial_{\tilde{\rho}} \notag \\& - 2(r-M)\tilde{\xi}^2Z_1+ 2\cotan\theta \tilde{\zeta}^2 Z_2 \label{tildeH} - 2\tilde{\eta}\tilde{\zeta}\cotan\theta Z_3
\end{align}
where $Z_1 := (1-\tilde{\xi}^2)\partial_{\tilde{\xi}}-\tilde{\eta}\tilde{\xi}\partial_{\tilde{\eta}}-\tilde{\zeta}\tilde{\xi}\partial_{\tilde{\zeta}}$, $Z_2 := (1-\tilde{\eta}^2)\partial_{\tilde{\eta}}-\tilde{\eta}\tilde{\zeta}\partial_{\tilde{\zeta}}-\tilde{\eta}\tilde{\xi}\partial_{\tilde{\xi}}$ and $Z_3 = (1-\tilde{\zeta}^2)\partial_{\tilde{\zeta}}-\tilde{\zeta}\tilde{\xi}\partial_{\tilde{\xi}}-\tilde{\zeta}\tilde{\eta}\partial_{\tilde{\eta}}$ are smooth vector fields on $\mathbb{S}^2$. Note that $\tilde{H}$ is smooth also on the rotation axis (see the computation in stereographic coordinates \eqref{HAxis} below).
\begin{remark}
Pay attention to the fact that $\partial_\theta$ in \eqref{tildeH} is different from $\partial_\theta$ in \eqref{H} since $\tilde{\rho}$ and $\tilde{\zeta}$ depend on $\theta$.
\end{remark}
Following \cite{vasy2013microlocal} (section 6.3), we define the sets $\Lambda_+ := \left\{\Delta_r = 0, \eta = 0, \zeta = 0, \xi>0\right\}$, $L_+ = \Lambda_+\cap \left\{\tilde{\rho} = 0\right\}$ and $\Lambda_-:=\left\{\Delta_r = 0, \eta = 0, \zeta = 0, \xi<0\right\}$, $L_- = \Lambda_-\cap \left\{\tilde{\rho} = 0 \right\}$.
The vector field $\tilde{H}$ extends smoothly to fiber infinity ($\tilde{\rho} = 0$).

\begin{definition}
Points of the characteristic set where the rescaled Hamitonian vector field vanishes are called radial points.
\end{definition}
\begin{prop}
Radial points are exactly $L_+\cup L_-$. 
\end{prop}

\begin{proof}
To cover the rotation axis, we use stereographic coordinates (of north pole) on the sphere so that:
$(x_N, y_N) = \left(cotan\frac{\theta}{2}\cos\phi_*, cotan\frac{\theta}{2}\sin\phi_*\right)$. We define similarly stereographic coordinates of south pole $(x_S, y_S)$. We know that the principal symbol is symmetric with respect to the reflection by the equatorial plan if we replace $a$ by $-a$. As a consequence, it is enough to do the analysis in the stereographic coordinates of north pole.
We write linear forms as $\xi \dd r + \gamma\dd x_N + \mu \dd y_N$.
We obtain:
\begin{align}\label{symbolAxis}
\mathfrak{s}(\hat{T}_s(\sigma)) = \Delta_r \xi^2 + 2a\xi(-y_N\gamma + x_N\mu) + \frac{1+x_N^2+y_N^2}{4}(\gamma^2+\mu^2)
\end{align}
We can compute the fiber infinity defining function $\tilde{\rho}$ (and see that it extends smoothly across the rotation axis):
\[
\tilde{\rho} = \frac{1}{\sqrt{\xi^2+\frac{(1+x_N^2+y_N^2)^2}{4}(\gamma^2+\mu^2)}}
\]
This expression prompts us to use the parametrization $\tilde{\rho}, (\tilde{\xi}, \tilde{\gamma}, \tilde{\mu})\in \mathbb{S}^2$ where $\tilde{\xi} = \tilde{\rho} \xi$, $\tilde{\gamma} = \tilde{\rho} \frac{(1+x_N^2+y_N^2)\gamma}{2}$ and $\tilde{\mu} = \tilde{\rho} \frac{(1+x_N^2+y_N^2)\mu}{2}$. Note that we have the relations 
\begin{align*}\tilde{\zeta} &= -\frac{y_N}{\sqrt{x_N^2+y_N^2}}\tilde{\gamma}+\frac{x_N}{\sqrt{x_N^2+y_N^2}}\tilde{\mu}\\
\tilde{\eta}&= -\frac{x_N}{\sqrt{x_N^2+y_N^2}}\tilde{\gamma} - \frac{y_N}{\sqrt{x_N^2+y_N^2}}\tilde{\mu}
\end{align*}
We can then rewrite the rescaled Hamiltonian vector field (and see that it extends smoothly across the rotation axis):
\begin{align}
\tilde{H} = &\left(2\Delta_r \tilde{\xi}+4a\frac{x_N\tilde{\mu}-y_N\tilde{\gamma}}{1+x_N^2+y_N^2}\right)\partial_r - 2a\tilde{\xi}y_N\partial_{x_N}+2a\tilde{\xi}x_N\partial_{y_N}+ (1+x_N^2+y_N^2)\left(\tilde{\gamma}\partial_{x_N} +\tilde{\mu}\partial_{y_N}\right)\notag\\ & +2(r-M)\tilde{\xi}^3\tilde{\rho}\partial_{\tilde{\rho}}-2(r-M)\tilde{\xi}^2Z_1 + 2x_N Z_3+2y_NZ_4 + 2a\tilde{\xi}Z_5 \label{HAxis}
\end{align}
where $Z_1$ is the same as previously and can be expressed in terms of $\tilde{\xi}, \tilde{\gamma}, \tilde{\mu}$:
$Z_1 = (1-\tilde{\xi}^2)\partial_{\tilde{\xi}}-\tilde{\xi}\tilde{\gamma}\partial_{\tilde{\gamma}} - \tilde{\xi}\tilde{\mu}\partial_{\tilde{\mu}}$. The vector fields $Z_4$ and $Z_5$ are smooth vector fields on $\mathbb{S}^2$ and are given by:
\begin{align*}
Z_3 &= -\tilde{\mu}^2\partial_{\tilde{\gamma}} + \tilde{\gamma}\tilde{\mu}\partial_{\tilde{\mu}}\\
Z_4 &= \tilde{\mu}\tilde{\gamma}\partial_{\tilde{\gamma}} - \tilde{\gamma}^2\partial_{\tilde{\mu}}\\
Z_5 &= -\tilde{\mu}\partial_{\tilde{\gamma}} + \tilde{\gamma}\partial_{\tilde{\mu}}
\end{align*}
With the expression of $\tilde{H}$, we see that the radial points (on the domain of the coordinates $(x_N, y_N, r)$) are exactly solutions of the following system:
\begin{numcases}{}
2\Delta_r \tilde{\xi} + 4a\frac{x_N\tilde{\mu}-y_N\tilde{\gamma}}{1+x_N^2+y_N^2} = 0 \label{line1}\\
-2a\tilde{\xi}y_N + (1+x_N^2+y_N^2)\tilde{\gamma} = 0 \label{line2}\\
2a\tilde{\xi}x_N +(1+x_N^2+y_N^2)\tilde{\mu} = 0 \label{line3}\\
2\tilde{\xi}^2(M-r)(1-\tilde{\xi}^2) = 0 \label{line4}\\
2(r-M)\tilde{\xi}^3\tilde{\gamma} + 2\tilde{\mu}\tilde{\gamma}y_N - 2 a\tilde{\mu}\tilde{\xi}-2\tilde{\mu}^2x_N = 0\\
2(r-M)\tilde{\mu}\tilde{\xi}^3 +2a\tilde{\gamma}\tilde{\xi} + 2 \tilde{\gamma}\tilde{\mu} x_N - 2 \tilde{\gamma}^2y_N = 0
\end{numcases}
Assume that we have a solution $(r,x_N,y_N, \tilde{\xi}, \tilde{\gamma},\tilde{\mu})$ of this system.
Using \eqref{line4}, (and the fact that $M-r>0$ on $\mathcal{M}_{\epsilon}$), we deduce $\tilde{\xi} \in \left\{-1, 0,1\right\}$. We exclude the case $\tilde{\xi} = 0$ since \eqref{line2} and \eqref{line3} give then $\tilde{\gamma} = \tilde{\mu}=0$ which is impossible since $\tilde{\xi}^2+\tilde{\gamma}^2+\tilde{\mu}^2 = 0$. We deduce that $\tilde{\xi}\in \left\{-1, 1\right\}$ and $\tilde{\gamma} = \tilde{\mu} = 0$. Using \eqref{line1}, we deduce that $\Delta_r = 0$ and therefore $r= r_+$. Finally, equations $\eqref{line3}$ and $\eqref{line4}$ give $x_N = y_N = 0$. We deduce that the set of radial points is exactly
 $\left\{x_N=y_N= \Delta_r = \tilde{\mu}=\tilde{\gamma} = 0, \tilde{\xi} = \pm 1\right\}\cup \left\{x_S=y_S=\Delta_r = \tilde{\mu} = \tilde{\gamma} = 0, \tilde{\xi} = \pm 1\right\}$ (the second set is obtained by symmetry).
\end{proof}

\begin{definition}
Let $A$ be a set of radial points which is a submanifold of fiber infinity or, in the scattering setting, of the boundary face $\left\{x=0\right\}$. Here we assume that $A$ does not intersect the other boundary face, but we could extend the definition to the case where $A$ is transversal to it.
We say that $A$ is a sink (resp. a source) for the Hamiltonian flow if there exists $\rho_0$ a non negative quadratic defining function of $A$ within the characteristic set\footnote{The restriction of $\rho_0$ to the characteristic set vanishes quadratically at $A$ and is non degenerate.} such that:
\begin{itemize}
\item $\tilde{H}\rho_0 = -\beta_1 \rho_0-F_2+F_3$ (resp. $\tilde{H}\rho_0 = \beta_1 \rho_0+F_2+F_3$) with $\beta_1, F_2, F_3$ are functions defined on a neighborhood of $A$ and $\beta_1$ is positive on $A$, $F_3$ vanishes cubically at $A$ and $F_2\geq 0$.
\item There exists $\mu$, a defining function of the boundary face containing $A$ (fiber infinity or $\left\{x=0\right\}$), such that $\tilde{H}\mu = -\beta_0 \mu$ (resp. $\tilde{H}\mu = \beta_0 \mu$) where $\beta_0$ is a function defined in a neighborhood of $A$ with $\beta_0>0$ on $A$.
\end{itemize}
\end{definition}

We check that $H \tilde{\rho} = 2(r-M)\tilde{\xi}^3$. So on $\Lambda_{\pm}$, $H\tilde{\rho} = \pm 2(r-M)$ (with the notation of \cite{vasy2013microlocal}, it means that $\beta_0 = 2(r_+-M)$). Moreover, the non negative homogeneous of degree zero function $\rho_0 := \tilde{\eta}^2+\tilde{\zeta}^2$ is a quadratic defining function of $\Lambda_{\pm}$ inside the characteristic set of $\hat{T}_s(\sigma)$ and we have $\tilde{\rho} H \rho_0 = 4(r-M)\tilde{\xi}^3\left(\tilde{\eta}^2 +\tilde{\zeta}^2\right) = 4(r-M)\tilde{\xi}^3\rho_0$. This shows that $L_+$ is a source for $\hat{T}_s(\sigma)$ and $L_-$ is a sink in the sense of \cite{vasy2013microlocal}. Note that there is a different sign convention in \cite{vasy2013microlocal} (the principal symbols given page 483 has a minus sign with respect to our choice).

Now we study more precisely the behavior of the bicharacteristics (integral curves of $\tilde{H}_{|\left\{\tilde{\rho} = 0\right\}}$ included in the characteristic set of $\mathfrak{s}(\hat{T}_s(\sigma))$). Let $\lambda \mapsto f(\lambda)=(r(\lambda),\theta(\lambda), \phi(\lambda), \tilde{\xi}(\lambda),\tilde{\zeta}(\lambda),\tilde{\eta}(\lambda))$ be a bicharacteristic defined on the maximally extended open interval $I$. By definition, we have for all $\lambda\in I$
\begin{align}
&\frac{\dd}{\dd \lambda}\tilde{\xi} = 2\tilde{\xi}^2(M-r)(1-\tilde{\xi}^2)\leq 0 \label{eqXi}\\
&\Delta_r\tilde{\xi}^2 + 2a\tilde{\xi}\tilde{\zeta}\sin\theta + \tilde{\zeta}^2+\tilde{\eta}^2 = 0
\end{align}
Note that if $\Delta_r > a^2$, we have
\begin{align*}
\Delta_r\tilde{\xi}^2 + 2a\tilde{\xi}\tilde{\zeta}\sin\theta + \tilde{\zeta}^2+\tilde{\eta}^2&>a^2\tilde{\xi}^2 + 2a\tilde{\xi}\tilde{\zeta}\sin\theta + \tilde{\zeta}^2 \\
&\geq a^2\tilde{\xi}^2-2\lvert a \rvert \lvert \tilde{\xi}\rvert \lvert \tilde{\zeta}\rvert + \tilde{\zeta}^2\\
&\geq (\lvert a\rvert \lvert\tilde{\xi}\rvert - \lvert\tilde{\zeta}\rvert)^2\\
&\geq 0
\end{align*}
This computation ensures that the bicharacteristic remains in $\left\{r\leq 2M\right\}$.

\begin{lemma}
\label{directionBeyondHorizon}
We have $\tilde{\xi}\tilde{H}r<0$ on the characteristic set when $\Delta_r<0$.
\end{lemma}
\begin{proof}
Since, $\tilde{H}r = 2\Delta_r \tilde{\xi} + 2a\tilde{\zeta}\sin\theta$, on the characteristic set we have:
\begin{align*}
0 =& \Delta_r\tilde{\xi}^2 + 2a\tilde{\xi}\tilde{\zeta}\sin\theta + \tilde{\zeta}^2+\tilde{\eta}^2\\
=& \tilde{\xi}\tilde{H}r-\Delta_r\tilde{\xi}^2 + \tilde{\zeta}^2+\tilde{\eta}^2\\
\tilde{\xi}\tilde{H}r =& \Delta_r\tilde{\xi}^2 - (\tilde{\zeta}^2+\tilde{\eta}^2) <0
\end{align*}
\end{proof}

We deduce from \eqref{eqXi} that $\tilde{\xi}$ is decreasing. Moreover, $\tilde{\xi}$ is bounded (since $\lvert\tilde{\xi}\rvert\leq 1$). Then 
\begin{itemize}

\item Case $1>\tilde{\xi}(0)> 0$: In this case, $\sup I = \lambda_+< +\infty$. Otherwise,   $\lim\limits_{\lambda \to +\infty}\tilde{\xi}= 0$ ($\xi$ has a limit because it is decreasing at least linearly when $\xi$ is in a fixed compact subset of $(0,1)$ ) but this contradicts the fact that for all $\lambda \in I$, $\Delta_r \tilde{\xi}^2 + 2a\tilde{\xi}\tilde{\zeta}\sin\theta + \tilde{\zeta}^2 + \tilde{\eta}^2 = 0$. As a consequence, $f$ leaves every compact set as $\lambda \to \lambda_{+}$. Since $\left\{r\geq r_+-\frac{\epsilon}{2}, \tilde{\rho} = 0\right\}\cap \left\{\tilde{\rho}^2 p = 0\right\}$ is compact, this means that there exists $\lambda_0 \in I$ such that for all $\lambda >\lambda_0$, $r(\lambda)<r_+-\frac{\epsilon}{2}$. By lemma \ref{directionBeyondHorizon}, $r$ cannot reach $r_+-\epsilon$ in the past therefore the lower bound is $-\infty$ and then $\lim\limits_{\lambda \to -\infty} \tilde{\xi} = 1$ and $\lim\limits_{\lambda \to -\infty} \tilde{\zeta} = \lim\limits_{\lambda \to -\infty} \tilde{\eta} = 0$ and using the fact that $\Delta_r \tilde{\xi}^2 + 2a\tilde{\xi}\tilde{\zeta}\sin\theta + \tilde{\zeta}^2 + \tilde{\eta}^2 = 0$, we deduce that $\lim\limits_{\lambda \to -\infty} r(\lambda) = r_+$. As a consequence, $f$ tends to $L_+$ when $\lambda\rightarrow -\infty$. 

\item Case $-1<\tilde{\xi}(0)<0$: We show as in the first case that $I$ is lower bounded by some $\lambda_-\in\R$ and there exists $\lambda_0 \in I$ such that for all $\lambda<\lambda_0$, $r(\lambda)<r_+-\frac{\epsilon}{2}$. Moreover $\sup(I)=+\infty$ and $f$ tends to $L_-$ when $\lambda \rightarrow +\infty$.

\item Case $|\tilde{\xi}(0)|=1$: In this case $\tilde{\eta} = \tilde{\zeta} = 0$. Therefore, since we are on the characteristic set, $\Delta_r = 0$. We conclude that the characteristic is included in a radial set.

\item Case $\tilde{\xi}(0) = 0$: This case is impossible since it would mean $\tilde{\zeta}^2+\tilde{\eta}^2 = 1$ which contradicts $\Delta_r \tilde{\xi}^2 + 2a\tilde{\xi}\tilde{\zeta}\sin\theta + \tilde{\zeta}^2 + \tilde{\eta}^2 = 0$.
\end{itemize}

Finally, possible behavior of bicharacteristics are summarized by the following proposition:
\begin{prop}\label{recapClassicalFlow}
Let $\lambda \mapsto f(\lambda)=(r(\lambda),\theta(\lambda), \phi(\lambda), \tilde{\xi}(\lambda),\tilde{\zeta}(\lambda),\tilde{\eta}(\lambda))$ be a bicharacteristic defined on the maximally extended open interval $I$ containing $0$. Then we are in one of the following cases:
\begin{enumerate}
\item $I = (-\infty, \lambda_+)$ for some $\lambda_+ \in (0, +\infty)$ and there exists a small non empty interval $J:=(\lambda_+-\eta, \lambda_+)$ such that $(r\circ f)_{|_J}$ is decreasing with values in $(r_+-\epsilon, r_+-\frac{\epsilon}{2})$ and $f$ tends to $L_+$ when $\lambda \to -\infty$. In this case $\tilde{\xi}$ is negative along the trajectory
\item $I = (\lambda_-, +\infty)$ for some $\lambda_- \in (-\infty, 0)$ and there exists a small non empty interval $J:=(\lambda_-, \lambda_-+\eta)$ such that $(r\circ f)_{|_J}$ is increasing with values in $(r_+-\epsilon, r_+-\frac{\epsilon}{2})$ and $f$ tends to $L_-$ when $\lambda \to +\infty$.
\item The curve is included in $L_+$ (resp. in $L_-$). In this case $\tilde{\xi}$ is positive (resp. negative) along the trajectory.
\end{enumerate}
\end{prop}

\begin{prop}
There exists an open neighborhood $U_+$ of $L_+$ such that if $f$ is a maximal bicharacteristic defined on the open interval $I$ containing $0$ and if $f(0)\in U_+\cap \Sigma$ ($\Sigma$ is the characteristic set of $p$), then $(-\infty, 0)\subset I$ and $f$ tends to $L_+$ as $\lambda \to -\infty$.
\end{prop}
\begin{proof}
We recall that by definition of $\epsilon$, $\Delta_{r_+-\frac{\epsilon}{2}} = a^2+\left(r_+-\frac{\epsilon}{2}\right)^2-2M\left(r_+-\frac{\epsilon}{2}\right)<0$.
We define $0<\alpha<1$ such that $\Delta_{r_+-\frac{\epsilon}{2}}(1-\alpha^2)+2a\alpha<0$. As a consequence, for all $1\geq \tilde{\xi}\geq\sqrt{1-\alpha^2}$ and all $\lvert\tilde{\zeta}\rvert\leq \alpha$ we have
\begin{equation} \Delta_{r_+-\frac{\epsilon}{2}}\tilde{\xi}^2 + 2a\tilde{\xi}\tilde{\zeta}\sin\theta <0 \label{rdot}\end{equation}
We define $U_+:=\left\{ \rho_0 < \alpha^2, \tilde{\xi} > \sqrt{1-\alpha^2}, r>r_+-\frac{\epsilon}{2}\right\}$. Let $f$ and $I$ be defined as in the statement of the proposition. Using \eqref{eqXi}, we have that $\tilde{\xi}$ is decreasing on $I$. Moreover, $\frac{\dd}{\dd \lambda}\rho_0\circ f = (\tilde H \rho_0)\circ f = 4(r-M)\tilde{\xi}^3\rho_0\circ f$ so $\rho_0$ is exponentially increasing along $f$. This shows that for every $\lambda \in I\cap (-\infty, 0]$, we have $\tilde{\xi}\geq \sqrt{1-\alpha^2}$ and $\lvert \tilde{\zeta} \rvert \leq \sqrt{\rho_0}\leq \alpha$. Then \eqref{rdot} ensures that for all $\lambda \in I\cap (-\infty, 0]$, $r(\lambda)>r_+-\frac{\epsilon}{2}$. As a consequence, $f$ remains in the compact set $\left\{r\geq r_+-\frac{\epsilon}{2}\right\}\cap \Sigma$ on $I\cap (-\infty, 0]$ and therefore, $(-\infty, 0]\subset I$. Then we deduce $\lim\limits_{\lambda \to -\infty}\rho_0\circ f (\lambda) = 0$ and the claimed convergence (since $\rho_0$ is a defining function of $L_+$ inside $\Sigma$)
\end{proof}
\begin{remark}
If we take the image of $U_+$ by the backward bicharacteristic flow, we get an open neighborhood of $L_+$ which is stable under the flow with the same property as $U_+$.
\end{remark}

We can prove similarly:
\begin{prop}
There exists an open neighborhood $U_-$ of $L_-$ such that if $f$ is a maximal bicharacteristic defined on the open interval $I$ containing $0$ and if $f(0)\in U_-\cap \Sigma$ ($\Sigma$ is the characteristic set of $p$), then $( 0,+\infty)\subset I$ and $f$ tends to $L_-$ as $\lambda \to +\infty$.
\end{prop}

\subsection{Analysis of the semiclassical flow} \label{semiClFlowSection}
We will use the semiclassical regime to study the behavior of the operator $\hat{T}_s(\sigma)$ when $\left|\sigma\right|\to +\infty$ and $0\leq \Im(\sigma) \leq C$ for some fixed constant $C>0$. Therefore, we introduce the semiclassical parameter $h =\left|\sigma\right|^{-1}$ and the rescaled operator $T_{s,h}(z):=h^2\hat{T}_s(h^{-1}z)$ for $z = \pm 1 + O(h)$. We denote by $z_0$ the element of $\left\{-1, 1\right\}$ such that $z-z_0 = O(h)$.
The semiclassical principal symbol is $p_{h}(\xi)=-\tilde{G}(\xi-z_0\dd \mathfrak{t})$  (it does not depend on the imaginary part of $z$). From this formula, we see that the semiclassical flow of $h^2\hat{T}_s(h^{-1}z)$ is closely related with the geodesic flow of $\mathcal{M}_{\epsilon}$.
Since the coordinate $\mathfrak{t}$ is not very convenient for computations, we study the effect of a change of time function (in the Fourier transform) of the form $t_2 = \mathfrak{t} + f(r)$ for a smooth function $f$ defined on an open interval $I\subset (r_+-\epsilon, +\infty)$. We define $U_I = I\times \mathbb{S}^2$ and $p_{2,h}(\xi) = -\tilde{G}(\xi-z_0\dd t_2)$ on $\mathcal{T}^*U_I$. In particular we have $p_{2,h}(\xi) = p_h(\xi-z_0 f'(r)\dd r)$.
\begin{lemma}
The map $\Psi: \xi \mapsto \xi-z_0 f'(r)\dd r$ is a symplectomorphism of $T^*U_I$
\end{lemma}
\begin{proof}
The map $\Psi$ is smooth with smooth inverse $\omega \mapsto \omega+z_0f'(r)\dd r$. Let $y \in U_I$ and $x:=(x_0,x_1, x_2)$ be smooth local coordinates around $y$ with $x_0 = r$. The symplectic form in local coordinates is therefore: $\omega =  \sum_{i=0}^2 \dd x_i \wedge \dd \xi_i$ (where $\xi_x:=(\xi_0, \xi_1, \xi_2)$ are the conjugated local coordinates) and $\Psi(x,\xi_r,\xi_1,\xi_2) = (x, \xi_0-z_0f'(r), \xi_1, \xi_2)$. Let $\xi \in T_y^*U_I$, $X = \sum_{i=0}^3 X_{x_i} \partial_{x_i} + \sum_{i=0}^3 X_{\xi_i}\partial_{\xi_i}$ and $Y = \sum_{i=0}^3 Y_{x_i} \partial_{x_i} + \sum_{i=0}^3 Y_{\xi_i}\partial_{\xi_i}$. We have:

\begin{align*}
\dd_{\xi} \Psi(X) &= X-X_0 z_0f'' \partial_{\xi_r} \in T_{\Psi(\xi)}T^*U_I\\
\omega(\dd_{\xi}\Psi X,\dd_{\xi}\Psi Y)&= \omega(X,Y)-\omega(X_0 z_0f'' \partial_{\xi_r}, Y)-\omega(X, Y_0 z_0f'' \partial_{\xi_r})\\
&= \omega(X,Y)+ Y_0X_0z_0f'' - Y_0X_0z_0f''\\
&= \omega(X,Y)
\end{align*}
\end{proof}
\begin{remark}
Note that $f$ extends to a diffeomorphism of $\overline{T}^*U_I$ (fiber radial compactification of $T^*U_I$) preserving fiber infinity. This is consistent with the fact that the classical principal symbol is not affected by a change of time function in the Fourier transform.
\end{remark}

\begin{lemma}
Let $\Psi$ be a symplectomorphism between the symplectic manifolds $(\mathcal{N}_1, \omega_1)$ and $(\mathcal{N}_2, \omega_2)$. Then for any $f \in C^{\infty}(\mathcal{N}_2)$, we have $H_{\Psi^* f} = \Psi^*H_f$. Therefore, if $\gamma$ is an integral curve of $H_f$, $\Psi^{-1}\circ\gamma$ is an integral curve of $H_{\Psi^*f}$.
\end{lemma}
\begin{proof}
By definition, for $x\in \mathcal{N}_1$ and $X \in T_x \mathcal{N}_1$, we have
\begin{align*}
\dd_x (f\circ \Psi) (X) &= \omega_1( H_{\Psi^*f}, X)
\end{align*}
But we also have:
\begin{align*}
\dd_x f\circ \Psi (X) &= \dd_{\Psi(x)} f (\dd_x \Psi(X))\\
&= \omega_2(H_f(\Psi(x)), \dd_x \Psi(X))\\
&= \omega_2(\dd_x \Psi(\dd_x \Psi)^{-1} (H_f(\Psi(x))), \dd_x \Psi(X))\\
&= \omega_1(\Psi^* H_f, X)
\end{align*}
By the non degeneracy of $\omega_1$, we have the equality.
\end{proof}

Using the two previous lemma, we see that we have an identification between the integral curves of $H_{p_h}$ and the integral curves of $H_{p_{2,h}}$ in $T^*U_I$ (and the projection on $U_I$ of two identified integral curves coincide). We use this identification implicitly in the remaining part of this section. In the following, we will make the three following choices in the concrete computations:
\begin{itemize}
\item $I = (r_+-\epsilon, 3M)$ and $f_I$ such that $t_* = \mathfrak{t}+f_I(r)$ (in particular, $f_I(r)$ is constant and the induced symplectomorphism is trivial)
\item $I = (6M, +\infty)$, we choose the time function ${}^*t = t-T(r)$ (in particular, $f_I(r)$ is constant and therefore the induced symplectomorphism is trivial).
\item $I = (r_{\text{min}}, r_{\text{max}})$ ($r_{\text{min}}$ and $r_{\text{max}}$ are constants which will be introduced in Lemmas \ref{convexityHorizon} and \ref{convexityInfty}), we choose $f_I$ such that $t = \mathfrak{t}+f_I(r)$
\end{itemize}
We also need a general elementary lemma about smooth vector fields.
\begin{lemma}
\label{flowNearNDLExtr}
Let $X$ be a smooth vector fields on a manifold $\mathcal{M}$ and $r$ be a smooth function on $\mathcal{M}$. If $\mathcal{M}$ has a boundary, we assume that $X$ is tangent to the boundary (in order to the integral curves to be locally defined near any point of $\mathcal{M}$).
Let $y \in \mathcal{M}$ be a point such that $Xr(y) = 0$ and $X^2r(y)>0$ (resp. $X^2r(y)<0$). Then there exists a neighborhood $V$ of $y$ such that for every integral curve $\gamma$ of $X$ starting at $V$ and maximally defined on an interval $J$ containing $0$, there exists $s >0$ in $J$ with $-s\in J$ such that $r(\gamma(s))>r(y)$ (resp. $r(\gamma(s))<r(y)$), $\dot{r}(s)>0$ (resp. $\dot{r}(s)<0$) and $r(\gamma(-s))>r(y)$ (resp. $r(\gamma(-s))<r(y)$), $\dot{r}(-s)<0$ (resp. $\dot{r}(-s)>0$). 
\end{lemma}
\begin{proof}
We do the proof in the case $X^2r>0$.
There exists an open neighborhood $U$ of $y$ on which $0< \eta <X^2r< 2\eta$. We fix $K\subset U$ a compact neighborhood of $y$. Since the existence time of integral curves in $U$ is lower semi-continuous, there exists $\delta>0$ such that every integral curve starting in $K$ exists and remains in $U$ on the interval $[-\delta, \delta]$ . Then for an integral curve starting at $K$ we have for every $-\delta\leq s\leq \delta$:
\begin{align*}
r(s) \geq r(0)+s\dot{r}(0)+\eta s^2\\
r(s)-r(y) \geq (r(0)-r(y)) + s(\dot{r}(0)+ \eta s)
\end{align*}
By shrinking $K$ until $|\dot{r}(x)|\leq \frac{\eta\delta}{2}$ and $r(x)-r(y) > - \frac{\eta\delta^2}{4}$  for all $x\in K$, we have $r(s)-r(y)>0$ at time $\delta$ and at time $-\delta$.
Moreover, for $s \in [-\delta, \delta]$, we have
\begin{align*}
\dot{r}(0)+s\eta \leq \dot{r}(s) \leq \dot{r}(0) + 2s\eta
\end{align*}
Therefore, using $|\dot{r}(0)|\leq \frac{\eta\delta}{2}$, we get
\begin{align*}
\dot{r}(-\delta)\leq& -\frac{3}{2}\eta\delta<0 \\
\dot{r}(\delta)\geq& \frac{1}{2}\eta\delta>0
\end{align*}
\end{proof}

We now compute the principal symbol near the horizon (using $t_*$ coordinates on the radial interval $I = (r_+-\epsilon, 3M)$):
We use coordinates $\xi = \xi_r \dd r+\zeta \dd \phi_* + \eta \dd \theta$ for cotangent vectors.
\begin{align*}
p_h(x,\xi) =& \Delta_r \xi_r^2 + \frac{1}{\sin^2\theta}\zeta^2 + \eta^2 +2a\zeta\xi_r - 2(a\zeta +(a^2+r^2)\xi_r)z_0 + a^2\sin^2\theta z_0^2 \\
=& \left(\frac{\zeta}{\sin(\theta)}-az_0\sin\theta\right)^2 + \eta^2 + \xi_r\left(-2(a^2+r^2)z_0+2a\zeta +\xi_r\Delta_r\right)\\
=& -\rho^2 G(\xi - z_0\dd t_*)\\
H_{p_h} = &2(-(a^2+r^2)z_0 + a\zeta + \xi_r \Delta_r)\partial_r - (2(r-M)\xi_r^2 - 4rz_0\xi)\partial_{\xi_r} + \left(-2az_0+\frac{2\zeta}{\sin^2\theta}+2a\xi_r\right)\partial_\phi\\
& + 2\eta\partial_\theta - 2\cos\theta\left(z_0^2a^2\sin\theta-\frac{\zeta^2}{\sin^3\theta}\right)\partial_\eta
\end{align*}

\begin{lemma}
The quantities $\zeta$, $p_h$ and $K:=\eta^2+\left(\frac{\zeta}{\sin(\theta)}-z_0a\sin\theta\right)^2$ are invariant under the flow. As a consequence $p_h-K = \xi_r(-2(a^2+r^2)z_0+2a\zeta+\xi_r\Delta_r)$ is also invariant.
\end{lemma}
\begin{proof}
It can be shown by a direct computation using the Hamiltonian vector field but it comes from the invariance of the energy $E$, the norm, the angular momentum $L$, and the Carter constant $K$ along geodesics. 
\end{proof}

\begin{lemma}
\label{AccumulationPointFiberInfinity}
If $y \in S^*U_I$ (at fiber infinity) is an accumulation point of a bicharacteristic $\gamma$ (maximally defined on some interval $J$). Then $\tilde{\xi}(y) = \pm 1$ and $\tilde{\eta}(y) = \tilde{\zeta}(y) = 0$. In particular we also have $\tilde{H}r(y) = 2\tilde{\xi}(y)\Delta_r = \pm 2\Delta_r$
\end{lemma}
\begin{proof}
We denote by $s_n$ an increasing sequence in $J$ such that $\lim\limits_{n\to +\infty} s_n = \sup(J)$ and $\lim\limits_{n\to +\infty}\gamma(s_n) =y$. By the conservation of $\eta^2+\left(\frac{\zeta}{\sin\theta}-az_0\sin\theta\right)^2$, $\lim\limits_{n\to +\infty}\tilde{\eta}(s_n) = \lim\limits_{n\to +\infty}\tilde{\zeta} = 0$ and $\lim\limits_{n\to+\infty}\tilde{\xi}(s_n) = \pm 1$. By continuity of $\tilde{\eta}$, $\tilde{\zeta}$ and $\tilde{\xi}$ we obtain the first part of the lemma. Then replacing in the explicit expression of the Hamiltonian flow at fiber infinity we get $\tilde{H}_{p_h}r(y) = \tilde{\xi}(y) 2\Delta_r = \pm 2 \Delta_r$
\end{proof}

\begin{lemma}
\label{connectedCompo}
If $z_0\in \left\{-1,1\right\}$, then:
$\left\{ -(r^2+a^2\cos^2\theta)z_0-Mr\xi_r = 0\right\}\cap p^{-1}(\left\{0\right\})=\emptyset$.
In particular, the characteristic set has at least two connected components corresponding to both signs of $-(r^2+a^2\cos^2\theta)z_0-Mr\xi_r$ (sign + correspond to the future light cone and - corresponds to the past light cone). Each bicharacteristic curve remains in one of the two connected components for all time.
\end{lemma}

\begin{proof}
We remark that the vector field $T = \partial_{t_*}-\frac{Mr}{\rho^2}\partial_{r_*}$ satisfies $g(T,T) = 1$.
Therefore, if $X$ is a non zero null vector $g(X,T)\neq 0$. Therefore, if $X^{\#}$ is the associated linear form, $X^{\#}(T)\neq 0$. In particular, we remark that for $p^{-1}_{h,z_0}(\xi) = 0$, we have $-\rho^2 G(\xi-z_0\dd {t_*}) = 0$ therefore $\xi-z_0\dd t_* = X^{\#}$ for some non zero null vector $X$. Therefore $-(r^2+a^2\cos^2\theta)z_0+Mr\xi_r = \rho^2(\xi-z_0\dd t_*)(T) = \rho^2g(X,T) \neq 0$.
\end{proof}

\begin{lemma}
\label{signLemma}
If $a^2\cos^2\theta + r^2 - 2Mr>0$ (i.e outside the ergoregion), then for all $\xi$ such that $p_{h,z_0}(\xi)=0$, $z_0$ and $(r^2+a^2\cos^2\theta)z_0+Mr\xi_r$ have the same sign. Therefore for $z_0= -1$, the component $(r^2+a^2\cos^2\theta)z_0+Mr\xi_r>0$ of the characteristic set lies inside the ergoregion. For $z_0=1$ the component$(r^2+a^2\cos^2\theta)z_0+Mr\xi_r<0$ of the characteristic set lies inside the ergoregion.
On $T^*U_I$, $\xi_r$ and $(r^2+a^2\cos^2\theta)z_0+Mr\xi_r$ have the same sign.
\end{lemma}
\begin{proof}
Outside the ergoregion, $T$ and $\partial_t$ are both timelike future oriented. Therefore, for all null vectors $X$, $g(T,X)$ and $g(\partial_t, X)$ have the same sign. Since $p_{h,z_0}(\omega)=0$, we have $\omega-z_0\dd t_* = X^{\#}$ with $X$ null and $(\xi-z_0\dd t_* )(\partial_t)=-z_0$ and $\rho^2(\xi-z_0\dd t_* )(T) = -(r^2+a^2\cos^2\theta)z_0-Mr\xi$ have the same sign.
Similarly, since $\partial_r$ is past oriented on $U_I$, we have $\rho^2(\xi-z_0\dd t_*)(T)$ and $(\xi-z_0\dd t_*)(\partial_r) = \xi_r$ have opposite signs. However, when $\partial_r$ is null, $\xi_r$ can vanish (contrarily to $-(r^2+a^2\cos^2\theta)z_0-Mr\xi_r$).
\end{proof}
Therefore for a fixed $z_0$, we can define $\Sigma_{\pm} := p_{h}^{-1}\left\{0\right\}\cap \left\{\pm((r^2+a^2\cos^2\theta)z_0+Mr\xi_r)>0\right\}$. With this definition $\Sigma_{\text{sgn}(-z_0)}$ lies in the ergoregion.
\begin{remark}
Note that $\tilde{\rho}(r^2+a^2\cos^2\theta)z_0 + Mr\tilde{\xi_r}$ has the same sign as $(r^2+a^2\cos^2\theta)z_0+Mr\xi_r$ on $T^*U_I$. Therefore, we can extend the definition of $\Sigma_{\pm}$ to include the component at fiber infinity given by $\Sigma_{\pm}\cap S^*U_I = \left\{ p_h = 0=\tilde{\rho}, \pm \tilde{\xi}_r > 0 \right\}$. As a consequence, we have $L_{\pm}\subset \Sigma_{\pm}$.
\end{remark}
\begin{lemma}
\label{convexityHorizon}
Let $z_0 = \pm 1$.
There exists $r_+<r_{\text{min}}< 3M$ such that:
$r_+\leq r\leq r_{\text{min}}$, $p_h = 0$ and $H_{p_h}r = 0$ $\Rightarrow$ $H_{p_h}^2r<0$.

For $r=r_+$ we even have that $H_{p_h}r$ does not vanish and has the sign of $\left(-(r^2+a^2\cos^2\theta)z_0-Mr\xi_r\right)$.
\end{lemma}
\begin{proof}
We have:
\begin{align*}
H_{p_h}r = 2(-(a^2+r^2)z_0 + a\zeta + \xi_r \Delta_r).
\end{align*}
First, we handle the particular case of $r = r_+$ (equivalently $\Delta_r=0$). 
Using the fact that $p_h = 0$, we get $\xi_r H_{p_h}r = -\left(\frac{\zeta}{\sin\theta}-az_0\sin\theta\right)^2-\eta^2$. If the right-hand side is equal to $0$, we get $\zeta = az_0\sin^2\theta$ and therefore $H_{p_h}r = 2\left(-(a^2+r^2)z_0+a^2z_0\sin^2\theta\right) = -z_0\rho^2\neq 0$, therefore we must have $\xi_r = 0$ and $H_{p_h}r$ has the same sign as $-(r^2+a^2\cos^2\theta)z_0-Mr\xi_r$. If the right-hand side is negative, we have $\xi_r\neq 0$ and $-(r^2+a^2\cos^2\theta)z_0-Mr\xi_r$ have opposite sign (by lemma \ref{signLemma}), therefore, $H_{p_h}r$ is non zero and has the sign of $\left(-(r^2+a^2\cos^2\theta)z_0-Mr\xi_r\right)$.

Now we assume $r>r_+$. If $H_{p_h}r = 0$, we have $\xi_r\left(-2(a^2+r^2)z_0+2a\zeta +\xi_r\Delta_r\right) = -\Delta_r \xi_r^2$. Using $p_h =0$ (for the first line) and $H_{p_h}r = 0$ for the second line, we get: 
\begin{align*}
\xi_r &= \pm \sqrt{\frac{\eta^2+\left(\frac{\zeta}{\sin\theta}-az_0\sin\theta\right)^2}{\Delta_r}}\\
\xi_r &= -\frac{(-(a^2+r^2)z_0 + a\zeta)}{\Delta_r}
\end{align*}
Let $\epsilon>0$:
We have 2 possibilities:
\begin{enumerate}
\item $\lvert\xi_r\rvert \geq \frac{\epsilon}{\sqrt{\Delta_r}}$\\
\item $\lvert\xi_r\rvert \geq \frac{r^2-\epsilon}{\Delta_r}$
\end{enumerate}
To prove this, assume that $\lvert\frac{\zeta}{\sin\theta}-az_0\sin\theta\rvert \leq \epsilon$ (otherwise, we are in the first case). Then, we have 
\begin{align*}\lvert a\zeta \rvert &\leq a^2\sin^2\theta + \lvert\sin\theta\rvert\epsilon \leq a^2\sin^2\theta + \epsilon\\
 \lvert-(a^2+r^2)z_0+a\zeta\rvert &\geq a^2+r^2 - \lvert a \zeta \rvert\\
& \geq a^2\cos^2\theta + r^2 - \epsilon \\
&\geq r^2-\epsilon
\end{align*}
Therefore, we are in the second case. We can use this with $\epsilon = \frac{r_+^2}{2}$ and we find $\lvert \xi_r\rvert\geq \frac{r_+^2}{2\sqrt{\Delta_r}}$ as soon as $\Delta_r\leq 1$.

Now we compute:
\begin{align*}
H_{p_h}^2r = -4\xi_r(-2rz_0 + \xi_r(M-r))\Delta_r + 4(-rz_0-\xi_r(M-r))H_{p_h}(r)
\end{align*}
Therefore, when $H_{p_h}r = 0$, we have $H_{p_h}^2r <0 $ as soon as $\lvert\xi_r\rvert \geq \frac{2r}{r-M}$ and this is true when $r>r_+$ is close to $r_+$ since $ \frac{2r}{r-M} = o(\frac{r_+^2}{2\sqrt{\Delta_r}})$ when $r\rightarrow r_+$.
\end{proof}

\begin{lemma}
\label{HrBeyondHorizon}
On $p_{h}^{-1}(\left\{0\right\})$ we have:
$\xi_r H_{p_h}r\leq 0$ when $\Delta \leq 0$. If $\Delta<0$ and $\xi_r\neq 0$, we have a strict inequality. If $\xi_r = 0$, $H_{p_h}r$ has the same sign as $-z_0$ (striclty) on $T^*U_I$.
\end{lemma}
\begin{proof}
We have:
$H_{p_h}r = 2(-(a^2+r^2)z_0 + a\zeta + \xi_r \Delta_r)$. Moreover, since $p_h = 0$ and $\Delta_r\leq 0$, we must have
$\xi_r\left(-2(a^2+r^2)z_0+2a\zeta +\xi_r\Delta_r\right)\leq 0 \leq -\Delta_r \xi_r^2$.
Therefore $\xi_r H_{p_h}r\leq 0$ with strict inequality when $\xi_r\neq 0$ and $\Delta_r <0$.
If $\xi_r = 0$, we have (on the characteristic set):
\begin{align*}
\left(\frac{\zeta}{\sin\theta} - az_0\sin\theta\right)^2+\eta^2=0
\end{align*}
and therefore
\begin{align*}
\eta =& 0\\
\zeta =& az_0\sin^2\theta
\end{align*}
If we replace in the expression of $H_{p_h}r$, we get:
\begin{align*}
H_{p_h}r =& 2(-(a^2+r^2)z_0 + a^2z_0\sin^2\theta)\\
=& -2\rho^2z_0
\end{align*}
\end{proof}

\begin{lemma}
\label{signXiBeyondHorizon}
Let $\gamma$ be a bicharacteristic defined on an interval $J$. Let $J_I = \left\{ s\in J: r(s)\in I\right\}$. On each connected component of $J_I$, we have the alternative $\xi_r(s) = 0$ or $\xi_r$ never vanishes.
\end{lemma}
\begin{proof}
Assume that $\xi_r(s_0)=0$ for some $s_0 \in J$. Then using that $p_h = 0$, we have $\eta^2 + \left(\frac{\zeta}{\sin\theta} - az_0\sin\theta\right)^2 = 0$. Therefore $\zeta = az_0\sin^2\theta$. This property is true for all $s\in J$ by the conservation law. Therefore for all $s\in J$ $\xi_r(-2(a^2\cos^2\theta+r^2)z_0 +\Delta_r \xi_r) = 0$. Since $\frac{2(a^2\cos^2\theta + r^2)z_0}{\Delta_r}>\delta>0$ for some constant $\delta$ independent of $r$ and $\theta$, the continuity of $\xi_r$ implies that $\xi_r(s) = 0$ for all $s$ in the connected component of $s_0$ in $J_I$.
\end{proof}

\begin{lemma}
\label{characteristicDecreasing}
Let $\gamma$ be a bicharacteristic defined on an interval $J$. We assume that there exists $s_0\in J$ such that $r(s_0)\leq r_{\text{min}}$ and $\dot{r}(s_0)<0$. Then we have two cases:
\begin{itemize}
\item
The curve $\gamma$ remains in $\left\{ r_+\leq r\right\}$ for all $s \in [s_0, \sup J)$ in which case $\sup J = +\infty$ and $\gamma$
tends to $L_+$ or to $L_-$.
\item 
$\gamma$ reach $\left\{r<r_+\right\}$ for some $s_1>s_0$ and in this case, $\sup J<+\infty$ and $\lim\limits_{s\to \sup J}r(\gamma(s))= r_+-\epsilon$.
\end{itemize}
\end{lemma}
\begin{proof}
Note that by lemma \ref{convexityHorizon} and lemma \ref{HrBeyondHorizon}, $r\circ \gamma$ cannot have a local minimum on $T^*U_I$ and therefore it must be strictly decreasing on $[s_0, \sup J)$. In particular, $\gamma$ remains in $\left\{r<r_{min}\right\}$ for $s>s_0$ where $r_{min}$ is defined in Lemma \ref{convexityHorizon}.

First assume that $\gamma$ does not reach $\left\{r<r_+\right\}$. The curve $\gamma$ remains in a compact set and therefore has an accumulation point $y \in \overline{T^*U_I}$ with $r(y) = r_\infty\geq r_+$ and $r_\infty = \inf\left\{r(\gamma(s)), s\geq s_0\right\}$. We denote by $s_n$ an increasing sequence in $J$ such that $\lim\limits_{n\to +\infty} s_n = \sup(J)$ and $\lim\limits_{n\to +\infty}\gamma(s_n) =y$. If $y$ is at fiber infinity, then by lemma \ref{AccumulationPointFiberInfinity}, we have $\tilde{H}_{p_h}r(y) = \pm 2\Delta_r$ and since $r_\infty \leq r(s)$ for all $s\in (s_0, \sup J)$, we conclude that $\tilde{H}_{p_h}r(y) = 0$ and therefore $r_\infty = r_+$ and $y \in L_+ \cup L_-$. Moreover, $y$ cannot be in $T^*U_I$, otherwise we have either $H_{p_h}r \neq 0$ or $H_{p_h}^2r <0$ by lemma \ref{convexityHorizon} and both cases contradicts $r_\infty = \inf_{s\in (s_0, \sup J)} r(s)$ (see lemma \ref{flowNearNDLExtr} for the second case). Therefore, the only possible accumulation points are $L_+$ and $L_-$ and by lemma \ref{connectedCompo}, $\gamma$ cannot have both accumulation points in $L_+$ and in $L_-$. Therefore, $\gamma$ tends to $L_+$ or to $L_-$ at $\sup J$ (and since $\gamma$ remains in a compact set $\sup J = +\infty$).

Now we assume that there exists $s_1>s_0$ such that $r(s_1)<r_+$. By Lemma \ref{HrBeyondHorizon}, $r\circ \gamma$ is then strictly decreasing on $(s_1, \sup J)$. The curve $\gamma$ cannot have any accumulation point in $\left\{r\geq r(s_1)\right\}$ at $\sup J$. Indeed, at such an accumulation point $y$, $\tilde{H}_{p_h}r(y)$ must vanish to be consistent with the fact that $r$ is decreasing but this never happens when $r<r_+$ (at fiber-infinity it can be ruled out since we show as previously $\tilde{H}_{p_h}(y) = \pm 2\Delta_r$ and elsewhere, it is a consequence of lemma \ref{HrBeyondHorizon}).
\end{proof}

We have an analog of lemma \ref{characteristicDecreasing} (with a similar proof) describing the behavior of some bicharacteritics in the past:
\begin{lemma}
\label{characteristicDecreasingPast}
Let $\gamma$ be a bicharacteristic defined on an interval $J$. We assume that there exists $s_0\in J$ such that $r(s_0)\leq r_{\text{min}}$ and $\dot{r}(s_0)>0$. Then we have two cases:
\begin{itemize}
\item
The curve $\gamma$ remains in $\left\{ r_+\leq r\right\}$ for all $s \in (\inf J, s_0]$ in which case $\inf J = -\infty$ and $\gamma$
tends to $L_+$ or to $L_-$.
\item 
$\gamma$ reach $\left\{r<r_+\right\}$ for some $s_1<s_0$ and in this case $Inf J>-\infty$ and $\lim\limits_{s\to \inf J}r(\gamma(s)) = r_+-\epsilon$.
\end{itemize}
\end{lemma}

\begin{lemma}
\label{characteristicIncreasing}
Let $\gamma$ be a bicharacteristic curve with at least one point in $\left\{r\leq r_{min,trapp}\right\}$. We assume that for all $s\in J$ such that $\gamma(s)\in \left\{ r\leq r_{\text{min}}\right\}$, we have $\dot{r}(s)\geq 0$. Then one of the following cases holds:
\begin{itemize}
\item Case 1: There exists $s_0 \in J$ such that $r(s_0)>r_+$. In this case, $\gamma$ reaches $\left\{r>r_{\text{min}}\right\}$ in finite time $s_1$ and remains in this set for all $s_1<s<\sup J$.
\item Case 2: For all $s\in J$, $r(s)\leq r_+$. In this case, $\sup J = +\infty$ and $\gamma$ tends to $L_+$ or to $L_-$ when $s \to +\infty$.
\end{itemize}
\end{lemma}

\begin{proof}
We use the notations of the lemma.
Assume that there exits $s_0\in J$ such that $r(s_0)>r_+$. Then we prove that $\gamma$ reaches $\left\{r>r_{\text{min}}\right\}$ by contradiction. Indeed, if it is not the case, there exists $y\in \overline{T^*U_I}$ with $r(y) = \sup_{s\in (s_0, \sup J)} r(\gamma(s)) > r_+$ and a sequence $s_n \in (s_0, \sup J)$ with $\lim\limits_{n\to +\infty} s_n = \sup J$ and $\lim\limits_{n\to +\infty}\gamma(s_n) = y$. If $y$ is at fiber infinity, by lemma \ref{AccumulationPointFiberInfinity} $\tilde{H}r(y) = \pm\Delta_r \neq 0$ which contradicts $r(y) = \sup_{s\in (s_0, \sup J)} r(\gamma(s))$. If $y$ is not at fiber infinity, for the same reason we must have $H_{p_h}r(y)=0$ and lemma \ref{convexityHorizon} implies that $H_{p_h}^2 r(y)<0$. But in this case, lemma \ref{flowNearNDLExtr} provides a point $s\in J$ with $\gamma(s)\in T^*U_I$ and $\dot{r}(s)<0$ which contradicts the hypothesis on $\gamma$. Therefore, such an accumulation point $y$ cannot exist and $\gamma$ reaches $\left\{r>r_{\text{min}}\right\}$ at some time $s_1\in J$. Moreover, for all $s\in J$ with $s>s_1$, $r(s) > r_{\text{min}}$ (otherwise we get a contradiction at $s_2 = \inf \left\{ s\in J: s>s_1, r(s)= r_{\text{min}}\right\}$ where we must have $\dot{r}(s_2) = 0$ since by hypthesis $\dot{r}(s_2)\geq 0$ and using lemma \ref{convexityHorizon}).

Assume that for all $s\in J$, $r(s) \leq r_+$. In particular $\gamma$ remains in a compact set and $\sup J = +\infty$. Then, by compactness, there exists $y\in \overline{T^*U_I}$ with $r(y) = \sup \left\{ r(s), s\in J\right\}\leq r_+$ (this implies that $\tilde{H}_{p_h}r(y) = 0$) such that $y$ is an accumulation point for $\gamma$ at $\sup J$. By lemmas \ref{HrBeyondHorizon} (for the case $r(y) <r_+$) and \ref{convexityHorizon} (for the case $r(y)=r_+$), $y$ cannot be in $T^*U_I$ therefore $y$ is at fiber infinity. By lemma \ref{AccumulationPointFiberInfinity}, the only possibility is $y\in L_\pm$. As usual we use lemma \ref{connectedCompo} to prove that $\gamma$ cannot have accumulation points both in $L_+$ and in $L_-$. We conclude that $\gamma$ tends to $L_+$ or to $L_-$.
\end{proof}

We also have an analog of lemma \ref{characteristicIncreasing} for past behavior:
\begin{lemma}
\label{characteristicIncreasingPast}
Let $\gamma$ be a bicharacteristic curve with at least one point in $\left\{r\leq r_{min,trapp}\right\}$. We assume that for all $s\in J$ such that $\gamma(s)\in \left\{ r\leq r_{\text{min}}\right\}$, we have $\dot{r}(s)\leq 0$. Then:
\begin{itemize}
\item Case 1: There exists $s_0 \in J$ such that $r(s_0)>r_+$. In this case, $\gamma$ reaches $\left\{r>r_{\text{min}}\right\}$ in finite time $s_1$ and remains in this set for all $\inf J<s<s_1$.
\item Case 2: For all $s\in J$, $r(s)\leq r_+$. In this case, $\inf J = -\infty$ and $\gamma$ tends to $L_+$ or to $L_-$ when $s \to -\infty$.
\end{itemize}
\end{lemma}

We now compute the principal symbol near infinity (using the time coordinate ${}^*t = t-T(r)$ on the radial interval range $(6M,+\infty)$):
We define $\xi = \xi_r\dd r + \zeta \dd \phi + \eta \dd \theta$.
\begin{align*}
p_h(\xi) =& -\rho^2 G(\xi - z_0\dd {}^*t)\\
=& a^2\sin^2\theta z_0^2-\left(\frac{a^2}{\Delta_r}-\frac{1}{\sin^2\theta}\right)\zeta^2+\frac{4Mar}{\Delta_r}\zeta z_0 + \eta^2+ \Delta_r \xi_r^2 + 2(a^2+r^2)z_0\xi_r \\
=& \eta^2+\left(\frac{\zeta}{\sin\theta}+az_0\sin\theta\right)^2 -\frac{a^2}{\Delta_r}\zeta^2-2a\zeta z_0+\Delta_r \xi_r^2+2(a^2+r^2)z_0\xi_r + \frac{4Mar}{\Delta_r}\zeta z_0\\
H_{p_h}(\xi) =& \left(2\xi_r\Delta_r + 2z_0(a^2+r^2)\right)\partial_r + 2\left(-a^2z_0^2\sin\theta+\frac{\zeta^2}{\sin^3\theta}\right)\cos\theta\partial_\eta\\
&+ 2\left(-\frac{2Maz_0\zeta}{\Delta_r}+\frac{a\zeta(M-r)(-4Mrz_0+a\zeta)}{\Delta_r^2}+\xi_r(-2rz_0+\xi_r(M-r))\right)\partial_{\xi_r}\\
& + 2\eta\partial_\theta  + 2\left(\frac{2Marz_0}{\Delta_r}-\zeta\left(\frac{a^2}{\Delta_r}-\frac{1}{\sin^2\theta}\right)\right)\partial_\phi
\end{align*}
\begin{remark}
The quantities $p_h$, $\zeta$ and $\eta^2+\left(\frac{\zeta}{\sin\theta}+az_0\sin\theta\right)^2$ are still conserved along the bicharacteristics.
\end{remark}

\begin{lemma}
\label{convexityInfty}
There exists $r_{\text{max}}>6M$ such that in $\left\{ r\geq r_{\text{max}} \right\}\cap p_h^{-1}(\left\{0\right\})$, the following implication holds:
$H_{p_h}r = 0 \Rightarrow H_{p_h}^2r>0$
\end{lemma}
\begin{proof}
On the set $\left\{H_{p_h}r=0\right\}$, we have:
\begin{align*}
\xi_r =& -\frac{a^2+r^2}{\Delta_r}z_0\\
=& -z_0+O(r^{-1}).
\end{align*}
On $\left\{H_{p_h}r=0\right\}\cap p_h^{-1}(\left\{0\right\})$, we have:
\begin{align*}
\eta^2+\left(\frac{\zeta}{\sin\theta}-az_0\sin\theta\right)^2+\left(-2az_0-\frac{a^2}{\Delta_r}\zeta\right)\zeta =& \Delta_r\xi_r^2\\
=& r^2+ O(r)
\end{align*}
Therefore, there exists a constant $B>0$ and $r_{c,1}>C'$ such that $|\zeta|\leq Br$ on this set for $r\geq r_{c,1}$. Otherwise, there would be a sequence $(\zeta_n, \eta_n, r_n, (\xi_r)_n)$ with $|\zeta_n|\geq nr_n$ and the right-hand side would be bigger than $\frac{n^2 r_n^2}{2}$ for $r$ large enough which is impossible. Similarly we can find $B'>0$ and $r_{c,2}>C'$ such that $|\eta|\leq B'r$ on this set for $r\geq r_{c,2}$. We can now compute, still on the set $\left\{H_{p_h}r=0\right\}\cap p_h^{-1}(\left\{0\right\})$ (on this set the expression simplifies):
\begin{align*}
H_{p_h}^2r =& -8Marz_0\zeta + 4\frac{a\zeta(M-r)}{\Delta_r}\left(-4Mrz_0+a\zeta\right)+4\xi_r\left(-2rz_0+\xi_r(M-r)\right)\Delta_r
\end{align*}
Using our previous estimates, we see that the dominant term is $-4\xi_r\Delta_r rz_0$ and $H_{p_h}^2 r = 4r^3+O(r^2)$ therefore, it is positive for $r$ larger than some $r_{c}>C'$. 
\end{proof}

\begin{remark}
Note that we are interested in the flow of the renormalized vector field $\tilde{H}_{p_h}:=r^{-1}H_{p_h}$ which extends to a continuous vector field on ${}^{sc}T^*(C',+\infty]\times\mathbb{S}^2$ (we do not need to compactify the fiber since the characteristic set of $p_h$ does not intersect fiber infinity) and using smooth (up to the boundary at infinity) coordinates $x:= \frac{1}{r}$, $\xi_{r}$, $\eta_{sc} = \frac{\eta}{r}$, $\zeta_{sc} = \frac{\zeta}{r}$ we find
\begin{align*}
\tilde{H}_{p_h} =& 2 x \left(- \xi_{r} \left(x^2\Delta_r\right) - z_0 \left(a^{2} x^{2} + 1\right)\right)\partial_{x}\\
& +2  \left(- \eta_{sc} \left(\xi_{r} \left(x^2\Delta_r\right) + z_0 \left(a^{2} x^{2} + 1\right)\right) + \left(- a^{2} x^{2} z_0^{2} \sin{\left (\theta \right )} + \frac{\zeta_{sc}^{2}}{\sin^3\theta}\right) \cos{\left (\theta \right )}\right)\partial_{\eta_{sc}}\\
& - 2 \zeta_{sc} \left(\xi_{r} \left(x^2\Delta_r\right) + z_0 \left(a^{2} x^{2} + 1\right)\right)\partial_{\zeta_{sc}}\\
& +\frac{2}{\left(x^2\Delta_r\right)^{2}} \left(-2 M a x^{2} z_0 \zeta_{sc} \left(x^2\Delta_r\right) + a x^{2} \zeta_{sc} \left(M x - 1\right) \left(-4 M z_0 + a \zeta_{sc}\right)\right.\\
&\hphantom{+\frac{2}{\left(x^2\Delta_r\right)^{2}} \left(\right.}\left. + \xi_{r} \left(\xi_{r} \left(M x - 1\right) - 2 z_0\right) \left(x^2\Delta_r\right)^{2}\right)\partial_{\xi_{r}}\\
& +2 \eta_{sc}\partial_{\theta} +\frac{2 \left(-2 M a x^{2} z_0 \sin^{2}{\left (\theta \right )} + \zeta_{sc} \left(2 M x - a^{2} x^{2} \cos^{2}{\left (\theta \right )} - 1\right)\right)}{\left(x^2\Delta_r\right) \sin^{2}{\left (\theta \right )}}\partial_{\phi} 
\end{align*}
Therefore, the restriction of $\tilde{H}_{p_h}$ to the boundary $\left\{x=0\right\}$ is
\begin{align*}
2 \xi_r \left(- \xi_r - 2 z_0\right)\partial_{\xi_r}+ \left(- 2 \eta_{sc} \xi_r - 2 \eta_{sc} z_0 + \frac{2 \zeta_{sc}^{2} \cos{\left (\theta \right )}}{\sin^{3}{\left (\theta \right )}}\right)\partial_{\eta_{sc}} +2 \zeta_{sc} \left(- \xi_r - z_0\right)\partial_{\zeta_{sc}}  +2 \eta_{sc}\partial_{\theta} - \frac{2 \zeta_{sc}}{\sin^{2}{\left (\theta \right )}}\partial_{\phi}
\end{align*}
As expected, it vanishes only for $(\xi_r, \zeta_{sc}, \eta_{sc})\in \left\{(0,0,0), (-2z_0, 0,0)\right\}$. We define $\mathcal{R}_{in}:=\left\{ \xi_r = \zeta_{sc} = \eta_{sc} = x= 0\right\}$ and $\mathcal{R}_{out}:=\left\{ \xi_r = -2z_0, \zeta_{sc} = \eta_{sc} = x = 0 \right\}$.
\end{remark}

\begin{lemma}
\label{characteristicInfIncreasing}
Let $\gamma$ be a bicharacteristic of the renormalized Hamiltonian vector field maximally defined on an interval $J$. If, there exists $s_0\in J$ such that $\dot{r}(s_0) >0$ and $r(s_0)\geq r_{\text{max}}$, then $\sup J = +\infty$ and $\gamma(s)$ converges to $\mathcal{R}_{in}$ or to $\mathcal{R}_{out}$ when $s\to +\infty$.
\end{lemma}
\begin{proof}
By the previous lemma, $\dot{r}(s)>0$ for all $s>s_0$ (and $\gamma( [s_0, \sup J))\subset \left\{r_{\text{max}}\leq r \leq +\infty\right\}$). Therefore $r(s)$ has a limit.  
First, we prove by contradiction that this limit is $+\infty$. Assume that $\lim\limits_{s\to \sup J}\gamma(s) = r_\infty<+\infty$. Then the curve remains in the compact set $\left\{r_{\text{max}}\leq r \leq +\infty\right\}\cap p_h^{-1}(\left\{0\right\})$ and has therefore an accumulation point $y$ with $r(y) = r_\infty$. However, $H_{p_h}r (y) = 0$ (otherwise, $\gamma$ should cross the hypersurface $\left\{r=r_\infty\right\}$ (lemma \ref{flowNearNDLExtr}) and therefore by lemma \ref{convexityInfty}, we have $H_{p_h}^2r>0$ but this contradicts the fact  $r_\infty = \sup_{s\in[s_0, \sup J)} \gamma(s)$ by lemma \ref{flowNearNDLExtr}.

Therefore, $r_\infty = +\infty$. Then by conservation of $\eta^2+\left(\frac{\zeta}{\sin\theta}-az_0\sin\theta\right)^2$ along the flow, we have $\lim\limits_{s\to +\infty}\eta_{sc}(s)=0$ and $\lim\limits_{s\to +\infty}\zeta_{sc}(s)=0$. Finally, using that $p_h = 0$ along the flow, we get that $\lim\limits_{s\to +\infty}\xi_r(s)(\xi_r(s) + 2z_0) = 0$ and using the continuity of $\xi_r(s)$, either $\lim\limits_{s\to +\infty}\xi_r(s) = 0$ or $\lim\limits_{s\to +\infty}\xi_r(s) = -2z_0$.
\end{proof}

The analog of lemma \ref{characteristicInfIncreasing} for past behavior is 
\begin{lemma}
\label{characteristicInfIncreasingPast}
Let $\gamma$ be a bicharacteristic maximally defined on an interval $J$. If, there exists $s_0\in J$ such that $\dot{r}(s_0) <0$ and $r(s_0)\geq r_{\text{max}}$, then $\inf J = -\infty$ and $\gamma(s)$ converges to $\mathcal{R}_{in}$ or to $\mathcal{R}_{out}$ when $s\to -\infty$.
\end{lemma}

To cover the behavior of all bicharacteristics in $T^* U_I$ we need the following lemma
\begin{lemma}
Let $\gamma$ be a bicharacteristic maximally defined on an interval $J$. We assume that for all $s\in J$ such that $r(s)\geq r_{\text{max}}$ we have $\dot{r}(s)\leq 0$. Then $\gamma$ reaches $\left\{r<r_{\text{max}}\right\}$ for some $s_0 \in J$ and stays in this set for $s_0\leq s <\sup J$.
\end{lemma}
\begin{proof}
As in the proof of lemma \ref{characteristicInfIncreasing}, we can prove using lemma \ref{convexityInfty} and lemma \ref{flowNearNDLExtr} that $\gamma$ cannot have an accumulation point at $\sup J$ in $\left\{r\geq r_{\text{max}}\right\}$ (in this case, accumulation points at $r= +\infty$ are impossible since $r\circ \gamma$ is decreasing). Therefore, $\gamma$ reaches $\left\{ r< r_{\text{max}}\right\}$ in finite time $s_0$. Moreover, if $r(s)\geq r_{\text{max}}$ for some $s>s_0$, then we could define $s_1 = \inf\left\{ s_0<s< \sup J: r(s) = r_{\text{max}} \right\}$. Since $\dot{r}(s_1)\leq 0$ by hypothesis and cannot be negative me must have $\dot{r}(s_1) = 0$ and by lemma \ref{convexityInfty} and lemma \ref{flowNearNDLExtr} we get a contradiction.
\end{proof}

Similarly we have
\begin{lemma}
Let $\gamma$ be a bicharacteristic maximally defined on an interval $J$. We assume that for all $s\in J$ such that $r(s)\geq r_{\text{max}}$ we have $\dot{r}(s)\geq 0$. Then $\gamma$ reaches $\left\{r<r_{\text{max}}\right\}$ in finite time and stays in this set.
\end{lemma}

We now compute the principal symbol near the trapped set (using the $t$ time coordinate on the radial interval $I=(r_{\text{min}}, r_{\text{max}})$):
We use coordinates $\xi = \xi_r\dd r + \zeta \dd \phi + \eta \dd \theta$ for cotangent vectors and we introduce the function $\alpha = \frac{-(r^2+a^2)z_0+a\zeta}{\Delta_r}$.
\begin{align*}
p_h(\xi) =& -\rho^2 G(\xi-z_0\dd t)\\
=& -\left(\frac{(a^2+r^2)^2}{\Delta_r}-a^2\sin^2\theta\right)z_0^2-\left(\frac{a^2}{\Delta_r}-\frac{1}{\sin^2\theta}\right)\zeta^2 + \left(\frac{2a(a^2+r^2)}{\Delta_r}-2a\right)\zeta z_0 + \eta^2+\Delta_r \xi_r^2\\
=& \Delta_r (\xi_r^2 - \alpha^2) + \eta^2 + \left(-az_0\sin\theta + \frac{\zeta}{\sin\theta}\right)^2 \\
H_{p_h} =& 2\xi_r\Delta_r \partial_r + \left(\alpha(-4rz_0-2(r-M)\alpha)-2\xi_r^2(r-M)\right)\partial_{\xi_r}+2\eta\partial_\theta + 2\cos\theta\left(\frac{\zeta^2}{\sin^3\theta}-a^2z_0^2\sin\theta\right)\partial_\eta\\
& + \frac{2}{\Delta_r}\left(-2Maz_0+a^2\zeta-\frac{\Delta_r}{\sin^2\theta}\zeta\right)\partial_{\phi}
\end{align*}
\begin{remark}
We still have the conserved quantities: $\zeta$, $p_h$ and $\eta^2+\left(-az_0\sin\theta  + \frac{\zeta}{\sin\theta}\right)^2$.
\end{remark}
Note that due to Lemmas \ref{characteristicDecreasing}, \ref{characteristicDecreasingPast}, \ref{characteristicInfIncreasing} and \ref{characteristicInfIncreasingPast}, we see that bicharacteristics leaving $(r_{\text{min}}, r_{\text{max}})$ exit any compact subset of $(r_+,+\infty)\times \mathbb{S}^2$ either in the past or in the future. Therefore, trapped geodesics remain in $T^*U_I$ for all time.
Proposition 3.2 of \cite{dyatlov2015asymptotics} translates\footnote{For more details about the translation from properties of the rescaled geodesic flow and the semiclassical flow see \cite{dyatlov2015asymptotics} section 2.3} in our context to:
\begin{prop}
Let $K_{z_0}$ be the union of bicharacteristics of $p_h(z_0)$ which are contained in $T^*U_I$ with $I = (r_{\text{min}}, r_{\text{max}})$. Then $K_{z_0} = \left\{ G(\xi-z_0\dd t) = \xi_r = ((r-M)\alpha + 2rz_0)\Delta_r = 0, (\sigma, \eta, \zeta)\neq (0,0,0)\right\}$. 
\end{prop}

Proposition 3.5 of \cite{dyatlov2015asymptotics} translates to
\begin{prop}
\label{ExprGammaPm}
Let $\Gamma_{\pm}$ be the union of bicharacteristics that are contained in $T^*U_I$ in the future (-) or in the past (+). We denote by $\hat{K}$ the projection of $K_{z_0}$ on $\mathbb{S}^2\times\R^2_{\zeta,\eta}$.  Then we have 
\begin{align*}
\Gamma_{\pm} =& \left\{(r,\hat{x},\xi_r, \hat{\xi}): (\hat{x},\hat{\xi})\in \hat{K}, \xi_r = \mp z_0\text{sgn}(r-r'_{\hat{x},\hat{\xi}})\sqrt{\frac{\Phi_{\hat{x},\hat{\xi}}(r)}{\Delta_r}}\right\}
\end{align*}
where
\begin{align*}
\Phi_{\hat{x},\hat{\xi}} = -\left(\eta^2+\left(-az_0\sin\theta  + \frac{\zeta}{\sin\theta}\right)^2\right) + \Delta_r \alpha^2
\end{align*}
and $r'_{\hat{x},\hat{\xi}}$ is the unique solution of $\Phi_{\hat{x},\hat{\xi}}(r)=0$.
\end{prop}

Finally, we need the following lemma which is implicit in \cite{dyatlov2015asymptotics}.
\begin{lemma}
\label{characteristicGammaPm}
Let $\gamma$ be a bicharacteristic included in $\Gamma_-$ (resp. $\Gamma_+$). We define $\hat{x}:=(\theta, \phi)(\gamma(0))$ and $\hat{\xi}:=(\zeta, \eta)(\gamma(0))$. Then, $\gamma$ is defined on an interval $J$ with $\sup(J) = +\infty$ (resp. $\inf(J)=-\infty$) and 
\begin{align*}
\forall s\in J, (\theta,\phi, \eta,\zeta)(\gamma(s))\in \hat{K}\\ 
\xi_r(\gamma(s))\xrightarrow[s\to \pm\infty]{}0\\
r(\gamma(s))\xrightarrow[s\to \pm\infty]{}r'_{\hat{x},\hat{\xi}}
\end{align*}
Therefore $\gamma$ is included in any neighborhood of $K_{z_0}$ for $s$ large enough (resp. for $-s$ large enough).
\end{lemma}
\begin{proof}
The fact that for all $s\in J$, $(\theta,\phi, \eta,\zeta)(\gamma(s))\in \hat{K}$ is a consequence of the fact that $\Gamma_{\pm}$ are stable by the flow and of the explicit expression given in proposition \ref{ExprGammaPm}. Moreover, we see that $\Phi_{\hat{x},\hat{\xi}}$ only depends on $\hat{x}$ and $\hat{\xi}$ through $\eta^2+\left(-az_0\sin\theta  + \frac{\zeta}{\sin\theta}\right)^2$ and $\zeta$ which are constant along the flow. Therefore, for all $s\in J$, $\xi_r(\gamma(s)) = \pm z_0\text{sgn}(r-r'_{\hat{x},\hat{\xi}})\sqrt{\frac{\Phi_{\hat{x},\hat{\xi}}(r)}{\Delta_r}}$.
Note that 
\begin{align*}\Phi_{\hat{x},\hat{\xi}}(\gamma(s)) =& \Delta_r\xi_r^2-p_h\\
=& \Delta_r\xi_r(\gamma(s))^2
\end{align*}
Therefore, it is enough to prove the convergence of $r(\gamma(s))$ to $r'_{\hat{x},\hat{\xi}}$, the convergence of $\xi_r(\gamma(s))$ will follow from the fact that $\Phi_{\hat{x},\hat{\xi}}(r'_{\hat{x},\hat{\xi}})=0$. We do the case of $\gamma\subset \Gamma_-$ and $z_0=-1$, the other cases are similar. To alleviate the notation, we define $r(s) := r(\gamma(s))$ and $\xi_r(s) := \xi_r(\gamma(s))$ and we denote by a dot the derivative with respect to $s$.
For all $s\in J$, we have 
\begin{align}
\dot{r}(s) =& H_{p_h}r (\gamma(s)) \notag\\
 =& 2\Delta_r\xi_r(s)\notag \\
=& -\text{sgn}(r-r'_{\hat{x},\hat{\xi}})\sqrt{\Delta_r\Phi_{\hat{x},\hat{\xi}}(r(s))} \label{equationOnR}
\end{align}
In particular $\dot{r}$ is negative when $r>r'_{\hat{x},\hat{\xi}}$ and positive when $r<r'_{\hat{x},\hat{\xi}}$. Moreover, at $r'_{\hat{x},\hat{\xi}}$ which is its only vanishing point, $\Phi_{\hat{x},\hat{\xi}}$ vanishes at order 2 (see the discussion before proposition 3.5 in \cite{dyatlov2015asymptotics}) therefore, the right-hand side is locally Lipschitz with respect to $r$. We conclude that we can solve \eqref{equationOnR} with initial condition $r(0)$ on an interval $\tilde{J}$ with $\sup \tilde{J} = +\infty$ and $\lim\limits_{s\to +\infty} r(s)= r'_{\hat{x},\hat{\xi}}$. But it means that $\gamma(s)$ remains in a compact set when $s\to \sup J$ (the characterstic set of $p$ intesected with $\left\{r'_{\hat{x},\hat{\xi}}-\epsilon\leq r \leq r'_{\hat{x},\hat{\xi}}+\epsilon\right\}$ is compact for $\epsilon$ small enough). Therefore $\sup{J} = +\infty$ and by uniqueness of the solution of \eqref{equationOnR}, $\lim\limits_{s\to+\infty} r(\gamma(s)) = r'_{\hat{x},\hat{\xi}}$.

Finally, the last statement of the lemma follows from the fact that \[F:=\left\{(r'_{\hat{x},\hat{\xi}},(\theta,\phi)(\gamma(s)),0,(\zeta,\eta)(\gamma(s))), s\in J\right\}\subset K_{z_0}\] and $\lim\limits_{s\to +\infty}d(\gamma(s), F)=0$ for any distance $d$ inducing the topology of $T^*U_I$.
\end{proof}

We can now completely describe the asymptotic behavior of bicharacteristics:
\begin{prop}
\label{sythesisScFlow}
We define the surface $B_{\epsilon} = \left\{ r=r_+-\epsilon\right\} \subset {}^{sc}\overline{T}^*(r_+-2\epsilon, +\infty]\times \mathbb{S}^2$. 
In this proposition, we say that a curve $\gamma$ defined on some interval $J$ is of type $(A,B)$ where $A$ and $B$ are two sets in $\left\{ L_-, L_+, B_{\epsilon}, \mathcal{R}_{in}, \mathcal{R}_{out}\right\}$ if $\gamma$ tends to $A$ at $\inf J$ and to $B$ at $\sup J$. Moreover, if $A  = B_{\epsilon}$ (resp. $B = B_{\epsilon}$), $\inf J > - \infty$ (resp. $\sup J < +\infty$), in all the other case the corresponding bound of $J$ is infinite. 

Let $\gamma$ be a bicharacteristic for the renormalized Hamiltonian flow on ${}^{sc}\overline{T^*(r_+-\epsilon, +\infty]\times \mathbb{S}^2}$ maximally defined on some interval $J$.

Let $z_0 = -1$. 
\begin{itemize}
\item If $\gamma \subset \Sigma_{+}$, then either $\gamma \subset L_+$ or $\gamma$ is of type $(L_+, B_{\epsilon})$.
\item If $\gamma \subset \Sigma_{-}$, then either $\gamma \subset L_- \cup K_{z_0} \cup \mathcal{R}_{in} \cup \mathcal{R}_{out}$ or $\gamma$ is of type $(\mathcal{R}_{in}, L_-)$, $(\mathcal{R}_{in}, K_{z_0})$, $(\mathcal{R}_{in},\mathcal{R}_{out})$ or $( B_{\epsilon}, L_-)$, $(B_{\epsilon}, K_{z_0})$, $(B_{\epsilon}, \mathcal{R}_{out})$, $(K_{z_0},L_-)$, $(K_{z_0},\mathcal{R}_{out})$
\end{itemize}

Let $z_0= 1$.
\begin{itemize}
\item If $\gamma\subset \Sigma_{+}$, then either $\gamma \subset L_+ \cup K_{z_0} \cup \mathcal{R}_{in} \cup \mathcal{R}_{out}$ or $\gamma$ is of type $(L_+, B_{\epsilon})$, $(L_+, K_{z_0})$, $(L_+, \mathcal{R}_{in})$, $(\mathcal{R}_{out}, K_{z_0})$, $(\mathcal{R}_{out}, B_{\epsilon})$, $(\mathcal{R}_{out},\mathcal{R}_{in})$,$(K_{z_0}, B_{\epsilon})$, $(K_{z_0},\mathcal{R}_{in})$.
\item If $\gamma \subset \Sigma_{-}$, then either $\gamma \subset L_{-}$ or $\gamma$ is of type $(B_{\epsilon}, L_-)$.
\end{itemize}
\end{prop} 
\begin{proof}
This proposition is a combination of lemmas \ref{characteristicDecreasing}, \ref{characteristicDecreasingPast}, \ref{characteristicIncreasing}, \ref{characteristicIncreasingPast}, \ref{characteristicGammaPm}, \ref{characteristicInfIncreasing}, \ref{characteristicInfIncreasingPast} and of the following observations about the hamiltonian vector field (which rule out some cases which are a priori compatible with the lemmas):
\begin{itemize}
\item If $z_0=-1$, $L_+$ is a source, $L_-$ is a sink, $\mathcal{R}_{in}$ is a source, $\mathcal{R}_{out}$ is a sink.
\item If $z_0=1$, $L_+$ is a source, $L_-$ is a sink, $\mathcal{R}_{in}$ is a sink, $\mathcal{R}_{out}$ is a source.
\item Bicharacteristic curves in $\Sigma_-$ only cross the horizon $\left\{ r= r_+\right\}$ towards the exterior of the black hole (see lemma \ref{convexityHorizon})
\item Bicharacteristic curves in $\Sigma_+$ only cross the horizon $\left\{ r=r_+ \right\}$ towards the interior of the black hole  (see lemma \ref{convexityHorizon})
\end{itemize}
\end{proof}
\begin{figure}[h]
\centering
\begin{tikzpicture}
\draw (-3, 2) -- (-3,-2) node[yshift = -0.2cm] (B) {$B_{\epsilon}$};
\draw (0,0) node[draw] (K) {$K_{1}$};
\draw (2,0) node[draw] (Rin) {$\mathcal{R}_{\text{in}}$};
\draw (2,-1) node[draw] (Rout) {$\mathcal{R}_{\text{out}}$};
\draw (-1, 2) node[draw] (Lp) {$L_+$};
\draw (-1, -2) node[draw] (Lm) {$L_-$};
\draw[decoration={markings, mark=at position 0.5 with {\arrow{>}}}, postaction={decorate}] (Lp)--(K);
\draw[decoration={markings, mark=at position 0.5 with {\arrow{>}}}, postaction={decorate}] (Lp)--(Rin);
\draw[decoration={markings, mark=at position 0.5 with {\arrow{>}}}, postaction={decorate}] (Lp)--(-3,0.5);
\draw[decoration={markings, mark=at position 0.5 with {\arrow{>}}}, postaction={decorate}] (K)--(Rin);
\draw[decoration={markings, mark=at position 0.5 with {\arrow{>}}}, postaction={decorate}] (K)--(-3,0);
\draw[decoration={markings, mark=at position 0.5 with {\arrow{>}}}, postaction={decorate}] (Rout)--(Rin);
\draw[decoration={markings, mark=at position 0.5 with {\arrow{>}}}, postaction={decorate}] (Rout)--(K);
\draw[decoration={markings, mark=at position 0.5 with {\arrow{>}}}, postaction={decorate}] (Rout)--(-3,-0.5);
\draw[decoration={markings, mark=at position 0.5 with {\arrow{>}}}, postaction={decorate}] (-3, -1.5)--(Lm);
\end{tikzpicture}\hspace{1cm}
\begin{tikzpicture}
\draw (-3, 2) -- (-3,-2) node[yshift = -0.2cm] (B) {$B_{\epsilon}$};
\draw (0,0) node[draw] (K) {$K_{-1}$};
\draw (2,0) node[draw] (Rin) {$\mathcal{R}_{\text{in}}$};
\draw (2,1) node[draw] (Rout) {$\mathcal{R}_{\text{out}}$};
\draw (-1, 2) node[draw] (Lp) {$L_+$};
\draw (-1, -2) node[draw] (Lm) {$L_-$};
\draw[decoration={markings, mark=at position 0.5 with {\arrow{>}}}, postaction={decorate}] (K)--(Lm);
\draw[decoration={markings, mark=at position 0.5 with {\arrow{>}}}, postaction={decorate}] (Rin)--(Lm);
\draw[decoration={markings, mark=at position 0.5 with {\arrow{>}}}, postaction={decorate}] (-3,-0.5)--(Lm);
\draw[decoration={markings, mark=at position 0.5 with {\arrow{>}}}, postaction={decorate}] (Rin)--(K);
\draw[decoration={markings, mark=at position 0.5 with {\arrow{>}}}, postaction={decorate}] (-3,0)--(K);
\draw[decoration={markings, mark=at position 0.5 with {\arrow{>}}}, postaction={decorate}] (Rin)--(Rout);
\draw[decoration={markings, mark=at position 0.5 with {\arrow{>}}}, postaction={decorate}] (K)--(Rout);
\draw[decoration={markings, mark=at position 0.5 with {\arrow{>}}}, postaction={decorate}] (-3,0.5)--(Rout);
\draw[decoration={markings, mark=at position 0.5 with {\arrow{>}}}, postaction={decorate}] (Lp)--(-3, 1.5);
\end{tikzpicture}
\caption{Structure of the semiclassical Hamiltonian flow (for $z_0=1$ on the left and $z_0 = -1$ on the right)}
\end{figure}
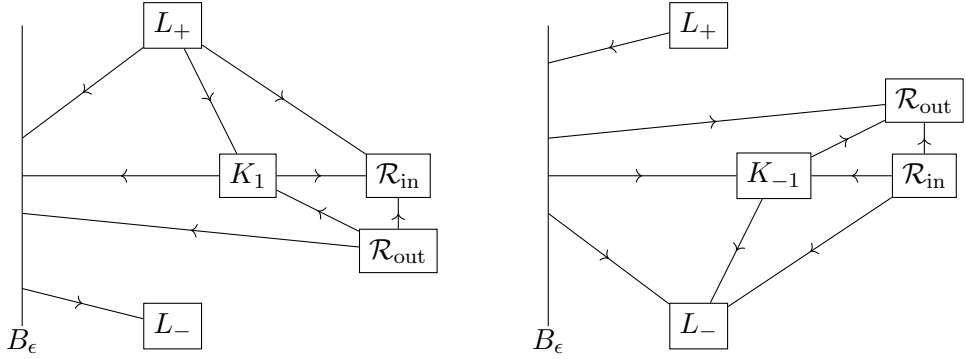
Finally, we state the following lemma for later use. It is a consequence of lemma 2.4 in \cite{dyatlov2016spectral} (the fact that normally hyperbolic trapping assumptions holds is proven in \cite{dyatlov2015asymptotics}). 
\begin{lemma}
\label{neighTrapping}
For all $\epsilon>0$ small enough, there exists a neighborhood $U$ of $K_{z_0}$ and $\phi_+$, $\phi_-$ smooth functions on $U$ such that:
\begin{itemize}
\item $\left\{\phi_+ = 0\right\} = \Gamma^+\cap U$
\item $\left\{\phi_- = 0\right\} = \Gamma^-\cap U$
\item There exists $\delta>0$ such that, $U_\delta:= \left\{ \lvert\phi_+\rvert<\delta, \lvert\phi_-\rvert<\delta, \lvert p\rvert<\delta\right\}$ is compactly contained in $U$.
\item $\left\{\phi_+, \phi_-\right\} >0$ on $U$
\item $H_p \phi_{\pm} = \mp c_{\pm} \phi_{\pm}$ with $c_\pm$ smooth positive bounded functions on $U$ with $c_\pm>\nu_{\min}-\epsilon$ where $\nu_{\text{min}}>0$ is defined in Proposition 3.7 of \cite{dyatlov2015asymptotics} (but we will not use its exact value).
\end{itemize}
\end{lemma}

\section{Fredholm property of \texorpdfstring{$\hat{T}_s(\sigma)$}{hat\{T\}\_ s( sigma)}}\label{secFred}
We have seen previously (see Corollary \ref{FourierTrans1} and Corollary \ref{FourierTrans2}) that the Cauchy problem can be reduced to a forcing problem and Fourier transformed into the problem:
\begin{align*}
\hat{T}_s(\sigma)\hat{v}(\sigma) = \hat{f}(\sigma)
\end{align*} 
where $\hat{f}$ is precisely characterized. To recover properties of $\hat{v}$, we want to write it as $\hat{v}(\sigma) = \hat{T}_s(\sigma)^{-1}\hat{f}(\sigma)$. In this section, we prove that the operator $\hat{T}_s(\sigma)$ is Fredholm between appropriate spaces. This together with a mode stability result (see \cite{whiting1989mode}, \cite{andersson2017mode} and \cite{andersson2022mode}) will provide the invertibility of $\hat{T}_s(\sigma)$ between well chosen spaces.

\subsection{Estimate near the horizon} \label{subSecEstHorizon}
To define adjoint operator, we use the volume form $\dd vol$ and the metric $\mathbb{m}$ (defined in Subsection \ref{SobolevSpaces}).
We now compute the subprincipal symbol, that is to say the principal symbol of $\frac{1}{2i}\left(\hat{T}_s(\sigma)-\hat{T}_s(\sigma)^*\right)$ near the horizon:
\begin{align*}
p_{sub} = \mathfrak{s}\left(\frac{1}{2i}\left(\hat{T}_s(\sigma)-\hat{T}_s(\sigma)^*\right)\right) = \left(-2a\Im(\sigma)+2\frac{a}{r}\right)\zeta + 2\left((r-M)s+\frac{\Delta_r}{r} - (a^2+r^2)\Im(\sigma)\right)\xi
\end{align*}

\begin{prop}
\label{radialPointEstInf1}
Let $k>m\geq \frac{1}{2} - \left(-s +  \frac{(a^2+r_+^2)}{r_+-M}\Im(\sigma)\right)$
If $A,B,G \in \Psi_b^{0,0}$ with compactly supported Schwarz kernels are such that $A$ and $G$ are elliptic on $L_{\pm}$ and every forward (or backward) classical bicharacteristic from $WF(B)$ tends to $L_{\pm}$ with closure in the elliptic set of $G$, 
\[Au \in \overline{H}^{m,l}_{(b)} \Rightarrow \left\|Bu\right\|_{\overline{H}^{k,l}_{(b)}}\leq C\left(\left\| G\hat{T}_s(\sigma)u \right\|_{\overline{H}^{k-1,l}_{(b)}} + \left\|u\right\|_{\overline{H}^{-N,l}_{(b)}}\right)\]
\end{prop}
\begin{prop}
\label{radialPointEstInf2}
Let $k < \frac{1}{2} + \left(-s + \frac{(a^2+r_+^2)}{r_+-M}\Im(\sigma)\right)$. Assume $A \in \Psi_b^{0,l}$ with compactly supported Schwarz kernel is elliptic on a neighborhood of $L_{\pm}$, $B\in \Psi_b^{0,0}$ and $G \in \Psi_b^{0,l}$ with compactly supported Schwarz kernels and $G$ is elliptic on $WF(A)\cup WF(B)$. Assume also that every forward (or backward) bicharacteristic of $P$ from a point of $WF(A)$ reach $Ell(B)$ while remaining in $Ell(G)$. Then
\[
\left\|Au\right\|_{\dot{H}^{k,l}_{(b)}}\leq C\left(\left\|G\hat{T}_s(\sigma)^*u\right\|_{\dot{H}^{k-1,l}_{(b)}} + \left\|Bu\right\|_{\dot{H}^{k,l}_{(b)}}+\left\|u\right\|_{\dot{H}^{-N,l}_{(b)}}\right)
\]
\end{prop}
\begin{remark}
These estimates are microlocalized away from $\partial \overline{X}$. As a consequence, the index $l$ does not play any role here and could be replaced in the right hand sides by $l-N$ for any $N \in \N$.
\end{remark}
\begin{proof}
These two propositions are a mild modification of propositions 2.3 and 2.4 in \cite{vasy2013microlocal}. There are only two differences: 
\begin{enumerate}
\item The proposition in \cite{vasy2013microlocal} is not stated for operators on sections of a complex vector bundles. However, in the case of a complex vector bundle of rank 1 (as it is the case here), the principal symbols of operators are scalar and the same proof can be applied (see remark 2.1 in \cite{vasy2013microlocal}). 
\item The subprincipal symbol $p_{sub}$ of the operator $\hat{T}_s(\sigma)$ has not exactly the form required \cite{vasy2013microlocal}. However, what is really needed in the proof is the positivity of $k- \frac{1}{2} + \frac{p_{sub}\rho}{\beta_0}$ in case of proposition 2.3 and the negativity of $k-\frac{1}{2} - \frac{p_{sub}\rho}{\beta_0}$ in case of proposition 2.4 where by definition $H\rho_{|L_{\pm}} = \pm \beta_0\rho$ (as mentioned previously, we have taken into account the different sign convention for the principal symbol in \cite{vasy2013microlocal}). As a consequence, the proof works if we adapt the hypothesis on $k$ as we did.
\end{enumerate}
\end{proof}

We also get the corresponding semiclassical estimates given by propositions 2.10 and 2.11 of \cite{vasy2013microlocal} which we state for the operator the operator $\hat{T}_{s,h}(z) = h^{2}\hat{T}_s(h^{-1}z)$ where $z = \frac{\sigma}{\lvert\sigma\rvert}$, $h= \lvert \sigma\rvert$. The estimates are uniform when $\sigma$ lies in a the region $\left\{ 0\leq\Im(\sigma)\leq \eta, \lvert\Re(\sigma)\rvert\geq \eta \right\}$ with $\eta>0$.
\begin{prop}
\label{radialPointMicrolocal1}
Let $k>m\geq \frac{1}{2} - \left(-s +  \frac{(a^2+r_+^2)}{r_+-M}\Im(h^{-1}z)\right)$.
If $A,B,G \in \Psi_{b,h}^{0,0}$ with compactly supported Schwartz kernels are such that $A$ and $G$ are elliptic on $L_{\pm}$ and every forward (or backward) bicharacteristic from $WF_h(B)$ tends to $L_{\pm}$ with closure in the elliptic set of $G$, 
\[Au \in \overline{H}^{m,l}_{(b),h} \Rightarrow \left\|Bu\right\|_{\overline{H}^{k,l}_{(b),h}}\leq C\left(h^{-1}\left\| G\hat{T}_{s,h}(z)u \right\|_{\overline{H}^{k-1,l}_{(b),h}} + h\left\|u\right\|_{\overline{H}^{-N,l}_{(b),h}}\right)\]
\end{prop}

\begin{prop}
\label{radialPointMicrolocal2}
Let $k < \frac{1}{2} + \left(-s + \frac{(a^2+r_+^2)}{r_+-M}\Im(\sigma)\right)$. Assume $A \in \Psi_{b,h}^{0,l}$ with compactly supported Schwartz kernel is elliptic on a neighborhood of $L_{\pm}$, $B\in \Psi_b^{0,0}$ and $G \in \Psi_b^{0,l}$ with compactly supported Schwartz kernels and $G$ is elliptic on $WF(A)\cup WF(B)$. Assume also that every forward (or backward) bicharacteristic of $P$ from a point of $WF(A)$ reach $Ell(B)$ while remaining in $Ell(G)$. Then
\[
\left\|Au\right\|_{\dot{H}^{k,l}_{(b),h}}\leq C\left(h^{-1}\left\|G\hat{T}_{s,h}(z)^*u\right\|_{\dot{H}^{k-1,l}_{(b),h}} + \left\|Bu\right\|_{\dot{H}^{k,l}_{(b),h}}+h\left\|u\right\|_{\dot{H}^{-N,l}_{(b),h}}\right)
\]
\end{prop}

We now look at the operator $\hat{T}_s(\sigma)$ on $\left\{r< r_+-\frac{\epsilon}{3}\right\}$. To estimate the solution in this region, we use classical hyperbolic estimates (following an idea presented in \cite[Section 3.2]{zworski2015resonances}). We need the two following propositions:
\begin{prop}
\label{hyperbolicEstimate1}
Let $u \in \dot{H}^{m,l}_{(b)}$ with $\supp(u)\subset \left\{r<r_+-\frac{\epsilon}{3}\right\}$. We denote by $\tilde{H}^{s}$ the space of distribution on $(-\infty, r_+-\frac{\epsilon}{3})$ which are supported in $(r_+-\epsilon, r_+-\frac{\epsilon}{3})$ and extendible at $r_+-\frac{\epsilon}{3}$ endowed with the corresponding norm. In particular we have $v:= u_{|_{r<r_+-\frac{\epsilon}{3}}} \in \tilde{H}^{s}$. With a slight abuse of notation, we write $\left\|u\right\|_{\tilde{H}^{s-1}}$ for $\left\|v\right\|_{\tilde{H}^{s-1}}$.
For all $s\in \R$, we have \[\left\|u\right\|_{\tilde{H}^{s}}\leq C\left\|\hat{T}_s(\sigma)^*u\right\|_{\tilde{H}^{s-1}}\] (with the convention that some terms may be infinite).
\end{prop}
The next proposition is similar but with the extendible and supported ends inverted.
\begin{prop}
\label{hyperbolicEstimate2}
Let $u \in \overline{H}^{m,l}_{(b)}$ with $\supp(u)\subset \left\{ r\leq r_+ - \frac{\epsilon}{3}\right\}$. Then for all $s\in \R$, we have \[\left\|u\right\|_{\overline{H}^{s,l}_{(b)}}\leq C\left\|\hat{T}_s(\sigma)u\right\|_{\overline{H}^{s-1,l}_{(b)}}\] (with the convention that some terms may be infinite).
\end{prop}

These two propositions are consequences of the standard hyperbolic theory for second order partial differential operators. Indeed, on $\left\{r_+-2\epsilon <r< r_+-\frac{\epsilon}{3}\right\}$, the operator $\hat{T}_s(\sigma)$ has (classical) principal symbol:
\begin{align*}
\mathfrak{s}(\hat{T}_s(\sigma)) =& \tilde{G}(\xi)
\end{align*}
where $\tilde{G}$ is a Lorentzian metric. Moreover $\tilde{G}(\dd r) = \Delta_r$ which is uniformly negative on $\left\{r_+-2\epsilon <r< r_+-\frac{\epsilon}{3}\right\}$. Therefore, $\hat{T}_s(\sigma)$ is strictly hyperbolic with respect to the level sets of $r$ in this region. We also need a semiclassical version of the hyperbolic estimate (see Proposition \ref{finalScHyperboEstimate}). Since it is less standard than classical hyperbolic estimate, we provide a detailed proof (in a more general setting) in the appendix, section \ref{semiclassicalHyperboEstSection}. This proof can easily be adapted to obtain a proof of Propositions \ref{hyperbolicEstimate1} and \ref{hyperbolicEstimate2}.

We can now combine the previous results to get a complete microlocal estimate near the horizon (and more precisely on the full region $\left\{r\leq 2M\right\}$)
\begin{prop}
\label{combinedNearHorizon}
Let $N$ and $N'$ be integers which will be taken large. Let $k'>\frac{1}{2} + s -\frac{a^2+r_+^2}{r_+-M}\Im(\sigma)$.
Let $u \in \overline{H}_{(b)}^{k', \infty}$ such that $\supp(u)\subset \left\{r\leq 2M + D\right\}$ (for some $D>0$) and $v \in \dot{H}_{(b)}^{-N, \infty}$ such that $\supp(v)\subset \left\{r\leq 2M+D\right\}$. Let $k$ be such that $k>k'$. We have
\begin{align*}
\left\| u \right\|_{\overline{H}_{(b)}^{k,l}}&\leq C\left(\left\| \hat{T}_s(\sigma)u\right\|_{\overline{H}_{(b)}^{k-1, -N'}} + \left\|u\right\|_{\overline{H}_{(b)}^{k', -N'}}\right)\\
\left\| v \right\|_{\dot{H}_{(b)}^{1-k, l}}&\leq C\left(\left\| \hat{T}_s(\sigma)^*v\right\|_{\dot{H}_{(b)}^{-k, -N'}} + \left\|v\right\|_{\dot{H}_{(b)}^{-N,-N'}}\right)
\end{align*}
\end{prop}

We first need a separation lemma (following from the analysis of the dynamical structure)
\begin{lemma}
\label{separationCharacteristics}
There exists disjoint open subsets $U_+$ and $U_-$ of $\overline{T}^*\mathcal{M}_{\epsilon}$ such that $\Sigma \cap \left\{r\geq r_+-\frac{\epsilon}{2}\right\} \subset U_+\cup U_-$, $\Lambda_+ \subset U_+$ and $\Lambda_-\subset U_-$ and 
\begin{itemize}
\item For every $x \in U_+\cap \Sigma$, the bicharacteristics through $x$ tends to $L_+$ at $-\infty$ and  is included in $\left\{r<r_+-\frac{\epsilon}{2}\right\}$ in the future.
\item For every $x \in U_-\cap \Sigma$, the bicharacteristics through $x$ tends to $L_-$ at $+\infty$ and  is included in $\left\{r<r_+-\frac{\epsilon}{2}\right\}$ in the past.
\end{itemize}
\end{lemma}
\begin{proof}
In view of the dynamical analysis obtained in Proposition \ref{recapClassicalFlow}, we can take $U_1 = \left\{\tilde{\xi}>0\right\}$ and $U_2 = \left\{\tilde{\xi}<0 \right\}$.
\end{proof}
\begin{proof}[Proof of proposition \ref{combinedNearHorizon}]
These estimates are proved by combining the previous estimates in this section. 
We first prove the first one:
Let $\Sigma_{\pm} := \Sigma \cap \left\{r\geq r_+-\frac{\epsilon}{2}\right\}\cap U_{\pm}$ (where $U_{\pm}$ are defined in lemma \ref{separationCharacteristics}) are compact subsets of $\overline{T}^*\mathcal{M}_{\epsilon}$. We cover them by a finite number of open sets $(U_i)_{i=1}^N$ included in $U_{\pm}$. Let $B_i\in \Psi^{0,0}_b$ with compactly supported Schwartz kernels such that $WF(B_i)\subset U_i$ and $WF(I-\sum_{i=1}^N B_i)\cap \Sigma \subset \left\{r<r_+-\frac{\epsilon}{2}\right\}$. Let $\chi \in C^{\infty}(\R)$ be such that $\chi = 1$ on $r<r_+-\frac{\epsilon}{2}$ and $\chi =0$ on $\left\{r>r_+-\frac{\epsilon}{3}\right\}$. We have:
\begin{align*}
\left\|u\right\|_{\overline{H}_{(b)}^{k,l} }&\leq \sum_{i=1}^N\left\|B_i u\right\|_{\overline{H}_{(b)}^{k,l}} + C\left\|\chi u\right\|_{\overline{H}_{(b)}^{k,l}}+ \left\| (1-\chi)(I-\sum_{i=1}^N B_i)u\right\|_{\overline{H}_{(b)}^{k,l}}\\
\end{align*}
We have $\left\|B_i u\right\|_{\overline{H}_{(b)}^{k,l}}\leq C\left\| \hat{T}_s(\sigma) u\right\|_{\overline{H}_{(b)}^{k-1,-N'}} + \left\|u\right\|_{\overline{H}_{(b)}^{-N,-N'}}$ using proposition \ref{radialPointEstInf1} for each $B_i$ (here we use the regularity assumption on $u$). Moreover $\left\| (1-\chi)(I-\sum_{i=1}^N B_i)u\right\|_{\overline{H}_{(b)}^{k,l}} \leq \left\|\hat{T}_s(\sigma)u\right\|_{H^{k-1,l}_{(b)}}+\left\|u\right\|_{\overline{H}_{(b)}^{-N,-N'}}$ using proposition \ref{ellipticEstimate} (elliptic estimate) since $WF((1-\chi)(I-\sum_{i=1}^N B_i))\cap \Sigma = \emptyset$.

Eventually, we have to bound the term $\left\|\chi u\right\|_{\overline{H}_{(b)}^{k,l}}$. Using proposition \ref{hyperbolicEstimate2}, we get $\left\|\chi u\right\|_{\overline{H}_{(b)}^{k,l}}\leq \left\|\hat{T}_s(\sigma) \chi u \right\|_{\overline{H}_{(b)}^{k-1,l}}$. Note that
\begin{align*}
\left\|\hat{T}_s(\sigma) \chi u \right\|_{\overline{H}_{(b)}^{k,l}} &\leq \left\|\hat{T}_s(\sigma) u \right\|_{\overline{H}_{(b)}^{k-1,-N'}} + \left\| [\hat{T}, \chi] u\right\|_{\overline{H}_{(b)}^{k-1,-N'}}\\
&\leq \left\|\hat{T}_s(\sigma) u \right\|_{\overline{H}_{(b)}^{k-1,-N'}} + C\left\| \tilde{\chi} u\right\|_{\overline{H}_{(b)}^{k,-N'}}
\end{align*}
where $\tilde{\chi}\in C^{\infty}_c$ and $\supp\tilde{\chi}\subset \left\{ r_+-\frac{\epsilon}{2}\leq r \leq r_+-\frac{\epsilon}{3}\right\}$. Therefore
\begin{align*}
\left\| \tilde{\chi} u\right\|_{\overline{H}_{(b)}^{k,-N'}} &\leq \sum_{i=1}^N \left\|B_i\tilde{\chi}u\right\|_{\overline{H}_{(b)}^{k,-N'}} + \left\| \left(I-\sum_{i=1}^NB_i\right)\tilde{\chi} u\right\|_{H_{(b)}^{k,-N'}} \\
&\leq C\left\| \hat{T}_s(\sigma) u\right\|_{\overline{H}_{(b)}^{k-1,-N'}} + \left\|u\right\|_{H_{(b)}^{-N,-N'}}\\
\end{align*}
where the last line is obtained using proposition \ref{radialPointEstInf1} for each $B_i\tilde{\chi}$ and proposition \ref{ellipticEstimate} for the last term.

We now prove the second estimate:
We define $B_i$ and $\chi$ in the same way but we arrange that $B_0$ and $B_1$ are such that $WF(I-B_0)\cap \Lambda_+ = \emptyset$ and $WF(I-B_1)\cap \Lambda_- = \emptyset$. 
We have
\begin{align*}
\left\|u\right\|_{\dot{H}_{(b)}^{1-k,l} }&\leq \sum_{i=1}^N\left\|B_i u\right\|_{\dot{H}_{(b)}^{1-k,l}} + C\left\|\chi u\right\|_{\dot{H}_{(b)}^{1-k,-N'}}+ \left\| (1-\chi)(I-\sum_{i=1}^N B_i)u\right\|_{\dot{H}_{(b)}^{1-k,-N'}}\\
\end{align*}
Using proposition \ref{radialPointEstInf2} for $B_0$ and $B_1$ and usual propagation of singularity for the other $B_i$, we get $\sum_{i=1}^N\left\|B_i u\right\|_{\dot{H}_{(b)}^{1-k,l}}\leq \left\|\hat{T}_s(\sigma)^* u\right\|_{H_{(b)}^{-k,-N'}} + \left\|\chi u\right\|_{H_{(b)}^{1-k,-N'}}$. As previously, using proposition \ref{ellipticEstimate}, we can write $\left\| (1-\chi)(I-\sum_{i=1}^N B_i)u\right\|_{\dot{H}_{(b)}^{1-k,-N'}} \leq \left\|\hat{T}_s(\sigma)^*u\right\|_{H_{(b)}^{-k, -N'}} + \left\|u\right\|_{-N,-N'}$. 
Finally, we have
\begin{align*}
\left\|\chi u \right\|_{H_{(b)}^{1-k,-N'}}&\leq C\left\| u\right\|_{\tilde{H}^{1-k}}\\
&\leq C\left\| \hat{T}_s(\sigma)^* u \right\|_{H_{(b)}^{-k,-N'}}
\end{align*}
where the last inequality comes from Proposition \ref{hyperbolicEstimate1}.
\end{proof}

\subsection{Estimate near \texorpdfstring{$x=0$}{x=0}}\label{subSecEstSpatInf}
Now we prove an estimate near the end $x=0$. As before, we also need a semiclassical version. Contrarily to the previous estimates which remained valid up to $\sigma = 0$, it is not automatic in this case and we need to prove a uniform version of the estimate down to $\sigma=0$. Unless otherwise indicated, computation in this section are localized on $\left\{0\leq x \leq \frac{1}{6M}\right\}$.
To be coherent with the notation in \cite{vasy2020limiting}, we introduce the operator:
\begin{align*}
P(\sigma) = e^{i\sigma x - 2iM\sigma \ln(x)}\frac{1}{\Delta_r}\hat{T}_s(\sigma)e^{-i\sigma x + 2iM\sigma \ln(x)}.
\end{align*}
Other choices are possible to fit in the framework of \cite{vasy2020limiting}, the choice made here imposes that the coefficient of $(xD_x)^2$ is equal to $1$.
Therefore, the operator after conjugation (defined in \cite[Section~3]{vasy2020limiting}) 
\begin{align}\label{conjugation}\hat{P}(\sigma) := e^{-i\sigma x} P(\sigma) e^{i\sigma x} = x^{-2iM\sigma}\frac{1}{\Delta_r}\hat{T}_s(\sigma)x^{2iM\sigma}
\end{align}
We will use this equality to translate estimates about $P(\sigma)$ into estimates about $\hat{T}_s(\sigma)$.

We have the following decomposition for $P(\sigma)$:
\begin{align*}
P(\sigma) = P(0) + \sigma Q -\sigma^2
\end{align*}
where $P(0) \in \text{Diff}^2_{sc}$ and $Q\in x\text{Diff}^1_{sc}$ and explicitly (using the trivialization $\mathcal{T}_m$)
\begin{align*}
P(0) = & (x^2D_x)^2 + \frac{2i(r-M)(1-s)}{\Delta_r}x^2D_x + \frac{1}{\Delta_r\sin\theta}D_\theta \sin\theta D_\theta + \frac{r^2+a^2\cos^2\theta - 2Mr}{\Delta_r^2\sin^2\theta}D_\phi^2\\
& + \frac{2s(\Delta_r\cos\theta + ia(M-r)\sin^2\theta)}{\Delta_r^2\sin^2\theta}D_\phi + \frac{s^2\text{cotan}^2\theta+s}{\Delta_r} \\
Q = & \frac{4Mar}{\Delta_r^2}D_\phi - \frac{4M(2Mr-a^2)}{r\Delta_r}x^2D_x - \frac{1}{r^2\Delta_r^2}q_0\\
q_0 =&  \left( 4 M \sigma + 2 i s\right)x^{-5} + x^{-4}C^{\infty}\left(\left[0, \frac{1}{6M}\right]_x\times \mathbb{S}^2\right)
\end{align*}

\begin{prop}\label{midFreqInfEstimate}
Let $\chi_1$ be a smooth cutoff compactly supported in $\left\{ x<\frac{1}{6M}\right\}$ with $\chi_1 = 1$ on $\left\{ x<\frac{1}{12M}\right\}$.

Assume that $l+\frac{1}{2}-2M\Im(\sigma)<0$, $\tilde{r}+l+\frac{1}{2}+ 2s+2M\Im(\sigma)>0$ and $\Im(\sigma)\geq 0$ and $\sigma\neq 0$. Then if $u\in H_{(b)}^{\tilde{r}',l}$ for some $\tilde{r}'$ such that $\tilde{r}'+l+\frac{1}{2}+2M\Im(\sigma)>0$, we have:
\begin{align*}
\left\|\chi_1 u\right\|_{H_{(b)}^{\tilde{r}, l}}&\leq C\left(\left\|\hat{P}(\sigma)u\right\|_{H_{(b)}^{\tilde{r},l+1}} + \left\|u\right\|_{H_{(b)}^{-N,l-1}}\right) \\
\left\|\chi_1 u\right\|_{H_{(b)}^{\tilde{r}, l}}&\leq C\left(\left\|\hat{P}(\sigma)u\right\|_{H_{(b)}^{\tilde{r}-1,l+2}} + \left\|u\right\|_{H_{(b)}^{-N,l-1}}\right) \\
\end{align*}
Moreover, if $v \in H_{(b)}^{\tilde{r}'',l''}$ for some $\tilde{r}'' \in \R$ (with no further condition on $\tilde{r}''$) and $l''+\frac{1}{2}+2M\Im(\sigma)>0$, we have:
\begin{align*}
\left\| \chi_1 v\right\|_{H_{(b)}^{-\tilde{r},-(l+1)}}\leq C\left(\left\|\hat{P}(\sigma)^*v\right\|_{H_{(b)}^{-\tilde{r},-l}} + \left\|v\right\|_{H_{(b)}^{-N,-l-2}}\right)\\
\end{align*}

The constants are uniforms with respect to $\sigma$ in a compact subset of $\left\{ \Im(\sigma)\geq 0\right\}\setminus\left\{0\right\}$
\end{prop}

\begin{proof}
First note that it is enough to prove that 
\begin{align}\label{reducedIneq}
\left\|\chi_1 u\right\|_{H_{(b)}^{\tilde{r}, l}}&\leq C\left(\left\|\hat{P}(\sigma)\chi_1 u\right\|_{H_{(b)}^{\tilde{r},l+1}} + \left\|\chi_1 u\right\|_{H_{(b)}^{-N,l-1}}\right)\end{align}
and similarly for the second inequality.
Indeed, we have that $[\hat{P}(\sigma), \chi_1]\in \text{Diff}^1_b$ is supported away from the boundary and on the elliptic set of $\hat{P}(\sigma)$. Therefore we have: 
\begin{align*}\left\|[\hat{P}(\sigma), \chi_1]u\right\|_{H_{(b)}^{\tilde{r},l+1}}&\leq \left\| u\right\|_{H_{(b)}^{\tilde{r}+1, -\infty}}\\
&\leq C\left\| \hat{P}(\sigma) u \right\|_{H_{(b)}^{\tilde{r}-1, -\infty}}\\
&\leq C\left\| \hat{P}(\sigma) u \right\|_{H_{(b)}^{\tilde{r}, l+1}}
\end{align*}
We can do the same for the second inequality. Therefore, we can consider $P(\sigma)$ as an operator acting on $\left[0,\frac{1}{6M}\right)_x\times \mathcal{B}_s$.

The estimate \eqref{reducedIneq}, is exactly the kind of estimate obtained in the proof of theorem 1.1 in \cite{vasy2020limiting}.

Note that the operator $P(\sigma)$ is not exactly of the form described in section 2 of \cite{vasy2020limiting} since it acts on a non trivial complex line bundle. However, as mentioned in remark 1.3 of \cite{vasy2020limiting}, the proof of the estimate is completely parallel in this case (and since in our case the bundle is of dimension 1, it requires even less adaptation).
The exact hypotheses that we use and which are sufficient to run the proof are the following (we do not intend to be as general as the adaptation of \cite{vasy2020limiting} to complex line bundles could be, in particular we do not consider conormal operator's coefficients and we restrict to a simpler form for $P(0)$ which is enough to treat our case):
Let $X = [0, \frac{1}{6M})_x\times \partial X$ be a smooth manifold with boundary of dimension $n$ with smooth boundary defining function $x$. Let $\omega_{\partial X}$ be a smooth volume form on the boundary, we denote by $\dd vol = x^{-n-1}|\dd x|\omega$. Let $p:E\rightarrow \partial X$ be a complex line bundle over $\partial X$ and $m$ be a smooth metric on $E$. We denote by $\tilde{E} = \pi_2^* E$ (the semitrivial bundle over $X$ associated to $E$) and $\tilde{m} = \pi_2^* m$ the associated metric. Let $P(\sigma) = P(0)+\sigma Q -\sigma^2$ be such that
\begin{enumerate}[label = (H\arabic*)]
\item \label{item:hyp1}
\begin{align*}
P(0) = (x^2D_x)^2 + x^2 A  + (i(n-1)+a_0)x(x^2D_x) 
\end{align*}
where $A \in C^{\infty}([0,M)_x,\text{Diff}^2(E))$ and $ a_0 \in C^{\infty}(X)$. We also require that the principal symbol of $x^{-2}P(0)$ is elliptic in the $b$ sense. Note that we imposed that the coefficients $a'$, $a_{0,j}$ are zero with respect to the local coordinate expression (3.4) in \cite{vasy2020limiting}, which is true in our case. 
\item
\begin{align*}
Q = b_0 x(x^2D_x) + x^2R + b'x
\end{align*}
where $R \in C^{\infty}([0,M)_x, \text{Diff}^1(E))$ and $b_0, b' \in C^{\infty}(X)$. Coefficients of $Q$ can depend smoothly on $\sigma$.

\item \label{item:hyp3} 
The operators $x^2A$ and $x^2R$ admit the following decomposition:
\begin{align*}
x^2A &= \sum_{j = 1}^N (A_j^*A_j + A^*_jA_j'+A^\dagger_j A_j) + A''\\
x^2R &= \sum_{j=1}^N (A^*_j R_j + R_j^\dagger A_j) + R''
\end{align*}
where $A_j \in xC^{\infty}([0,M)_x,\text{Diff}^1(E)) \subset x\text{Diff}^1(\tilde{E})$, $A_j',A^\dagger_j,R_j,R_j^\dagger \in xC^{\infty}(X)$ and $A'', R'' \in x^2 C^{\infty}(X)$.
This hypothesis is needed to enforce the conclusion of lemma 3.3 in \cite{vasy2020limiting} (which is automatically true when we have the precise local form given in (3.4) and (3.5), we avoid giving local expression here since it would also depend on the choice of a local trivialization of $E$ in addition to the choice of coordinates on $\partial X$  but it could be done). 
\end{enumerate}

Under the previous hypotheses, we have:
\begin{align*}
\hat{P}(\sigma) &:= e^{-i\frac{\sigma}{x}}P(\sigma)e^{i\frac{\sigma}{x}}\\
&= P(0) + \sigma \hat{Q} - 2\sigma(x^2D_x + i\frac{n-1}{2}x + x\tilde{\alpha}_+(\sigma))
\end{align*}
with $\hat{Q} = Q-xb'$ and $\tilde{\alpha}_+(\sigma) = \frac{a_0 - b'+b_0 \sigma}{2}$. Therefore the normal operator is formally identical to the one in \cite{vasy2020limiting} (the only difference is that we consider it as a differential operator on $\tilde{E}$). The threshold value for $\tilde{r}$ is then $-\frac{1}{2}+\Im(\alpha_+(\sigma))$ where $\alpha_{+}(\sigma) = \lim\limits_{x\to 0}\tilde{\alpha}_+(\sigma)$. Similarly, the threshold value for $\tilde{r}+l$ is $-\frac{1}{2} + \Im(\alpha_-(\sigma))$ where $\alpha_-(\sigma) := \lim\limits_{x\to 0}\tilde{\alpha}_-(\sigma)$ and $\tilde{\alpha}_-(\sigma) := \frac{a_0+b'+b_0\sigma}{2}$.

Let $(y_i)$ be local coordinates on $\partial X$. Let $\tau$ be the variable associated to $x$ in ${}^{sc}T^*X$ and $\mu_i$ be the variable associated to $y_i$. Hypotheses \ref{item:hyp1} and \ref{item:hyp3} imply that the principal symbol of $P(0)$ in the scattering decay sense is $\tau^2 + q(y)(\mu)$ where $q(y)$ is a positive definite quadratic form (positivity comes from \ref{hyp3} and definiteness from the fact that $x^{-2}P(0)$ is elliptic in the $b$ sense). This is the correct form to run commutator estimates with the same method as in \cite{vasy2020limiting} (and we have the conclusion of Lemma 3.3 in \cite{vasy2020limiting} thanks to hypothesis \ref{item:hyp3} to deal with the case of non real $\sigma$).

We just have to check that the hypotheses are satisfied our case.
\begin{enumerate}[label = (H\arabic*)]
\item The hypothesis is true with $a_0 = \frac{2i(Mr(1+s)-a^2-r^2s)}{\Delta_r} = -2is + O(x)$, 
\begin{align} \label{defA}
A=& \frac{r^2}{\Delta_r\sin\theta}D_\theta \sin\theta D_\theta + \frac{r^2\left(r^2+a^2\cos^2\theta - 2Mr\right)}{\Delta_r^2\sin^2\theta}D_\phi^2 + \frac{2sr^2(\Delta_r\cos\theta + ia(M-r)\sin^2\theta)}{\Delta_r^2\sin^2\theta}D_\phi \notag\\ 
&+ \frac{r^2\left(s^2\text{cotan}^2\theta+s\right)}{\Delta_r}.
\end{align}
\item The hypothesis is true with $b' = -\frac{1}{r\Delta_r^2}q_0 = -4M\sigma-2is + O(x)$, $b_0 = -\frac{4M(2Mr-a^2)}{\Delta_r}$ and $R = \frac{4Mar^3}{\Delta_r^2}D_\phi$ (which is smooth as a differential operator on $\tilde{E}$).
\item We define the following smooth operators on $\mathcal{B}_3(s)$ by their coordinates in local trivialization $A_m$:
\begin{align*}
(\tilde{Z}_1)_m &= -\sin\phi \partial_{\theta} + \cos\phi \left(-\frac{is}{\sin\theta} - \frac{\cos\theta}{\sin\theta}\partial_\phi\right) \\
(\tilde{Z}_2)_m &= -\cos\phi \partial_{\theta} -\sin\phi \left(-\frac{is}{\sin\theta} - \frac{\cos\theta}{\sin\theta}\partial_\phi\right) \\
(\tilde{Z}_3)_m &= \sqrt{1-\frac{a^2}{\Delta_r}}\partial_\phi. \\
\end{align*}
We have 
\begin{align*}
x^2A =& A_1^*A_1 + A_2^*A_2+ A_3^*A_3 + A_3^\dagger A_3 + A''\\
x^2R =& R_3^{\dagger}A_3 + R''
\end{align*}
with
\begin{align*}
A_i &= \frac{1}{\sqrt{\Delta_r}}\tilde{Z}_i \\
A_3^\dagger &= \frac{a(M-r)}{\Delta_r^2\sqrt{1-\frac{a^2}{\Delta_r}}}\\
A_3'' &= \frac{s(1-s)}{\Delta_r}\\
R_3 &= \frac{4Mar}{\Delta_r^2\sqrt{1-\frac{a^2}{\Delta_r}}}. \\
\end{align*}
\end{enumerate}
\end{proof}

If we translate the previous estimate into an estimate on $\frac{1}{\Delta_r}\hat{T}_s(\sigma)$ using \eqref{conjugation}, we get:
\begin{prop}
\label{estimateNearInfMedFreq}
Let $\chi_1$ be a smooth cutoff compactly supported in $\left\{ x<\frac{1}{C}\right\}$ with $\chi_1 = 1$ on $\left\{ x<\frac{1}{2C}\right\}$.

Assume that $l+\frac{1}{2}<0$, $\tilde{r}+l+\frac{1}{2}+2s + 4M\Im(\sigma)>0$ and $\Im(\sigma)\geq 0$ and $\sigma\neq 0$. Then if $u\in H_{(b)}^{\tilde{r}',l}$ for some $\tilde{r}'$ such that $\tilde{r}'+l+\frac{1}{2}+2s+4M\Im(\sigma)>0$, we have:
\begin{align*}
\left\|\chi_1 u\right\|_{H_{(b)}^{\tilde{r}, l}}&\leq C\left(\left\|\hat{T}_s(\sigma)u\right\|_{H_{(b)}^{\tilde{r},l-1}} + \left\|u\right\|_{H_{(b)}^{-N,l-1}}\right)\\
\left\|\chi_1 u\right\|_{H_{(b)}^{\tilde{r}, l}}&\leq C\left(\left\|\hat{T}_s(\sigma)u\right\|_{H_{(b)}^{\tilde{r}-1,l}} + \left\|u\right\|_{H_{(b)}^{-N,l-1}}\right).
\end{align*}
Moreover, if $v \in H_{(b)}^{\tilde{r}'',l''}$ for some $\tilde{r}'' \in \R$ (with no further condition on $\tilde{r}''$) and $l''>-\frac{1}{2}$, we have:
\begin{align*}
\left\| \chi_1 v\right\|_{H_{(b)}^{-\tilde{r},1-l}}\leq C\left(\left\|\hat{T}_s(\sigma)^*v\right\|_{H_{(b)}^{-\tilde{r},-l}} + \left\|v\right\|_{H_{(b)}^{-N,-l}}\right).\\
\end{align*}

The constants are uniform with respect to $\sigma$ in a compact subset of $\left\{ \Im(\sigma)\geq 0\right\}\setminus\left\{0\right\}$
\end{prop}

We can also get the semiclassical version of this estimate (meaning uniform with respect to $\sigma$ in a strip $\left\{ 0\leq \Im(\sigma)\leq \eta, \lvert\Re(\sigma)\rvert>A\right\}$). This is done in \cite{vasy2020limiting}, in section 5. 
We give the version for $\hat{T}_h(\tilde{z}) = h^2\hat{T}_s(h^{-1}\tilde{z})$ with $0\leq \Im(\tilde{z})\leq \eta h$ and $\left|\Re(\tilde{z})\right|\geq Ah$. Note that the non trapping assumption mentioned in \cite{vasy2020limiting} plays no role for the microlocal version of the estimate. Note that we use $b$ Sobolev spaces in the statement instead of second microlocal Sobolev spaces, thus we lose some precision with respect to \cite{vasy2020limiting}.
\begin{prop}
\label{scRinEstimate}
Let $U$ be a neighborhood of $\mathcal{R}_{in}$ separated from fiber infinity. Let $B_0, B_1, G \in \Psi^{0,0}_{sc,h}$ with $WF(B_0)\cup WF(B_1)\cup WF(G)\subset U$, $WF(B_0)\cup WF(B_1)\subset Ell(G)$ and every bicharacteristic from $WF(B_0)$ reaches $Ell(B_1)$ in finite time (with the time having the same sign as $\Re(\tilde{z})$) while remaining in $Ell(G)$.
If $0\leq \Im(\tilde{z})\leq \eta h$, $\lvert \Re(\tilde{z})\rvert \geq Ah$, $l+\frac{1}{2}<0$ and $\tilde{r}\in \R$. For all $u\in H^{\tilde{r}',l'}_{(b),h}$ (with no conditions on $\tilde{r}'$ and $l'$) and all $N\in \N$, we have a constant $C>0$ uniform with respect to $h$, $\tilde{z}$ and $u$ such that:
\begin{align*}
\left\|B_0 u\right\|_{H_{(b),h}^{\tilde{r},l}}\leq& C\left(h^{-1}\left\|G\hat{T}_h(\tilde{z})u\right\|_{H_{(b),h}^{\tilde{r},l-1}}+\left\|B_1 u\right\|_{H^{\tilde{r},l}_{(b),h}} + h^N\left\|u\right\|_{H^{-N, l}_{(b),h}}\right)\\
\left\|B_0 u\right\|_{H_{(b),h}^{\tilde{r},l}}\leq& C\left(h^{-1}\left\|G\hat{T}_h(\tilde{z})u\right\|_{H_{(b),h}^{\tilde{r}-1,l}}+\left\|B_1 u\right\|_{H^{\tilde{r},l}_{(b),h}} + h^N\left\|u\right\|_{H^{-N, l}_{(b),h}}\right)
\end{align*}
Moreover, if $u\in H^{\tilde{r}',l'}_{(b),h}$ with $l'>-\frac{1}{2}$:
\begin{align*}
\left\|B_0 u\right\|_{H_{(b),h}^{-\tilde{r}, 1-l}}\leq& C\left(h^{-1}\left\|G\hat{T}_h(\tilde{z})^*u\right\|_{H_{(b),h}^{-\tilde{r},-l}}+ h^{N}\left\|u\right\|_{H^{-N,1-l}_{(b),h}}\right)\\
\end{align*}
\end{prop}
We have the corresponding estimate near $\mathcal{R}_{out}$:
\begin{prop}
\label{scRoutEstimate}
Let $U$ be a neighborhood of $\mathcal{R}_{out}$ separated from fiber infinity. Let $B_0, B_1, G \in \Psi^0_{sc,h}$ with $WF(B_0)\cup WF(B_1)\cup WF(G)\subset U$, $WF(B_0)\cup WF(B_1)\subset Ell(G)$ and every bicharacteristic from $WF(B_0)$ reaches $Ell(B_1)$ in finite time (with the time having the same sign as $-\Re(\tilde{z})$) while remaining in $Ell(G)$.
If $0\leq Im(z)\leq \eta h$, $\lvert \Re(z)\rvert \geq Ah$, $\tilde{r}+l+\frac{1}{2}+2s + 4M\Im(\sigma)>0$ and $u\in H_{(b)}^{\tilde{r}',l'}$ with $\tilde{r}'+l'+\frac{1}{2}+2s+ 4M\Im(\sigma)>0$, then for all $N\in \N$, there exists $C>0$ independent of $u$, $h$ and $z$ such that:
\begin{align*}
\left\|B_0 u\right\|_{H_{(b),h}^{\tilde{r}, l}}\leq C\left(h^{-1}\left\|G\hat{T}_h(z)u\right\|_{H_{(b),h}^{\tilde{r},l-1}}+ h^{N}\left\|u\right\|_{H^{-N,l}_{(b),h}}\right) \\
\left\|B_0 u\right\|_{H_{(b),h}^{\tilde{r}, l}}\leq C\left(h^{-1}\left\|G\hat{T}_h(z)u\right\|_{H_{(b),h}^{\tilde{r}-1,l}}+ h^{N}\left\|u\right\|_{H^{-N,l}_{(b),h}}\right) 
\end{align*}
and if $u\in H_{(b)}^{\tilde{r}',l'}$ with no condition on $\tilde{r}', l'$, we have:
\begin{align*}
\left\|B_0 u\right\|_{H_{(b),h}^{-\tilde{r},1-l}}\leq C\left(h^{-1}\left\|G\hat{T}_h^*(z)u\right\|_{H_{(b),h}^{-\tilde{r},-l}}+\left\|B_1 u\right\|_{H^{-\tilde{r},1-l}_{(b),h}} + h^N\left\|u\right\|_{H^{-N, 1-l}_{(b),h}}\right)
\end{align*}
\end{prop}

We also need to state an estimate which is true uniformly up to $\sigma = 0$. We begin by recalling the definition of the effective normal operator (which is compatible with Definition 2.4 in \cite{vasy2020resolvent} although phrased in a slightly different setting).
\begin{definition}\label{defNormalOp}
The effective normal operator of $\hat{P}(\sigma)$ denoted by $N_{\eff}(\hat{P}(\sigma))$ is $\hat{P}(\sigma)$ modulo $x(x+\sigma)^2\text{Diff}^2_b([0,\frac{1}{r_+-\epsilon})_x\times \mathcal{B}_s)$. Similarly, $N_{\eff}(\hat{T}_s(\sigma))$ is $N_{\eff}(\hat{T}_s(\sigma))$ modulo $x^{-1}(x+\sigma)^2\text{Diff}^2_b([0,\frac{1}{r_+-\epsilon})_x\times \mathcal{B}_s)$.
\end{definition}

\begin{prop}
Let $\chi_1$ be a smooth cutoff compactly supported in $\left\{ x<\frac{1}{6M}\right\}$ with $\chi_1 = 1$ on $\left\{ x<\frac{1}{12M}\right\}$.
Assume that $\alpha \in \left(l+\frac{1}{2}+s-\left|s\right|, l+\frac{3}{2}+s+\left|s\right|\right)$, $l+\tilde{r} >-\frac{1}{2}-2s$ and $l<-\frac{1}{2}$. Then, there exists $\sigma_0>0$ such that for all $\sigma\in \C$ with $\Im(\sigma)\geq 0$ and $\lvert\sigma\rvert \leq \sigma_0$, we have:
\begin{align*}
\left\| (x+\lvert \sigma\rvert)^\alpha \chi_1 u \right\|_{H^{\tilde{r},l}_{(b)}}& \leq C \left(\left\| (x+\lvert \sigma \rvert)^\alpha \hat{P}(\sigma) u\right\|_{H^{\tilde{r}-1, l+2}_{(b)}} + \left\|u\right\|_{H^{-N,l-1}_{(b)}}\right)\\
\end{align*}
\end{prop}

\begin{proof}
As for Proposition \ref{midFreqInfEstimate}, using that $[\hat{P}(\sigma), \chi_1]\in \text{Diff}^1_b$ is supported away from the boundary and on the elliptic set of $\hat{P}(\sigma)$, it is enough to prove
\begin{align}\label{reducedIneq2}
\left\|(x+\lvert \sigma\rvert)^\alpha \chi_1 u\right\|_{H_{(b)}^{\tilde{r}, l}}&\leq C\left(\left\|(x+\lvert \sigma\rvert)^\alpha \hat{P}(\sigma)\chi_1 u\right\|_{H_{(b)}^{\tilde{r}-1,l+2}} + \left\|\chi_1 u\right\|_{H_{(b)}^{-N,l-1}}\right)
\end{align}
and similarly for the second inequality.
The estimate is therefore localized in a neighborhood of $\left\{x=0\right\}$ where the microlocal estimates obtained in the proof of Theorem 2.5 in \cite{vasy2020resolvent} applies. As for the previous proposition, because we consider an operator acting on a complex line bundle, we are not exactly in the setting of \cite{vasy2020resolvent}. However, the proof is completely parallel under the following hypotheses (which are not the most general ones):
\begin{align*}
P(\sigma) = P(0) + \sigma Q + (1-R)\sigma^2
\end{align*}
\begin{enumerate}[label = (H\arabic*)]
\item
$R \in xC^{\infty}(X)$ and $\Im(R) \in x^2C^{\infty}(X)$
\item
\label{hypFormPZero}
\begin{align*}
P(0) = (x^2D_x)^2 + x^2 A  + (i(n-1)+\beta)x(x^2D_x) + x^2a'
\end{align*}
where $A \in C^{\infty}([0,M)_x,\text{Diff}^2(E))$ such that $\frac{x^2}{2i}(A-A^*)\in x^3C^{\infty}([0,M)_x, \text{Diff}^1(E))$ and $ a',\beta \in C^{\infty}(X)$. We also require that the principal symbol of $x^{-2}P(0)$ is elliptic in the $b$ sense. Note that we imposed that the coefficients $a_{0,j}$ in the local coordinate expression in the proof of proposition 2.1 in \cite{vasy2020resolvent} 	are zero, which is true in our case. 
\item
\begin{align*}
Q = b_0 x(x^2D_x) + x^2S + \gamma x
\end{align*}
where $S \in C^{\infty}([0,M)_x, \text{Diff}^1(E))$ is principally self-adjoint and $b_0, \gamma \in C^{\infty}(X)$ is such that $\Im(b_0)\in xC^{\infty}(X)$. The coefficients of $Q$ are supposed to be independent of $\sigma$

\item \label{hyp3} The operators $x^2A$ and $x^2R$ admit the following decomposition:
\begin{align*}
x^2A &= \sum_{j = 1}^N (A_j^*A_j + A^*_jA_j'+A^\dagger_j A_j) + A''\\
x^2S &= \sum_{j=1}^N (A^*_j R_j + R_j^\dagger A_j) + R''
\end{align*}
where $A_j \in xC^{\infty}([0,M)_x,\text{Diff}^1(E)) \subset x\text{Diff}^1(\tilde{E})$, $A_j',A^\dagger_j,R_j,R_j^\dagger \in xC^{\infty}(X)$ and $A'', R'' \in x^2 C^{\infty}(X)$.
\end{enumerate}

With the previous hypotheses, we see that (2.2) and (2.3) in \cite{vasy2020resolvent} are satisfied with $\beta_I = \Im(\beta)$ and $\beta_I' = -\Re(\beta)\frac{n-2}{2} + \Im(a')$ and $\gamma_I = \Im(\gamma)$. We have the following form for the conjugated operators:
\begin{align*}
\hat{P}(\sigma) = P(0) + \sigma \hat{Q} + \sigma^2 \hat{R} - 2\sigma\left(x^2D_x + i\frac{n-1}{2}x+\frac{\beta-\gamma}{2}x\right)
\end{align*}
where
\begin{align*}
\hat{Q} &= Q-x\gamma\\
\hat{R} &= R-xb_0\\
\end{align*}

Under these hypotheses, the effective normal operator (see Definition \ref{defNormalOp}) is
\begin{align*}
N_{\eff}(\hat{P}(\sigma)) =& (x^2D_x)^2+ i(n-1)x(x^2D_x)+ x^2A(0) + \beta_{|_{\partial_x}} x^2\left(xD_x + i\frac{n-2}{2}\right)\\
& + x^2(\beta')_{|_{\partial X}} -2\sigma\left(x^2D_x + i\frac{n-1}{2}x + \frac{\beta_{|_{\partial X}}-\gamma_{|_{\partial X}}}{2}x\right)
\end{align*}
where $\beta' = a'-i\beta\frac{n-2}{2}$ (therefore $\Im(\beta') = \beta_I'$).
Note that $(x^2D_x)^2 + i(n-1)x(x^2D_x) + x^2A(0)$ is similar to the term $\Delta_{g_0}$ in \cite{vasy2020resolvent}.

We add two last hypotheses:
\begin{enumerate}[resume, label = (H\arabic*)]
\item \label{hyp4} $\beta_{|_{\partial X}} \in i\R$ and $\Re(\beta'_{|_{\partial X}})+\lambda_0>\frac{\left(\beta_{|_{\partial X}}\right)^2}{4}-\left(\frac{n-2}{2}\right)^2$ where $\lambda_0$ is the smallest eigenvalue of $A(0)$. This hypothesis is present in theorem 2.5 in \cite{vasy2020resolvent} (but with $\lambda_0 = 0$).
\item \label{hyp5} For every $\lambda_k$ eigenvalue of $A(0)$:
\begin{align*}\frac{1}{2}+\frac{i\gamma}{2}+\sqrt{\left(\frac{n-2}{2}\right)^2-\frac{\beta^2}{4}+\lambda_k+\beta'}\notin -\N \\
\frac{1}{2}-\frac{i\gamma}{2}+\sqrt{\left(\frac{n-2}{2}\right)^2-\frac{\beta^2}{4}+\lambda_k+\beta'} \notin -\N
\end{align*}
This hypothesis is not mentioned in \cite{vasy2020resolvent} but it seems necessary (it is at least sufficient) to ensure that Proposition 5.4 in \cite{vasy2020resolvent} applies to the effective normal operator with no error term. Indeed, under this hypothesis, we can prove that the effective normal operator has no kernel in the space of interest (see appendix \ref{absenceOfKernel}, note that by hypothesis \ref{hypFormPZero}, $A(0)$ is an elliptic formally selfadjoint operator and by hypothesis \ref{hyp4} there exists some constant $C>0$ such that $A(0)+C$ is positive. Therefore, even replacing $\beta'$ by $\beta'-C$ and taking $L:= A(0)+C$, we can apply the result of appendix \ref{absenceOfKernel}).
\end{enumerate}

We now check the hypothesis in our case:
We write: $P(\sigma) = P(0) + \sigma Q -(1-R)\sigma^2$ with $P(0)$ already defined and 
\begin{align*}
Q =& \frac{4 M a r}{\Delta_r^{2}}D_{\phi} -\frac{4 M \left(2 M r - a^{2}\right)}{r \Delta_r}x^2D_x -\frac{1}{r^{2} \Delta_r^{2}} q_1\\
q_1 = & 2 i r^{5} s + r^{4} \left(- 6 i M s - 2 a s \cos{\left (\theta \right )}\right) + r^{3} \left(- 8 i M^{2} s + 4 M a s \cos{\left (\theta \right )} + 2 i a^{2} s\right)\\
& + r^{2} \left(8 i M^{3} s + 6 i M a^{2} s + 2 i M a^{2} - 2 a^{3} s \cos{\left (\theta \right )}\right) + r \left(- 4 i M^{2} a^{2} s - 4 i M^{2} a^{2}\right)+ 2 i M a^{4} \\
r^2\Delta_r^2R = & - 16 M^{3} a^{2} r + 4 M^{2} a^{4} - 4 M r^{5} + r^{4} \left(4 M^{2} + a^{2} \sin^{2}{\left (\theta \right )}\right)\\
& + r^{3} \left(- 2 M a^{2} \sin^{2}{\left (\theta \right )} - 4 M a^{2}\right) + r^{2} \left(16 M^{4} + a^{4} \sin^{2}{\left (\theta \right )}\right)
\end{align*}

\begin{enumerate}[label = (H\arabic*)]
\item
The first hypothesis is obvious with the definition of $R$.
\item
We have the correct form for $P(0)$ if we define $A$ as in \eqref{defA} and $\beta := \frac{2i(Mr(1+s)-a^2-r^2s)}{\Delta_r} = -2is + O(x)$ and $a' = 0$. We check that
\begin{align*}
\frac{1}{2i}(A-A^*) &= \frac{2s a(M-r)}{\Delta_r^2}D_\phi
\end{align*}
\item We have the correct form for $Q$ if we define $b_0 := -\frac{4 M \left(2 M r - a^{2}\right)}{\Delta_r} = O(x)$, $S := \frac{4 M a r^3}{\Delta_r^{2}}D_{\phi}$ (self adjoint) and $\gamma := -\frac{1}{r \Delta_r^{2}} q_1 = -2is + O(x)$
\item
The decomposition of $x^2A$ has already been checked and the decomposition of $x^2S$ is obtained by taking $R^\dagger_3:=\frac{4iMar}{\sqrt{1-\frac{a^2}{\Delta_r}}\Delta_r^2}$ (and all the other terms equal to zero).
\end{enumerate}
Therefore in our case, the effective normal operator is
\begin{align*}
N_{\eff}(\hat{P}(\sigma)) =& (x^2D_x)^2+ 2ix(x^2D_x)\\ 
&+ x^2\left(\frac{1}{\sin\theta}D_\theta \sin\theta D_\theta + \frac{1}{\sin^2\theta}D_\phi^2 + \frac{2s\cos\theta}{\sin^2\theta}D_\phi + s^2\text{cotan}^2\theta+s\right)\\
& -2is x^2\left(xD_x + \frac{i}{2}\right) -sx^2 -2\sigma\left(x^2D_x + ix\right)
\end{align*}
To be coherent with the notation in \cite{vasy2020resolvent}, we study the conjugated and renormalized operator:
\begin{align*}
x^{\frac{-5}{2}}N(\hat{P}(0))x^{\frac{1}{2}} =& (xD_x)^2+A(0)+\frac{1}{4} -2isxD_x -s
\end{align*}
After a Mellin transform in $x$, we obtain:
\begin{align*}
\tau^2_b -2is\tau_b -s + A(0) + \frac{1}{4}
\end{align*}
To deduce the central weight interval (the analog of (2.12) in \cite{vasy2020resolvent}, we have to compute the eigenvalues of $A(0)-s$ (which is the spin-s-weighted Laplacian). The eigenvalues of $A(0)-s$ are $(l+s)(l-s)+l$ for $l \in \N$ such that $l\geq \lvert s\rvert$. Therefore, the central interval for weights for the scattering end (the equivalent of (2.12) in \cite{vasy2020resolvent}) is $\left(-\frac{3}{2}-s-\left|s\right|, -\frac{1}{2}-s+\left|s\right|\right)$.
\end{proof}
\begin{remark}\label{estimateSecondMic}
For $\alpha = 0$, $l\in (-\frac{3}{2}-s-\left|s\right|,-\frac{1}{2})$, $l+r>-\frac{1}{2}$ and $\Im(\sigma)\geq 0$ the second resolved microlocal spaces $H^{r,m,l}_{sc,b,res}$ used in \cite{vasy2020resolvent} enables to state a more precise estimate:
\begin{align*}
\left\| \chi_1 u \right\|_{H^{r,r+l,l}_{sc,b,res}}& \leq C \left(\left\| (x+\lvert \sigma \rvert)^{-1}\hat{P}(\sigma) u\right\|_{H^{r-2,l+r+1 ,l+1}_{sc,b,res}} + \left\|u\right\|_{H^{-N,l-1}_{(b)}}\right)
\end{align*}
We can then use the bound (for $\epsilon\in[0,1]$): $(x+\lvert \sigma \rvert)^{-1}\leq \lvert\sigma\rvert^{\epsilon-1}(x+\lvert \sigma \rvert)^{-\epsilon}$ and the fact that $(x+\lvert \sigma \rvert)^{\epsilon}H_{sc,b,res}^{r-2,l+r+1,l+1} \supset H_{sc,b,res}^{r-2,l+r+1, l+1+\epsilon} \supset H_{(b)}^{r-\epsilon, l+1+\epsilon}$ to get:
\begin{align*}
\left\| \chi_1 u \right\|_{H^{r,l}_{(b)}}& \leq C \left(\lvert\sigma\rvert^{\epsilon-1}\left\|\hat{P}(\sigma) u\right\|_{H^{r-\epsilon ,l+1+\epsilon}_{(b)}} + \left\|u\right\|_{H^{-N,l-1}_{(b)}}\right)
\end{align*}
\end{remark}
 
\subsection{High frequency estimate at the trapped set}\label{subSecTrapp}
We use notations introduced in section \ref{semiClFlowSection}. We recall that the operator $\hat{T}_{s,h}$ was defined at the end of Section \ref{analyticFramework}. We begin by some preliminary computations near the trapped set (on the region $U_I$ with $I = (r_{min},r_{max})$).

\begin{lemma}
\label{subSymbTrapp}
Let $z \in \C$.
On $U_I$, we have:
\begin{align*}
\frac{1}{2hi}(\hat{T}_{s,h}(z)-\hat{T}_{s,h}^*(z)) =& \frac{2(\Delta_r+rs(r-M)}{r}hD_r + \frac{2a(s(M-r)+2Mr\Im(z))}{\Delta_r}hD_\phi\\& - 2s\Re(z)\frac{M(a^2-r^2)+r\Delta_r}{\Delta_r} +\frac{i}{r^2}(2Mr(s+2)-a^2-3r^2(s+1))\\
&+2\Im(z)as\cos\theta
\end{align*}
(where as usual, the adjoint is computed with respect to the volume form $r^2\sin\theta\dd t \dd r\dd \phi \dd \theta$).
If $z = z_0 + O(h)$ with $z_0\in \left\{-1,1\right\}$, the principal symbol of $\frac{1}{2i}(T_s-T_s^*)$ is $0$ on the trapped set $K_{z_0}$.
\end{lemma}

\begin{remark}\label{remarkVolumeForm}
Note that since the principal symbol of $\hat{T}_s$ is $\tilde{G}$, the most natural choice for the volume form is the one associated to the metric $\tilde{g}$ which is $\sin(\theta)\dd r\dd\theta\dd\phi$. However, the factor $r^2$ is harmless here since it amounts to replace $\hat{T}_{s,h}^*$ by $r^{-2} \hat{T}_{s,h} r^{2}$ and the difference $r^{-2}[\hat{T}_{s,h}, r^2]$ has principal symbol $r^{-2}H_{p_{h}} r^2= 0$ on the trapped set. However, since $H_{p_{h}}$ does not vanish on the trapped set, the choice of the volume form matters.
\end{remark}
\begin{proof}
By a direct computation we find
\begin{align*}
p_{sub} := \frac{2(\Delta_r+rs(r-M))}{r}\xi - \frac{2as(r-M)}{\Delta_r}\zeta - 2z_0s\frac{M(a^2-r^2)+r\Delta_r}{\Delta_r}
\end{align*}
Then, using that $\xi = 0$ on $K_{z_0}$, we find:
\begin{align*}
p_{sub} &= - \frac{2as(r-M)}{\Delta_r}\zeta - 2z_0s\frac{M(a^2-r^2)+r\Delta_r}{\Delta_r}
\end{align*}
We introduce the function $\alpha := \frac{-(r^2+a^2)z_0+a\zeta}{\Delta_r}$. We have $a\zeta = \Delta_r\alpha + (r^2+a^2)z_0$ and, on $K_{z_0}$, we have $(r-M)\alpha \Delta_r =- 2r\Delta_r z_0$. Therefore on $K_{z_0}$, we have:
\begin{align*}
p_{sub} &= -\frac{2s(r-M)}{\Delta_r}(\alpha\Delta_r +(r^2+a^2)z_0) - 2sz_0\frac{M(a^2-r^2)+r\Delta_r}{\Delta_r}\\
&= 4srz_0 - \frac{2sz_0(r-M)(r^2+a^2)}{\Delta_r} - 2z_0s\frac{M(a^2-r^2)}{\Delta_r}-2sz_0r\\
&= 2srz_0 -\frac{2sz_0}{\Delta_r}r\Delta_r\\
&= 0
\end{align*}
\end{proof}

We need the following modified version of theorem 1 in \cite{dyatlov2016spectral}:
\begin{theorem}
\label{trappingDyatlov}
Let $P_h\in \Psi^m_h$ be a principally scalar and principally real operator on a smooth (complex) vector bundle (with a fixed smooth hermitian inner product) over some orientable smooth manifold $X$ (with a fixed volume form). We denote by $p_h$ its semiclassical principal symbol, by $\Sigma$ the semiclassical characteristic set $p_h^{-1}\left\{0\right\}$ and by $H_b$ the semiclassical Hamiltonian flow. We assume that there exists $\phi_{\pm}$ smooth functions defined on a bounded open set $U$ of $T^*X$ such that for a fixed small $\epsilon>0$:
\begin{enumerate}
\item There exists $\delta>0$ such that, $U_\delta:= \left\{ \lvert\phi_+\rvert<\delta, \lvert\phi_-\rvert<\delta, \lvert p\rvert<\delta\right\}$ is compactly contained in $U$.
\item $H_p \phi_\pm = \mp c_\pm \phi_\pm$ with $c_\pm$ smooth bounded positive functions on $U$ satisfying $\inf c_{\pm} >\nu_{min}-\epsilon>0$. 
\item $\left\{\phi_+, \phi_-\right\}>0$
\end{enumerate}
Let $V\Subset U$ be a neighborhood of $K:= \left\{\phi_+=\phi_- = 0\right\}$. We assume that $\left\|\frac{1}{2ih}\mathfrak{s}_h(P_h-P^*_h)\right\|_{End(E)}<\frac{\nu_{\min}-\epsilon}{2}$ on $K$. There exist $B_0, B_1, G \in \Psi^0_h$ with:
\begin{itemize}
\item $WF_h(B_0)\cup WF_h(B_1) \cup WF_h(G)\subset U$
\item $\mathfrak{s}_h(B_0) = 1$ on $V$
\item $WF_h(B_1)\cap \left\{\phi_+=0\right\} = \emptyset$
\end{itemize}
such that:
\begin{align*}
\left\|B_0 u\right\|_{L^2}\leq C\left(\left\|B_1 u\right\|_{L^2} + o(h^{-2})\left\|GP_h(z)u\right\|_{L^2}\right)
\end{align*}
\end{theorem}
\begin{remark}
We use Lemma \ref{neighTrapping} for the existence of $\phi_{\pm}$. Note that in particular, $WF_h(B_0)\cap \Gamma_+ = \emptyset$. 
\end{remark}
\begin{remark}
The principal differences with respect to \cite{dyatlov2016spectral}, theorem 1 are the fact that $P_h$ is not of the form $\tilde{P}_h - \lambda$ with $\tilde{P}_h$ selfadjoint (but the hypothesis $\left\|\frac{1}{2ih}\mathfrak{s}_h(P_h-P^*_h)\right\|_{End(E)}<\frac{\nu_{\min}}{2}$ on $K$ replaces self adjointness in the proof) and the fact that we do not use an absorbing potential (and therefore we keep the estimate microlocal). However, the proof given in \cite{dyatlov2016spectral} adapts with minor changes.
\end{remark}

We can now apply Theorem \ref{trappingDyatlov} to the operators $\hat{T}_{s,h}(z)$, $-\hat{T}_{s,h}(z)$, $\hat{T}_{s,h}(z)^*$ and $-\hat{T}_{s,h}(z)^*$ (see Lemma \ref{neighTrapping} and Lemma \ref{subSymbTrapp} to see that the requirements are met) to get the semiclassical estimate near the trapped set:
\begin{prop}\label{estimateTrappedSet}
Let $z = \pm 1 +O(h)$. There exists $U$ a bounded neighborhood of the trapped set $K$ and $B_{K}$, $B_{0}$ in $\Psi_{b,h}^{0,0}$ with $B_{K}$ elliptic on a neighborhood of $K$,  $WF(B_0)\cap \Sigma\subset \Sigma_{\text{sgn}\Re(z)}$ with either $WF(B_0)\cap \Gamma_{+} = \emptyset$ or $WF(B_0)\cap \Gamma_{-} = \emptyset$ (we can chose between the two) such that for all $u\in \overline{H}^{\tilde{r},l}_b$ and all $v\in \dot{H}^{\tilde{r},l}_{b}$:
\begin{align}
\left\|B_K u \right\|_{\overline{H}^{\tilde{r},l}_{b,h}}\leq C\left(\left\|B_0 u \right\|_{\overline{H}^{\tilde{r},l}}+h^{-2}\left\|\hat{T}_{s,h}(z)u\right\|_{\overline{H}^{\tilde{r},l-1}_{b,h}} + h\left\|u\right\|_{\overline{H}^{\tilde{r},l}_{b,h}}\right) \label{EstimateK1} \\
\left\|B_K v \right\|_{\dot{H}^{\tilde{r},l}_{b,h}}\leq C\left(\left\|B_0 v \right\|_{\overline{H}^{\tilde{r},l}}+h^{-2}\left\|\hat{T}_{s,h}(z)^*v\right\|_{\dot{H}^{\tilde{r},l-1}_{b,h}} + h\left\|v\right\|_{\dot{H}^{\tilde{r},l}_{b,h}}\right) \label{EstimateK2}
\end{align}
\end{prop}

\subsection{Global estimates}\label{subSecGlob}
We recall the following elementary lemma.
\begin{lemma}
\label{FredholmEstimate}
Let $X,Y$ be Banach spaces. Let $\mathcal{X}$ and $Z$ be Banach spaces such that we have compact inclusions $\mathcal{X}\subset X$, $Y^*\subset Z$.
If $P:\mathcal{X}\rightarrow Y $ is a bounded operator satisfying the following estimates:
\begin{align}
\left\|u\right\|_{\mathcal{X}}&\leq C(\left\|Pu\right\|_{Y} + \left\|u\right\|_{X}) \label{ineqP} \\
\left\|v\right\|_{Y^*}&\leq C(\left\|P^*u\right\|_{\mathcal{X}^*} + \left\|u\right\|_{Z}) \notag
\end{align}
Then $P$ is Fredholm.
\end{lemma}

\begin{proof}
Let $A:=\left\{u\in B(0,1)_{\mathcal{X}}: \left\|u\right\|_{\mathcal{X}}\leq C\left\|u\right\|_{X}\right\}$. This set is compact in $X$ (closed and included in $B(0,1)_{\mathcal{X}}$). The topology induced by $\mathcal{X}$ is the same by definition. So it is also compact in $\mathcal{X}$.
$Ker(P)\cap B(0,1)_{\mathcal{X}}\subset A $ so the $\mathcal{X}$-unit ball of $ker(P)$ is compact for the $\mathcal{X}$ induced topology. Then $Ker(P)$ is finite dimensional.
Similarly $Ker(P^*) = Ran(P)^{\perp}$ is finite dimensional.
We now prove that $Ran(P)$ is closed. We take $y_1,...,y_k$ a normed basis of $Ker(P)$ and we denote by $y_1^*,...,y_k^*$ extensions (of norm less than 1) of the dual basis (obtained by Hahn-Banach theorem).
Assume that $Pu_n \to y$ in $Y$. Without loss of generality, we can assume that $y_i^*(u_n) = 0$ (even replacing $u_n$ by $u_n - \sum_{i\leq k} y_i^*(u_n)y_i$).
By contradiction assume that $u_n$ is unbounded in $\mathcal{X}$. Extracting a subsequence, we can assume $\left\|u_n\right\|_{\mathcal{X}} \to +\infty$. Then $P\frac{u_n}{\left\|u_n\right\|_{\mathcal{X}}}\to 0$. By compactness, we can assume (even extracting a subsequence) that $\frac{u_n}{\left\|u_n\right\|_{\mathcal{X}}}$ converges to $z\in X$ for the topology of $X$. Using inequality \eqref{ineqP}, we deduce that $\frac{u_n}{\left\|u_n\right\|_{\mathcal{X}}}$ is Cauchy in $\mathcal{X}$. We deduce $z\in Ker(P)$ and we have the convergence in $\mathcal{X}$. Then by continuity, for $i=1,...,k$, we have $y_i^*(z) = 0$ and because $z\in Ker(P)$ we deduce $z=0$. But it is a contradiction since $\left\|\frac{u_n}{\left\|u_n\right\|_{\mathcal{X}}}\right\|_{\mathcal{X}}=1$. So $u_n$ is bounded in $\mathcal{X}$.
Even extracting a subsequence, we can assume that $u_n \to x$ in $X$. Then we use \eqref{ineqP} an deduce that $u_n$ is Cauchy in $\mathcal{X}$. Finally, $u_n \to x$ in $\mathcal{X}$ and $Pu_n \to Px$ which proves that $y= Px \in Ran(P)$.
\end{proof}

\begin{definition}
We define the following spaces $\mathcal{X}^{\tilde{r},l}_{\sigma} := \left\{u\in \overline{H}_b^{\tilde{r},l}: \hat{T}_s(\sigma)u \in \overline{H}_{(b)}^{\tilde{r},l-1} \right\}$ endowed with the norm $\left\|u\right\|_{\mathcal{X}^{\tilde{r},l}_{\sigma}} := \left\|u\right\|_{\overline{H}_{(b)}^{\tilde{r},l}} + \left\|\hat{T}_s(\sigma)u\right\|_{\overline{H}_{(b)}^{\tilde{r}, l-1}}$.
We also define $\mathcal{W}^{\tilde{r},l}_{\sigma} := \left\{ u\in \overline{H}^{\tilde{r},l}_{(b)}: \hat{T}_s(\sigma) u \in \overline{H}^{\tilde{r}, l}_{(b)}\right\}$ endowed with the norm $\left\|u\right\|_{\mathcal{W}^{\tilde{r},l}_{\sigma}} := \left\|u\right\|_{\overline{H}_{(b)}^{\tilde{r},l}} + \left\|\hat{T}_s(\sigma)u\right\|_{\overline{H}_{(b)}^{\tilde{r}, l}}$.
\end{definition}
The Fredholm estimates that we want to prove are the following
\begin{theorem}
Let $K\subset \left\{\sigma \in \C : \Im(\sigma)\geq 0\right\}\setminus\left\{0\right\}$ be compact.
If $l < -\frac{1}{2}$, $\tilde{r}+l> -\frac{1}{2}-2s - 4M\Im(\sigma)$ and $\tilde{r}>\frac{1}{2}+s-\frac{a^2+r_+^2}{r_+-M}\Im(\sigma)$, then for all $u\in H^{\tilde{r},l}_{(b)}$ and $v\in H^{-\tilde{r},1-l}_{(b)}$ we have:
\begin{align*}
\left\| u \right\|_{\overline{H}_{(b)}^{\tilde{r},l}} &\leq C\left(\left\| \hat{T}_s(\sigma)u\right\|_{\overline{H}_{(b)}^{\tilde{r},l-1}}+ \left\|u\right\|_{\overline{H}_{(b)}^{\tilde{r}-1, l-1}}\right)\\
\left\| u \right\|_{\overline{H}_{(b)}^{\tilde{r},l}} &\leq C\left(\left\| \hat{T}_s(\sigma)u\right\|_{\overline{H}_{(b)}^{\tilde{r}-1,l}}+ \left\|u\right\|_{\overline{H}_{(b)}^{\tilde{r}-1, l-1}}\right)\\
\left\| v \right\|_{\dot{H}_{(b)}^{-\tilde{r}, 1-l}} &\leq C\left(\left\|\hat{T}_s(\sigma)^*v\right\|_{\dot{H}_{(b)}^{-\tilde{r}, -l}} + \left\|v\right\|_{H_{(b)}^{-\tilde{r}-1, -l}} \right)\\
\end{align*}
\end{theorem}
\begin{proof}
The key point is that the junction between the estimates take place in an elliptic region for $\hat{T}_s(\sigma)$ and $\hat{T}_s(\sigma)^*$, therefore, we can use the elliptic estimate to bound the commutator terms. We prove the first estimate as an example.
Let $u\in \mathcal{X}^{\tilde{r},l}_{\sigma}$. Then we define three smooth cut-off functions $\chi_1$, $\chi_2$ and $\chi_3$ such that $\chi_1+\chi_2+\chi_3=1$ and such that $\supp(\chi_1)\subset (r_+-\epsilon, 2M+D)$ (where $D$ is defined in Proposition \ref{combinedNearHorizon}), $\chi_1= 1$ on $(r_+-\epsilon, 2M+D)$, $\supp(\chi_3)\subset (6M,+\infty)$ and $\chi_3 = 1$ in a neighborhood of $+\infty$. We have:
\begin{align*}
\left\|u\right\|_{H_{(b)}^{\tilde{r},l}}&\leq \left\|\chi_1 u \right\|_{H_{(b)}^{\tilde{r},l}} + \left\|\chi_2 u \right\|_{H_{(b)}^{\tilde{r},l}}+ \left\|\chi_3 u \right\|_{H_{(b)}^{\tilde{r},l}}.\\
\end{align*}
We use proposition \ref{combinedNearHorizon} to bound $\left\|\chi_1 u \right\|_{H_{(b)}^{\tilde{r},l}}$, proposition \ref{ellipticEstimate} (elliptic estimate) to bound $\left\|\chi_2 u \right\|_{H_{(b)}^{\tilde{r},l}}$ and proposition \ref{estimateNearInfMedFreq} to bound $\left\|\chi_3 u \right\|_{H_{(b)}^{\tilde{r},l}}$. We get:
\begin{align*}
\left\|u\right\|_{H_{(b)}^{\tilde{r},l}}&\leq C\left(\left\| \hat{T}_s(\sigma) \chi_1 u \right\|_{H_{(b)}^{\tilde{r},l-1}} + \left\| \hat{T}_s(\sigma) \chi_2 u \right\|_{H_{(b)}^{\tilde{r},l-1}} + \left\| \hat{T}_s(\sigma) \chi_3 u \right\|_{H_{(b)}^{\tilde{r},l-1}} + \left\|u\right\|_{\overline{H}_{(b)}^{\tilde{r}-1, l-1}}\right) \\
&\leq C\left(\left\|\hat{T}_s(\sigma) u \right\|_{H_{(b)}^{\tilde{r}, l-1}} + \sum_{i=1}^3\left\|[\hat{T}_s(\sigma), \chi_i] u\right\|_{H_{(b)}^{\tilde{r}, l-1}} + \left\|u\right\|_{\overline{H}_{(b)}^{\tilde{r}-1, l-1}} \right)
\end{align*}
The terms $\left\|[\hat{T}_s(\sigma), \chi_i] u\right\|_{H_{(b)}^{\tilde{r}, l-1}}$ can be estimated by an elliptic estimate (proposition \ref{ellipticEstimate}).
We prove the bound for $\hat{T}_s(\sigma)^*$ exactly in the same way.
\end{proof}

\begin{remark}
We see in the proof that we loose some regularity away from $\partial X$, we could obtain more precise bounds by using second microlocal spaces (see \cite{vasy2020limiting}).
\end{remark}

\begin{coro}\label{coroFredholmMidFreq}
Let $\sigma \in \C\setminus \left\{0\right\}$ with $\Im(\sigma)\geq 0$.
If $l < -\frac{1}{2}$, $\tilde{r}+l> -\frac{1}{2}-2s - 4M\Im(\sigma)$ and $\tilde{r}>\frac{1}{2}+s-\frac{a^2+r_+^2}{r_+-M}\Im(\sigma)$, then $\hat{T}_s(\sigma)$ is Fredholm as an operator from $\mathcal{X}^{\tilde{r},l}_{\sigma}$ to $\overline{H}^{\tilde{r},l-1}_{(b)}$.
\end{coro}

\subsubsection{Global estimates near zero energy}
\begin{prop}
\label{estimateNearZero}
Let $l < -\frac{1}{2}$, $\tilde{r}+l>-\frac{1}{2} -2s$, $\tilde{r}>\frac{1}{2}+s$ and \[\alpha \in \left(l+\frac{1}{2}+s-\left|s\right|, l+\frac{3}{2}+s+\left|s\right|\right)\]. Then, there exists $\sigma_0>0$ such that for $\lvert \sigma \rvert \leq \sigma_0$ and $\Im(\sigma)\geq 0$
\begin{align*}
\left\| (x+\lvert \sigma\rvert)^\alpha u\right\|_{\overline{H}_{(b)}^{\tilde{r},l}} \leq C\left(\left\| (x+\lvert \sigma \rvert)^\alpha \hat{T}_s(\sigma)u\right\|_{\overline{H}_{(b)}^{\tilde{r}-1, l}}+\left\|u\right\|_{\overline{H}_{(b)}^{\tilde{r}-1, l-1}}\right)
\end{align*}
 and the constant $C$ is independent of $\sigma$.
We have also the version coming from the second microlocalized resolved space (see Remark \ref{estimateSecondMic}) with $l\in (-\frac{3}{2}-s-\left|s\right|,-\frac{1}{2})$, $l+r>-\frac{1}{2}$, $\Im(\sigma)\geq 0$ and $\epsilon \in (0,1)$:
\begin{align*}
\left\|u\right\|_{\overline{H}_{(b)}^{\tilde{r},l}}\leq C\left(\lvert \sigma \rvert^{\epsilon-1} \left\|\hat{T}_s(\sigma)u\right\|_{\overline{H}_{(b)}^{\tilde{r}-\epsilon, l-1+\epsilon}}+ \left\|u\right\|_{H^{-N,l-1}_{(b)}}\right)
\end{align*}

\end{prop}
\begin{proof}
The proof is very similar to the previous one since the gluing is made on elliptic regions.
\end{proof}


\begin{lemma}
\label{restrictedFredholm}
Let $X$ and $Y$ be Banach spaces.
Let $P: X\rightarrow Y$ be a Fredholm operator. Let $B$ be a Banach space which is continuously and densely included in $Y$, then $P$ is a Fredholm operator from $A = P^{-1}(B)$ to $B$
($A$ is endowed with the graph norm). Morever, $\text{ind}(P)\leq \text{ind}(P_{|_A})$.
\end{lemma}
\begin{proof}
First, $ker(P) = ker(P_{|_A})$ is still finite dimensional.
We have an inclusion from $B/(B\cap P(X))$ into $Y/P(X)$. Therefore, $\dim(B/(B\cap P(X)))\leq \dim(Y/P(X))<+\infty$. Moreover, $B\cap P(X)$ is closed into $B$ since $P(X)$ is closed in $Y$.
\end{proof}

\begin{coro}\label{coroFredholmNearZero}
Assume:
\begin{align*}
-\frac{3}{2}-s-\left|s\right|&< l < -\frac{1}{2} \\
\tilde{r}+l &> -\frac{1}{2}-2s\\
\tilde{r} &> \frac{1}{2}+s
\end{align*}
Then, for every $\sigma \in \C \setminus\left\{0\right\}$ such that $\Im(\sigma)\geq 0$, we have that $\hat{T}_s(\sigma)$ is a Fredholm operator between $\mathcal{W}_{\sigma}^{\tilde{r},l}$ and $\overline{H}_{(b)}^{\tilde{r},l}$. Moreover, the index of $\hat{T}_s(\sigma)$ as an operator between $\mathcal{W}_{\sigma}^{\tilde{r},l}$ and $\overline{H}_{(b)}^{\tilde{r},l}$ is larger than the index as an operator between $\mathcal{X}_{\sigma}^{\tilde{r},l}$ and $\overline{H}_{(b)}^{\tilde{r}-1,l}$.
\end{coro}
\begin{proof}
We use Lemma \ref{restrictedFredholm} with the space $B:= \overline{H}_{(b)}^{\tilde{r},l}$, $Y:=\overline{H}_{(b)}^{\tilde{r}-1,l}$, $X := \mathcal{X}_{\sigma}^{\tilde{r},l}$ since we know by Corollary \ref{coroFredholmMidFreq} that $\hat{T}_s(\sigma)$ is Fredholm from $X$ to $Y$ under the hypotheses of Corollary \ref{coroFredholmNearZero}.
\end{proof}
\begin{remark}
The advantage of the previous corollary is that, once we have proved the invertibility, estimate \eqref{estimateNearZero} translates into a uniform bound for the inverse up to $\sigma = 0$. This is not the case if $\mathcal{W}_{\sigma}^{\tilde{r},l}$ is replaced by $\mathcal{X}_{\sigma}^{\tilde{r},l}$.
\end{remark}

\subsubsection{Global semiclassical estimates}
\begin{prop}
\label{globalSemiclassicalEstimate}
Let $\eta>0$. Let  $l < -\frac{1}{2}$, $\tilde{r}+l> -\frac{1}{2}-2s - 4M\Im(\sigma)$ and $\tilde{r}>\frac{1}{2}+s-\frac{a^2+r_+^2}{r_+-M}\Im(\sigma)$.
There exits $A>0$ such that, for all $\sigma \in \C$ with $0\leq \Im(\sigma)\leq \eta$ and $\lvert \sigma \rvert\geq A$
and for all $u\in \overline{H}^{\tilde{r},l}_{(b),\lvert\sigma\rvert^{-1}}$, $v\in \dot{H}^{-\tilde{r}, 1-l}_{(b),\lvert \sigma\rvert^{-1}}$ we have:
\begin{align*}
\left\| u \right\|_{\overline{H}_{(b), \lvert\sigma\rvert^{-1}}^{\tilde{r},l}} &\leq C\left\| \hat{T}_s(\sigma)u\right\|_{\overline{H}_{(b), \lvert\sigma\rvert^{-1}}^{\tilde{r},l-1}}\\
\left\| u \right\|_{\overline{H}_{(b), \lvert\sigma\rvert^{-1}}^{\tilde{r},l}} &\leq C\left\| \hat{T}_s(\sigma)u\right\|_{\overline{H}_{(b), \lvert\sigma\rvert^{-1}}^{\tilde{r}-1,l}}\\
\left\| v \right\|_{\dot{H}_{(b), \lvert\sigma\rvert^{-1}}^{-\tilde{r}, 1-l}} &\leq C\left\| \hat{T}_s(\sigma)^*v\right\|_{\dot{H}_{(b), \lvert\sigma\rvert^{-1}}^{-\tilde{r}, -l}}\\
\end{align*}
\end{prop}
\begin{remark} \label{invertibilityHF}
In particular, we have that $\hat{T}_s(\sigma)$ is invertible in this range of $\sigma$ between $\mathcal{X}^{\tilde{r},l}_{\sigma}$ and $\overline{H}^{\tilde{r},l-1}_{(b)}$ (it is injective with closed range by the first estimate and has dense range by the third).
\end{remark}
\begin{proof}
We introduce the semiclassical parameter $h = \lvert \sigma \rvert^{-1}$, $z = h\sigma$ and the operator $\hat{T}_{s,h}(z) = h^2\hat{T}_s(h^{-1}z)$. Note that $z = z_0+O(h)$ with $z_0\in \left\{\pm 1\right\}$. We denote by $p_{h,z}$ the (semiclassical) principal symbol of $\hat{T}_{s,h}(z)$. To be concrete, we present the case $z_0 = 1$ for the first estimate (the other cases are similar).

The estimate is proved by using microlocal propagation of singularities and elliptic estimates in semiclassical second microlocal spaces on the model of \cite[Section~5]{vasy2020limiting}. The difference in our case is that the non trapping hypothesis is not satisfied. We recall the structure of the proof in \cite{vasy2020limiting} to see where this assumption comes into place and why we can replace it by our analysis of the flow (Proposition \ref{sythesisScFlow}) and the semiclassical estimates obtained previously (Propositions \ref{radialPointMicrolocal1}, \ref{radialPointMicrolocal2}, \ref{scRinEstimate}, \ref{scRoutEstimate}, \ref{estimateTrappedSet}) and the hyperbolic semiclassical estimate\footnote{Note that $\hat{T}_{s,h}(z)$ satisfies the hypotheses of proposition \ref{finalScHyperboEstimate} on $[r_+-\frac{\epsilon}{2}, +\infty)$ and $\chi u \in \mathcal{H}^{\tilde{r}}_h$ if $r$ plays the role of the time variable (and is affinely reparametrized so that $r_+-\frac{\epsilon}{2}$ correspond to time $0$ and $r_+-\epsilon$ correspond to time $T$).} (Proposition \ref{finalScHyperboEstimate}). By compactness of the semiclassical characteristic set on ${}^{sc}T^*X\cap \left\{x\leq \frac{1}{r_+-\frac{\epsilon}{2}}\right\}$, it is enough to prove a microlocal estimate on a neighborhood of each point of the semiclassical characteristic set (near a point outside of the characteristic set, the estimate is obtained by ellipticity). Note that, with the hypotheses on $\tilde{r}$ and $l$, as proved in \cite{vasy2020limiting} and restated in Proposition \ref{scRoutEstimate}, we have a source estimate near $\mathcal{R}_{\text{out}}$. Under the non trapping assumption, this can be used to initialize all the propagation of singularities estimates and the sink estimate at $\mathcal{R}_{\text{in}}$. In our case, by proposition \ref{sythesisScFlow}, for every $x\in \Sigma_{p_h}\setminus \mathcal{R}_{\text{in}}\cup K_1$, there exists a neighborhood of $V_x$ and $t\in \R$ such that $e^{tH_{p_h}}V_x \subset \mathcal{V}(\mathcal{R}_{\text{in}})$ or $\mathcal{V}(K_1)$ or $\mathcal{V}(L_+)$ or $\mathcal{V}(L_-)$ where $\mathcal{V}(A)$ is an arbitrary (but fixed) neighborhood of $A$. Therefore, a microlocal estimate of $u$ in the neighborhood of all these sets is enough to initialize propagation of singularities, sink estimate at $\mathcal{R}_{\text{in}}$ (as proved in \cite{vasy2020limiting} and restated in Proposition \ref{scRinEstimate}) and semiclassical hyperbolic estimate on $\left\{x>\frac{1}{r_+-\frac{\epsilon}{2}}\right\}$ (Proposition \ref{finalScHyperboEstimate}). With our assumptions on $\tilde{r}$ and $l$, we already have unconditional\footnote{In the sense which does not require an estimate in an other region of space time.} estimates on a neighborhood of $\mathcal{R}_{\text{out}}$ (obtained in \cite{vasy2020limiting} and restated in Proposition \ref{scRoutEstimate}), and on a neighborhood of $L_+$ and $L_-$ (see Proposition \ref{radialPointMicrolocal1}). The only remaining estimate to get is on the neighborhood of $K_1$. By proposition \ref{estimateTrappedSet}, we have an estimate which require control on $WF_h(B_0)$ (where $B_0$ is as defined in Proposition \ref{estimateTrappedSet} with $WF_h(B_0)\cap \Gamma_+=\emptyset$). Since $WF_h(B_0)\cap \Gamma_+=\emptyset$, we can get the estimate on $WF_h(B_0)$ by propagation of singularities only from the estimates on $\mathcal{R}_{\text{in}}$, $L_+$ and $L_-$ which conclude the proof.
\end{proof}

\subsection{Index zero property}\label{subSecIndexZero}
Note that the spaces $\mathcal{X}_{\sigma}^{\tilde{r},l}$ and $\mathcal{W}_{\sigma}^{\tilde{r},l}$ depend on $\sigma$, which prevent us to use directly the stability of the index for $\hat{T}_s(\sigma)$. However, we can adapt the deformation argument in the proof of Theorem 6.1 in \cite{hafner2021linear}.
We begin by the following lemma which is a generalization of the stability of invertibility to the case of operators. The proof relies on an argument given in \cite{vasy2020resolvent} (in the proof of (5.8)).
\begin{lemma}
\label{coerciveEstimateLemma}
Let $E, E_w, F, F_w$ be Banach spaces with continuous and compact inclusion $E\subset E_w$ and continuous inclusion $F\subset F_w$.  and let $(P_j)_{j\in\N}$ be a sequence of bounded operators from $E_w$ to $F_w$ such that $\lim\limits_{j\to +\infty}P_j = P_{\infty}$ in the operator norm topology. We also assume that we have uniform half Fredholm estimates for all $j\in \N \cup \left\{\infty\right\}$:
\begin{align*}
\left\|u\right\|_E\leq C\left(\left\|P_j u\right\|_F + \left\|u\right\|_{E_w}\right)
\end{align*}
in the strong sense that if the right-hand side is finite for some $u\in E_w$, then $u\in E$ and the inequality holds.
Then, if $ker(P_{\infty})\cap E = \left\{0\right\}$, there exists $N\in \N$ and $C'>0$ such that for all $n\geq N$ (and $n=\infty$),
\begin{align}
\label{coerciveEstimate}
\left\|u\right\|_E\leq C'\left\|P_n u\right\|_F
\end{align}
in the strong sense that if the right-hand side is finite for some $u\in E_w$, then $u\in E$ and the inequality holds.
\end{lemma}

\begin{remark}
In view of the half Fredholm estimate, $ker(P_{\infty})\cap E= ker(P_{\infty})$. Moreover, if $u\in E_w$ satisfies $P_n u \in F$, we get $u\in E$. Therefore, the strong character of the second estimate is a consequence of the strong character of the first.
\end{remark}

\begin{proof}
We argue by contradiction. If \eqref{coerciveEstimate} is false, there exists a sequence $(j_n)_{n\in \N}$ in $\N\cup\left\{\infty\right\}$ with $\lim\limits_{n\to +\infty}j_n = +\infty$ and $(u_n)_{n\in N}$ with $\left\|u_n\right\|_{E} = 1$ such that $\lim\limits_{n \to +\infty} P_{j_n}u_n = 0$ in $F$. There exists $v\in E$ such that $u_n \rightharpoonup v$ and by compactness of the inclusion of $E$ into $E_w$, we have strong convergence in $E_w$. By the half Fredholm estimate, we have $\left\|v\right\|_{E_w}\geq C^{-1}$ and therefore $v\neq 0$. We now show that $P_{\infty}v=0$. To do that we estimate:
\begin{align*}
\left\| P_{j_n}u_n - Pv\right\|_{F_w}\leq& \left\|(P_{j_n}-P)u_j\right\|_{F_w} + \left\|P(u_j-v)\right\|_{F_w}\\
\leq& \left\|P_{j_n}-P\right\|_{\mathcal{L}(E_w,F_w)}\left\|u_j\right\|_{E_w}+ \left\|P\right\|_{\mathcal{L}(E_w,F_w)}\left\|u_j-v\right\|_{E_w}
\end{align*}
Finally, using that $u_j$ converges to $v$ in $E_w$ (and in particular is bounded) and the convergence of $P_{j_n}$ towards $P$ in the operator norm topology, we find: $\lim\limits_{n\to +\infty} P_{j_n}u_n = Pv$ in $F_w$ since the limit is 0 in $F$, we find $Pv = 0$. Using the half Fredholm estimate, we find that $v\in E$ and this contradicts $ker(P_{\infty})\cap E = \left\{0\right\}$.
\end{proof}

As a consequence of the previous lemma, we have the following proposition:
\begin{prop}
\label{continuityIndex}
Let $E, E_w, E_s, F, F_w, F_s$ be reflexive Banach spaces such that $E_s\subset E\subset E_w$, $F_s\subset F\subset F_w$  are continuous and dense and $E\subset E_w$ and $F^*\subset F_s^*$ are compact.  Let $(P_j)_{j\in\N}$ be a sequence of bounded operators in $\mathcal{L}(E_w, F_w)$ with $(P_j)_{|_{E_s}}\in \mathcal{L}(E_s, F_s)$ such that $\lim\limits_{j\to +\infty}P_j = P_{\infty}$ in both operator norm topologies (therefore we have the same convergence property for the adjoint operators). We also assume that we have uniform Fredholm estimates for all $j\in \N \cup \left\{\infty\right\}$:
\begin{align*}
\left\|u\right\|_E\leq C\left(\left\|P_j u\right\|_F + \left\|u\right\|_{E_w}\right)\\
\left\|u\right\|_{F^*}\leq C\left(\left\|P_j^* u\right\|_{E^*} + \left\|u\right\|_{F_s^*}\right)
\end{align*}
in the strong sense that if the right-hand side is finite, so is the left hand side and the inequality holds.
Then, there exists $N\in \N$ such that for all $n\geq N$, $P_n: \left\{u\in E: P_n u \in F\right\}\rightarrow F$ and $P_\infty: \left\{u\in E: P_\infty u\in F\right\}\rightarrow F$ are Fredholm and have the same index.
\end{prop}

\begin{proof}
The fact that $P_n$ is Fredholm for any $n\in \N\cup \left\{\infty\right\}$ comes from lemma \ref{FredholmEstimate}. Then we can adapt the standard proof of stability of the index\footnote{See for example \cite[Theorem C.5]{dyatlov2019mathematical}} using lemma \ref{coerciveEstimateLemma} on the operators and on their adjoint. 
\end{proof}

As a consequence, we have that
\begin{prop}\label{indexZeroMidFreq}
Let $\sigma \in \C\setminus \left\{0\right\}$ with $\Im(\sigma)\geq 0$, $l < -\frac{1}{2}$, $\tilde{r}+l> -\frac{1}{2}-2s$ and $\tilde{r}>\frac{1}{2}+s$.
The index of $\hat{T}_s(\sigma)$ as a Fredholm operator from $\mathcal{X}_\sigma^{\tilde{r},l} \rightarrow \overline{H}^{\tilde{r},l-1}_{(b)}$ is zero.
\end{prop}
\begin{proof}
As stated in Remark \ref{invertibilityHF}, we know that $\hat{T}_s(\sigma)$ is invertible for $\Re(\sigma)$ large enough. Moreover, the index is locally constant by proposition \ref{continuityIndex} (with $E := \overline{H}_{(b)}^{\tilde{r},l}$ and $F = \overline{H}_{(b)}^{\tilde{r},l-1}$, $E_w := \overline{H}_{(b)}^{\tilde{r}-1,l-1}$, $F_s := \overline{H}_{(b)}^{\tilde{r}+1,l}$, $E_s = \overline{H}_{(b)}^{\tilde{r}+3,l+1}$ and $F_w := \overline{H}_{(b)}^{\tilde{r}-3, l-2}$).
\end{proof}
\begin{remark}\label{nonNegIndex}
By Lemma \ref{coroFredholmNearZero}, the index of $\hat{T}_s(\sigma)$ as a Fredholm operator from $\mathcal{W}_\sigma^{\tilde{r},l} \rightarrow \overline{H}^{\tilde{r},l-1}_b$ is non negative for all $\sigma \in \left\{\Im(\sigma)\geq 0, \sigma \neq 0\right\}$.
\end{remark}

\subsection{Fredholm property for \texorpdfstring{$\hat{T}_s(0)$}{hat\{T\}\_ s(0)}}\label{subSecHatTZero}
Near infinity, we have
\begin{align*}
\hat{T}_s(0) =& \left(\frac{a^2}{\Delta_r}-\frac{1}{\sin^2\theta}\right)\partial_{\phi}^2 - \frac{1}{\sin\theta}\partial_{\theta}\sin\theta\partial_\theta-\Delta_r^{-s}\partial_r \Delta_r^{s+1}\partial_r\\& + 4s(r-M)\partial_r-2s\left(\frac{a(r-M)}{\Delta_r}+i\frac{\cos\theta}{\sin^2\theta}\right)\partial_\phi + s^2\cotan^2\theta + s\\
\sigma_b(\hat{T}_s(0)) =& -\frac{a^2}{\Delta_r}\zeta_b^2 + \lvert (\zeta_b,\eta_b)\rvert^2_{\mathbb{S}^2} + \xi_b^2\\
N(\hat{T}_s(0)) =& -\Delta_{\mathbb{S}^2}+(xD_x)^2+i(2s+1)xD_x - 4isxD_x + 2s\frac{\cos\theta}{\sin^2\theta}D_\phi + s^2\cotan^2\theta +s 
\end{align*}
We rewrite it:
\begin{align*}
N(\hat{T}_s(0)) =& (xD_x)^2 +(i-2is)xD_x + \Delta^{[s]} + s
\end{align*}
where we have used the non negative spin weighted Laplacian (acting on $\mathcal{B}_s$ with spectrum \\$\left\{(l+s)(l-s)+l, l\in \lvert s\rvert + \N\right\}$):
\begin{align*}
\Delta^{[s]} := -\Delta_{\mathbb{S}^2} + \frac{2s\cos\theta}{\sin^2\theta}D_\phi + s^2\cotan^2\theta.
\end{align*}

\label{sectionNotationLapSpin}
For $j\in \lvert s\rvert+\N$, we denote by $Y_{j}$ the eigenspace of $\Delta^{[s]}$ associated with the eigenvalue $(j-s)(j+s)+j$. For any $\tilde{r}\in \R$, we denote by $\Pi_j$ the operator on $H^{\tilde{r}}(\mathcal{B}_s)$ which is the identity on $Y_{j}$ and which vanishes on all $Y_{j'}$ for $j'\neq j$ (it is well defined since $\oplus_{j'}Y_{j'}$ is dense in $H^{\tilde{s}}(\mathcal{B}_s)$ and it has norm smaller than $1$). We denote by $\mathcal{Y}_{j}^{\tilde{r},l}$ the completion of $C^\infty_c((0,\frac{1}{r_+-\epsilon}])\otimes Y_j$ in the norm $\overline{H}_b^{\tilde{r},l}$. Note that $\text{Id}\otimes\Pi_j$ (initially defined on $\overline{C}^\infty_c((0,\frac{1}{r_+}])\otimes H^{\tilde{r}}(\mathcal{B}_s))$ extends by continuity to a bounded operator from $\overline{H}_{(b)}^{\tilde{r},l}$ to $\overline{H}_{(b)}^{\tilde{r},l}((0,\frac{1}{r_+}],Y_j)$ which is the identity on $\mathcal{Y}^{\tilde{r},l}_j$ and zero on $\mathcal{Y}^{\tilde{r},l}_p$ for $p\neq j$. With a slight abuse of notation, we also denote this operator by $\Pi_j$. We define $\mathcal{Y}^{\tilde{r},l}_{\leq j} := \oplus_{p\leq j}\mathcal{Y}_p^{\tilde{r},l}$ and $\mathcal{Y}^{\tilde{r},l}_{>j}:= \oplus_{p>j}\mathcal{Y}^{\tilde{r},l}_{p}$ (Hilbert sum). For $j\in \N$, we also introduce $\mathcal{S}^{\tilde{r},l}_{\geq \left|s\right|+j}$ as the completion of $C^\infty_c((0,\frac{1}{r_+-\epsilon}))\otimes \left(\Gamma(\mathcal{B}_s)\cap \oplus_{\substack{\left|m\right|\geq \left|s\right|+j\\ m-s\in \Z}}\text{Ker}(\partial_\phi-im)\right)$ in the norm of $\overline{H}^{\tilde{r},l}_{(b)}$.

\begin{prop}
\label{FredholmPropertyAtZero}
Let $l\in (-\frac{3}{2}-s-\left|s\right|, -\frac{1}{2}-s+\left|s\right|)$ and $\tilde{r}>\frac{1}{2}+s$.
$\hat{T}_s(0)$ is Fredholm from $\left\{u\in H_{(b)}^{\tilde{r},l}: \hat{T}_s(0)u\in \overline{H}_{(b)}^{\tilde{r}-1,l}\right\}$ to $\overline{H}^{\tilde{r}-1,l}_{(b)}$. By lemma \ref{restrictedFredholm}, it is also Fredholm from $\mathcal{W}_0^{\tilde{r},l}$ to $\overline{H}^{\tilde{r},l}_{(b)}$.
\end{prop}

\begin{proof}
We fix $\frac{1}{2}+s<\tilde{r}'<\tilde{r}$ and $N$ a large integer such that $-N<1-\tilde{r}$.
First note that the estimates of Proposition \ref{combinedNearHorizon} still apply since we have not assumed $\sigma \neq 0$ in the proof. Moreover, using the principal symbol computation, we have that $\hat{T}_s(0)\in \Psi^{2,0}_b$ is elliptic in the b sense on $\left\{x\leq \frac{1}{2M}\right\}$. Therefore, for all $u\in \overline{H}^{\tilde{r}',l}_{(b)}$ and $v\in \dot{H}^{-N,l}_{(b)}$ we get the global estimates:
\begin{align*}
\left\|u\right\|_{\overline{H}_{(b)}^{\tilde{r},l}}\leq C\left(\left\|\hat{T}_s(0)u\right\|_{\overline{H}_{(b)}^{\tilde{r}-1,l}} + \left\|u\right\|_{\overline{H}_{(b)}^{\tilde{r}',l}}\right)\\
\left\|v\right\|_{\dot{H}_{(b)}^{1-\tilde{r},-l}}\leq C\left(\left\|\hat{T}_s(0)^*v\right\|_{\overline{H}_{(b)}^{-\tilde{r},-l}} + \left\|v\right\|_{\overline{H}_{(b)}^{-N,-l}}\right)
\end{align*} 
To deduce a full Fredholm estimate, we improve the decay of the error term by a normal operator argument.
To perform normal operator argument, we consider the slightly simpler conjugated operator:
\begin{align}
\tilde{N}(\hat{T}_s(0)):=& x^{-\frac{1}{2}}N(\hat{T}_s(0))x^{\frac{1}{2}} \label{directFredZero}\\
=& (xD_x)^2 -2is xD_x + \Delta^{[s]}+ \frac{1}{4}.
\end{align}
The Mellin-transform of this operator is 
\begin{align*}
M(\tau):=\tau^2 -2is \tau + \Delta^{[s]}+\frac{1}{4}.
\end{align*}
For $\Im(\tau)\in (-\frac{1}{2}+s-\left| s\right|, \frac{1}{2}+s+\left|s\right|)$, we have $\Re(\tau^2-2is\tau + \frac{1}{4})>-\left|s\right|$. In particular, the operator $M(\tau)$ is invertible on each $Y_j$ for $j\in \left|s\right|+\N$ with inverse $\left(\tau^2-2is\tau + (j-s)(j+s)+j+\frac{1}{4}\right)^{-1}$ which is bounded by $C\Re(\tau)^{-2}$ uniformly with respect to $j$ when the imaginary part of $\tau$ is fixed in the interval $(-\frac{1}{2}+s-\left| s\right|, \frac{1}{2}+s+\left|s\right|)$. 
Let $w\in \overline{H}^{-N',l'}_{(b)}$ with support in $\left\{x\leq 1\right\}$, $l'\in (-\frac{3}{2}-s-\left|s\right|, -\frac{1}{2} -s +\left|s\right|)$ and $N'\in \N$. Assume that $N(\hat{T}_s(0))w \in \overline{H}^{\tilde{r}-1, l''}_{(b)}$ with $l''\in (-\frac{3}{2}-s-\left|s\right|, -\frac{1}{2} -s +\left|s\right|)$, we have the equality for $\Im(\tau) = -l''+1$:
\begin{align*}
\mathcal{M}x^{-\frac{1}{2}}w(\tau) = \sum_{j\in \left|s\right|+\N}\left(\tau^2-2is\tau + (j+s)(j-s)+j+\frac{1}{4}\right)^{-1}\Pi_j\mathcal{M}\left(\tilde{N}(\hat{T}_s(0))x^{-\frac{1}{2}}w\right)
\end{align*}
where the sum converges in $\left<\Re(\tau)\right>^{\tilde{r}+1}L^2(\R_{\Re(\tau)}\times \mathbb{S}^2)$ and is bounded by \begin{align*}C\left<\Re(\tau)\right>^{\tilde{r}-1}\left\|\left(\tilde{N}(\hat{T}_s(0))x^{-\frac{1}{2}}\mathcal{M}w(\tau)\right)\right\|_{L^2(\R_{\Re(\tau)}\times \mathbb{S}^2)}.\end{align*} In particular, we get\footnote{See Lemma \ref{sumMode} and \ref{contourDeformationMellin} for more details about this step} that $w\in \overline{H}_{(b)}^{\tilde{r}+1,l''}$ and:
\begin{align*}
\left\|w\right\|_{\overline{H}_{(b)}^{\tilde{r}+1,l''}} \leq C \left\|N(\hat{T}_s(0))w\right\|_{\overline{H}_{(b)}^{\tilde{r}-1,l''}}.
\end{align*}

We get back to $u$ already defined. Let $\chi$ be a smooth cutoff with support in $\left\{x\leq 1\right\}$ equal to $1$ in a neighborhood of $\left\{x=0\right\}$. Since $N(\hat{T}_s(0))-\hat{T}_s(0)\in x\text{Diff}^2_b$ and $\hat{T}_s(0)u \in \overline{H}^{\tilde{r}-1,l}$, we deduce that $N(\hat{T}_s(0))\chi u \in \overline{H}^{\tilde{r}'-2,l}$. Therefore, we have:
\begin{align*}
\left\|u\right\|_{\overline{H}_{(b)}^{\tilde{r}',l}} \leq& \left\|u\right\|_{\overline{H}_{(b)}^{\tilde{r}',l-1}} +\left\|\chi u\right\|_{\overline{H}_{(b)}^{\tilde{r}',l}}\\
\leq& \left\|u\right\|_{\overline{H}_{(b)}^{\tilde{r}',l-1}}+ C\left\|N(\hat{T}_s(0))\chi u\right\|_{\overline{H}_{(b)}^{\tilde{r}'-2,l}}\\ 
\leq& C'\left(\left\|u\right\|_{\overline{H}_{(b)}^{\tilde{r}',l-1}}+ \left\|\hat{T}_s(0)\chi u\right\|_{\overline{H}_{(b)}^{\tilde{r}'-2,l}}\right)\\
\leq& C''\left(\left\|u\right\|_{\overline{H}_{(b)}^{\tilde{r}',l-1}}+\left\|\hat{T}_s(0)\chi u\right\|_{\overline{H}_{(b)}^{\tilde{r}-1,l}}\right).
\end{align*}
Using this in estimate \eqref{directFredZero}, we get
\begin{align*}
\left\|u\right\|_{\overline{H}^{\tilde{r},l}_{(b)}}\leq C\left(\left\|\hat{T}_s(0)u\right\|_{\overline{H}^{\tilde{r}-1,l}_{(b)}}+ \left\|u\right\|_{\overline{H}^{\tilde{r}',l-1}_{(b)}}\right).
\end{align*}
We can do a similar normal operator argument for the adjoint $\hat{T}_s(0)^*$ (but this time, we find that the decay index has to belong to $(\frac{1}{2}+s-\left|s\right|, \frac{3}{2}+s+\left|s\right|)$ ) and we get an estimate of the form:
\begin{align*}
\left\|v\right\|_{\dot{H}^{1-\tilde{r},-l}_{(b)}}\leq C\left(\left\|\hat{T}^*(0)v\right\|_{\overline{H}^{-\tilde{r},-l}_{(b)}}+ \left\|v\right\|_{\dot{H}^{-\tilde{r},-l-1}_{(b)}}\right).
\end{align*}
Finally, by Lemma \ref{FredholmEstimate}, we get the proposition.
\end{proof}

We now describe precisely the kernel of $\hat{T}_s(0)$ and $\hat{T}_s(0)^*$. To do so, it is easier to write $\hat{T}_s(0)$ in $(r,\theta,\phi_*)$ coordinates (which correspond to the expression of the operator near the horizon):
\begin{align*}
\hat{T}_s(0) = -\Delta_r\partial_r^2 + 2((r-M)(s-1)-a\partial_{\phi_*})\partial_r+\Delta^{[s]} + s\\
r^2\hat{T}_s(0)^*r^{-2}= -\Delta_r\partial_r^2 + 2((r-M)(-s-1)-a\partial_{\phi_*})\partial_r+\Delta^{[s]} - s
\end{align*}
Where the adjoint is taken with respect to the density $r^2\dd r\dd \phi^* \dd \theta$.
We for $l\in \lvert s\rvert+\N$ and $-l\leq m\leq l$ (with $m-s\in \Z$), we denote by $f_{l,m}$ the normalized spin weighted spherical harmonics associated to the eigenvalue $(l+s)(l-s)+l$ and such that $\partial_{\phi_*} f_{l,m} = imf_{l,m}$. These functions form a Hilbert basis of $L^2(\mathcal{B}_s)$.

Following \cite[Section 9]{olver1997asymptotics}, we introduce the hypergeometric function 
\begin{align*}\mathbf{F}(a,b,c, z):= \sum_{k=0}^{+\infty}\frac{(a)_k(b)_k}{\Gamma(c+k)}\frac{z^k}{k!}
\end{align*} where $(a)_k := \prod_{j=0}^{k-1}(a+j)$.

\begin{prop}
\label{preciseComputationOfTheKernel}
We have:
\begin{align*}
\text{ker}(\hat{T}_s(0))\cap \overline{H}^{s+\frac{1}{2},q}_{(b)} =& \bigoplus_{\substack{\lvert s\rvert \leq l <-\frac{3}{2}-s-q \\ \lvert m\rvert \leq l} }\C u_{l,m}\\
\text{ker}(\hat{T}_s(0)^*)\cap \dot{H}^{\left(-s+\frac{1}{2}\right)-,-q}_{(b)} r^2\dd r\dd \phi^*\dd \theta =& \bigoplus_{\substack{\lvert s\rvert \leq l < \frac{1}{2}+s+q\\ \lvert m\rvert \leq l }}\C \overline{u^*_{l,m}}
\end{align*}
where 
\[
\begin{array}{|c|c|}
\hline
& \text{ If } m\neq 0 \text{ or } s\notin \Z \\
\hline
u_{l,m} & \mathbf{F}\left(-l-s, 1+l-s, 1-s+\frac{2iam}{r_+-r_-}, \frac{r_+ - r}{r_+ -r_-}\right)f_{l,m}(\theta,\phi_*) \\
\hline
u^*_{l,m} & \frac{(r-r_+)^{-s-\frac{2iam}{r_+-r_-}}_+}{(r_+-r_-)^{-s-\frac{2iam}{r_+-r_-}}}\mathbf{F}\left(-l-\frac{2iam}{r_+-r_-}, l-\frac{2iam}{r_+-r_-}, 1-s-\frac{2iam}{r_+-r_-}, \frac{r_+-r}{r_+-r_-}\right)f_{l,m}(\theta,\phi_*)\dd r\dd \phi^* \dd \theta \\
\hline
\end{array}\]
\[
\begin{array}{|c|c|}
\hline
& \text{ If } m=0 \text{ and } s\in -1-\N\\
\hline
u_{l,m} & \mathbf{F}\left(-l-s, 1+l-s, 1-s, \frac{r_+ - r}{r_+ -r_-}\right)f_{l,m}(\theta,\phi_*)\\
\hline
u_{l,m}^* & \mathbb{1}_{(r_+, +\infty)} \left(\frac{r-r_+}{r_+-r_-}\right)^{-s}\left(\frac{r-r_-}{r_+-r_-}\right)^{-s}\mathbf{F}\left(-l-s, 1+l-s, 1-s, \frac{r_+-r}{r_+-r_-}\right)f_{l,m}(\theta,\phi_*) \dd r\dd \phi_*\dd \theta \\
\hline
\end{array}\]
\[ \begin{array}{|c|c|}
\hline
& \text{ If } m=0 \text{ and } s\in \N \\
\hline
u_{l,m} & \left(\frac{r-r_+}{r_+-r_-}\right)^{s}\left(\frac{r-r_-}{r_+-r_-}\right)^s\mathbf{F}\left(-l+s, 1+l+s, 1+s, \frac{r_+-r}{r_+-r_-}\right)f_{l,m}(\theta,\phi_*)\\
\hline
u_{l,m}^* & \left(\frac{(l+s)!}{(l-s)!}(r_+-r_-)^{-s}\mathbb{1}_{(r_+,+\infty)}\mathbf{F}\left(-l+s, 1+l+s, 1+s, \frac{r_+ - r}{r_+ -r_-}\right)\right.\\
& \left.+ \sum_{j=0}^{s-1}\frac{(l+s-1-j)!}{(l-s+1+j)!(s-1-j)!}(r_+-r_-)^{-s+1+j}\delta_{r_+}^{(j)}\right)f_{l,m}(\theta,\phi_*)\dd r\dd \phi^* \dd \theta\\
\hline
\end{array}
\]
\end{prop}

\begin{proof}
First, we decompose $u\in \overline{H}^{s+\frac{1}{2},q}_{(b)}$ on $\oplus_{l\geq \lvert s\rvert, -l\leq m\leq l}\overline{H}^{s+\frac{1}{2},q}_{(b)}([r_+-\epsilon, +\infty))\otimes f_{l,m}$ (we denote by $u_{l,m}$ the element of $\overline{H}^{s+\frac{1}{2},q}_{(b)}([r_+-\epsilon, +\infty))$ in this decomposition). 
If $u\in \text{ker}(\hat{T}_s(0))$, we have 
\begin{align*}
\hat{T}(0)_{l,m,s} u_{l,m} :=& (-\Delta_r\partial_r^2 + 2((r-M)(s-1)-aim)\partial_r+ (l-s+1)(l+s)) u_{l,m}\\
 =& 0.
\end{align*}
 This is a hypergeometric equation and we can put it in canonical form (see \cite{olver1997asymptotics} chapter 5 section 8.1) by the change of variable: $z = \frac{r-r_-}{r_+-r_-}$. With this new variable, we have:
\begin{align}
\label{normalHypergeoEq}
\hat{T}(0)_{l,m,s} = z(1-z)\partial_z^2 + \left(1-s-\frac{2iam}{r_+-r_-}-(2-2s)z\right)\partial_z + (l-s+1)(l+s)
\end{align}
We define $\alpha := -l-s$, $\beta := l-s+1$ and $\gamma := 1-s-\frac{2iam}{r_+-r_-}$ which are the coefficient of the canonical hypergeometric equation (correspond to coefficients $a,b,c$ in \cite{olver1997asymptotics}.
Note that we can in the same way decompose $v\in \dot{H}^{\tilde{r}, -q-2}$ and if we assume that $v\dd r\dd \phi^*\dd \theta\in \text{Ker}(T(0)^*)$, we get that for all $l\in \lvert s\rvert +\N$ and $m \in \Z+l$ with $\lvert m\rvert \leq l$:
\begin{align*}
\hat{T}(0)_{l,m,-s}v_{l,m} = 0
\end{align*}
Therefore, we are reduced to find for all $s\in \frac{1}{2}\Z$ and $l,m$ as above $\text{Ker}(\hat{T}(0)_{l,m,s})\cap \overline{H}_{(b)}^{\frac{1}{2}+s, q}$ and $\text{Ker}(\hat{T}(0)_{l,m,s})\cap \dot{H}_{(b)}^{\tilde{r}, -q-2}$. Therefore, we fix $s\in \frac{1}{2}\Z$,  $l\in \lvert s\rvert +\N$ and $\lvert m\rvert \leq l$, $u_{m,l} \in \text{Ker}(\hat{T}(0)_{l,m,s})\cap \overline{H}_{(b)}^{\frac{1}{2}+s, q}$ and $v_{l,m}\in \text{Ker}(\hat{T}(0)_{l,m,s})\cap \dot{H}_{(b)}^{\tilde{r}, -q-2}$.

We have three cases (because of possible integer coincidence):
\begin{itemize}
\item If $s\notin \Z$ or $m\neq 0$. In this case $1-\gamma \notin \Z$ and $1+\gamma-\alpha-\beta \notin \Z$.
We have two solutions of the equation $\hat{T}(0)_{l,m,s}u = 0$: 
\begin{align*}u_1:= \mathbf{F}(\alpha, \beta, \alpha+\beta+1-\gamma, 1-z)
\end{align*} which is an analytic solution on $\C\setminus (-\infty, 0]$ and equal to $\frac{1}{\Gamma(\alpha+\beta+1-\gamma)}$ at $z=1$ and 
\begin{align*}u_2 := (z-1)^{\gamma-\alpha-\beta}\mathbf{F}(\gamma-\beta, \gamma-\alpha, 1+\gamma-\alpha-\beta, 1-z)
\end{align*} (where we use the continuous branch of the logarithm defined on $\C\setminus (-\infty, 0]$ which is the usual logarithm on $(0,+\infty)$) which is an analytic solution on $\C\setminus (-\infty, 1]$. Note that $u_1$ and $u_2$ are linearly independent. By the study of the regular singularity at infinity (by Frobenius method), we see that the indicial roots are $\alpha$ and $\beta$ (see for example \cite{olver1997asymptotics} chapter 5, section 5). Since $\beta-\alpha \in 1+\N$, we have an integer coincidence and we can find a connection formula analogous to (10.13) in \cite{olver1997asymptotics}. It gives the following form for $u_1$ and $u_2$ when $z\in \C\setminus (-\infty,1]$:
\begin{align}
u_1(z) =&\frac{(\beta-\alpha-1)!}{\Gamma(1-\gamma+\beta)\Gamma(\beta)}z^{-\alpha}H(z^{-1})+ c_1z^{-\beta}\ln(z)G(z^{-1}) \label{asymptoticsUOne}\\
\label{asymptoticsUTwo}
u_2(z) =& e^{i\pi(\gamma-\alpha-\beta)}\frac{(\beta-\alpha-1)!}{\Gamma(\gamma-\alpha)\Gamma(1-\alpha)}z^{-\alpha}H(z^{-1}) + c_2z^{-\beta}\ln(z)G(z^{-1})
\end{align}
where $H$ is holomorphic on $\C\setminus [1, +\infty)$ and $H(0) = 1$ and $G$ is holomorphic on $\C\setminus [1,+\infty)$ and $c_1,c_2$ are constants. In particular since $\frac{(\beta-\alpha-1)!}{\Gamma(1-\gamma+\beta)\Gamma(\beta)}z^{-\alpha}\neq 0$ we have $u_1(z)\sim \frac{(\beta-\alpha-1)!}{\Gamma(1-\gamma+\beta)\Gamma(\beta)}z^{-\alpha}$ when $z\to +\infty$.
We know that the restriction of $(u_{l,m})_{|_{(r_+,+\infty)}} = \lambda_1 u_1\left(\frac{r-r-}{r_+-r_-}\right) + \lambda_2 u_2\left(\frac{r-r_-}{r_+-r_-}\right)$. In particular, $u_{l,m}-\lambda_1 u_1\left(\frac{r-r_-}{r_+-r_-}\right)$ is an extension of $\lambda_2 u_2\left(\frac{r-r_-}{r_+-r_-}\right)$ in $\overline{H}_{(b)}^{s+\frac{1}{2}, loc}(r_+-\epsilon, +\infty)$. But this is impossible except if $\lambda_2 = 0$. Indeed, $u_2\left(\frac{r-r_-}{r_+-r_-}\right) = \left(\frac{r-r_+}{r_+-r_-}\right)^{\gamma-\alpha-\beta}\Gamma(1+\gamma-\alpha-\beta)^{-1}+  \left(\frac{r-r_+}{r_+-r_-}\right)^{\gamma-\alpha-\beta+1}F(r-r_+)$ with $F$ extending smoothly on $[r_+-\epsilon, +\infty)$. In particular, since $\left(\frac{r-r_++i0}{r_+-r_-}\right)^{\gamma-\alpha-\beta+1}F(r-r_+) \in \overline{H}_{(b)}^{(\frac{3}{2}+s)-, loc}$, $u_2$ has an extension in $\overline{H}_{(b)}^{\frac{1}{2}+s, loc}$ if and only if $\left(\frac{r-r_+}{r_+-r_-}\right)^{\gamma-\alpha-\beta}\Gamma(1+\gamma-\alpha-\beta)^{-1}$ has a such an extension which is not the case.
Therefore, $(u_{l,m})_{|_{(r_+,+\infty)}} = \lambda_1 u_1$ and a similar argument shows that $(u_{l,m})_{|_{(r_+-\epsilon, r_+)}} = \lambda_1' u_1$. We deduce that there exists $a_0, ..., a_k$ such that $v_1:=(\lambda_1-\lambda_1')\mathbb{1}_{(r_+,+\infty)}u_1 + \sum_{j=0}^k a_j\delta_{r_+}^{(k)}$ is a solution in $\overline{H}_{(b)}^{s+\frac{1}{2}}$. Using the equation, we have:
\begin{align*}
\hat{T}(0)_{m,l,s} v_1 =& 2(\lambda_1-\lambda_1')u_1(r_+)\left((r_+-M)s-aim\right)\delta_{r_+}\\
& + \sum_{j=0}^k \left(2a_j((r_+-M)(s+j+1)-aim)\delta_{r_+}^{(j+1)}-a_j((j+2)(j+1)\right.\\
&\hphantom{+\sum_{j=0}^k}\left.+2(s-1)(j+1)-(l-s+1)(l+s))\delta_{r_+}^{(j)}\right)
\end{align*}
Therefore, the coefficient in front of each $\delta$ term has to vanish, in particular the coefficient of $\delta_{r_+}^{(k+1)}$:
\begin{align*}
2a_k((r_+-M)(s+k+1)-aim) = 0
\end{align*}
Since we assumed that $m \neq 0$ or $s\notin \Z$, we get $a_k=0$. Repeating the argument we have $a_0 = ... = a_k = 0$. Finally, the vanishing of the coefficient of $\delta_{r_+}$ gives $\lambda_1 = \lambda_1'$. We conclude that $u_{l,m} \in \C u_1\left(\frac{r-r_-}{r_+-r_-}\right)$. However, in view of \eqref{asymptoticsUOne}, we have $u_1 \in \overline{H}_{(b)}^{\frac{1}{2}+s, (-\frac{3}{2}-l -s)-}$ (and no better in term of decay). Therefore, $u_{l,m}$ can be non zero if and only if $q<-\frac{3}{2}-l-s$.

We now determine the form of $v_{l,m}$. Since $v_{l,m}$ is a distribution supported on $[r_+-\epsilon, +\infty)$ and is the solution to a non degenerate differential equation on $(r_-,r_+)$, it has to be supported on $[r_+, +\infty)$. We also have that $(v_{l,m})_{|_{(r_+,+\infty)}} = \lambda_1 u_1 + \lambda_2 u_2$. Since $u_2\left(\frac{r-r_+}{r_+-r_-}\right) = (r-r_+)^{\gamma-\alpha-\beta}F(z)$ where $F$ is smooth near $r-r_+$ and $\gamma-\alpha-\beta\notin \Z$, we have that $\tilde{u}_2(r-r_+)_+^{\gamma-\alpha-\beta}F(z)$ (where $x_+^{\gamma-\alpha-\beta}$ denote the only homogeneous distribution supported in $[0,+\infty)$ whose restriction to $(0,+\infty)$ is $x^{\gamma-\alpha-\beta}$)  is in $ker(\hat{T}(0)_{l,m,s})\cap \dot{H}_{(b)}^{(\frac{1}{2}+s)-, loc}$. Therefore, $v_{l,m} - \lambda_2 \tilde{u}_2  = \lambda_1 \mathbb{1}_{[0,+\infty)}u_1\left(\frac{r-r_+}{r_+-r_-}\right) + \sum_{j=0}^k c_k\delta^{(j)}_{r_+} \in ker(\hat{T}(0)_{l,m,s})$. But we have seen that in this case $\lambda_1 = a_0 = ...=a_k =0$. Therefore, 
\begin{align*}
v_{l,m} = \lambda_2 \tilde{u}_2\left(\frac{r-r_+}{r_+-r_-}\right).
\end{align*} Using \eqref{asymptoticsUOne}, we see that $\tilde{u}_2 \in \dot{H}_{(b)}^{(\frac{1}{2}+s)-, -\frac{3}{2} -s -l}$ and therefore $\lambda_2$ can be non zero if and only if $-q-2<-\frac{3}{2}-s-l$.

\item In the case of $m=0$ and $s\in -\N$, we have $\gamma-\alpha-\beta = s\in -\N$. 
A basis of solutions for \eqref{normalHypergeoEq} is given by:
\begin{align*}
u_1(z):=&\mathbf{F}(\alpha, \beta, 1-s, 1-z)\\
u_2(z):=&\tilde{\mathbf{H}}(\alpha, \beta, 1-s, 1-z) \\
:=& (-1)^{s}(\partial_c)_{|_{c=1-s}}\left( \vphantom{\prod_{j=1}^{-s-1}}(z-1)^{1-c}\mathbf{F}(\alpha+1-c, \beta+1-c, 2-c,1-z)\right. \\
&\hphantom{(-1)^{s}(\partial_c)_{|_{c=1-s}}}\left.- e^{i\pi(1-c)}\prod_{j=1}^{-s-1}(j-\alpha)(j-\beta)\mathbf{F}(\alpha,\beta, c,1-z)\right)
\end{align*} 
$u_1$ is holomorphic on $\C\setminus(-\infty, 0]$ and $u_1(1) = \frac{1}{(-s)!}$. $u_2$ is holomorphic on $\C\setminus(-\infty, 1]$ and near $z=1$, we have: 
\begin{align*}
u_2 = \ln(z-1) F(z-1) + (1-z)^{s}G(z-1)
\end{align*} with $F$ and $G$ holomorphic near $0$, $F(0) = -\frac{(\alpha+s)_{-s}(\beta+s)_{-s}}{(-s)!}$ (where $(x)_{k} := \prod_{j=0}^{k-1}(x+j)$) and if $-s-1\geq 0$, $G(0) = (-1)^{-s}(-s-1)!$.
We still have that for $z\in \C\setminus (-\infty,1]$:

\begin{align}
\label{devUOneInteg}
u_1(z) =& \frac{z^{l+s}(2l)!}{l!(l-s)!}G(z^{-1}) + z^{-1-l+s}\ln(z)F(z^{-1})\\
\label{devUTwoInteg}
u_2(z) =& \frac{z^{l+s}(2l)!\left(\Psi(1+l+s)+\Psi(1+l)\right)}{l!(l+s)!}G(z^{-1})+z^{-1-l+s}\ln(z)F(z^{-1})
\end{align}
with $F$ and $G$ holomorphic near zero with $G(0)=1$ and $\Psi$ denote the digamma function (see \cite[Chapter~2, Section~2]{olver1997asymptotics}).

As before, we use the fact that $u_2\left(\frac{r-r_-}{r_+-r_-}\right)$ has no extension in $\overline{H}_{(b)}^{\frac{1}{2}+s, loc}((r_+-\epsilon, +\infty))$ to conclude that there exists $\lambda \in \C$ such that $v_1:= u_{l,m} - \lambda u_1\left(\frac{r-r_-}{r_+-r_-}\right) = \frac{a_{-1}}{u_1(r_+)}\mathbb{1}_{(r_+-\epsilon,+\infty)}u_1 + \sum_{j=0}^k a_j\delta_{r_+}^{(j)}$.
First note that since $v_1 = u_{l,m} - \lambda u_1\in \overline{H}_{(b)}^{\frac{1}{2}+s,loc}$, only the $a_j$ with $-1\leq j<-1-s$ can be non zero therefore we can assume that $k\leq-1-s$. Moreover we have:
\begin{align}\label{eq:v1}
\hat{T}(0)_{m,l,s} v_1 = & \sum_{j=0}^k\left( 2a_j(r_+-M)(s+j+1)\delta_{r_+}^{(j+1)}\right.\\
&\hphantom{\sum_{j=0}^k}\left.-a_j((j+2)(j+1)+2(s-1)(j+1)-(l-s+1)(l+s))\delta_{r_+}^{(j)}\right)\\
& + 2a_{-1}(r_+-M)s\delta_{r_+}
\end{align} and since $\hat{T}(0)_{l,m,s}v_1 = 0$ we deduce that:
\begin{align*}
2a_k(r_+-M)(s+k+1) = 0
\end{align*}
and therefore, $a_k=0$. By induction, we get that $v_1= 0$ and $u_{l,m} \in \C u_1\left(\frac{r-r_-}{r_+-r_-}\right)$. Finally, using \eqref{devUOneInteg} we deduce that $u_{l,m}$ can be non zero if and only if $q<-\frac{3}{2}-l-s$.

We now consider $v_{l,m}$. As before, $(v_{l,m})_{|_{(r_+-\epsilon, r_+)}} = 0$ and $(v_{l,m})_{|_{(r_+, +\infty)}} = \lambda_1 u_1 + \lambda_2 u_2$. We define $a_{-s-1} :=1$ and for all $-1\leq j\leq -s-2$:
\begin{align*}
a_{j} = \frac{a_{j+1}((j+3)(j+2)+2(s-1)(j+2)-(l-s+1)(l+s))}{2(r_+-M)(s+j+1)}.
\end{align*}
Note that if $s=0$, we have only one term $a_{-1} = 1$ in the sequence and if $s<0$, using that $((j+3)(j+2)+2(s-1)(j+2)-(l-s+1)(l+s)) = (j+2 + l+s)(j+1 -l+s)$ we see that $a_j$ never vanishes for $-1\leq j\leq -s-2$ and therefore $a_{-1}\neq 0$. We can solve explicitly the recurrence relation and we find:
\begin{align*}
a_{j} = \frac{(l-s-1-j)!}{(l+s+1+j)(-s-1-j)!}(r_+-r_-)^{s+1+j}
\end{align*}
We have (according to \eqref{eq:v1}) $\tilde{u}_1 :=  \frac{a_{-1}}{u_1(r_+)}\mathbb{1}_{(r_+,+\infty)}u_1 + \sum_{j=0}^{-s-1} a_j\delta_{r_+}^{(k)}$ is a solution in $\dot{H}^{(\frac{1}{2}+s)-, loc}_{(b)}$. As a consequence $v_1:= v_{l,m} - \frac{\lambda_1 u_1(r_+)}{a_{-1}}\tilde{u}_1$ is of the form $\lambda_2 \tilde{u}_2 + \sum_{j=0}^k d_j \delta^{(j)}_{r_+}$ where $\tilde{u}_2$ is an extension of $u_2$. We introduce the distribution $\ln(x)_+$ which is the $L^1_{loc}$ function equal to zero on $(-\infty,0)$ and to $\ln(x)$ on $(0,+\infty)$.
Since $u_2(z) = \ln(z-1)F(z-1) + (z-1)^{s}G(z-1)$ with $F$ and $G$ holomorphic near zero, we can choose $\tilde{u}_2$ of the form $\ln\left(\frac{r-r_+}{r_+-r_-}\right)_+ F\left(\frac{r-r_+}{r_+-r_-}\right)+ \frac{(r-r_+)^{s}_+}{(r_+-r_-)^s}G\left(\frac{r-r_+}{r_+-r_-}\right)$ (where $(r-r_+)^{s}_+$ is defined for example in \cite{hormanderI} equation (3.2.5)) if $-s\geq 1$ and of the form $\ln\left(\frac{r-r_+}{r_+-r_-}\right)_+ F\left(\frac{r-r_+}{r_+-r_-}\right)+\mathbb{1}_{(r_+, +\infty)}G\left(\frac{r-r_+}{r_+-r_-}\right)$ if $s=0$. When $s=0$ we can compute that $\hat{T}(0)_{l,m,s} \tilde{u}_2 = -(r_+-r_-)F(0)\delta_{r_+}$ (using for example that in the sense of distributions $\tilde{u}_2 = \lim\limits_{\epsilon\to 0} \mathbb{1}_{(r_++\epsilon, +\infty)} u_2$) we show as before that every coefficient $d_j = 0$, then $\lambda_2 = 0$ and finally $v_1=0$. If $s<0$, we consider $(r-r_+)^{-s}u_2$ which is of the form $\lim\limits_{\epsilon \to 0}\mathbb{1}_{(r_++\epsilon, +\infty)}u_3$ where $u_3$ is a solution of $(r-r_+)^{-s}\hat{T}(0)_{l,m,s}(r-r_+)^{s} u_3 = 0$.
We compute 
\begin{align*}
&(r-r_+)^{-s}\hat{T}(0)_{l,m,s}(r-r_+)^s = -\Delta_r \partial_r^2 + 2((r_--M)s+M-r)\partial_r + (l-s)(l+s)+l+s^2\\
&\lim\limits_{\epsilon \to 0}[(r-r_+)^{-s}\hat{T}(0)_{l,m,s}(r-r_+)^s, \mathbb{1}_{(r_++\epsilon, +\infty)}]u_3 = -(r_+-r_-)^{1-s}G(0)\delta_{r_+}
\end{align*}
and therefore $(r-r_+)^{-s}\hat{T}(0)_{l,m,s} \tilde{u}_2 = -(r_+-r_-)^{1-s}G(0)\delta_{r_+}$.
On the other hand, we know that $\hat{T}(0)_{l,m,s} \tilde{u}_2$ is of the form $\sum_{j=0}^{N} c_j \delta_{r_+}^{(j)}$. We deduce that $N= -s$ and $c_{-s} = \frac{(-1)^{1-s}(r_+-r_-)^{1-s}}{(-s)!}G(0)$. Finally, using that $G(0)\neq 0$ we deduce that $v_1 = 0$. Therefore, we conclude in both cases that $v_{l,m} \in \C \tilde{u}_1$. Using \eqref{devUOneInteg}, we see that $v_{l,m}$ can be a non trivial multiple of $\tilde{u}_1$ if and only if $-q-2<-\frac{3}{2}-l-s$.

\item In the case $m=0$, $s\in 1+\N$:
We have the two independent solutions:
\begin{align*}
u_1 = (z-1)^sz^s \mathbf{F}(-l+s, 1+l+s, 1+s, 1-z)\\
u_2 = (z-1)^sz^s\tilde{\mathbf{H}}(-l+s, 1+l+s, 1+s, 1-z)
\end{align*}
Where $\tilde{\mathbf{H}}$ is the function defined in the previous case. Note that 
\begin{align*}
u_2\left(\frac{r-r_-}{r_+-r_-}\right) = -\frac{(-l)_s(1+l)_s}{s!}\left(\frac{r-r_+}{r_+-r_-}\right)^{s}\ln\left(\frac{r-r_+}{r_+-r_-}\right) + \overline{H}_{(b)}^{\frac{1}{2}+s, loc}
\end{align*} and therefore has no extension in $\overline{H}_{(b)}^{\frac{1}{2}+s, loc}$. It is even easier to conclude that $u_{l,m} \in \C u_1\left(\frac{r-r_-}{r_+-r_t}\right)$ since $u_{l,m}$ is in particular continuous by Sobolev embedding. Since by \eqref{devUOneInteg}, $u_1\left(\frac{r-r_-}{r_+-r_-}\right)\in \overline{H}_{(b)}^{\infty, (-\frac{3}{2}-l-s)-}$, we deduce that $u_{l,m}$ can be a non zero multiple of $u_1$ if and only if $q<-\frac{3}{2}-l-s$.

We now consider $v_{l,m}$. As before, $v_{l,m} = 0$ on $(r_+-\epsilon, r_+)$ and $v_{l,m} = \lambda_1 u_1 + \lambda_2 u_2$ on $(r_+, +\infty)$. Noting that $\tilde{u}_1 :=\mathbb{1}_{(r_+, +\infty)} u_1$ is a distribution solution of the equation, we deduce that $v_1 := v_{l,m}-\lambda_1 \tilde{u}_1$ is also a solution of the equation. It is of the form $\lambda{u}_2 +\sum_{j=0}^N c_j \delta_{r_+}^{(j)}$ where $\tilde{u}_2 := \mathbb{1}_{(r_+,+\infty)} u_2$ (which is in $L^1_{loc}$). Using that $\tilde{u}_2 = \lim\limits_{\epsilon \to 0} \mathbb{1}_{(r_++\epsilon,+\infty)} u_2$ and the expression of $\tilde{\mathbf{H}}(-l+s, 1+l+s, 1+s, 1-z)$ near $z=1$, we compute: 
\begin{align*}
\hat{T}(0)_{l,m,s} \tilde{u}_2 =& (-1)^ss!(r_+-r_-)\delta_{r_+}\\
\hat{T}(0)_{l,m,s} \left(\lambda_2\tilde{u}_2 + \sum_{j=0}^N c_j \delta_{r_+}^{(j)}\right) =& \lambda_2 (-1)^ss!(r_+-r_-)\delta_{r_+}+ \sum_{j=0}^N 2c_j(r_+-M)(s+j+1)\delta_{r_+}^{(j+1)}\\
&-\sum_{j=0}^N c_j(j+2 + l+s)(j+1 -l+s)\delta_{r_+}^{(j)}
\end{align*}
and we prove by induction that $c_N = ... = c_0 = 0$ and $\lambda_2 = 0$. As a consequence, $v_{l,m} \in \C \tilde{u}_1$ and it can be non zero if and only if $-q-2 < -\frac{3}{2}-l-s$.

\end{itemize}
\end{proof}

\section{Existence, boundedness and regularity of the resolvent} \label{secAnalRes}
\begin{prop}
\label{invertibility}
If $l < -\frac{1}{2}$, $\tilde{r}+l> -\frac{1}{2}-2s - 4M\Im(\sigma)$ and $\tilde{r}>\frac{1}{2}+s-\frac{a^2+r_+^2}{r_+-M}\Im(\sigma)$, then the operator $\hat{T}_s(\sigma)$ is invertible from $\mathcal{X}^{\tilde{r},l}_\sigma$ to $\overline{H}^{\tilde{r},l-1}$ when $\Im(\sigma)\geq 0$ and  $\sigma\neq 0$.
If 
\begin{align*}
-\frac{3}{2}-s-\left|s\right|&< l < -\frac{1}{2}, \\
\tilde{r}+l &> -\frac{1}{2}-2s,\\
\tilde{r} &> \frac{1}{2}+s,
\end{align*}
the operator $\hat{T}_s(\sigma)$ is invertible from $\mathcal{W}^{\tilde{r},l}_\sigma$ to $\overline{H}^{\tilde{r},l}$ when $\Im(\sigma)\geq 0$.
\end{prop}
\begin{proof}
We have already proved that the operators are Fredholm with non negative index (see Proposition \ref{indexZeroMidFreq} and Remark \ref{nonNegIndex}). Therefore, it remains to prove that their kernel is trivial. 
This is a consequence of \cite{whiting1989mode} for $\Im(\sigma)>0$, of the generalization \cite{andersson2017mode} for $\Im(\sigma)=0$, $\sigma\neq 0$. For the precise proof that the absence of mode for the radial equation imply the triviality of the kernel in our case, see the proof of case (1) in Theorem 4.5 in \cite{andersson2022mode}. The case $\sigma=0$ is a consequence of the explicit computation of the kernel and cokernel of $\hat{T}_s(0)$ given in Proposition \ref{preciseComputationOfTheKernel}, see also \cite{andersson2022mode}.
\end{proof}

As a consequence of the previous proposition, we can define the resolvent operator $R(\sigma) = \hat{T}_s(\sigma)^{-1}$ for $\Im(\sigma)\geq 0$. The following property summarize its basic properties.
\begin{prop}
\label{defResolvent}
For every $\eta\in [0,1]$ $R(\sigma)$ is a bounded operator from $\overline{H}^{\tilde{r},l}_{(b)}$ to $\overline{H}^{\tilde{r},l+1-\eta}_{(b)}$ for $\Im(\sigma)\geq 0$, $\sigma\neq 0$, $\lvert \sigma \rvert \leq c$ and
\begin{align*}
-\frac{3}{2}-s-\lvert s\rvert&< l+1-\eta < -\frac{1}{2}\\
\tilde{r}+l+1-\eta &> -\frac{1}{2}-2s\\
\tilde{r} &> \frac{1}{2}+s.
\end{align*}
Moreover, in this case, we have the following bound (uniform in $\lvert \sigma\rvert \leq c$ for $c$ small enough):
\begin{align}
\label{lowFreqUniform}
\left\|R(\sigma)\right\|_{\mathcal{L}(\overline{H}_{(b)}^{\tilde{r},l}, \overline{H}_{(b)}^{\tilde{r},l+1-\eta})}\leq C\vert\sigma\rvert^{\eta-1}
\end{align}
It is also a bounded operator from $\overline{H}^{\tilde{r},l}_{(b)}$ to $\overline{H}^{\tilde{r},l+1}_{(b)}$ for $\sigma\neq 0$ and
 $l+1 < -\frac{1}{2}$, $\tilde{r}+l+1> -\frac{1}{2}-2s - 4M\Im(\sigma)$ and $\tilde{r}>\frac{1}{2}+s-\frac{a^2+r_+^2}{r_+-M}\Im(\sigma)$ and in this case we have the bound (uniform for $\sigma$ in a strip $\left\{ 0\leq\Im(\sigma)\leq A, \lvert\sigma\rvert>\frac{1}{A}\right\}$:
\begin{align}
\label{highFreqInStrip}
\left\|R(\sigma)\right\|_{\mathcal{L}(\overline{H}_{(b),\lvert\sigma\rvert^{-1}}^{\tilde{r},l}, \overline{H}_{(b), \lvert\sigma\rvert^{-1}}^{\tilde{r},l+1})}\leq C
\end{align}
\end{prop}

\begin{proof}
It is a consequence of proposition \ref{invertibility} and of the estimates of proposition \ref{estimateNearZero} (where we can drop the error term since the kernel is trivial) and proposition \ref{globalSemiclassicalEstimate}.
\end{proof}

\begin{remark}
Note that for $\sigma,\sigma' \in \C\setminus\left\{0\right\}$ with $\Im(\sigma),\Im(\sigma')\geq 0$, we have $\mathcal{X}^{\tilde{r}+1,l}_{\sigma}\subset \mathcal{X}^{\tilde{r},l}_{\sigma'}$. This is due to the fact that $\hat{T}_s(\sigma)-\hat{T}_s(\sigma') \in x^{-1}\text{Diff}^1_b$.
The following identities are useful in the computations:
\begin{enumerate}[label = (I\arabic*)]
\item \label{firstId} If
\begin{align*}
l+1 &< -\frac{1}{2}\\
\tilde{r}+l+1 &> -\frac{1}{2}-2s\\
\tilde{r} &> \frac{1}{2}+s
\end{align*}
then $\hat{T}_s(\sigma)R(\sigma) = \text{Id}_{\overline{H}_{(b)}^{\tilde{r},l}}$.
\item \label{secondId} If
\begin{align*}
l&< -\frac{1}{2}\\
\tilde{r}+l &> -\frac{1}{2}-2s\\
\tilde{r} &> \frac{1}{2}+s
\end{align*}
then $R(\sigma)\hat{T}_s(\sigma) = \text{Id}_{\mathcal{X}^{\tilde{r},l}_{\sigma}}$.
\item \label{thirdId} For any $\sigma'\neq 0$ with $\Im(\sigma')\geq 0$, if
\begin{align*}
l+1 &< -\frac{1}{2}\\
\tilde{r}+l &> -\frac{1}{2}-2s\\
\tilde{r}-1 &> \frac{1}{2}+s
\end{align*}
then we have $R(\sigma)\hat{T}_s(\sigma)R(\sigma') = R(\sigma')$ on $\overline{H}_{(b)}^{\tilde{r},l}$.
Since $R(\sigma')\overline{H}_{(b)}^{\tilde{r},l} = \mathcal{X}^{\tilde{r},l+1}_{\sigma'}$ (using the inequalities on $(\tilde{r},l)$) and since $\mathcal{X}^{\tilde{r},l+1}_{\sigma'}\subset \mathcal{X}^{\tilde{r}-1,l+1}_{\sigma}$, the identity is a consequence of \ref{secondId}.
\end{enumerate}
Under the hypotheses:
\begin{align*}
l+1 < -\frac{1}{2}\\
\tilde{r}+l &> -\frac{1}{2}-2s\\
\tilde{r}-1 &> \frac{1}{2}+s,
\end{align*}
identities \ref{firstId} and \ref{thirdId} enable to deduce the identity of the resolvent on $\overline{H}^{\tilde{r},l}_{(b)}$:
\begin{align*}
R(\sigma)-R(\sigma')= R(\sigma)\left(\hat{T}_s(\sigma')-\hat{T}_s(\sigma)\right)R(\sigma').
\end{align*}
Similarly, under the hypotheses:
\begin{align*}
-\frac{3}{2}-s-\left|s\right|<&l<-\frac{1}{2}\\
\tilde{r}+l>&-\frac{1}{2}-2s\\
\tilde{r}>&\frac{1}{2}+s
\end{align*}
 the following resolvent identity holds on $\overline{H}_{(b)}^{\tilde{r},l}$:
 \begin{align*}
 R(\sigma)-R(0) = R(\sigma)\left(\hat{T}_s(0)-\hat{T}_s(\sigma)\right)R(0).
 \end{align*}
\end{remark}

\subsection{Regularity of the resolvent}
\begin{prop}
\label{regHighFreq}
Let $l$ and $\tilde{r}$ be such that $l+1 < -\frac{1}{2}$, $\tilde{r}+l> -\frac{1}{2}-2s$ and $\tilde{r}-1>\frac{1}{2}+s$.
The family $R(\sigma)$ is holomorphic in $\Im(\sigma)>0$ as a family of operators in $\mathcal{L}(\overline{H}_{(b)}^{\tilde{r}+1,l}, \overline{H}_{(b)}^{\tilde{r},l+1})$. Moreover, it is locally Lipschitz on $D:=\left\{\sigma\in \C, \Im(\sigma)\geq 0, \sigma\neq 0\right\}$ as a family in $\mathcal{L}(\overline{H}_{(b)}^{\tilde{r},l}, \overline{H}_{(b)}^{\tilde{r}-1,l+1})$.
\end{prop}
\begin{proof}
For $\sigma, \sigma'\in D$, we have the resolvent identity on $\overline{H}_{(b)}^{\tilde{r},l}$:
\begin{align*}
R(\sigma)-R(\sigma') =& R(\sigma)(\hat{T}_s(\sigma)-\hat{T}_s(\sigma'))R(\sigma')\\
\end{align*} 
Note that $\lim\limits_{\sigma'\to \sigma}\frac{\hat{T}_s(\sigma)-\hat{T}_s(\sigma')}{\sigma-\sigma'} = -2\sigma a_{\mathfrak{t},\mathfrak{t}} + ia_{\mathfrak{t},\phi}\partial_\phi + ia_{\mathfrak{t},r}\partial_r + ia_\mathfrak{t} \in x^{-1}\text{Diff}_b^1$ the limit being in the norm topology of $\mathcal{L}(\overline{H}_{(b)}^{\tilde{r},l+1}, \overline{H}_{(b)}^{\tilde{r}-1,l})$. Therefore $\hat{T}_s(\sigma)$ is holomorphic and in particular locally Lipschitz continuous in this space. Since $R(\sigma')$ is uniformly bounded in $\mathcal{L}(\overline{H}_{(b)}^{\tilde{r},l}, \overline{H}_{(b)}^{\tilde{r}, l+1})$ (locally in $\sigma'\in D$), we obtain the local Lipschitz continuity property for $R(\sigma)$. 
If we perform the same computation on $\overline{H}_{(b)}^{\tilde{r}+1,l}$, using the continuity we just proved, we have that $R(\sigma')$ is continuous in $\mathcal{L}(\overline{H}_{(b)}^{\tilde{r}+1,l}, \overline{H}_{(b)}^{\tilde{r},l+1})$. We then deduce the following limit in the norm topology of $\mathcal{L}(\overline{H}_{(b)}^{\tilde{r}+1, l},\overline{H}_{(b)}^{\tilde{r}-1,l})$:
\begin{align*}
\lim\limits_{\sigma'\to \sigma}\frac{R(\sigma)-R(\sigma')}{\sigma-\sigma'} =& R(\sigma)\partial_{\sigma}\hat{T}_s(\sigma)R(\sigma)
\end{align*}
and as a consequence, $R(\sigma)$ is holomorphic on the interior of $D$ as a family in $\mathcal{L}(\overline{H}_{(b)}^{\tilde{r}+1,l}, \overline{H}_{(b)}^{\tilde{r}-1,l+1})$. Finally, using the Cauchy integral formula (see \cite{hille1996functional}, theorem 3.3.11 for a proof of the formula in the case of vector valued holomorphic functions) and the continuity of $R(\sigma)$ in $\mathcal{L}(\overline{H}_{(b)}^{\tilde{r}+1,l}, \overline{H}_{(b)}^{\tilde{r},l+1})$, we deduce that the holomorphic property holds in $\mathcal{L}(\overline{H}_{(b)}^{\tilde{r}+1,l}, \overline{H}_{(b)}^{\tilde{r},l+1})$.
\end{proof}
\begin{remark}
Note that if we assume $-\frac{3}{2}-s-\left|s\right|<l+1$ and use \eqref{lowFreqUniform} (with $\eta=0$), then we get $\left\|R(\sigma)\right\|_{\mathcal{L}\left(\overline{H}_{(b)}^{\tilde{r},l},\overline{H}_{(b)}^{\tilde{r}-1,l+1}\right)}\leq C\left|\sigma\right|^{-1}$ near zero. Since for $\eta \in [0,1]$, $\mathcal{L}\left(\overline{H}_{(b)}^{\tilde{r},l},\overline{H}_{(b)}^{\tilde{r}-1,l+1}\right)$ is continuously included in $\mathcal{L}\left(\overline{H}_{(b)}^{\tilde{r},l+\eta}, \overline{H}^{\tilde{r}-1, l+1}_{(b)}\right)$, we also have the continuity in this space and thanks to \eqref{lowFreqUniform} we see that the behavior near $\sigma = 0$ improves when $\eta$ increases (up to locally bounded for $\eta = 1$). We can get a local Hölder continuity statement by weakening a bit more the operator norm (see the proof of Proposition \ref{globalContinuity}).
\end{remark}

The following proposition provides a global (operator norm) continuity statement up to the real axis including at $\sigma = 0$. It also provides a rough uniform bound when $\left|\sigma\right|\to +\infty$. It will be used in the contour deformation argument.
\begin{prop}
\label{globalContinuity}
Let $-\frac{3}{2}-s-\lvert s\rvert<l<-\frac{1}{2}$, $\tilde{r}-1>\frac{1}{2}+s$ and $\tilde{r}+l-1>-\frac{1}{2}-2s$.
The family $R(\sigma)$ is continuous on $\left\{\Im(\sigma)\geq 0\right\}$ for the norm topology of $\mathcal{L}(\overline{H}_{(b)}^{\tilde{r},l}, \overline{H}_{(b)}^{\tilde{r}-1,l-})$ and for every $C>0$, there exists $D>0$ such that for every $\sigma$ with $0\leq \Im(\sigma)\leq C$:
\begin{align*}
\left\|R(\sigma)\right\|_{\mathcal{L}(\overline{H}_{(b)}^{\tilde{r},l}, \overline{H}_{(b)}^{\tilde{r}-1,l-})}\leq D\left<\sigma\right>^{\tilde{r}}
\end{align*}
\end{prop}
\begin{proof}
For the continuity at 0, we use the resolvent identity:
\begin{align*}
R(\sigma)-R(0) = R(\sigma)\left(\hat{T}_s(0)-\hat{T}_s(\sigma)\right)R(0)
\end{align*}
and under the hypotheses on $l$ and $\tilde{r}$, there exists $C'>0$ such that:
\begin{itemize}
\item $R(0)$ is well defined and bounded from $\overline{H}_{(b)}^{\tilde{r},l}$ to $\overline{H}_{(b)}^{\tilde{r},l}$ by propositions \ref{FredholmPropertyAtZero} and \ref{preciseComputationOfTheKernel}.
\item $\left\|\hat{T}_s(0)-\hat{T}_s(\sigma)\right\|_{\mathcal{L}(\overline{H}_{(b)}^{\tilde{r},l},\overline{H}_{(b)}^{\tilde{r}-1,l-1})}\leq C'\lvert\sigma\rvert$ since $\hat{T}_s(0)-\hat{T}_s(\sigma)\in x^{-1}\sigma \text{Diff}^1_b$
\item By \eqref{lowFreqUniform}, uniformly near zero we have $\left\|R(\sigma)\right\|_{\mathcal{L}(\overline{H}_{(b)}^{\tilde{r}-1,l-1},\overline{H}_{(b)}^{\tilde{r}-1,l-\delta})}\leq C'\lvert\sigma\rvert^{\delta-1}$ for $\delta>0$ such that $l-\delta>-\frac{3}{2}-s-\lvert s\rvert$ and $\tilde{r}+l-1-\delta>-\frac{1}{2}-2s$.
\end{itemize}
We get:
\begin{align*}
\left\|R(\sigma)-R(0)\right\|_{\mathcal{L}(\overline{H}_{(b)}^{\tilde{r},l},\overline{H}_{(b)}^{\tilde{r}-1,l-\delta})} \leq C''\lvert\sigma\rvert^{\delta}
\end{align*}
and therefore the continuity at zero.

The continuity at $\sigma\neq 0$ has already been proved in the stronger space $\mathcal{L}\left(\overline{H}_{(b)}^{\tilde{r},l-1}, \overline{H}_{(b)}^{\tilde{r}-1,l}\right)$ (see Proposition \ref{regHighFreq}).
 
To get the estimate 
\begin{align*}
\left\|R(\sigma)\right\|_{\mathcal{L}(\overline{H}_{(b)}^{\tilde{r},l}, \overline{H}_{(b)}^{\tilde{r}-1,l-})}\leq D\left<\sigma\right>^{\tilde{r}}
\end{align*}
we combine \eqref{highFreqInStrip} and \eqref{lowFreqUniform} (with $\eta = 1$) which provides uniform bounds in stronger norms than $\mathcal{L}(\overline{H}^{\tilde{r},l}_{(b)}, \overline{H}^{\tilde{r}-1,l-}_{(b)})$ (note that we compare semiclassical norms with usual norms using $\left\|u\right\|_{\overline{H}_{(b)}^{\tilde{r},l}}\left<\sigma\right>^{-\tilde{r}}\leq \left\|u\right\|_{\overline{H}^{\tilde{r},l}_{(b), \lvert \sigma\rvert^{-1}}}\leq \left\|u\right\|_{\overline{H}_{(b)}^{\tilde{r},l}}$ true for $\lvert \sigma\rvert \geq 1$).
\end{proof}

The following proposition establishes regularity of the resolvent on the real axis.
\begin{prop}
\label{differentiabilityRes}
Let $m\in\N$.
Let $l$ and $r$ be such that $l+1 < -\frac{1}{2}$, $\tilde{r}+l-2m> -\frac{1}{2}-2s$ and $\tilde{r}-2m-1>\frac{1}{2}+s$.
We have $R\in C^m\left(\R_\sigma \setminus\left\{0\right\}, \mathcal{L}\left(\overline{H}_{(b)}^{\tilde{r},l}, \overline{H}_{(b)}^{\tilde{r}-1-2m, l+1}\right)\right)$ and $\lvert\sigma\rvert^{-m}\partial_\sigma^m R(\sigma)$ is uniformly bounded in $\mathcal{L}\left(\overline{H}_{(b),\lvert\sigma\rvert^{-1}}^{\tilde{r},l}, \overline{H}_{(b),\lvert\sigma\rvert^{-1}}^{\tilde{r}-2m-1,l+1}\right)$ for $\sigma\in\R\setminus (-\alpha, \alpha)$ (where $\alpha>0$ is arbitrary).

Moreover, if we assume in addition that $-\frac{3}{2}-s-\left|s\right|<l+1$, there exists $\sigma_0>0$ such that for all $\eta\in [0,1]$, $\sigma^{m+1-\eta}\partial_{\sigma}^m R(\sigma)$ is uniformly bounded in $\mathcal{L}\left(\overline{H}_{(b)}^{\tilde{r},l+\eta}, \overline{H}_{(b)}^{\tilde{r}-1-2m, l+1}\right)$ for $\sigma\in(-\sigma_0,\sigma_0)\setminus\left\{0\right\}$.
\end{prop}
\begin{proof}
We first prove by induction that for all $m\in \N$, $R(\sigma)$, if $l+1 < -\frac{1}{2}$, $\tilde{r}+l-2m> -\frac{1}{2}-2s$ and $\tilde{r}-2m-1>\frac{1}{2}+s$ then $R(\sigma)\in C^{m}\left(\R_\sigma \setminus\left\{0\right\}, \mathcal{L}\left(\overline{H}_{(b)}^{\tilde{r},l}, \overline{H}_{(b)}^{\tilde{r}-1-2m, l+1}\right)\right)$. For $m=0$, it follows from Proposition \ref{regHighFreq}. Assume that the property is true for some $m\geq 0$ and assume $l+1 < -\frac{1}{2}$, $\tilde{r}+l-2(m+1)> -\frac{1}{2}-2s$ and $\tilde{r}-2(m+1)-1>\frac{1}{2}+s$. Then using the resolvent identity, we get that $R(\sigma)\in C^1\left(\R_{\sigma}\setminus\left\{0\right\},  \mathcal{L}\left(\overline{H}_{(b)}^{\tilde{r},l}, \overline{H}_{(b)}^{\tilde{r}-3, l+1}\right)\right)$ and  $\partial_\sigma R(\sigma) = R(\sigma)\partial_\sigma \hat{T}_s(\sigma) R(\sigma)$. By the induction hypothesis for all $m'\leq m$, $R(\sigma)\in C^{m'}\left(\R_{\sigma}\setminus\left\{0\right\}, \mathcal{L}\left(\overline{H}_{(b)}^{\tilde{r}-2(m-m')-2,l}, \overline{H}_{(b)}^{\tilde{r}-1-2(m+1), l+1}\right)\right)$ and $R(\sigma) \in C^{m'}\left(\R_{\sigma}\setminus\left\{0\right\}, \mathcal{L}\left(\overline{H}_{(b)}^{\tilde{r},l}, \overline{H}_{(b)}^{\tilde{r}-1-2m', l+1}\right)\right)$. We deduce that 
\begin{align*}
\partial_\sigma R(\sigma) \in C^{m}\left(\R_{\sigma}\setminus\left\{0\right\}, \mathcal{L}\left(\overline{H}_{(b)}^{\tilde{r},l}, \overline{H}_{(b)}^{\tilde{r}-1-2(m+1), l+1}\right)\right)
\end{align*} and by Leibniz rule:
\begin{align}\label{formulaDerivative}
\partial_{\sigma}^{m+1}R(\sigma) = \sum_{i_1+i_2+i_3 = m} a_{i_1,i_2,i_3}\partial_{\sigma}^{i_1}R(\sigma)\partial_\sigma^{i_2+1}\hat{T}_s(\sigma)\partial_{\sigma}^{i_3}R(\sigma)
\end{align}where $a_{i_1,i_2,i_3}$ are absolute combinatorial constants.
To prove that $\left|\sigma\right|^{-m}\partial_\sigma^m R(\sigma)$ is bounded on $\R\setminus (-\alpha,\alpha)$, we use an induction. The case $m = 0$ correspond to estimate \eqref{highFreqInStrip} and we conclude by using the induction hypothesis in \eqref{formulaDerivative} (together wit the fact that $\left|\sigma\right|^{-1}\partial_{\sigma}^{i_2+1}\hat{T}_s(\sigma)$ is bounded in $x^{-1}\text{Diff}^1_b$ for $\sigma \in \R\setminus (-\alpha,\alpha)$.
To prove that for all $\eta\in[0,1]$, $\left|\sigma\right|^{-m+1-\eta}\partial_\sigma^m R(\sigma)$ is bounded on $(-\sigma_0, \sigma_0)$ in the norm of $\mathcal{L}\left(\overline{H}_{(b)}^{\tilde{r},l+\eta}, \overline{H}_{(b)}^{\tilde{r}-1-2m, l+1}\right)$, we again use an induction on $m$. The case $m=0$ follows from estimate \eqref{lowFreqUniform} and we conclude by using the induction hypothesis with $\eta \in [0,1]$ for $\partial_{\sigma}^{i_3}R(\sigma)$  and with $\eta = 0$ for $\partial_{\sigma}^{i_1}R(\sigma)$ in \eqref{formulaDerivative} (noting that $\partial_{\sigma}^{i_2+1}$ is uniformly bounded on $(-\sigma_0,\sigma_0)$).
\end{proof}
\begin{remark}
Using the resolvent identity at $\sigma = 0$, we have:
\begin{align*}
R(\sigma)-R(0) = R(\sigma)\left(\hat{T}_s(0)-\hat{T}_s(\sigma)\right)R(0)
\end{align*}
Let $m\in \N$ and $l\in (-\frac{3}{2}-s-\left|s\right|,-\frac{1}{2}-s+\left|s\right|)$,$l_c\in (-\frac{3}{2}-s-\left|s\right|,-\frac{1}{2})$, $l_c\leq l\leq l_c+1$, $l_c+\tilde{r}-2m-1>-\frac{1}{2}-2s$ and $\tilde{r}-2m-1>\frac{1}{2}-s$.
Using the fact that $\sigma^{-1}\left(\hat{T}_s(0)-\hat{T}_s(\sigma)\right)R(0)$ maps $\overline{H}^{\tilde{r},l}$ to $\overline{H}_{(b)}^{\tilde{r},l-1}$ (with uniform bound with respect to $\sigma\in (-\sigma_0, \sigma_0)$) and the previous proposition, we deduce that for $\sigma \in (-\sigma_0, \sigma_0)$:
\begin{align*}
\left\|\left|\sigma\right|^{m}\partial_m (R(\sigma)-R(0))\right\|_{\mathcal{L}(\overline{H}_{(b)}^{\tilde{r},l},\overline{H}_{(b)}^{\tilde{r}-2m-1,l_c})}\leq C\left|\sigma\right|^{l-l_c}
\end{align*}
where the constant $C>0$ does not depend on $\sigma$.
\end{remark}

\subsection{More precise regularity near \texorpdfstring{$\sigma = 0$}{sigma = 0}}
First we establish more precisely the mapping properties of $R(0)$ which will be used throughout the section using a normal operator argument. Then we describe precisely the singular part of $R(\sigma)f$ near zero. We use the notations $Y_j$, $\mathcal{Y}_j^{\tilde{r},l}$, $\mathcal{Y}_{>j}^{\tilde{r},l}$, $\mathcal{Y}_{\leq j}^{\tilde{r},l}$ and $\Pi_j$ already introduced in Section \ref{sectionNotationLapSpin}.

\begin{lemma}\label{sumMode}
Let $l\geq l_0$, $\tilde{r}\geq \tilde{r}_0$. Let $u\in \overline{H}_{(b)}^{\tilde{r}_0,l_0}$ with $u=0$ in a neighborhood of $\left\{x\geq\frac{1}{r_+}\right\}$. We assume that for all $j\in \lvert s\rvert +\N$, $u_j:=\Pi_j u \in \mathcal{Y}^{\tilde{r},l}_j$ and $\sum_{p=0}^{+\infty}\int_\R (1+\lvert \tau\rvert^2+p^2)^{\tilde{r}}\lvert \hat{u}_{\left|s\right|+p}(\tau-i(l+\frac{3}{2}))\rvert^2\dd \tau<+\infty$ (where $\hat{u}_j$ is the Mellin transform of $u_j$). Then, $u\in \overline{H}^{r,l}_{(b)}$ and \begin{align*}\left\|u\right\|_{\overline{H}_{(b)}^{\tilde{r},l}} \leq C\sum_{p=0}^{+\infty} \int_\R (1+\lvert \tau\rvert^2+p^2)^{\tilde{r}}\left| \hat{u}_{\left|s\right|+p}\left(\tau-i\left(l+\frac{3}{2}\right)\right)\right|^2\dd \tau\end{align*} for some constant $C>0$ independent of $u$.
\end{lemma}
\begin{proof}
Let $u$ be defined as in the statement of the proposition. First note that the algebraic sum $\oplus_{j\in \lvert s\rvert + \N} \mathcal{Y}_j^{\tilde{r},l}$ is dense in $\overline{H}_{(b)}^{\tilde{r},l}$ and therefore also in $\overline{H}_{(b)}^{\tilde{r}_0,l_0}$. This density property can be obtained using first the density of $C^{\infty}_c((0,\frac{1}{r_+}])\otimes H^{\tilde{r}}(\mathcal{B}_s)$ in $\overline{H}_{(b)}^{\tilde{r},l}$ and then using the density of the algebraic sum $\oplus_{j\in\left|s\right|+\N}Y_j$ in $H^s(\mathcal{B}_s)$. This density and the orthogonality of $\mathcal{Y}_j^{\tilde{r}_0,l_0}$ and $\mathcal{Y}_j^{\tilde{r}_0,l_0}$ when $j\neq j'$ for the norm inducing scalar product $\left<v,w\right>_{\tilde{r}_0,l_0}:= \left<x^{-l_0}(1-(x\partial_x)^2+\Delta^{[s]})^{\frac{\tilde{r}_0}{2}})v, x^{-l_0}(1-(x\partial_x)^2+\Delta^{[s]})^{\frac{\tilde{r}_0}{2}})w\right>_{L^2_{(b)}}$ imply that the sequence $(\sum_{p=0}^N\Pi_{\left|s\right|+p}(u))_{N\in \N}$ converges towards $u$ in $\overline{H}_{(b)}^{\tilde{r}_0,l_0}$. Since $\mathcal{Y}_j^{\tilde{r},l}$ and $\mathcal{Y}_{j'}^{\tilde{r},l}$ are orthogonal when $j\neq j'$ for the norm-inducing scalar product $\left<v,w\right>_{\tilde{r},l}$, we get:
\begin{align*}
\left\|\sum_{p=0}^N u_{\left|s\right|+p}\right\|_{\overline{H}_{(b)}^{\tilde{r},l}}^2 =& \sum_{p=0}^N \left\|u_{\left|s\right|+p}\right|^2_{\overline{H}_{(b)}^{\tilde{r},l}}\\
=& \sum_{p=0}^N \left\|(1-(x\partial_x)^2 + p(p+2\left|s\right|)+p+\left|s\right|)^\frac{r}{2}u_{\left|s\right|+p}\right|^2_{\overline{H}_{(b)}^{\tilde{r},l}}\\
=& C\sum_{p=0}^{N} \int_\R (1+\lvert \tau\rvert^2+p^2)^{\tilde{r}}\lvert \hat{u}_{\left|s\right|+p}(\tau-i(l+\frac{3}{2}))\rvert^2\dd \tau.
\end{align*}
We deduce that $\sum_{p=0}^{+\infty} u_{\left|s\right|+p}$ is absolutely convergent in $\overline{H}_{(b)}^{\tilde{r},l}$ and therefore, the limit (which is equal to $u$ by uniqueness of the limit in $\overline{H}^{\tilde{r}_0,l_0}_{(b)}$) belongs to $\overline{H}_{(b)}^{\tilde{r},l}$. Moreover the computation shows that $\left\|u\right\|_{\overline{H}_{(b)}^{\tilde{r},l}} \leq C\sum_{p=0}^{+\infty} \int_\R (1+\lvert \tau\rvert^2+p^2)^{\tilde{r}}\lvert \hat{u}_{\left|s\right|+p}(\tau-i(l+\frac{3}{2}))\rvert^2\dd \tau$.
\end{proof}

We also record the following useful lemma about the Mellin transform.
\begin{lemma}\label{contourDeformationMellin}
Let $u\in \overline{H}^{\tilde{r}_0,l_0}_{(b)}$ with $\supp(u)$ contained in $\left\{x\leq D\right\}$ for some $D>0$. The Mellin transform $M u$ is holomorphic on $\left\{\Im(\tau)>-(l_0+\frac{3}{2})\right\}$ and $\alpha\mapsto Mu(.+i\alpha)$ is continuous from $[-\left(l_0+\frac{3}{2}\right), +\infty)$ to $\left< 1+\tau^2+\Delta^{[s]}\right>^{-\frac{\tilde{r}_0}{2}}L^2(\R_{\tau}, L^2(\mathcal{B}_s))$. Moreover, if for $l>l_0$, $Mu$ admits a meromorphic extension to $\left\{\Im(\tau)>-(l+\frac{3}{2})\right\}$ with a finite number of poles $a_0, ... a_N$ included in some compact subset $K$ of $\left\{\Im(\tau)>-(l+\frac{3}{2})\right\}$ and with bounds of the form:
\begin{align}\label{boundMellin}
\int \left|\left(1+\tau^2+\Delta^{[s]}\right)^{\frac{\tilde{r}_0}{2}}1_{\C\setminus K}Mu (\tau +i\alpha)\right|^2\dd \tau \leq C_0
\end{align} for some $C_0>0$ independent of $\alpha$.
Then there exists a family of non negative integers $(c_j)_{j=0}^N$ and a family of complex numbers $(b_{j,k})_{0\leq j\leq N\\ 0\leq k\leq c_j}$ such that, $u \in \sum_{j=0}^N \sum_{k=0}^{c_j}b_{j,k}\ln(x)^{k}x^{ia_j}\chi_0(x) + \overline{H}^{\tilde{r}_0, l}_{(b)}$ where $\chi_0$ is any smooth cutoff equal to $1$ near zero and compactly supported.
\end{lemma}
\begin{proof}
By definition, for $\Im(\tau)>-\left(l_0+\frac{3}{2}\right)$, $Mu (\tau) = \int x^{-i\tau}u(x)\frac{\dd x}{x}$. Since for every $\epsilon>0$, $x^{-(l_0+\frac{3}{2})+\epsilon}u(x) \in L^1_b([0,+\infty))$, we get that $Mu$ is holomorphic on $\left\{ \Im(\tau)>-(l_0+\frac{3}{2})\right\}$. If $\alpha, \alpha' \leq l_0$, $x^{-\alpha}u \in H^{\tilde{r}_0,0}_{(b)}$ and 
\begin{align*}\left\|(x^{-\alpha}-x^{-\alpha'})u\right\|_{H^{\tilde{r}_0,0}_{(b)}}\leq C\sum_{j=0}^{N'}\left\|(xD_x)^j(x^{l_0-\alpha}-x^{l_0-\alpha'})\chi_1\right\|_{L^\infty}\left\|x^{-l_0}u\right\|_{\overline{H}^{\tilde{r}_0,0}_{(b)}}
\end{align*}
where $\chi_1 = 1$ on $x(\supp(u))$ and $\chi_1$ is smooth and compactly supported and $N'$ is a large integer. Using that $\lim\limits_{\alpha\to \alpha'}\left\|(xD_x)^j(x^{l_0-\alpha}-x^{l_0-\alpha'})\chi_1\right\|_{L^\infty} = 0$ for any $j\in \N$, we deduce that $\alpha \mapsto x^{\alpha}u$ is continuous in $\dot{H}^{\tilde{r}_0,l_0}$ and therefore, by continuity of the Mellin transform from $\dot{H}^{\tilde{r}_0,l_0}$ to $\left< 1+\tau^2+\Delta^{[s]}\right>^{-\frac{\tilde{r}_0}{2}}L^2(\R_{\tau}, L^2(\mathcal{B}_s))$ we conclude the first part of the proof. Now we assume that $Mu$ admits an extension as in the statement.
Then by the residue theorem, we have for all $A>0$ and $\epsilon>0$ small enough such that $K\subset \left\{ z: -(l+\frac{3}{2})+\epsilon\leq\Im(z)\leq -(l_0 + \frac{3}{2})\right\}$ and for all $x\in (0,+\infty)$:
\begin{align*}
\int_{-A}^A x^{i(\tau -i(l_0+\frac{3}{2})}Mu(\tau -i(l_0+\frac{3}{2}))\dd \tau =& \sum_{j=0}^N 2i\pi\text{Res}_{|_{\tau = a_j}}x^{i\tau}Mu(\tau)\\
& + \int_{-A}^A x^{i(\tau -i(l+\frac{3}{2})+i\epsilon}Mu(\tau -i(l+\frac{3}{2})+i\epsilon)\dd \tau \\
&+ i\int_{-(l+\frac{3}{2})+\epsilon}^{-(l_0+\frac{3}{2})}x^{i(A +iy)}Mu(A+iy)\dd y\\
&-i\int_{-(l+\frac{3}{2})+\epsilon}^{-(l_0+\frac{3}{2})}x^{i(-A+iy)}Mu(-A+iy)\dd y
\end{align*}
We have the following limit in the sense of distributions:
\begin{align*}
\lim\limits_{A\to +\infty}\int_{-A}^A x^{i(\tau -i(l_0+\frac{3}{2})}Mu(\tau -i(l_0+\frac{3}{2}))\dd \tau = u(x)
\end{align*}
Moreover, (using the bound \eqref{boundMellin}) we get that:
\begin{align*}
\int_{-A}^A x^{i(\tau -i(l+\frac{3}{2})+i\epsilon}Mu(\tau -i(l+\frac{3}{2})+i\epsilon)\dd \tau
\end{align*}
has a limit $v\in x^{l+\frac{3}{2}-\epsilon}H^{\tilde{r}_0,0}_{b}((0,+\infty)\times \mathcal{B}_s)$ in the sense of distributions when $A\to +\infty$.
We deduce that
\begin{align*}
f_A := i\int_{-(l+\frac{3}{2})+\epsilon}^{-(l_0+\frac{3}{2})}x^{i(A +iy)}Mu(A+iy)\dd y-i\int_{-(l+\frac{3}{2})+\epsilon}^{-(l_0+\frac{3}{2})}x^{i(-A+iy)}Mu(-A+iy)\dd y
\end{align*}
also has a limit in the sense of distributions. We prove that this limit is zero. Let $\phi \in C^{\infty}_c((0,+\infty)\times \mathcal{B}_s)$. There exists a constant $C>0$ (independent of $A$) such that for a large integer $K$:
\begin{align*}
\left|\left<f_A, \phi\right>\right|\leq C\left|x^{C}D_x^K \phi\right|_{L^1((0,+\infty),H^{-\tilde{r}_0}(\mathcal{B}_s))}&\left(\int_{-(l+\frac{3}{2})+\epsilon}^{-(l_0+\frac{3}{2})} A^{-K}\left|(1+\Delta^{[s]})^{\frac{\tilde{r}_0}{2}}Mu(A+iy)\right|^2_{L^2(\mathcal{B}_s)}\right.\\
&\hphantom{(\int_{-(l+\frac{3}{2})+\epsilon}^{-(l_0+\frac{3}{2})}}\left.+ A^{-K}\left|(1+\Delta^{[s]})^{\tilde{r}_0{2}}Mu(A+iy)\right|^2_{L^2(\mathcal{B}_s)}\dd y\right)
\end{align*}
Since for $A_0$ and $K$ large enough\footnote{More explicitly, if $\tilde{r}_0\geq 0$, we can take any $K\geq 0$ and if $\tilde{r}_0<0$, we take $K\geq -\frac{\tilde{r}_0}{2}$}:
\begin{align*}\int_{A_0}^{+\infty}\left(\int_{-(l+\frac{3}{2})+\epsilon}^{-(l_0+\frac{3}{2})} A^{-K}\left|(1+\Delta^{[s]})^{\frac{\tilde{r}_0}{2}}Mu(A+iy)\right|^2_{L^2(\mathcal{B}_s)}\right.&\\
\left.+A^{-K}\left|(1+\Delta^{[s]})^{\tilde{r}_0{2}}Mu(A+iy)\right|^2_{L^2(\mathcal{B}_s)}\dd y\right)\dd A &\leq CC_0
\end{align*}
where $C_0$ is the constant in \eqref{boundMellin}. We deduce that there exists a sequence $A_n$ with $\lim\limits_{n\to +\infty}A_n = +\infty$ such that 
\begin{align*}
\lim\limits_{n\to +\infty}\left(\int_{-(l+\frac{3}{2})+\epsilon}^{-(l_0+\frac{3}{2})} A_n^{-K}\left|(1+\Delta^{[s]})^{\frac{\tilde{r}_0}{2}}Mu(A_n+iy)\right|^2_{L^2(\mathcal{B}_s)}\right.&\\
\left.+ A_n^{-K}\left|(1+\Delta^{[s]})^{\tilde{r}_0{2}}Mu(A_n+iy)\right|^2_{L^2(\mathcal{B}_s)}\dd y\right) &= 0
\end{align*}
Using this particular sequence, we deduce that $\lim\limits_{A\to +\infty}f_{A} = 0$ in the sense of distributions. Therefore writing $\sum_{k=0}^{c_j}b_{j,k}\ln(x)^kx^{ia_j}:= 2i\pi \text{Res}_{|_{\tau = a_j}}x^{i\tau}Mu(\tau)$ we get:
\begin{align*}
u = v + \sum_{j=0}^N \sum_{k=0}^{c_j}b_{j,k}\ln(x)^{k}x^{ia_j}
\end{align*}
We use this equality with different values of $\epsilon$ and we get for all $0<\epsilon < \epsilon_0$:
\begin{align*}
\left\|v\right\|_{x^{l+\frac{3}{2}-\epsilon}H_b^{\tilde{r}_0,0}}\leq C_1C_0
\end{align*}
where $C_1$ is independent of $\epsilon$ and $u$. Let $\phi \in C^{\infty}_c((0,+\infty)\times \mathcal{B}_s)$, we have:
\begin{align}
\left|\left<v,\phi\right>\right|\leq& \left\|v\right\|_{x^{l+\frac{3}{2}-\epsilon}H_b^{\tilde{r}_0,0}}\left\|\phi\right\|_{x^{-l-\frac{3}{2}+\epsilon}H_b^{-\tilde{r}_0,0}} \notag\\
\leq& C_1C_0\left\|\phi\right\|_{x^{-l-\frac{3}{2}+\epsilon}H_b^{-\tilde{r}_0,0}}\label{estimateV}
\end{align}
But for all $N\in \N$, $\lim\limits_{\epsilon \to 0}\left\|\left(x^{l+\frac{3}{2}-\epsilon}-x^{l+\frac{3}{2}}\right)\phi\right\|_{H^{\tilde{r}_0,0}_b} = 0$ (all $b$ derivatives converge uniformly on the support of $\phi$). Taking the limit $\epsilon\to 0$ in $\eqref{estimateV}$, we get $v\in H_b^{\tilde{r}_0,l+\frac{3}{2}}$ with a norm smaller than $C_1C_0$.
Using the support condition of $u$, we have:
\begin{align*}
u = \chi_0v + \sum_{j=0}^N \sum_{k=0}^{c_j}b_{j,k}\ln(x)^{k}x^{ia_j}\chi_0
\end{align*}
where $\left\|\chi_0v\right\|_{\overline{H}^{\tilde{r}_0,l}_{(b)}}\leq C_2C_0$ with $C_2$ independent of $u$.
\end{proof}
\begin{remark}\label{contourDeformationMellinRk}
The estimate \eqref{boundMellin} can be replaced by an other similar estimate. For example, we can assume that $Mu$ extends meromorphically to $\left\{\Im(\tau)>D_0\right\}$ (where $D_0\in \R$) with a finite number of poles and that we have bounds of the form:
\begin{align*}
\left\|Mu(\tau)\right\|^2_{H^{\tilde{r}}(\mathcal{B}_s)}\leq C(1+\Re(\tau)^2)^{-\tilde{r}}e^{C\left|\Im(\tau)\right|}
\end{align*}
for $\tau \in \left\{\Im(\tau)\geq D\right\}\setminus K$ where $D>D_0$ and $\Re(K)$ compact. We also assume that there is no pole with imaginary part equal to $D$. In, this case we have that $u \in  \sum_{j=0}^N \sum_{k=0}^{c_j}b_{j,k}\ln(x)^{k}x^{ia_j}\chi_0+ \overline{H}^{\tilde{r},-D-\frac{3}{2}}_{(b)}$.
\end{remark}

\begin{lemma}
\label{mappingHighHarmonics}
Let $j\in \N$, $-\frac{3}{2}-s-\lvert s \rvert<l<-\frac{1}{2}+j-s+\lvert s\rvert$ and $\tilde{r}+1>\frac{1}{2}+s$, then $\hat{T}_s(0)^{-1}(\mathcal{Y}_{>j+\left|s\right|}^{\tilde{r},l})\subset \mathcal{Y}_{>j+\left|s\right|}^{\tilde{r}+1,l}$.
\end{lemma}
\begin{proof}
Since $\hat{T}_s(0)$ commutes with $\Delta^{[s]}$, for $\epsilon>0$ small enough, $\hat{T}_s(0)^{-1}(\mathcal{Y}^{\tilde{r},l}_{>j+\left|s\right|})\subset \mathcal{Y}^{\tilde{r}+1, -\frac{3}{2}-s-\lvert s\rvert+\epsilon}_{>j+\left|s\right|}$. Let $u\in \hat{T}_s(0)^{-1}(\mathcal{Y}^{\tilde{r},l}_{>j+\left|s\right|})$. We have $\hat{T}_s(0)\chi_0 u = f^0\in \mathcal{Y}^{\tilde{r},l}_{>j+\left|s\right|}$ (where $\chi_0$ is a cutoff depending only on $x$ and localizing near zero) and $N(\hat{T}_s(0))\chi_0 u = f^0-x\text{Diff}^2_b u =: f^1$. Applying $\Pi_{p+\left|s\right|}$ for $p\geq j$ and the Mellin transform, we get that $(\tau^2 + (i-2is)\tau + p(p+2\lvert s\rvert)+\lvert s\rvert +p + s)\widehat{\chi_0u_{p+\left|s\right|}} = \hat{f}^1_{p+\left|s\right|}$ where $\widehat{\chi_0u}_{p+\left|s\right|}$, is holomorphic on $\left\{\Im(\tau)> s+\lvert s\rvert-\epsilon\right\}$. Since $\hat{f}^0_{p+\left|s\right|}$ is holomorphic on $\left\{\Im(\tau)>-(l+\frac{3}{2})\right\}$ we deduce that $\hat{f}^1_{p+\left|s\right|}$ is holomorphic on $\left\{\Im(\tau)>\max(-(l+\frac{3}{2}), -1+s+\lvert s\rvert-\epsilon)\right\}$. Since $(\tau^2 + (i-2is)\tau + p(p+2\lvert s\rvert)+\lvert s\rvert +p + s)$ has no zero on $\left\{-1-j+s-\left|s\right|<\Im(\tau)<s+\left|s\right|+p\right\}$ we deduce that $\widehat{\chi_0u}_{p+\left|s\right|}$ extends holomorphically to $\left\{\Im(\tau)>\max(-(l+\frac{3}{2}),-1+s+\lvert s\rvert-\epsilon)\right\}$. An iteration of this argument proves that $\widehat{\chi_0u}_{p+\left|s\right|}$ and $\hat{f}^1_{p+\left|s\right|}$ extend holomorphically to $\left\{\Im(\tau)>-(l+\frac{3}{2})\right\}$. 
Moreover, for any $\alpha \geq s+\lvert s\rvert -\epsilon$, there exists $C>0$ independent of $p$ such that 
\begin{align*}\int_{\R}\left|\widehat{\chi_0u}_{p+\left|s\right|}(\tau+i\alpha)\right|^2\left(\left|\tau\right|^2+p^2+1\right)^{\tilde{r}+1}\dd \tau& \leq C\left\|x^{\alpha}\chi_0u_{p+\left|s\right|}\right\|^2_{\overline{H}^{\tilde{r}+1,0}_b}\\
& \leq C\left\|x^{\alpha}\chi_0u_{p+\left|s\right|}\right\|^2_{\overline{H}^{\tilde{r}+1,-s-\lvert s\rvert +\epsilon}_b}
\end{align*} and for any $\beta \geq -(l+\frac{3}{2})$, 
\begin{align*}
\int_{\R} (1+p^2 + \left|\tau\right|^2)^{\tilde{r}}\left|\hat{f}_{p+\left|s\right|}^0(\tau+i\beta)\right|^2\dd \tau \leq& C\left\|x^{\beta}f^0_{p+\left|s\right|}\right\|^2_{\overline{H}^{\tilde{r},0}_b}\\
\leq& C\left\|f^0_{p+\left|s\right|}\right\|^2_{\overline{H}^{\tilde{r},l}_{(b)}}
\end{align*}
Using that $\left|\tau^2 + (i-2is)\tau + p(p+2\lvert s\rvert)+\lvert s\rvert +p + s\right|^{-1}\leq C(\Re(\tau)^2 + p^2 + 1)^{-1}$ where $C$ is uniform with respect to $p\geq j$ when $-(l+\frac{3}{2})\leq \Im(\tau)\leq s+\left|s\right|-\epsilon$, we deduce that for any $\max(-(l+\frac{3}{2}), -1+s+\left|s\right|-\epsilon)\leq \alpha\leq s+\left\|s\right|-\epsilon $, there exists $C'$ independent of $p\geq j$ such that:
\begin{align*}
\int_{\R}\left|\widehat{\chi_0u}_{p+\left|s\right|}(\tau+i\alpha)\right|^2\left(\left|\tau\right|^2+p^2+1\right)^{\tilde{r}+1}\dd \tau\leq C'\left(\left\|\chi_0u_{p+\left|s\right|}\right\|^2_{\overline{H}^{\tilde{r}+1,-\frac{3}{2}-s-\left|s\right|+\epsilon}_{(b)}} + \left\|x^{\alpha}f^0_{p+\left|s\right|}\right\|^2_{\overline{H}^{\tilde{r},l}_{(b)}}\right)
\end{align*}
Once again we can iterate the argument and deduce that there exists $C''$ independent of $p$ such that for $\alpha = -(l+\frac{3}{2})$:
\begin{align*}
\int_{\R}\left|\widehat{\chi_0u}_{p+\left|s\right|}(\tau+i\alpha)\right|^2\left(\left|\tau\right|^2+p^2+1\right)^{\tilde{r}+1}\dd \tau \leq C''\left(\left\|\chi_0 u_{p+\left|s\right|}\right\|_{\overline{H}^{\tilde{r}+1,-\frac{3}{2}-s-\left|s\right|+\epsilon}_{(b)}}+\left\|f^0\right\|_{\overline{H}^{\tilde{r},l}_{(b)}}\right)
\end{align*}
Using Lemma \ref{contourDeformationMellin}, we get that for each $p\geq j$, $\chi_0u_{p} \in \overline{H}^{\tilde{r},l}_{(b)}$ with 
\begin{align*}\left\|\chi_0 u_{p+\left|s\right|}\right\|^2_{\overline{H}^{\tilde{r},l}_{(b)}}\leq C''\left(\left\|\chi_0 u_{p+\left|s\right|}\right\|^2_{\overline{H}^{\tilde{r}+1,-\frac{3}{2}-s-\left|s\right|+\epsilon}_{(b)}}+\left\|f^0_{p+\left|s\right|}\right\|^2_{\overline{H}^{\tilde{r},l}_{(b)}}\right).
\end{align*}
Since the constant is uniform with respect to $p\geq j$, we deduce:
\begin{align*}
\sum_{p=j}^{+\infty} \int_{\R}\left|\widehat{\chi_0u}_{p+\left|s\right|}(\tau-i(l+\frac{3}{2}))\right|^2\left(\left|\tau\right|^2+p^2+1\right)^{\tilde{r}+1}\dd \tau \leq C''\left(\left\|\chi_0 u\right\|^2_{\overline{H}^{\tilde{r}+1,-\frac{3}{2}-s-\left|s\right|+\epsilon}_{(b)}}+\left\|f^0\right\|^2_{\overline{H}^{\tilde{r},l}_{(b)}}\right)
\end{align*}
We conclude by Lemma \ref{sumMode} that $u\in \mathcal{Y}^{\tilde{r}+1, l}_{>j+\left|s\right|}$.
\end{proof}

\begin{lemma}
\label{MeroExtension}
Let $k\in \N$, $\alpha\in\C$, $\epsilon>0$, $\chi$ a smooth cutoff localizing near zero and $c\in E\setminus\left\{0\right\}$ (where $E$ is some Banach space) and $u = \chi(x) x^{\alpha}\ln(x)^k c$. Then the Mellin transform of $u$ has a meromorphic extension to $\C$ with a pole of order $k+1$ at $-i\alpha$.
\end{lemma}
\begin{proof}
We have that $(xD_x+i\alpha)^{k}u \in C^{\infty}_c((0,+\infty),E)$. Therefore, $(\tau + i\alpha)^k \hat{u}$ has a holomorphic extension to $\C$ which prove the lemma.
\end{proof}

In the remaining part of the section, $f_{l,m}$ is used with variables $(\theta,\varphi)$. Note that we have 
\begin{align*}
f_{l,m}(\theta,\varphi) = e^{im(\varphi-\phi_*)}f_{l,m}(\theta, \phi^*)
\end{align*} where $\varphi-\phi_*$ is a function of $r$. In the following lemma, we record some mapping properties of $\hat{T}_s(0)^{-1}$. We first give an expansion of $\hat{T}_s(0)^{-1}u$ for a general $u\in \overline{H}_{(b)}^{\tilde{r},l}$ (see the statement for the precise conditions on $l$ and $\tilde{r}$). Then, we give an expansion for $u$ of the form $x^{\alpha}\ln(x)^p$ near $x=0$. It will be useful since such terms naturally arise in the iteration in the proof of Proposition \ref{highDecay}.
\begin{lemma}
\label{normalOperatorArgument}
The indices $\lvert s\rvert +j$ denote an element of $Y_{\lvert s\rvert + j}$.
Let $-\frac{1}{2}-s+\lvert s\rvert+k<l< -\frac{1}{2}-s+\lvert s \rvert + k+1$ for $k\in \N$, $\tilde{r}+1>\frac{1}{2}+s$ and $u \in \overline{H}_{(b)}^{\tilde{r},l}$, then we have:
\begin{align*}
\hat{T}_s(0)^{-1}u = \chi(x)x^{1-s+\lvert s\rvert}(\sum_{j=0}^{k}\sum_{i=j}^{k}x^{i} v_{i,\lvert s\rvert+j}) + \overline{H}_{(b)}^{\tilde{r}+1, l}
\end{align*}
Moreover, we have for $0\leq j\leq k$:
\begin{align*}
v_{j,\lvert s\rvert+j} = \frac{1}{1+2(\lvert s\rvert+j)}\sum_{\lvert m\rvert \leq \lvert s\rvert+j}\left<\overline{U^{*}_{\lvert s\rvert+j,m}},u\right>f_{\lvert s\rvert+j,m}
\end{align*}
and for $0\leq j\leq k-1$:
\begin{align*}
v_{j+1, \lvert s\rvert +j} = \sum_{\lvert m\rvert\leq \lvert s\rvert + j}\frac{\left(M(j^2+(2\lvert s\rvert-s)j+s(s-\lvert s\rvert)+2\lvert s\rvert-s+1) -iams\right)}{(1+2(\lvert s\rvert+j))(j+\lvert s\rvert + 1)}\left<\overline{U^{*}_{\lvert s\rvert+j,m}},u\right>f_{\lvert s\rvert+j,m}
\end{align*}
Where $U^*_{\lvert s\rvert+j,m}$ is the unique multiple of $u^*_{\lvert s\rvert+j,m}$ equivalent to $r^{\lvert s\rvert+j-s}f_{\lvert s\rvert +j,m}(\theta,\varphi)\dd r\dd \varphi\dd \theta$ when $r\to +\infty$.
Let $p\in \N$ and $\chi$ be a smooth cutoff localizing near zero.
If $\alpha = 1-s+\lvert s\rvert+j-k>-s-\lvert s \rvert-j$ for some $k,j\in \N$:
\begin{align*}
\hat{T}_s(0)^{-1}\chi(x)x^{\alpha}\ln(x)^p u_{\lvert s\rvert+j} =& \chi(x)x^{\alpha}\sum_{j'=0}^{p}\sum_{i=0}^{n}x^i\ln(x)^{j'}v_{i,j',\lvert s \rvert+j} + \chi(x)x^{\alpha}\sum_{i=k}^{n}x^{i}\ln(x)^{p+1}w_{i,\lvert s\rvert+j}\\& + \overline{H}_{(b)}^{\infty, \alpha+n-\frac{1}{2}-}
\end{align*}
Moreover, if $\alpha < 1-s+\lvert s\rvert +j$, we have 
$v_{0,p,\lvert s\rvert +j} = \frac{u_{\lvert s\rvert+j}}{k(1+2\lvert s\rvert +2j-k)}$, if $\alpha=-s+\lvert s\rvert +j$, $j+\lvert s\rvert\neq 0$ and $p=0$, we have 
\begin{align*}w_{1,\lvert s\rvert+j}=\sum_{\lvert m\rvert \leq \lvert s\rvert +j}\frac{\left<f_{\lvert s\rvert+j,m},u_{\lvert s\rvert +j}\right>f_{\lvert s\rvert+j,m}}{2(j+\lvert s\rvert)(1+2j+2\lvert s\rvert)}&\left( 2M(j^2+(2\lvert s\rvert-s+2a-2)j +2(a-1)\lvert s\rvert+s(s-\lvert s\rvert))\right. \\ &\left.-2iams\right)
\end{align*} and if $\alpha = 1-s+\lvert s\rvert +j$, we have $w_{0,\lvert s\rvert+j} = -\frac{u_{\lvert s\rvert+j}}{(p+1)(1+2\lvert s\rvert +2j)}$.
\end{lemma}

\begin{proof}
Let $-\frac{1}{2}-s+\lvert s\rvert+k<l\leq -\frac{1}{2}-s+\lvert s\rvert+k+1$ for $k\in \N$ and $u \in \overline{H}_{(b)}^{\tilde{r},l}$.
Let $\chi_0$ be a smooth cutoff localizing near 0. We define $v := \hat{T}_s(0)^{-1}u \in \overline{H}_{(b)}^{\tilde{r}+1,-\frac{1}{2}-s+\lvert s\rvert -}$.
By definition, $\hat{T}_s(0)v = u$.
\begin{align*}
N(\hat{T}_s(0))\chi_0 v =& u-x\text{Diff}^{2}_b v =: f\\
\end{align*}
Since $N(\hat{T}_s(0))$ and $\hat{T}_s(0)-N(\hat{T}_s(0))$ commutes with $\Delta^{[s]}$, we can project the equality on each eigenspace. We get for $j\in \N$:
\begin{align*}
((xD_x)^2+(i-2is)xD_x + j(j+2\lvert s\rvert) + \lvert s\rvert +j+s) \chi_0 v_{\lvert s\rvert+j} =& f_{\lvert s\rvert+j}
\end{align*}
We compute the Mellin transform:
\begin{align}\label{relMellin}
(\tau^2 + (i-2is)\tau + j(j+2\lvert s\rvert)+\lvert s\rvert +j + s) \widehat{\chi_0v}_{\lvert s\rvert+j} =& \hat{f}_{\lvert s\rvert+j}
\end{align}
We define $k$ as in the statement of the Proposition (i.e. the smallest integer such that $-1-(k+1)+s-\left|s\right|<-l-\frac{3}{2}$). We fix $K$ a compact subset of $\C$ such that $K\cap \left\{\Im(\tau) = -l-\frac{3}{2}\right\} = \emptyset$ and $K$ contains a neighborhood of $\mathcal{P}:=\left\{ i(-1-j+s-\left|s\right|), 0\leq j \leq k\right\}$ and such that if $\Im(\tau)\geq -l-\frac{3}{2}$ and $\tau \in K$, then $\tau-i \in K$.
Note that for all $j \in \N$, $(\tau^2 + (i-2is)\tau + j(j+2\lvert s\rvert)+\lvert s\rvert +j + s)^{-1}$ is meromorphic on $\left\{-1-(k+1)+s-\left|s\right|<\Im(\tau)<-\frac{1}{2}+s-\left|s\right| \right\}$ with poles of order 1 included in $\mathcal{P}$. Moreover we have a constant $C>0$ independent of $j$ such that for every $\tau \in \left\{ -l-\frac{3}{2}\leq \Im(\tau)\leq -\frac{1}{2}+s-\left|s\right| \right\}\setminus K$:
\begin{align}\left|\tau^2 + (i-2is)\tau + j(j+2\lvert s\rvert)+\lvert s\rvert +j + s\right|^{-1}\leq C(\Re(\tau)^2 + j^2+1)^{-1}. \label{boundMellinMultiplier}
\end{align}
Using \eqref{relMellin}, we have that:
\begin{itemize}
 \item $\widehat{\chi_0v}_{\lvert s\rvert+j}$ is holomorphic on $\mathcal{D}_0:=\left\{\Im(\tau)>\max(-1-j+s-\lvert s\rvert,-l-\frac{3}{2})\right\}$
\item There exists $C>0$ independent of $j$ such that for all $b>\max(-1-j+s-\lvert s\rvert, -l-\frac{3}{2})$, $\int_{\R}\left|1_K(\tau)(1+\tau^2+j^2)^{\frac{r+1}{2}}\widehat{\chi_0v}_{\lvert s\rvert+j}(\tau+ib)\right|^2\dd \tau\leq C\left(\left\|u_{j+\left|s\right|}\right\|^2_{\overline{H}^{\tilde{r},l}_{(b)}} + \left\|v_{j+\left|s\right|}\right\|^2_{\overline{H}^{\tilde{r}+1,-\frac{1}{2}-s+\left|s\right|-\epsilon}_{(b)}}\right)$
\end{itemize}
We prove by induction on $m\in \N$ that 
\begin{itemize}
\item $\widehat{\chi_0v}_{\lvert s\rvert+j}$ has a meromorphic extension (with at most simple poles at $i(-1-j-\N+s-\lvert s\rvert)$) to $\mathcal{D}_m:=\left\{\Im(\tau)>\max(-1-j-m+s-\lvert s\rvert,-l-\frac{3}{2})\right\}$ for $m\in \N$.
\item There exists $C>0$ independent of $j$ such that for all $b>\max(-1-j-m+s-\lvert s\rvert, -l-\frac{3}{2})$, $\int_{\R} \left|1_{\C\setminus K}(\tau)(1+\lvert \tau\rvert^2+j^2)^{\frac{r+1}{2}}\widehat{\chi_0v}_{\lvert s\rvert+j}(\tau+ib)\right|^2\dd \tau\leq C\left(\left\|u_{j+\left|s\right|}\right\|^2_{\overline{H}^{\tilde{r},l}_{(b)}} + \left\|v_{j+\left|s\right|}\right\|^2_{\overline{H}^{\tilde{r}+1,-\frac{1}{2}-s+\left|s\right|-\epsilon}_{(b)}}\right)$ 
\end{itemize}
We assume the induction hypothesis for some $m\in \N$. We have $f_{\lvert s\rvert+j} = \chi_1 u_{\lvert s\rvert+j} - \chi_1 x\text{Diff}^2_b v_{\lvert s\rvert +j}$ (where $\chi_1$ is a cutoff localizing near zero such that $\chi_1 = 1$ on the support of $\chi_0$). We have $\chi_1 x\text{Diff}^2_b v_{\lvert s\rvert+j} = x\text{Diff}^2_b \chi_0 v_{\lvert s\rvert +j} + w$ where $w \in \chi_1 x\text{Diff}^2_b (1-\chi_0) v_{\left|s\right|+j}$ has compact support with respect to $x$ in $(0,+\infty)$ and is therefore in $\overline{H}^{\tilde{r}-1, \infty}$ and for any $N\in \R$, there exists $C_N>0$ independent of $v$ such that 
\begin{align*}
\left\|w\right\|_{\overline{H}^{\tilde{r}-1,N}_{(b)}}\leq C_N\left\|v_{\left|s\right|+j}\right\|_{\overline{H}^{\tilde{r}+1, -\frac{1}{2}-s+\left|s\right|-\epsilon}_{(b)}}.
\end{align*} We deduce that $\hat{f}_{\lvert s\rvert+j}$ has a meromorphic extension to $\left\{\Im(\tau)>\max(-1-j-m-1+s-\lvert s\rvert, -l-\frac{3}{2})\right\}$ with at most simple poles at $\left\{i(-1-j-p+s-\lvert s\rvert), 1\leq p\leq m+1\right\}$. Moreover, we have for $b > \max\left(-1-j-(m+1)+s-\left|s\right|, -l-\frac{3}{2}\right)$: 
\begin{align*}
\int_{\R}\left|1_{\C\setminus K}(\tau+ib)(1+\tau^2+ j^2)^{\frac{\tilde{r}-1}{2}}\hat{f}_{\left|s\right|+j}(\tau +ib)\right|^2\dd \tau \leq& C\left(\int_{\R} \left|(1+\tau^2+j^2)^{\frac{\tilde{r}}{2}}\hat{u}\right|^2\dd \tau\right.\\
&\hspace{-2.3cm}+ \int_{\R}1_{\C\setminus K}(\tau+ib)\left|(1+\tau^2+j^2)^{\frac{\tilde{r}+1}{2}}\widehat{\chi_0v}_{\left|s\right|+j}(\tau+i(b+1))\right|^2\dd\tau\\
&\left. + \int_{\R}\left|(1+\tau^2+j^2)^{\frac{\tilde{r}-1}{2}}\hat{w}(\tau+ib)\right|^2\dd \tau\right)\\
\leq& C\left(\left\|u_{\left|s\right|+j}\right\|^2_{\overline{H}^{\tilde{r},l}_{(b)}}+ \left\|v_{\left|s\right|+j}\right\|^2_{\overline{H}^{\tilde{r}+1,-\frac{1}{2}-s+\left|s\right|-\epsilon}_{(b)}}\right)
\end{align*}
where we used the induction hypothesis to bound the term involving $\widehat{\chi_0v}_{\left|s\right|+j}(\tau+i(b+1))$.
Using the identity $(\tau^2 + (i-2is)\tau + j(j+2\lvert s\rvert)+\lvert s\rvert +j + s) = (\tau -i(-1+s-\lvert s\rvert-j))(\tau -i(s+\lvert s\rvert+j))$, we conclude that $\widehat{\chi_0 v}_{\lvert s\rvert+j}$ has a meromorphic extension to $\mathcal{D}_{m+1}$ with (at most) simple poles at $\left\{i(-1-j-p+s-\lvert s\rvert), 0\leq p\leq m+1\right\}$. Using \eqref{boundMellinMultiplier} we get 
\begin{align*}\int_{\R}\left|1_{\C\setminus K}(1+\tau^2+j^2)^{\frac{\tilde{r}+1}{2}}\widehat{\chi_0 v_{\left|s\right|+j}}(\tau+ib)\right|^2\dd \tau \leq& C\int_{\R}\left|1_{\C\setminus K}(\tau+ib)(1+\tau^2+ j^2)^{\frac{\tilde{r}-1}{2}}\hat{f}_{\left|s\right|+j}(\tau +ib)\right|^2\dd \tau.
\end{align*} We can therefore apply Lemma \ref{contourDeformationMellin} to get :
\begin{align*}
v_{\lvert s\rvert +j} = x^{1-s+\lvert s\rvert+j}\sum_{m = 0}^{k-j}x^{m}v_{m,\lvert s\rvert + j} + \overline{H}^{\tilde{r}+1,l}_{(b)}
\end{align*}
Using this decomposition for $j\leq k$ and lemma \ref{mappingHighHarmonics} to prove that $v-\sum_{p=0}^k v_{\lvert s\rvert+p} \in \overline{H}_{(b)}^{\tilde{r}+1,l}$, we get the first claim of the lemma.

To determine exactly the form of $v_{j,\lvert s\rvert+j}$, we fix $0\leq j\leq k$ and compute for every $\epsilon>0$:
\begin{align*}
\left<\overline{U^{*}_{\lvert s\rvert+j,m}}, u\right> =& \left<\overline{U^{*}_{\lvert s\rvert+j,m}}, \hat{T}_s(0)\hat{T}_s(0)^{-1}u\right> \\
=& \left<\overline{U^*_{\lvert s\rvert+j,m}}, \hat{T}_s(0)\left(\chi(\frac{x}{\epsilon})x^{1-s+\lvert s\rvert}(\sum_{j'=0}^{k}\sum_{i=j'}^{k}x^{i} v_{i,\lvert s\rvert+j'}) + \overline{H}_{(b)}^{\tilde{r}+1, l}\right) \right>\\
\end{align*}
Note that $U^*_{\left|s\right|+j,m} \in \dot{H}_{(b)}^{(-s+\frac{1}{2})-, (\frac{1}{2}+s-\left|s\right|-j)-}r^{2}\dd r \dd \phi \dd \theta$. In particular, it belongs to the dual space of $\overline{H}_{(b)}^{\tilde{r}-1, l}$ under the hypothesis on $\tilde{r}$ and $l$. Therefore, for all $g \in \overline{H}_{(b)}^{\tilde{r}+1, l}$, 
\begin{align*}\left<\overline{U^*_{\lvert s\rvert+j,m}}, \hat{T}_s(0) g\right> = \left< \overline{\hat{T}_s(0)^* U^*_{\left|s\right|+j,m}}, g\right>. 
\end{align*}
Using the fact that $\hat{T}_s(0)^* U^*_{\lvert s\rvert+j,m} = 0$ and the orthogonality property of $f_{\lvert s\rvert+j,-m}$, we get:
\begin{align*}
\left<\overline{U^{*}_{\lvert s\rvert+j,m}}, u\right> =& \left<\overline{U^{*}_{\lvert s\rvert+j,m}}, \hat{T}_s(0)\chi\left(\frac{x}{\epsilon}\right)x^{1-s+\lvert s\rvert+j}v_{j,\lvert s\rvert+j}\right>\\
\end{align*}
We use the fact that $\chi\left(\frac{x}{\epsilon}\right)x^{1-s+\lvert s\rvert+j}v_{j,\lvert s\rvert+j}$ tends to zero in $\overline{H}_{(b)}^{r,(-\frac{1}{2}-s+\lvert s\rvert +j)-}$ when $\epsilon$ tends to zero, the fact that $N(\hat{T}_s(0))x^{1-s+\lvert s\rvert+j}v_{j,\lvert s\rvert+j} = 0$ and the fact that $U^*_{\lvert s\rvert+j,-m}-x^{-2+s-\lvert s\rvert-j}\dd x \in \overline{H}_{(b)}^{\infty,(\frac{3}{2}+s-\lvert s\rvert-j)-}x^{-4}\dd x$ to see that:
\begin{align*}
\left<\overline{U^{*}_{\lvert s\rvert+j,m}}, u\right> =& \lim\limits_{\epsilon\to 0}\int_0^{+\infty}x^{-2+s-\lvert s\rvert -j}[N(\hat{T}_s(0),\chi\left(\frac{x}{\epsilon}\right)]x^{1-s+\lvert s\rvert +j}\dd x\left<f_{\lvert s\rvert+j,m},v_{j,\lvert s\rvert+j}\right>\\
[N(\hat{T}_s(0)), \chi\left(\frac{x}{\epsilon}\right)]x^{1-s+\lvert s\rvert +j} =& -2(1-s+\lvert s\rvert+j)\frac{x^{2-s+\lvert s\rvert +j}}{\epsilon}\chi'\left(\frac{x}{\epsilon}\right)\\
&-x^{2-s+\lvert s\rvert +j}\partial_x\left(x\partial_x\chi\left(\frac{x}{\epsilon}\right)\right)+\left(1-2s\right)\frac{x^{2-s+\lvert s\rvert +j}}{\epsilon}\chi'\left(\frac{x}{\epsilon}\right)
\end{align*}
We conclude that:
\begin{align*}
\left<\overline{U^{*}_{\lvert s\rvert+j,m}}, u\right> =& \left(1+2\left(\lvert s\rvert +j\right)\right)\left<f_{\lvert s\rvert+j,m},v_{j,\lvert s\rvert+j}\right>
\end{align*}
Since $\left<f_{\lvert s\rvert+j,m},v_{j,\lvert s\rvert+j}\right>$ is the coefficient of $f_{\lvert s\rvert+j,m}$ in the decomposition of $v_{j,\lvert s\rvert+j}$, we get the claimed decomposition.
In the case $0\leq j\leq k-1$, we get (after projection on $\mathcal{Y}_{\lvert s\rvert +j}$):
\begin{align*}
\hat{T}_s(0)(v_{j,\lvert s\rvert +j}x^{1-s+\lvert s\rvert +j}+v_{j+1,\lvert s\rvert +j}x^{2-s+\lvert s\rvert +j} + \overline{H}^{\tilde{r}, \min(\frac{3}{2}-s+\lvert s\rvert+j-, l)}) = u_{\lvert s\rvert +j}
\end{align*}
We deduce that the coefficients of $x^{1-s+\lvert s\rvert+j}$ and $x^{2-s+\lvert s\rvert +j}$ in $\hat{T}_s(0)(v_{j,\lvert s\rvert +j}x^{1-s+\lvert s\rvert +j}+v_{j+1,\lvert s\rvert +j}x^{2-s+\lvert s\rvert +j})$ have to vanish. The coefficient of $x^{1-s+\lvert s\rvert +j}$ vanish since $1-s+\lvert s\rvert$ is an indicial root of $N(\hat{T}_s(0))$ and the vanishing or the coefficient of $x^{2-s+\lvert s\rvert+j}$ gives:
\begin{align*}
v_{j+1, \lvert s\rvert +j} = \sum_{\lvert m\rvert\leq \lvert s\rvert + j}\frac{\left(2M(j^2+(2\lvert s\rvert-s)j+s(s-\lvert s\rvert)+2\lvert s\rvert-s+1) -2iams\right)}{2(1+2(\lvert s\rvert+j))(j+\lvert s\rvert + 1)}\left<\overline{U^{*}_{\lvert s\rvert+j,m}},u\right>f_{\lvert s\rvert+j,m}
\end{align*}
We can iteratively compute the other terms.

We now prove the second claim. Let $\alpha = 1-s+\lvert s\rvert+j-k>-s-\lvert s \rvert$ for some $k\in \N$. Let $v = \hat{T}_s(0)^{-1}x^{\alpha}\ln(x)^pu_{\lvert s\rvert+j}$. A priori for all $\epsilon>0$, $v\in \overline{H}^{\infty, \alpha-\frac{3}{2}-\epsilon}_{(b)}$. As previously, we have $(\tau^2 + (i-2is)\tau + j(j+2\lvert s\rvert)+\lvert s\rvert +j + s)\widehat{\chi_0v} = \hat{f}$ with $f = x^{\alpha}\ln(x)^{p}u_{\lvert s\rvert+j} - x\text{Diff}^2_b v$. We prove by induction on $m\in \N$ as previously (and using lemma \ref{MeroExtension}) that we have:
\begin{itemize}
\item $\widehat{\chi_0v}$ has a meromorphic extension to $\left\{\Im(\tau)>-\alpha-m\right\}$ with poles of order at most $p+1$ at $-i\alpha$, ..., $-i(\alpha+k-1)$ and poles of order at most $p+2$ at 
\begin{align*}-i(\alpha + k) = -i(1-s+\lvert s\rvert+j), ..., -i(\alpha + m-1).
\end{align*}
\item $\widehat{f}$ has a meromorphic extension to $\left\{\Im(\tau)>-\alpha-m-1\right\}$ with poles of order at most $p+1$ at $-i\alpha$, ..., $-i(\alpha+k)$ and poles of order at most $p+2$ at 
\begin{align*}-i(\alpha + k+1) = -i(2-s+\lvert s\rvert), ..., -i(\alpha + m).
\end{align*}
\item For every $N\in\N$, there exists $C>0$ such that for all $\tau \in \C$ with $\left|\Re(\tau)\right|\geq 1$: 
\begin{align*}\lvert \hat{f}\rvert \leq C(1+\Re(\tau)^2)^{-N}e^{C\lvert \Im\tau\rvert}.
\end{align*}
\end{itemize}
It is enough to use Remark \ref{contourDeformationMellinRk}. At a pole $\theta$ of order at most $q$ of $\widehat{\chi_0v}$, the residu of $x^{i\tau} \widehat{\chi_0 v} = x^{i\theta}(e^{i\ln(x)(\tau-\theta)}\widehat{\chi_0 v}$ is of the form $ x^{i\theta} \sum_{j'=0}^{q-1} \frac{(i\ln(x))^{j'}}{j'!}\text{Res}_{|_{\tau = \theta}}((\tau-\theta)^{j'}\widehat{\chi_0 v})$. Therefore, we get the claimed expansion for $v$.
To compute the exact expression of the principal term, we use the fact that 
\begin{align*}&\hat{T}_s(0)\left(\chi(x)x^{\alpha}\sum_{j'=0}^{p}\sum_{i=0}^{n}x^i\ln(x)^{j'}v_{i,j',\lvert s \rvert+j}\right.\\
&\hphantom{\hat{T}_s(0)\left(\right.}\left. + \chi(x)x^{\alpha}\sum_{i=k}^{n}x^{i}\ln(x)^{p+1}w_{i,\lvert s\rvert+j} + \overline{H}_{(b)}^{\infty, \alpha+n-\frac{1}{2}-}\right) = \chi(x)x^{\alpha}\ln(x)^p u_{\lvert s\rvert+j}
\end{align*} together with the computation:
\begin{align*}
N(\hat{T}_s(0))x^\alpha\ln(x)^pv_{\lvert s\rvert +j} =& -(\alpha+\lvert s\rvert +j+s)(\alpha-\lvert s\rvert-j+s-1)x^\alpha\ln(x)^pv_{\left|s\right|+j}\\
&+ p(1-2s-2\alpha)x^\alpha\ln(x)^{p-1}v_{\left|s\right|+j} - p(p-1)x^\alpha\ln(x)^{p-2}v_{\left|s\right|+j}
\end{align*}
Remark that we could compute the other terms iteratively and in particular, in the case $p=0$ and $\alpha = -s+\lvert s\rvert$, we get the claimed expression for $w_{1,\lvert s\rvert +j}$.
\end{proof}

For $k\in\N$ denote by $W_b^{k,\infty}(\R_\sigma, H_{(b)}^{\tilde{r},l})$ the set of elements of functions $f$ from $\R$ to $H_{(b)}^{\tilde{r},l}$ which are $k$ times weakly differentiable in $\sigma\neq 0$ and such that for all $0\leq j\leq k$, $(\sigma\partial_{\sigma})^jf$ is bounded on every compact neighborhood of zero.
In this part, we fix $f\in C^\infty(\R_{\sigma}, H_{(b)}^{\tilde{r},l})$ with compact support and for different values of the decay rate $l$ (we also impose condition on $\tilde{r}$ which should be considered as large). Our goal is to study the regularity of $R(\sigma)f(\sigma)$ with respect to $\sigma$ in the various cases.

\begin{prop}

\label{moderateDecay}
Let $f\in C_c^\infty(\R_\sigma, \overline{H}_{(b)}^{\tilde{r}, l})$ with $l\in (-\frac{5}{2}-s-\lvert s\rvert, \frac{1}{2}-s+\lvert s\rvert)\setminus \left\{-\frac{1}{2}-s+\lvert s\rvert \right\}$, $\tilde{r}+l_c-2k-1>-\frac{1}{2}-2s$ and $\tilde{r}-2k-1>\frac{1}{2}+s$. Let $l_c \in (-\frac{3}{2}-s-\lvert s\rvert, -\frac{1}{2})$ such that $l_c\leq l+1$. Then, we have that $R(\sigma)f \in C^{\infty}(\R_\sigma, \overline{H}_{(b)}^{\tilde{r}, l_c})+\lvert \sigma \rvert^{l-l_c}W_b^{k,\infty}(\R_\sigma, \overline{H}_{(b)}^{\tilde{r}-2k-1,l_c})$.
\end{prop}

\begin{remark}
If $l = \frac{1}{2}-s+\lvert s\rvert$, using the result with $l-\epsilon$ for every $\epsilon>0$, we have that $R(\sigma)f \in C^{\infty}(\R_\sigma, \overline{H}_{(b)}^{\tilde{r}, l_c})+\lvert \sigma \rvert^{(l-l_c)-}W_b^{k,\infty}(\R_\sigma, \overline{H}_{(b)}^{\tilde{r}-2k-1, l_c})$.
\end{remark}
\begin{remark}
\label{polynomialTerm}
If $f$ is constant with respect to $\sigma$ near the origin we can be more precise about the smooth term. It is a polynomial near the origin and if in addition $l-1<l_c$, then it is constant with respect to $\sigma$ equal to $R(0)f$.
\end{remark}

\begin{proof}
We proceed in three steps.
\begin{itemize}
\item \emph{If $l\leq l_c$:} In this case, we use directly proposition \ref{differentiabilityRes} to conclude.
\item \emph{If $l\in (-\frac{3}{2}-s-\lvert s\rvert, -\frac{1}{2}-s+\lvert s\rvert)$ and $l>l_c$:} In this case, we use the resolvent identity to write:
\begin{align*}
R(\sigma)f = R(0)f - R(\sigma)(\hat{T}_s(\sigma)-\hat{T}_s(0))R(0)f
\end{align*}
The term $R(0)f$ is in $C^\infty\left(\R_\sigma, \overline{H}_{(b)}^{\tilde{r}+1, l}\right)$ and $f_1:=(\hat{T}_s(\sigma)-\hat{T}_s(0))R(0)f \in \sigma C^{\infty}\left(\R_\sigma, \overline{H}_{(b)}^{\tilde{r}, l-1}\right)$. We can iterate this procedure until $l-k\leq l_c$.
\item \emph{If $l\in (-\frac{1}{2}-s+\lvert s\rvert, \frac{1}{2}-s+\lvert s\rvert)$:} In this case we also write $R(\sigma) f = R(0)f+ R(\sigma)(\hat{T}_s(\sigma)-\hat{T}_s(0))R(0)f$ and we use a normal operator argument to see that $R(0)f = cx^{1-s+\lvert s\rvert}\chi(x) + f_1$ where $c$ is a constant depending smoothly on $\sigma$, $\chi$ is a cutoff near zero and $f_1 \in \overline{H}_{(b)}^{\tilde{r},l}$ depend smoothly on $\sigma$. Therefore, we have:
\begin{align*}
R(\sigma)f = R(0)f + R(\sigma)v_1
\end{align*}
where $v_1 = C^{\infty}(\R_\sigma)(\hat{T}_s(\sigma)-\hat{T}_s(0))x^{1-s+\lvert s\rvert}\chi(x)+ \sigma C^{\infty}\left(R_\sigma, \overline{H}_{(b)}^{\tilde{r}, l-1}\right)$. Note that if $-s+\lvert s\rvert = 0$, we have $v_1 \in \sigma C^{\infty}\left(R_\sigma, \overline{H}_{(b)}^{\tilde{r}, l-1}\right)$ (since $1$ is the indicial root of the normal operator of $\hat{T}_s(\sigma)-\hat{T}_s(0)$) and we can conclude immediately. If $-s+\lvert s \rvert$ is a positive integer, we can define $v_{k+1} = (\hat{T}_s(\sigma)-\hat{T}_s(0))R(0)v_k$ for $1\leq k\leq -s+\lvert s\rvert$ and we have $R(\sigma)v_1 \in C^{\infty}(\R_\sigma, \overline{H}_{(b)}^{\tilde{r},l_c})+R(\sigma)v_{-s+\lvert s\rvert}$. By induction of the normal operator argument, we get that for $k\leq -s+\lvert s\rvert$, $v_{k} = \sigma^k C^\infty(\R_\sigma) x^{1-k-s+\lvert s\rvert}\chi(x) + \sigma^{k}C^\infty(\R_\sigma, \overline{H}_{(b)}^{\tilde{r}, l-k})$. For $k = -s+\lvert s\rvert + 1$, we have the indicial root cancellation which gives $v_{k} \in \sigma^{k}C^\infty(\R_\sigma, \overline{H}_{(b)}^{\tilde{r}, l-k})$ and $l-1+s-\lvert s \rvert\in (-\frac{3}{2}-s-\lvert s\rvert, -\frac{1}{2}-s+\lvert s\rvert)$, therefore we can conclude by the second case.
\end{itemize}
\end{proof}

We now describe the behavior of $R(\sigma)f$ when $f$ has a higher spatial decay.
\begin{prop}
\label{highDecay}
Let $f\in C^\infty(\R_\sigma, \overline{H}_{(b)}^{\tilde{r}, l})$ with $l>\frac{1}{2}-s+\lvert s\rvert$, $\tilde{r}-\frac{3}{2}-s-\left|s\right|+\epsilon-2k-1>-\frac{1}{2}-2s$, $\tilde{r}-2k-1>\frac{1}{2}+s$. We have the following equality:
\begin{align*}
R(\sigma)f = C^\infty(\R_\sigma, \overline{H}_{(b)}^{\tilde{r}, -\frac{1}{2}-s+\lvert s\rvert-}) + \sigma^{2\lvert s\rvert + 2}R(\sigma)v + \lvert \sigma\rvert^{2\lvert s\rvert +2 + \epsilon}W_b^{k,\infty}(\R_\sigma, \overline{H}_{(b)}^{\tilde{r}-2k-1,-\frac{3}{2}-s-\lvert s\rvert +\epsilon})
\end{align*}
where $v = x^{-s-\lvert s\rvert}\chi(x)c_f + \overline{H}_{(b)}^{\tilde{r},-\frac{3}{2}-s-\lvert s\rvert + \epsilon}$
for $\chi$ a smooth cutoff localizing near zero and $c_f$ an element of $Y_{\lvert s\rvert}$.
More precisely, we have:
\begin{align*}
c_f = \sum_{\lvert m\rvert \leq \lvert s\rvert, m-s\in \Z} (-1)^{\lvert s\vert -s + 1}(2i)^{1+2\lvert s\rvert}\frac{\mathfrak{c}_m}{(2\lvert s\rvert)!}
\end{align*}
where
\begin{align*}
\mathfrak{c}_m :=& \frac{4iMs+\mathfrak{c}^{(2)}_m + 2i(1+\lvert s\rvert -s)\mathfrak{c}^{(3)}_m + \mathfrak{c}^{(4)}_m}{(1+2\lvert s\rvert )^{2}}\left<\overline{U^*_{\lvert s\rvert, m}}, f(0)\right>f_{\lvert s\rvert, m}\\
\mathfrak{c}^{(2)}_m :=& 2as\beta_{\lvert s\rvert, m}\\
\mathfrak{c}^{(3)}_m :=& \frac{-iams + M(s(s-\lvert s\rvert) + 2\lvert s\rvert -s + 1)}{(\lvert s\rvert + 1)}\\
\mathfrak{c}^{(4)}_m :=& \begin{cases}-2i(-2iams+2M(2(a-1)\lvert s\rvert + s(s-\lvert s\rvert)) \text{ if $ s <0$}\\
0 \text{ if $s\geq 0$}
\end{cases}\\
\beta_{\lvert s\rvert,m} :=& \left<f_{\lvert s\rvert ,m},\cos(\theta)f_{\lvert s\rvert ,m}\right>
\end{align*}
\end{prop}

\begin{proof}
First note that since $f(\sigma) = \sum_{i=0}^{2\lvert s\rvert +2}\sigma^{i}f_i + \sigma^{2\lvert s\rvert +3}C^{\infty}(\R_\sigma, \overline{H}_{(b)}^{\tilde{r},l})$ and $R(\sigma)\sigma^{2\left|s\right|+3}C^{\infty}(\R_\sigma, \overline{H}_{(b)}^{\tilde{r},l})\subset \lvert \sigma\rvert^{2\lvert s\rvert +2+\epsilon}W_b^{k,\infty}(\R_{\sigma},\overline{H}_{(b)}^{\tilde{r}-2k-1, -\frac{1}{2}-s-\lvert s\rvert-})$ (see Proposition \ref{differentiabilityRes}), we are reduce to the case of $f$ independent of $\sigma$.
In this proof $c_{\left|s\right|}$ is an element of $Y_{\lvert s\rvert}$ and $c_{\left|s\right|+1}$ is an element of $Y_{\lvert s\rvert +1}$ but each instance (even in the same line) can be different.
For this proof, we also record the expressions of $\hat{T}_s(\sigma)-\hat{T}_s(0)$ and $\hat{T}_s(0)$ near $x= 0$:
\begin{align*}
\hat{T}_s(\sigma)-\hat{T}_s(0) =& a^2\sin^2\theta \sigma^2 + \frac{4Mar}{\Delta_r}\sigma D_\phi - 2(a^2+r^2)\sigma x^2D_x \\
&- i\sigma\left(\frac{4Msa^2-2sr(a^2+r^2)}{\Delta_r}+2(s+1)r+2ias\cos\theta\right)\\
=& -2\sigma(D_x + ix^{-1})+a^2\sin^2\theta \sigma^2 + i\sigma 4Ms +2\sigma as\cos\theta+\sigma x \text{Diff}^1_b\\
\hat{T}_s(0)=& \frac{a^2}{\Delta_r}\partial_\phi^2 + \Delta^{[s]}-\Delta_r^{-s}\partial_r\Delta_r^{s+1}\partial_r + 4s(r-M)\partial_r - 2s\frac{a(r-M)}{\Delta_r}\partial_\phi + s
\end{align*}
Note that $\hat{T}_s(0)$ commutes with $\Delta^{[s]}$ but it is not the case of $\hat{T}_s(\sigma)-\hat{T}_s(0)$.
We define recursively $u_1 = (\hat{T}_s(\sigma)-\hat{T}_s(0))\hat{T}_s(0)^{-1}f$ and $u_{j+1} = (\hat{T}_s(\sigma)-\hat{T}_s(0))\hat{T}_s(0)^{-1}u_j$ for $1\leq j\leq 2\lvert s\rvert + 2$. It is not clear yet that this sequence is well defined since $\hat{T}_s(0)^{-1}$ is only defined on $\overline{H}_{(b)}^{\tilde{r},l}$ for $l>-\frac{3}{2}-s-\lvert s\rvert$ but we will see in the proof that $u_j$ remains in this space when $j\leq 2\lvert s \rvert + 1$.
We prove recursively using lemma \ref{normalOperatorArgument}, the explicit form of $\hat{T}_s(\sigma)-\hat{T}_s(0)$ and the fact that $\cos\theta Y_{\left|s\right|}\subset Y_{\left|s\right|}+Y_{\left|s\right|+1}$ (see for example (2.33) in \cite{ma2021sharp}) that for $j \leq \lvert s \rvert-s$:
\begin{align*}
u_j =& \chi(x)\left(c_{\left|s\right|} \sigma^{j}x^{1-s+\lvert s\rvert -j} + \sigma^j x^{2-s+\lvert s \rvert -j}(c_{\left|s\right|}(j-1)\ln(x) + c_{\left|s\right|} + c_{\left|s\right|+1})\right)\\
& + \sigma^{j+1}P(\sigma)\overline{H}_{(b)}^{\tilde{r}, \frac{1}{2}-s+\lvert s\rvert-j-} + \sigma^j Q(\sigma)\overline{H}_{(b)}^{\tilde{r},\frac{1}{2}-s+\lvert s\rvert - j + \epsilon}
\end{align*}
where $P$ and $Q$ are polynomials in $\sigma$ and $\chi$ is any smooth cutoff localizing near $x=0$.
For $j = \lvert s\rvert - s +1$:
\begin{align*}
u_j =& \chi(x)\left(\sigma^j x^{2-s+\lvert s \rvert -j}((j-1)c_{\left|s\right|}\ln(x) + c_{\left|s\right|} + c_{\left|s\right|+1}) + \sigma^{j+1}P(\sigma)\overline{H}_{(b)}^{\tilde{r}, \frac{1}{2}-s+\lvert s\rvert-j-}\right)\\
& + \sigma^j Q(\sigma)\overline{H}_{(b)}^{\tilde{r},\frac{1}{2}-s+\lvert s\rvert - j + \epsilon}
\end{align*}
For $\lvert s\rvert -s+2 \leq j \leq 2\lvert s\rvert + 2$:
\begin{align*}
u_j = \chi(x)\left(\sigma^j x^{2-s+\lvert s \rvert -j}c_{\left|s\right|}  + \sigma^{j+1}P(\sigma)\overline{H}_{(b)}^{\tilde{r}, \frac{1}{2}-s+\lvert s\rvert-j-}\right) + \sigma^jQ(\sigma) \overline{H}_{(b)}^{\tilde{r},\frac{1}{2}-s+\lvert s\rvert - j + \epsilon}
\end{align*}
And finally, for all $0\leq j\leq 2\lvert s\rvert+1$,
\begin{align*}
R(\sigma)f = \sum_{k=0}^{j}(-1)^k R(0)u_k + (-1)^{j+1}R(\sigma)u_{j+1}
\end{align*}

To compute precisely the $c_{\left|s\right|}$ term in the expression of $u_{2\left|s\right|+2}$, we see that we have two cases:
\begin{itemize}
\item Case $s\geq 0$: The sequence of terms $c_{\left|s\right|}$ in the expressions of $u_j$ (call $c^{j}_{\left|s\right|}$ the term appearing in $u_j$) can be computed from the recurrence relation (obtained using Lemma \ref{normalOperatorArgument}):
\begin{align*}
c^{1}_{\left|s\right|} =& \sum_{\left|m\right|\leq \left|s\right|}\frac{2i(1-s+\left|s\right|)\mathfrak{c}^{(3)}_m + 4iMs+2as\beta_{\left|s\right|,m}}{1+2\left|s\right|} \left<\overline{U}^*_{\left|s\right|,m},f(0)\right>f_{\left|s\right|,m} \\
c^{2}_{\left|s\right|} =& -2i\frac{c^{1}_{\left|s\right|}}{1+2\left|s\right|}\\
c^{j+1}_{\left|s\right|} = &2i(1-s+\left|s\right|-j)\frac{c^{j}_{\left|s\right|}}{(j-1)(2+2\left|s\right|-j)} \text{, for $j\geq 2$}
\end{align*}
\item Case $s<0$: The logarithmic term in the expression of $u_j$ for $2\leq j\leq \left|s\right|-s+1$ and the terms $c_{\left|s\right|}$ in the expression of $u_j$ for $\lvert s\rvert -s+2 \leq j \leq 2\lvert s\rvert + 2$ can be computed recursively (we call $c^{j}_{\left|s\right|}$ the coefficients appearing in these terms): 
First, the term $c^{2}_{\left|s\right|}$ is obtained by computing the logarithmic term in $u_2$.
\begin{align*}
\alpha_2:=& \sum_{\left|m\right|\leq \left|s\right|}\frac{2i(1-s+\left|s\right|)\mathfrak{c}^{(3)}_m + 4iMs+2as\beta_{\left|s\right|,m}}{1+2\left|s\right|} \left<\overline{U}^*_{\left|s\right|,m},f(0)\right>f_{\left|s\right|,m}\\
c^{2}_{\left|s\right|} =& 2i(1-s+\left|s\right|)\left(-\frac{\alpha_2}{1+2\left|s\right|}\right.\\
&\left. + \sum_{\left|m\right|\leq \left|s\right|} \frac{2i}{(1+2\left|s\right|)^2}\left(2M(2(a-1)\left|s\right|+s(s-\left|s\right|))-2iams\right)\left<\overline{U}^*_{\left|s\right|,m},f(0)\right>f_{\left|s\right|,m}\right)\\\
\end{align*}
Then we have the following recursive relation (on the logarithmic terms) for $2\leq j\leq -s+\left|s\right|$:
\begin{align*}
c^{j+1}_{\left|s\right|} =& 2i(1-s+\left|s\right|-j)\frac{c^{j}_{\left|s\right|}}{(j-1)(2+2\left|s\right|-j)}
\end{align*}
Then we have:
\begin{align*}
c^{2+2\left|s\right|}_{\left|s\right|} =& \frac{2ic^{1+2\left|s\right|}_{\left|s\right|}}{(\left|s\right|-s)(1+2\left|s\right|)}
\end{align*}
\end{itemize}
\end{proof}

Now we want a more precise description of the term $R(\sigma)v$ in Proposition \ref{highDecay}. 
We need to introduce the effective normal operator (see Definition \ref{defNormalOp}). We have:
\begin{align*}
N_{\eff}(\hat{P}(\sigma)) =& (x^2D_x)^2+ 2ix(x^2D_x)\\ 
&+ x^2\left(\frac{1}{\sin\theta}D_\theta \sin\theta D_\theta + \frac{1}{\sin^2\theta}D_\phi^2 + \frac{2s\cos\theta}{\sin^2\theta}D_\phi + s^2\text{cotan}^2\theta+s\right)\\
& -2is x^2\left(xD_x + \frac{i}{2}\right) -sx^2 -2\sigma\left(x^2D_x + ix\right)
\end{align*}
we deduce that
\begin{align*}
N_{\eff}(\hat{T}_s(\sigma)) =& x^{-2}N_{\eff}(\hat{P}(\sigma))\\
=& (xD_x)^2 + ixD_x + \Delta^{[s]} -2is\left(xD_x + \frac{i}{2}\right)-2\frac{\sigma}{x} \left(xD_x +i\right)
\end{align*}
For $\sigma >0$, the change of variable $X = \frac{x}{\lvert \sigma\rvert}$ gives:
\begin{align*}
N_{\eff}^{+} :=& (XD_X)^2 + iXD_X + \Delta^{[s]} -2is\left(XD_X + \frac{i}{2}\right)-2X^{-1} \left(XD_X +i\right)
\end{align*}
and for $\sigma<0$, it gives:
\begin{align*}
N_{\eff}^{-} :=& (XD_X)^2 + iXD_X + \Delta^{[s]} -2is\left(XD_X + \frac{i}{2}\right)+2X^{-1} \left(XD_X +i\right)
\end{align*}
We use Definition \ref{defbSpacesTwoEnds} (but here the variable $X$ plays the role of $x$ in the definition) to introduce the spaces $H_b^{\tilde{r},l,\nu}$.
Using corollary \ref{invertibilityNormalEff}, we get that $N_{\eff}^{\pm}$ are invertible from $\left\{u\in H^{\tilde{r},l,\nu}_b: N_{\eff}^{\pm}(\hat{T}_s(\sigma))u\in H^{\tilde{r}, l-1, \nu}_b\right\}$ to $H^{\tilde{r},l-1, \nu}_b$ where $l<-\frac{1}{2}$, $\tilde{r}+l>-\frac{1}{2}-2s$ and $\nu \in \left(\frac{1}{2}+s-\lvert s\rvert, \frac{3}{2}+s+\lvert s\rvert\right)$. 

Let $B$ be a Banach space. For $\alpha \in \R$, we denote by $\mathcal{A}([0,1)_{\sigma}, \sigma^{\alpha}B)$ the space of smooth functions $u$ from $(0,1)$ to $B$ such that for all $k\in \N$, $\sup_{\sigma\in (0,1)}\left\|\sigma^{-\alpha}\left(\sigma\partial_{\sigma}\right)^k u \right\|_{B}<+\infty$. 
\begin{lemma}
\label{traductionEffectiveSpace}
Let $u \in H^{\infty,l,\nu}_b$, we have 
\begin{align*}\chi(x)\chi_1(\sigma)u\left(\frac{x}{\sigma}\right)\in \mathcal{A}([0,1)_\sigma, \sigma^{\min(-l-\frac{3}{2},\nu-\frac{3}{2})}\overline{H}_{(b)}^{\infty, l} )\cap \mathcal{A}([0,1)_\sigma, \sigma^{\nu-\frac{3}{2}}\overline{H}_{(b)}^{\infty, \min(l, -\nu)})
\end{align*} (here $\chi,\chi_1$ are smooth cutoffs with $\chi = \chi_1 = 1$ in a neighborhood of $0$ and $\chi_1$ has compact support in $[0,1)$ while $\chi$ has compact support in $[0, \frac{1}{r_+-\epsilon}$).
\end{lemma}
\begin{proof}
The strategy is to bound a family of quantities $(B_{k,k'}(u))_{k,k'\in \N}$ which dominates the seminorms of $\mathcal{A}([0,1)_\sigma, \sigma^{\min(-l-\frac{3}{2},\nu-\frac{3}{2})}\overline{H}_{(b)}^{\infty, l} )$ and $\mathcal{A}([0,1)_\sigma, \sigma^{\nu-\frac{3}{2}}\overline{H}_{(b)}^{\infty, \min(l, -\nu)})$.
Let $k\in \N$ and $k'\in \N$, we define:
\begin{align*}
B_{k,k'}(u):=\sup_{\sigma\in (0,1)}\int_{0}^{(r_+-\epsilon)^{-1}}\left\|\sigma^{-\nu+\frac{3}{2}}x^{-l}\left(x+\sigma\right)^{\nu+l}(x\partial_x)^{k'}(\sigma\partial_{\sigma})^{k} \chi(x)\chi_1(\sigma)u\left(\frac{x}{\sigma}\right)\right\|^2_{H^{k'}(\mathcal{B}_s)}\frac{\dd x}{x^4}
\end{align*}
Note that for all $N,N'\in \N$, there exists $C_{N,N'}>0$ a constant independent of $u$ such that:
\begin{align*}
\sup_{\sigma \in (0,1)}\left\|\sigma^{-\min(-l-\frac{3}{2},\nu-\frac{3}{2})}(\sigma \partial_{\sigma})^{N}\chi(x)\chi_1(\sigma)u\left(\frac{x}{\sigma}\right)\right\|_{\overline{H}_{(b)}^{N',l}} \leq C_{N,N'} \sum_{\substack{0\leq k\leq N\\ 0\leq k'\leq N'}} B_{k,k'}(u)\\
\sup_{\sigma \in (0,1)}\left\|\sigma^{-\nu+\frac{3}{2}}(\sigma \partial_{\sigma})^{N}\chi(x)\chi_1(\sigma)u\left(\frac{x}{\sigma}\right)\right\|_{\overline{H}_{(b)}^{N',\min(l,-\nu)}} \leq C_{N,N'} \sum_{\substack{0\leq k\leq N\\ 0\leq k'\leq N'}} B_{k,k'}(u)\\
\end{align*}
Moreover, we have (writing $X = \frac{x}{\sigma}$):
\begin{align*}
(x\partial_x)^k(\sigma\partial_\sigma)^{k'}\left(\chi(x)\chi_1(\sigma)u\left(\frac{x}{\sigma}\right)\right) = & \sum_{j=0}^{k}\sum_{m=0}^{k'}\begin{pmatrix}k\\ j\end{pmatrix}\begin{pmatrix}k' \\ m\end{pmatrix}(x\partial_x)^{m}\chi (\sigma\partial_\sigma)^j\chi_1 (-1)^{k-j}(X\partial_X)^{k+k'-m-j}u\left(X\right)
\end{align*}
There exists $C_{k,k'}>0$ such that for all $j\leq k$ and $m\leq k'$:
\begin{align*}
\sup_{\substack{\sigma\in (0,1)\\ x\in \left(0, \frac{1}{r_+-\epsilon}\right)}}\left| \begin{pmatrix}k\\ j\end{pmatrix}\begin{pmatrix}k' \\ m\end{pmatrix}(x\partial_x)^{m}\chi (\sigma\partial_\sigma)^j\chi_1 \right| \leq C_{k,k'}
\end{align*} 
Making the change of variable $X = \frac{x}{\sigma}$ in the definition of $B_{k,k'}(u)$, we deduce that there exists $C'_{k,k'}>0$ such that :
\begin{align*}
B_{k,k'}(u) \leq& \sum_{j=0}^{k}\sum_{m=0}^{k'}\int_{0}^{+\infty} C_{j,m}\left\|\left(\frac{X}{X+1}\right)^{-l}(X+1)^{\nu}(X\partial_X)^{k+k'-j-m}u \right\|_{H^{k'}(\mathcal{B}_s)}^2\frac{\dd X}{X^4}\\
\leq& C'_{k,k'} \left\|u\right\|_{H_b^{k+2k', l, \nu}}^2
\end{align*}
\end{proof}

We define 
\begin{align*}
N_{\eff, \lvert s\rvert}^{\pm} := (XD_X)^2 + iXD_X + \lvert s\rvert -2is\left(XD_X + \frac{i}{2}\right)\mp 2X^{-1} \left(XD_X +i\right)
\end{align*}
since for all $\lvert m\rvert \leq \lvert s\rvert$ and $u \in H_b^{\tilde{r},l,\nu}([0,+\infty]_X)$, we have $N_{\eff}^{\pm}uf_{\lvert s\rvert, m}=(N_{\eff, \lvert s\rvert}^{\pm}u)f_{\lvert s\rvert, m}$ we deduce that $N_{\eff, \lvert s\rvert}^{\pm}$ are invertible between $\left\{u\in H^{r,l,\nu}_b([0,+\infty]_X): N_{\eff,\lvert s\rvert}^{\pm}u\in H^{r, l-1, \nu}([0,+\infty]_X)\right\}$ to $H^{r,l-1, \nu}_b([0,+\infty]_X)$ where $l<-\frac{1}{2}$ and  $r+l>-\frac{1}{2}-2s$ and $\nu \in \left(\frac{1}{2}+s-\lvert s\rvert, \frac{3}{2}+s+\lvert s\rvert\right)$ (we can also prove it directly in the spirit of the proof of corollary \ref{invertibilityNormalEff}).

We define $\tilde{u}^{\pm} = (\pm 1)^{-s-\lvert s\rvert}(N_{\eff,\lvert s\rvert}^{\pm})^{-1}X^{-s-\lvert s\rvert}$. since $X^{-s-\left|s\right|}\in H_b^{\infty, (-s-\left|s\right|-\frac{3}{2})-, (s+\left|s\right|+\frac{3}{2})-}([0,+\infty]_X)$, we get $\tilde{u}^{\pm} \in H_b^{\infty, (-s-\left|s\right|-\frac{1}{2})-, (s+\left|s\right|+\frac{3}{2})-}$. Moreover, we have $\tilde{u}^{-} = (-1)^{-s-\vert s\rvert}\overline{\tilde{u}^{+}}$.

\begin{remark}\label{devUtildePlus}
Using a normal operator argument on $\chi \tilde{u}^+$ and $(1-\chi)\tilde{u}^+$ (where $\chi$ is a smooth cutoff localizing near $X = 0$) on the model of what is done in the proof of Proposition \ref{normalOperatorArgument}, we deduce the following asymptotic expansions:
\begin{align*}
(1-\chi)\tilde{u}^+(X) =& (1-\chi)X^{-s-\left|s\right|}\left(\frac{\ln(X)}{1+2\left|s\right|}+b\right) + H_b^{\infty, \infty, (s+\left|s\right|+\frac{5}{2})-}\\
\chi\tilde{u}^+(X) =& \chi \left(\sum_{k=1}^{s+\left|s\right|}\frac{(-1)^{k+1}i^k(k-1)!}{2^k\left(s+\left|s\right|-k+1\right)}X^{k-s-\left|s\right|} + X\left(\frac{(-i)^{s+\left|s\right|+1}(s+\left|s\right|)!}{2^{s+\left|s\right|+1}}\ln(X)+b'\right)\right)\\
&+H_b^{\infty, (-\frac{3}{2})-,\infty}
\end{align*}
Where $b,b'$ are complex constants (depending on $s$). Note that the sum in the second line is empty when $s\leq 0$.
\end{remark}

\begin{prop}
\label{expressionRSigmav}
Let $c_f \in Y_{\lvert s\rvert}$ be defined as in Proposition \ref{highDecay}. Let $u^{(0)}(c_f)$ be the unique element of $\text{Ker}(\hat{T}_s(0))\cap \overline{H}^{\infty, -\frac{3}{2}-s-\lvert s\rvert-}_{(b)}$ such that $u^{(0)}(c_f) - x^{-s-\lvert s\rvert}c_f$ is of order $x^{1-s-\lvert s\rvert}$ at $x=0$.
Let $v = x^{-s-\lvert s\rvert}\chi(x)c_f$ with $\chi$ a smooth cutoff localizing near $x=0$. Then for $\sigma$ in a small punctured real neighborhood of zero:
\begin{align*}R(\sigma)v =& \sigma^{-s-\lvert s\rvert}\chi(x) \left(\mathbb{1}_{\sigma>0}\tilde{u}^{+}\left(\frac{x}{\sigma}\right)+(-1)^{-s-\lvert s\rvert}\mathbb{1}_{\sigma<0}\overline{\tilde{u}^{+}}\left(-\frac{x}{\sigma}\right)\right)c_f\\
& + \frac{\ln\lvert\sigma\rvert}{1+2\lvert s\rvert}(\chi(x)x^{-s-\lvert s\rvert}c_f - u^{(0)}(c_f))- H(\sigma)2i\Im(b)\left(\chi(x) x^{-s-\lvert s\rvert}c_f - u^{(0)}(c_f)\right)\\
& + C^{\infty}\left((-1,1)_\sigma, \overline{H}_{(b)}^{\infty, -\frac{1}{2}-s-\lvert s\rvert-\epsilon}\right)+  \sigma^{\epsilon-}W^{\infty,\infty}_b(\overline{H}_{(b)}^{\infty,-\frac{1}{2}-s-\lvert s\rvert-\epsilon -})
\end{align*} for some $\epsilon \in (0,1)$. In the previous expression, $H$ is the Heaviside function and $b$ is a complex constant which appears in the development of $\tilde{u}^+$ (we compute exactly $\Im(b)$ later).
\end{prop}
\begin{remark}
In the proposition, $c_f$ could be replaced by any element of $Y_{\left|s\right|}$.
\end{remark}

\begin{proof}
Note that for every $\sigma>0$, $\chi(x) \tilde{u}^+\left(\frac{x}{\sigma}\right)\in H_{(b)}^{\infty, (-\frac{1}{2}-s-\left|s\right|)-}$. In particular, we have $R(\sigma)\hat{T}_s(\sigma)\chi(x) \tilde{u}^+\left(\frac{x}{\sigma}\right)  = \chi(x) \tilde{u}^+\left(\frac{x}{\sigma}\right)$.
For $\sigma>0$, we deduce $R(\sigma)v = \sigma^{-s-\lvert s\rvert}\chi(x)\tilde{u}^{+}c_f+ R(\sigma)(v-\sigma^{-s-\lvert s\rvert}\hat{T}_s(\sigma)\chi \tilde{u}^{+}c_f)$ and
\begin{align*}
v-\sigma^{-s-\lvert s\rvert}\hat{T}_s(\sigma)\chi \tilde{u}^{+}c_f =& v-\sigma^{-s-\lvert s\rvert}N_{\eff}^{+}\chi\tilde{u}^{+}c_f + \sigma^{-s-\lvert s\rvert}(N_{\eff}^{+}-\hat{T}_s(\sigma))\chi\tilde{u}^{+}c_f\\
=& -\sigma^{-s-\lvert s\rvert}[N(\hat{T}_s(0)), \chi]\tilde{\chi}\tilde{u}^{+}c_f + \sigma^{-s-\lvert s\rvert}(N(\hat{T}_s(0))- \hat{T}_s(0)) \chi \tilde{u}^{+}c_f + \sigma \text{Diff}^1_b \tilde{\chi} \tilde{u}^{+}c_f\\
\end{align*}
where $\tilde{\chi} = 1$ on $\text{supp}\chi$.
Note that by the mapping properties of $N_{\eff}^{+}$, 
\begin{align*}\tilde{u}^{+} \in H_b^{\infty,(-\frac{1}{2}-s-\lvert s\rvert)-, (s+\lvert s\rvert + \frac{3}{2})-}([0,+\infty]_X).
\end{align*} By lemma \ref{traductionEffectiveSpace}, we have 
\begin{align*}\sigma^{-s-\lvert s\rvert}\chi_1(\sigma)\tilde{\chi}(x)\tilde{u}^{+}c_f \in \mathcal{A}\left([0,1)_\sigma, \lvert \sigma\rvert^{-1-}\overline{H}_{(b)}^{\infty,(-\frac{1}{2}-s-\lvert s\rvert)-}\right)\cap \mathcal{A}\left([0,1)_\sigma, \lvert \sigma\rvert^{0-}\overline{H}_{(b)}^{\infty,(-\frac{3}{2}-s-\lvert s\rvert)-}\right).
\end{align*} In particular, 
\begin{align*}\sigma^{-s-\lvert s\rvert}\chi_1(\sigma)\sigma\text{Diff}^1_b\tilde{\chi}(x)\tilde{u}^{+}c_f \in \mathcal{A}\left([0,1)_\sigma, \lvert \sigma\rvert^{1-}\overline{H}_{(b)}^{\infty,(-\frac{3}{2}-s-\lvert s\rvert)-}\right)
\end{align*} and we can use Proposition \ref{differentiabilityRes} to conclude that 
\begin{align*}\sigma^{-s-\lvert s\rvert}\chi_1(\sigma)R(\sigma)\sigma\text{Diff}^1_b\tilde{\chi}(x)\tilde{u}^{+}c_f \in \lvert \sigma\rvert^{(1-\eta)-}W^{\infty,\infty}_b\left(\overline{H}_{(b)}^{\infty, -\frac{3}{2}-s-\lvert s\rvert+\eta-}\right)
\end{align*} for $\eta \in (0,1)$.
It is therefore part of the error term.

For $\epsilon \in (0,1)$, we have that 
\begin{align}\label{expansion_u_plusResolved}
\tilde{u}^{+}\left(\frac{x}{\sigma}\right)c_f = \chi_2\left(\frac{\lvert \sigma\rvert}{x}\right)\left(\frac{x}{\sigma}\right)^{-s-\lvert s\rvert}\left(\frac{1}{1+2\lvert s\rvert}\ln\left(\frac{x}{\lvert \sigma\rvert}\right)+b\right)c_f + H_b^{\infty, -\frac{1}{2}-s-\lvert s\rvert-, s+\lvert s\rvert +\frac{3}{2}+\epsilon}
\end{align} and therefore by Lemma \ref{traductionEffectiveSpace} we have 
\begin{align}
\sigma^{-s-\lvert s\rvert}\chi_1(\sigma)\tilde{\chi}(x)\tilde{u}^{+}c_f =& \chi_1(\sigma)\tilde{\chi}(x)\chi_2\left(\frac{\lvert \sigma\rvert}{x}\right)x^{-s-\lvert s\rvert}\left(\frac{1}{1+2\lvert s\rvert}\ln\left(\frac{x}{\lvert \sigma\rvert}\right)+b\right)c_f \notag\\ \label{expansion_u_plus}
& + \mathcal{A}\left([0,1)_\sigma, \sigma^{\epsilon}\overline{H}^{\infty, -\frac{3}{2}-s-\lvert s\vert -\epsilon}_{(b)}\right)
\end{align} (and similarly with $\tilde{\chi}$ replaced by $\chi$). Since $[N(\hat{T}_s(0)), \chi] \in x^{\infty}\text{Diff}^1_b$ and $N(\hat{T}_s(0))-\hat{T}_s(0) \in x\text{Diff}^2_b$, we get
\begin{align*}
A:=&-\sigma^{-s-\lvert s\rvert}\chi_1(\sigma)[N(\hat{T}_s(0)), \chi]\tilde{\chi}\tilde{u}^{+}c_f + \sigma^{-s-\lvert s\rvert}\chi_1(\sigma)(N(\hat{T}_s(0))- \hat{T}_s(0)) \chi \tilde{u}^{+}c_f\\
  =& -[N(\hat{T}_s(0)), \chi]\chi_1(\sigma)\tilde{\chi}(x)\chi_2\left(\frac{\lvert \sigma\rvert}{x}\right)x^{-s-\lvert s\rvert}\left(\frac{1}{1+2\lvert s\rvert}\ln\left(\frac{x}{\lvert \sigma\rvert}\right)+b\right)c_f\\
& + (N(\hat{T}_s(0))- \hat{T}_s(0))\chi_1(\sigma)\chi(x)\chi_2\left(\frac{\lvert \sigma\rvert}{x}\right)x^{-s-\lvert s\rvert}\left(\frac{1}{1+2\lvert s\rvert}\ln\left(\frac{x}{\lvert \sigma\rvert}\right)+b\right)c_f\\
& +  \mathcal{A}\left([0,1)_\sigma, \sigma^{\epsilon-}\overline{H}^{\infty, -\frac{1}{2}-s-\lvert s\vert -\epsilon -}_{(b)}\right)
\end{align*}
Using Proposition \ref{differentiabilityRes}, we get $R(\sigma)\mathcal{A}\left([0,1)_\sigma, \sigma^{\epsilon-}\overline{H}^{\infty, -\frac{1}{2}-s-\lvert s\vert -\epsilon -}_{(b)}\right) \subset \sigma^{\epsilon-}W_b^{\infty,\infty}(\overline{H}^{\infty, -\frac{1}{2}-s-\lvert s\vert -\epsilon -}_{(b)})$ and is therefore also an error term.

We write $\chi_2\left(\frac{\sigma}{x}\right) = 1 + \frac{\sigma}{x}\Psi\left(\frac{\sigma}{x}\right)$ with $(z\partial_z)^N\Psi (z)$ smooth and bounded for all $N$. We deduce:
\begin{align*}
\chi_1(\sigma)\chi(x)\chi_2\left(\frac{\lvert \sigma\rvert}{x}\right)x^{-s-\lvert s\rvert}\left(\frac{1}{1+2\lvert s\rvert}\ln\left(\frac{x}{\sigma}\right)+b\right)c_f =& \chi_1(\sigma)\chi(x)x^{-s-\lvert s\rvert}\left(\frac{1}{1+2\lvert s\rvert}\ln\left(\frac{x}{\sigma}\right)+b\right)c_f \\&+ \sigma \mathcal{A}([0,1)_\sigma, \overline{H}^{\infty, -\frac{5}{2}-s-\lvert s\rvert-}_{(b)})
\end{align*}
We can use Proposition \ref{differentiabilityRes} to prove that 
\begin{align*}R(\sigma)(N(\hat{T}_s(0))- \hat{T}_s(0))\sigma \mathcal{A}([0,1)_\sigma, \overline{H}^{\infty, -\frac{5}{2}-s-\lvert s\rvert-}_{(b)}) \subset \sigma^{\epsilon}W_b^{\infty,\infty}(\overline{H}_{(b)}^{\infty, -\frac{1}{2}-s-\lvert s\rvert-\epsilon -})
\end{align*} and similarly 
\begin{align*}R(\sigma)[N(\hat{T}_s(0)), \chi]\sigma \mathcal{A}([0,1)_\sigma, \overline{H}^{\infty, -\frac{5}{2}-s-\lvert s\rvert-}_{(b)}) \subset \sigma^{\epsilon}W_b^{\infty,\infty}(\overline{H}_{(b)}^{\infty, -\frac{1}{2}-s-\lvert s\rvert-\epsilon -}).
\end{align*}

We have 
\begin{align*}(N(\hat{T}_s(0))- \hat{T}_s(0))\chi_1(\sigma)\chi(x)x^{-s-\lvert s\rvert}\left(\frac{1}{1+2\lvert s\rvert}\ln(x)+b\right)c_f \in& \overline{H}_{(b)}^{\infty, -\frac{1}{2}-s-\lvert s\rvert-}\\
&+ C^{\infty}_c((0,1]_\sigma, \overline{H}_{(b)}^{\infty, -\frac{1}{2}-s-\lvert s\rvert-})
\end{align*} and 
\begin{align*}[N(\hat{T}_s(0)), \chi]\chi_1(\sigma)\tilde{\chi}(x)x^{-s-\lvert s\rvert}\left(\frac{1}{1+2\lvert s\rvert}\ln(x)+b\right)\in \overline{H}_{(b)}^{\infty, -\frac{1}{2}-s-\lvert s\rvert-}+ C^{\infty}_c((0,1]_\sigma, \overline{H}_{(b)}^{\infty, -\frac{1}{2}-s-\lvert s\rvert-})
\end{align*} and therefore, by remark \ref{polynomialTerm}, we have that their image by $R(\sigma)$ is in 
\begin{align*}R(0)w^{+}+ C^{\infty}_c((0,1]_\sigma, \overline{H}_{(b)}^{\infty, -\frac{1}{2}-s-\lvert s\rvert-\epsilon-})+ \sigma^{\epsilon}W_b^{\infty,\infty}(\overline{H}_{(b)}^{\infty, -\frac{1}{2}-s-\lvert s\rvert -\epsilon-})
\end{align*} for $\epsilon \in (0,1)$ where 
\begin{align*}w^{+} := -[N(\hat{T}_s(0)),\chi]\tilde{\chi} x^{-s-\lvert s\rvert}\left(\frac{\ln(x)}{1+2\lvert s\rvert}+b\right)c_f+\left(N(\hat{T}_s(0))-\hat{T}_s(0)\right)\chi(x)x^{-s-\lvert s\rvert}\left(\frac{\ln(x)}{1+2\lvert s\rvert}+b\right)c_f.
\end{align*} Note that $w^{+} \in \overline{H}_{(b)}^{\infty,-\frac{1}{2}-s-\lvert s\rvert}$. 

So far, we are reduced to calculating the image by $R(\sigma)$ of
\begin{align*}
-\chi_1(\sigma)\frac{1}{1+2\lvert s\rvert}\ln(\sigma)\left(-[N(\hat{T}_s(0)), \chi]\tilde{\chi}(x)x^{-s-\lvert s\rvert}c_f + (N(\hat{T}_s(0))- \hat{T}_s(0))\chi(x)x^{-s-\lvert s\rvert}c_f\right)\\
= \chi_1(\sigma)\frac{\ln(\sigma)}{1+2\lvert s\rvert}\hat{T}_s(0)\chi(x)x^{-s-\lvert s\rvert}c_f
\end{align*}
where we used $\chi(x)N(\hat{T}_s(0))x^{-s-\lvert s\rvert}c_f = 0$ to get the right-hand side.
Since $\hat{T}_s(0)\chi(x)x^{-s-\lvert s\rvert}c_f\in \overline{H}^{\infty, -\frac{1}{2}-s-\lvert s\rvert-}_{(b)}$ ($x^{-s-\lvert s \rvert}c_f$ is in the kernel of the normal operator of $\hat{T}_s(0)$), we can write:
\begin{align*}
\chi_1(\sigma)\frac{\ln(\sigma)}{1+2\lvert s\rvert}R(\sigma)\hat{T}_s(0)\chi(x)x^{-s-\lvert s\rvert}c_f =& \chi_1(\sigma)\frac{\ln(\sigma)}{1+2\lvert s\rvert}\tilde{v} + \chi_1(\sigma)\ln(\sigma)R(\sigma)\left(\hat{T}_s(\sigma)-\hat{T}_s(0)\right)\tilde{v}
\end{align*}
where $\tilde{v} := R(0)\hat{T}_s(0)\chi(x)x^{-s-\lvert s\rvert}c_f$. By definition, $\tilde{v}$ is the only element in  $\overline{H}^{\infty, -\frac{1}{2}-s-\lvert s\rvert-}_{(b)}$ such that $\hat{T}_s(0)\tilde{v} = \hat{T}_s(0)\chi(x)x^{-s-\lvert s\rvert}c_f$. We write $u^{(0)}(c_f)$ the element of $\text{Ker}(\hat{T}_s(0))\cap \overline{H}^{\infty, -\frac{3}{2}-s-\lvert s\rvert-}_{(b)}$ such that $u^{(0)}(c_f) - x^{-s-\lvert s\rvert}c_f$ is of order $x^{1-s-\lvert s\rvert}$ at $x=0$. With this definition, we have $\tilde{v} =\chi(x)x^{-s-\lvert s\rvert}c_f-u^{(0)}(c_f)$. By proposition \ref{moderateDecay}, we get that $R(\sigma)\left(\hat{T}_s(\sigma)-\hat{T}_s(0)\right)\tilde{v} \in \sigma^{\epsilon}W_b^{\infty,\infty}(\overline{H}^{\infty, -\frac{1}{2}-s-\lvert s\rvert-\epsilon-})$ and is therefore an error term. Finally, modulo error terms we get (for $\sigma>0$ small enough):
\begin{align*}
R(\sigma)v = \sigma^{-s-\lvert s\rvert}\chi(x)\tilde{u}^{+}\left(\frac{x}{\sigma}\right)c_f + \frac{\ln(\sigma)}{1+2\lvert s\rvert}(\chi(x)x^{-s-\lvert s\rvert}c_f - u^{(0)}(c_f))+ R(0)w^{+}
\end{align*}

We perform the same computation to determine the form of $R(\sigma)v$ in a left neighborhood of $0$ (it amounts to replacing $N_{\eff}^{+}$ by $N_{\eff}^{-}$ and $\tilde{u}^{+}$ by $\tilde{u}^{-}$). We find that modulo error terms (for $\sigma<0$ small enough):
\begin{align*}
R(\sigma)v = (-\sigma)^{-s-\lvert s\rvert}\chi(x)\overline{\tilde{u}^{+}}\left(-\frac{x}{\sigma}\right)c_f + \frac{\ln(-\sigma)}{1+2\lvert s\rvert}(\chi(x)x^{-s-\lvert s\rvert}c_f - u^{(0)}(c_f))+ R(0)w^{-}
\end{align*}
where $w^{-}:= -[N(\hat{T}_s(0)),\chi]\tilde{\chi} x^{-s-\lvert s\rvert}\left(\frac{\ln(x)}{1+2\lvert s\rvert}+\overline{b}\right)c_f+\left(N_{\eff}-\hat{T}_s(0)\right)\chi(x)x^{-s-\lvert s\rvert}\left(\frac{\ln(x)}{1+2\lvert s\rvert}+\overline{b}\right)c_f$.

We see that $R(0)w^{+}\mathbb{1}_{\sigma>0} + R(0)w^{-}\mathbb{1}_{\sigma<0}$ has no smooth extension in a neighborhood of $\sigma = 0$ in general. We rewrite it (to get the second line we use that $N(\hat{T}_s(0))x^{-s-\lvert s\rvert}c_f = 0$):
\begin{align*}
R(0)w^{+}\mathbb{1}_{\sigma>0} + R(0)w^{-}\mathbb{1}_{\sigma<0} =& R(0)w^{-}+ H(\sigma)R(0)(w^{+}-w^{-}) \\
=& R(0)w^{-}- H(\sigma)R(0)\hat{T}_s(0)2i\Im(b)\chi(x)x^{-s-\lvert s\rvert}c_f \\
=& R(0)w^{-} - H(\sigma)2i\Im(b)\left(\chi(x) x^{-s-\lvert s\rvert}c_f - u^{(0)}(c_f)\right)
\end{align*}
The term $R(0)w^{-}$ is now constant in $\overline{H}_{(b)}^{\tilde{r},-\frac{1}{2}-s-\lvert s\rvert-}$ and can be put in the error term  $C^{\infty}\left((-1,1)_\sigma, \overline{H}_{(b)}^{\infty, -\frac{1}{2}-s-\lvert s\rvert-\epsilon}\right)$.
\end{proof}

\begin{remark}
Contrary to what it seems, the principal term (i.e. $R(\sigma)v$ modulo $C^{\infty}\left((-1,1)_\sigma, \overline{H}_{(b)}^{\infty, -\frac{1}{2}-s-\lvert s\rvert-\epsilon}\right)+  \sigma^{\epsilon-}W^{\infty,\infty}_b(\overline{H}_{(b)}^{\infty,-\frac{1}{2}-s-\lvert s\rvert-\epsilon -})$) does not depend on the choice of $\chi$. Indeed, if $\chi_1$ and $\chi_2$ are two smooth cutoffs wit compact support in $\left[0,\frac{1}{r_+}\right)_x$ and equal to $1$ in a neighborhood of $0$, we define the difference:
\begin{align*}
D_{\chi_1,\chi_2}:=&\sigma^{-s-\lvert s\rvert}(\chi_1(x)-\chi_2(x)) \left(\mathbb{1}_{\sigma>0}\tilde{u}^{+}\left(\frac{x}{\sigma}\right)+(-1)^{-s-\lvert s\rvert}\mathbb{1}_{\sigma<0}\overline{\tilde{u}^{+}}\left(-\frac{x}{\sigma}\right)\right)c_f\\
& + \frac{\ln\lvert\sigma\rvert}{1+2\lvert s\rvert}\left(\chi_1(x)-\chi_2(x)\right)x^{-s-\lvert s\rvert}c_f - H(\sigma)2i\Im(b)\left(\chi_1(x)-\chi_2(x)\right) x^{-s-\lvert s\rvert}c_f 
\end{align*}
We have:
\begin{align*}
D_{\chi_1, \chi_2} \in C^{\infty}\left((-1,1)_\sigma, \overline{H}_{(b)}^{\infty, -\frac{1}{2}-s-\lvert s\rvert-\epsilon}\right)+  \sigma^{\epsilon-}W^{\infty,\infty}_b(\overline{H}_{(b)}^{\infty,-\frac{1}{2}-s-\lvert s\rvert-\epsilon -})
\end{align*}
This can be seen using \eqref{expansion_u_plusResolved} combined with Lemma \ref{traductionEffectiveSpace} which provides that for $\sigma$ in a small interval $(0,\eta)$ with $\eta<1$:
\begin{align*}
(\chi_1(x)-\chi_2(x))\tilde{u}^+\left(\frac{x}{\sigma}\right) =& \left(\chi_1(x)-\chi_2(x)\right)\left(\frac{\ln(x)-\ln(\left|\sigma\right|)}{1+2\left|s\right|}+b\right)\left(\frac{x}{\left|\sigma\right|}\right)^{-s-\left|s\right|}c_f\\
& + \left|\sigma\right|^{s+\left|s\right|+\epsilon}\mathcal{A}([0,1)_{\sigma}, H^{\tilde{r},-\frac{1}{2}-s-\left|s\right|-}_{(b)})
\end{align*}
\end{remark}

\begin{definition}
The previous proposition leads us to define
\begin{align*}
\tilde{u}(\sigma,x) := \left(\mathbb{1}_{\sigma>0}\tilde{u}^{+}\left(\frac{x}{\sigma}\right)+(-1)^{-s-\lvert s\rvert}\mathbb{1}_{\sigma<0}\overline{\tilde{u}^{+}}\left(-\frac{x}{\sigma}\right)\right)c_f
\end{align*}
\end{definition}

Combining proposition \ref{highDecay}, proposition \ref{expressionRSigmav}, we get
\begin{coro}
\label{lowEnergyRecap}
Let $f\in C^\infty(\R_\sigma, \overline{H}_{(b)}^{\tilde{r}, l})$ with $l>\frac{1}{2}-s+\lvert s\rvert$, $\tilde{r}-\frac{3}{2}-s-\left|s\right|+\epsilon-2k-1>-\frac{1}{2}-2s$ and $\tilde{r}-2k-1>\frac{1}{2}+s$. There exists $\epsilon>0$ such that we have the following equality for $\sigma$ in a punctured neighborhood of zero:
\begin{align*}
R(\sigma)f =& \sigma^{2-s+\lvert s\rvert}\chi(x) \tilde{u}\left(\sigma,x\right) + \sigma^{2+2\lvert s\rvert}\frac{\ln\lvert\sigma\rvert}{1+2\lvert s\rvert}(\chi(x)x^{-s-\lvert s\rvert}c_f - u^{(0)}(c_f))\\
& -\sigma^{2+2\lvert s\rvert} H(\sigma)2i\Im(b)\left(\chi(x) x^{-s-\lvert s\rvert}c_f - u^{(0)}(c_f)\right)\\
&+\lvert \sigma\rvert^{2\lvert s\rvert +2 + \epsilon}W_b^{k,\infty}(\R_\sigma, \overline{H}_b^{\tilde{r}-2k-1,-\frac{3}{2}-s-\lvert s\rvert+\epsilon})+C^\infty(\R_\sigma, \overline{H}_{(b)}^{\tilde{r}, -\frac{1}{2}-s-\lvert s\rvert-\epsilon})
\end{align*}
for $\chi$ a smooth cutoff localizing near zero and $c_f$ an element of $Y_{\lvert s\rvert}$.
\end{coro}

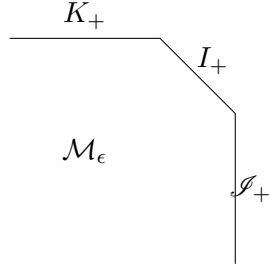
\begin{figure}
\begin{center}
\begin{tikzpicture}
\draw (-1,1)--(1,1)node[midway, yshift = 0.3cm] {$K_+$}--(2,0)node[midway, xshift = 0.2cm,yshift = 0.2cm] {$I_+$}--(2,-2)node[midway, xshift=0.2cm] {$\mathscr{I}_+$};
\draw (0,-0.5)node {$\mathcal{M}_\epsilon$};
\end{tikzpicture}
\end{center}
\caption{Blow up of the corner $x=0$, $\mathfrak{t}^{-1}=0$}\label{FigAsympto}
\end{figure}

We add boundaries to $\overline{\mathcal{M}}_{\epsilon}$ to describe the asymptotic behavior of the the solution for large time.
\begin{definition}
We first add the boundary corresponding to the boundary defining function $\mathfrak{t}^{-1}$ to $\overline{\mathcal{M}}_{\epsilon}$. We can then blow up the corner $\left\{\mathfrak{t}^{-1} = 0, x=0\right\}$. By a slight abuse of notation, we still call $\overline{\mathcal{M}}_{\epsilon}$ the manifold with corners resulting from this procedure. Concretely, the manifold $\overline{\mathcal{M}}_{\epsilon}\cap \left\{\mathfrak{t}>1\right\}$ has now three boundary faces (see Figure \ref{FigAsympto}):
\begin{itemize}
\item $K_+$ with boundary defining function $(x\mathfrak{t}+1)^{-1}$.
\item The front face of the blow up $I_+$ with boundary defining function $\frac{x\mathfrak{t}+1}{\mathfrak{t}}$.
\item The face $\mathscr{I}^+$ (or rather its closure in the blow up space) with boundary defining function $\frac{x\mathfrak{t}}{x\mathfrak{t}+1}$.
\end{itemize}
\end{definition}

We also introduce a new function space to measure decay with respect to $\mathfrak{t}$ and $x$.
\begin{definition}
Let $k\in \N$.
$Z^{k,\tilde{r}}_{\alpha,\beta,\gamma}$ is the space of bounded functions $u$ from $\R_{\mathfrak{t}}$ to $\overline{H}_{(b)}^{\tilde{r},\gamma}$ which have $k$ bounded derivatives and such that for all $j\leq k$:
\begin{align*}
\sup_{\mathfrak{t} \in \R}\left\| (x\mathfrak{t}+1)^{\alpha}\left(\frac{x\mathfrak{t}+1}{\mathfrak{t}}\right)^{-\beta}\left(\frac{x\mathfrak{t}}{x\mathfrak{t}+1}\right)^{-\gamma}(\mathfrak{t}\partial_{\mathfrak{t}})^{j}u\right\|_{\overline{H}_{(b)}^{\tilde{r}, -\frac{3}{2}}}<+\infty
\end{align*}
\end{definition}
\begin{remark}
Note that by Sobolev embedding, $\left\|u\right\|_{C^0}\leq \left\|u\right\|_{\overline{H}_{(b)}^{\frac{3}{2}+,-\frac{3}{2}-}}$ (note that the $-\frac{3}{2}$ offset comes from the fact that we use $sc$ volume form instead of a $b$ volume form). Therefore, if $\tilde{r}>\frac{3}{2}+m$ for some $m\in \N$, then for $u\in Z^{k,\tilde{r}}_{\alpha,\beta,\gamma}$, $j\leq k$ and $p\leq m$ we have the uniform bound
\begin{align*}
\lvert (\mathfrak{t}\partial_{\mathfrak{t}})^{j}(xD_x)^p u\rvert \leq C(x\mathfrak{t}+1)^{-\alpha}\left(\frac{x\mathfrak{t}+1}{\mathfrak{t}}\right)^{\beta}\left(\frac{x\mathfrak{t}}{x\mathfrak{t}+1}\right)^{\gamma}
\end{align*} 
\end{remark}

\begin{remark}
The weights $(x\mathfrak{t}+1)^{-1}$, $\frac{x\mathfrak{t}+1}{\mathfrak{t}}$ and $\frac{x\mathfrak{t}}{x\mathfrak{t}+1}$ measure decay at $K_+$, $I_+$ and $\mathscr{I}^+$ (see Figure \ref{FigAsympto}).
\end{remark}

We will need a slightly more precise version of Lemma 3.6 in \cite{hintz2022sharp}.
\begin{lemma}
\label{boundedFT}
Let $X$ be a Banach space.
Let $\beta>-1$ and $k>\beta+1$, let $\hat{\Phi} \in L^1_{c}((-1,1)_\sigma, X)\cap\lvert\sigma\rvert^{\beta} W_b^{k,\infty}((-1,1)_\sigma\setminus\left\{0\right\}, X)$. For all $j\in \N$ such that $j<k-(\beta+1)$:
\begin{align*}
\lvert (\mathfrak{t}\partial_{\mathfrak{t}})^j\Phi(\mathfrak{t})\rvert \leq C_j\left<\mathfrak{t}\right>^{-1-\beta}
\end{align*} where $\Phi$ is the inverse Fourier transform of $\hat{\Phi}$.

Moreover, if $\hat{\Phi} \in C^p(\R_{\sigma}, X)\cap L^1(\R_{\sigma},X)$ is supported on $\R\setminus [-\epsilon, \epsilon]$ and for $0\leq q\leq p$, $(\sigma\partial_\sigma)^q \hat{\Phi} \in \left<\sigma\right>^{-k+p}L^1(\R_{\sigma},X)$ (for some $k\geq 0$) we have, for all $j\leq k$:
\begin{align*}
\lvert (\partial_{\mathfrak{t}})^j\Phi(\mathfrak{t})\rvert \leq C_j\left<\mathfrak{t}\right>^{-p}
\end{align*}
\end{lemma}
\begin{proof}
The proof of the first claim is formally identical to the proof of Lemma 3.6 in \cite{hintz2022sharp}. The difference is that we only require the lowest conormal regularity permitted by the proof ($\hat{\Phi} \in W_b^{k,\infty}((-1,1)_\sigma\setminus\left\{0\right\}, X)$ instead of $\hat{\Phi}\in W_b^{\infty,\infty}((-1,1)_\sigma\setminus\left\{0\right\}, X)$). We now prove the second claim. Since $\hat{\Phi}$ is integrable, $\Phi$ is bounded with respect to $\mathfrak{t}$. Therefore, it is enough prove the estimate for $\left|t\right|\geq 1$. We have:
\begin{align*}
\Phi(\mathfrak{t}) =& \frac{1}{2\pi}\int_{\R\setminus[-\epsilon, \epsilon]} e^{-i\sigma\mathfrak{t}}\hat{\Phi}(\sigma)\dd \sigma\\
D_{\mathfrak{t}}^j\Phi(\mathfrak{t}) =& \frac{1}{2\pi}\int_{\R\setminus[-\epsilon, \epsilon]} \sigma^{j} e^{-i\sigma\mathfrak{t}}\hat{\Phi}(\sigma)\dd \sigma\\
=& \frac{1}{2\pi}\int_{\R\setminus[-\epsilon, \epsilon]} \sigma^{j} \mathfrak{t}^{-p}(D_\sigma)^{p}e^{-i\sigma\mathfrak{t}}\hat{\Phi}(\sigma)\dd \sigma\\
=& \frac{1}{2\pi}\mathfrak{t}^{-p}\int_{\R\setminus[-\epsilon, \epsilon]}e^{-i\sigma\mathfrak{t}}\sum_{l+m = p}\begin{pmatrix}p\\l\end{pmatrix}D_\sigma^l \sigma^j D_\sigma^m \hat{\Phi}(\sigma)\dd \sigma\\
\end{align*}
Since $D_\sigma^l \sigma^j D_\sigma^m \hat{\Phi}(\sigma)$ is integrable on $\R\setminus [-\epsilon, \epsilon]$ for all $l,m$ such that $l+m=p$ we get for all $\mathfrak{t}\geq 1$:
\begin{align*}
\left\|D_{\mathfrak{t}}^j\hat{\Phi}(\mathfrak{t})\right\|_{X} \leq C\left|\mathfrak{t}\right|^{-p}.
\end{align*}
\end{proof}

For $p,m\in \R$, we define the space $\mathcal{A}^{p,m}(\R_{\sigma}\setminus \left\{0\right\})$ as the set of complex valued tempered distributions $f$ on $\R$ such that for all $k\in \N$:
\begin{align*}
\left<\sigma\right>^{m}\left(\frac{\sigma}{\left<\sigma\right>}\right)^{-p}(\sigma D_{\sigma})^k f \in L^{\infty}(\R)
\end{align*}
We also define the space $S^{m}(\R_{\sigma})$ (symbol of order $m$) as the set of complex valued smooth functions $f$ on $\R$ such that for all $k\in \N$:
\begin{align*}
(D_{\sigma})^k f \leq \left<\sigma\right>^{m-k}
\end{align*}
We will need the following lemma:
\begin{lemma}\label{FTconormal}
Let $m\leq 0$ and $p>-1$.
Let $\hat{\Phi} \in \mathcal{A}^{p,m}(\R\setminus \left\{0\right\})$. Then we have:
$\Phi \in\mathcal{A}^{m-1,p+1}(\R\setminus \left\{0\right\})$.
\end{lemma}
\begin{proof}
We fix $\chi$ a smooth cutoff with compact support equal to 1 near zero. By lemma \ref{boundedFT}, we have that $\chi u$ has Fourier transform in $S^{-(p+1)}(\R)$. Using  $S^{-(p+1)}+\mathcal{A}^{m-1,\infty}\subset \mathcal{A}^{m-1,p+1}$ and $(1-\chi)u\in S^{-m}$, we are reduced to proving that the Fourier transform maps $S^{-m}$ to $\mathcal{A}^{m-1,\infty}$.
This follows for example from Lemma 2.3 (and the estimate (2.7) just after) in \cite{taylor2013partial} for the estimate near zero and from Proposition 8.2 in \cite{taylor2011partial} for the estimate at infinity.
\end{proof}

\begin{lemma}
\label{errorLowEnergy}
We denote by $\mathcal{F}^{-1}$ the inverse Fourier transform.
We have:
$\mathcal{F}^{-1}C^{\infty}_c(\R_\sigma, \overline{H}_{(b)}^{\tilde{r},-\frac{1}{2}-s-\lvert s\rvert-})\subset Z^{\infty, \tilde{r}}_{\infty, \infty, -\frac{3}{2}-s-\lvert s\rvert-}$ and, for $k\geq 3+2\lvert s\rvert +\epsilon$,
\[\mathcal{F}^{-1}\lvert \sigma\rvert^{2\lvert s\rvert +2 + \epsilon}W_b^{k,\infty}((-1,1)\setminus\left\{0\right\}, \overline{H}_{(b)}^{\tilde{r},-\frac{3}{2}-s-\lvert s\rvert +\epsilon})\subset Z^{E(k-3-2\lvert s\rvert -\epsilon), \tilde{r}}_{ 3+2\lvert s\rvert + \epsilon, 3-s+\lvert s\rvert +2\epsilon, -s-\lvert s\rvert +\epsilon}\] (where $E(x)$ denotes the integer part of $x$).
More generally for $\beta>-1$ and $k> \beta+1$, we have:
\begin{align*}
\mathcal{F}^{-1}\lvert \sigma\rvert^{\beta}W_b^{k,\infty}((-1,1)\setminus\left\{0\right\}, \overline{H}_{(b)}^{\tilde{r},l})\subset Z^{E(k-(\beta+1)), \tilde{r}}_{ \beta+1, \beta+l+\frac{5}{2}, l+\frac{3}{2}}
\end{align*}
\end{lemma}

\begin{proof}
These results are obtained by application of lemma \ref{boundedFT}.
\end{proof}

\begin{definition}
We define:
\begin{align*}
u^0_{I_+}(v) := \left(\frac{v}{v+1}\right)^{3+\left|s\right|-s}\mathcal{F}_Y^{-1}\left(Y^{2+\lvert s\rvert -s}\tilde{u}(Y,1)\right)(v)
\end{align*}
Using the fact that $\overline{\mathcal{F}(g)} = \mathcal{F}(\overline{\check{g}})$ (where $\check{g}(Y)=g(-Y)$), we equivalently get:
If $s\in \Z$, we have:
\[u^0_{I_+}(v):=2\left(\frac{v}{v+1}\right)^{3+\lvert s\rvert - s}\Re\left(\mathcal{F}_Y^{-1}H(Y)Y^{2+\lvert s\rvert -s}\tilde{u}^{+}(Y^{-1})\right)c_f\]
If $s\notin \Z$, we have:
\[u^0_{I_+}(v):=2i\left(\frac{v}{v+1}\right)^{3+\lvert s\rvert - s}\Im\left(\mathcal{F}_Y^{-1}H(Y)Y^{2+\lvert s\rvert -s}\tilde{u}^{+}(Y^{-1})\right)c_f\]
The function $u^0_{I_+}$ is smooth on $(0,+\infty)$.
\end{definition}

\begin{remark}\label{expansionUZeroIPlus}
We can compute the asymptotic expansion of $u^{0}_{I_+}$ at $v=0$ and $v=+\infty$. Indeed for $\chi$ a smooth compactly supported cutoff localizing near $0$, using remark \ref{devUtildePlus}, we have: 
\begin{align*}
H(Y)Y^{2+\lvert s\rvert -s}\chi(Y) \tilde{u}^+(Y^{-1}) =& H(Y)Y^{2+2\left|s\right|}\left(-\frac{\ln(Y)}{1+2\left|s\right|}+b\right)+\mathcal{A}^{3+2\left|s\right|-, -2-2\left|s\right|-}(\R\setminus \left\{0\right\})\\
H(Y)Y^{2+\lvert s\rvert -s}(1-\chi(Y)) \tilde{u}^+(Y^{-1}) =& H(Y)\sum_{k=1}^{s+\left|s\right|}\frac{(-1)^{k+1}i^k(k-1)!}{s+\left|s\right|-k+1}Y^{2\left|s\right|+2-k}\\
&+ H(Y)Y^{1+\left|s\right|-s}\left(-\frac{(-i)^{s+\left|s\right|+1}(s+\left|s\right|)!}{2^{s+\left|s\right|+1}}\ln(Y) + b'\right)\\
& + \mathcal{A}^{(1+\left|s\right|-s)-, -(\left|s\right|-s)-}(\R\setminus \left\{0\right\})
\end{align*}
Using Lemma \ref{FTconormal} (and the proof of Proposition \ref{lowEnergyPartFinal} for the computation of $\Im(b)$), we get that:
\begin{itemize}
\item if $s\in \Z$:
\begin{align*}
u^0_{I_+}(v) =& (1-\chi(v))\frac{(-1)^{1+\left|s\right|}(2+2\left|s\right|)!}{1+2\left|s\right|}v^{-3-2\left|s\right|}c_f\\
 &+ \chi(v)(1+\left|s\right|+s)!\left((-1)^{1+\left|s\right|}\frac{(s+\left|s\right|)!}{2^{s+\left|s\right|+2}}+(-1)^{1+\frac{\left|s\right|-s}{2}}\frac{\Re(b')}{\pi}\right)vc_f\\
 &+\mathcal{A}^{2-,(4+2\left|s\right|)-}(\R\setminus \left\{0\right\})c_f
\end{align*}
\item if $s\in \frac{1}{2}\Z$:
\begin{align*}
u^0_{I_+}(v) =& -i(1-\chi(v))\frac{(-1)^{\frac{1+2\left|s\right|}{2}}(2+2\left|s\right|)!}{1+2\left|s\right|}v^{-3-2\left|s\right|}c_f\\
 &+ \chi(v)i(1+\left|s\right|-s)!\left(-\frac{(-1)^{\frac{1+2\left|s\right|}{2}}(s+\left|s\right|)!}{2^{s+\left|s\right|+2}}+(-1)^{\left|s\right|-s}\frac{\Im(i^{2+\left|s\right|-s}b')}{\pi}\right)vc_f\\
&+\mathcal{A}^{2-,(4+2\left|s\right|)-}(\R\setminus \left\{0\right\})c_f
\end{align*}
\end{itemize}
\end{remark}

\begin{lemma}
\label{remainderFTI_+}
We have:
\begin{align*}
2\chi(\mathfrak{t}^{-1})\mathcal{F}^{-1}\left(\chi(x)\chi(\sigma)\sigma^{2+\lvert s\rvert -s}\tilde{u}\left(\sigma,x\right)\right)(\mathfrak{t}) = \chi(x)\chi(\mathfrak{t}^{-1})\left(\frac{x\mathfrak{t} + 1}{\mathfrak{t}}\right)^{3+\lvert s\rvert -s}u^0_{I_+}(x\mathfrak{t}) + R.
\end{align*}
with $R$ smooth with bounds (for all $N,M,K \in N$):
\begin{align*}
\left| \partial_{\mathfrak{t}}^{N}(x\partial_x)^{M} R\right|\leq C_{N,M,K} \left<\mathfrak{t}\right>^{-K}x^{1-}
\end{align*}
In particular, $R\in Z^{\infty,\infty}_{\infty,\infty,1-}$.
\end{lemma}

\begin{proof}
Let $\chi$ be a smooth compactly supported cutoff equal to $1$ near zero.
\begin{align*}
\chi(\mathfrak{t}^{-1})\mathcal{F}^{-1}\left(\chi(x)\chi(\sigma)\sigma^{2+\lvert s\rvert -s}\tilde{u}\left(\sigma,x\right)\right) =& I+I'+R\\
I :=& \chi(\mathfrak{t}^{-1})\chi(x)\mathcal{F}^{-1}\left(H(\sigma)\sigma^{2+\lvert s\rvert -s}\tilde{u}^+\left(\frac{x}{\sigma}\right)\right)c_f\\
I' :=& \chi(x)\chi(\mathfrak{t}^{-1})\mathcal{F}^{-1}\left((-1)^{-s-\left|s\right|}H(-\sigma)\sigma^{2+\lvert s\rvert -s}\overline{\tilde{u}^+\left(\frac{x}{\sigma}\right)}\right)c_f\\
R :=& \chi(\mathfrak{t}^{-1})\mathcal{F}^{-1}\left(\chi(x)(1-\chi(\sigma))\sigma^{2+\lvert s\rvert -s}\tilde{u}\left(\sigma, x\right)\right)\\
\end{align*}
We see in particular that $I' = \overline{I}$ if $s\in \Z$ and $I' = -\overline{I}$ if $s\in \frac{1}{2}\Z$. Moreover, by a change of variable $Y=\frac{\sigma}{x}$ in the Fourier transform, $I$ is exactly $\chi(\mathfrak{t}^{-1})\chi(x)x^{3+\lvert s\rvert -s}\mathcal{F}^{-1}\left(H(Y)Y^{2+\lvert s\rvert -s}\tilde{u}^+(Y)\right)(x\mathfrak{t})$. Therefore, it remains to prove the claimed bound on $\partial_{\mathfrak{t}}^N (x\partial_x)^M R$. It is enough to prove that for all $N,M,K\in \N$ with $K$ large enough, there exists $C_{N,M,K}>0$ such that for all $\mathfrak{t}>0$ and $x>0$:
\begin{align*}
\chi(x)\left|\mathcal{F}^{-1}\left(H(\sigma)(1-\chi(\sigma))\sigma^{N+2+\lvert s\rvert -s}(x\partial_x)^{M}\tilde{u}^+\left(\frac{x}{\sigma}\right)\right)(\mathfrak{t})\right|\leq C_{N,M,K}\mathfrak{t}^{-K}x
\end{align*}
Using the fact that for $0<X<C$, we have $\lvert(X\partial_X)^{N}\tilde{u}^+\rvert\leq C_N X^{1-}$, we deduce that 
\begin{align*}\left|\chi(x)(1-\chi(\sigma))\partial_\sigma^N (x\partial_x)^{M}\tilde{u}^+\left(\frac{x}{\sigma}\right)\right|\leq C_N\sigma^{-N-1+}x^{1-} 
\end{align*} (note that for $x\in \supp(\chi)$ and $\sigma \in \supp(1-\chi)$, we have $0<\frac{x}{\sigma}<C$).
Using these properties for $M\in \N$ and $K>N+2+\lvert s\rvert -s$, we deduce:
\begin{align*}
A:=&\chi(x)\mathcal{F}^{-1}\left(H(\sigma)(1-\chi(\sigma))\sigma^{N+2+\lvert s\rvert -s}(x\partial_x)^{M}\tilde{u}^+\left(\frac{x}{\sigma}\right)\right)(\mathfrak{t})\\
 =& \left(\frac{i}{\mathfrak{t}}\right)^{K} \chi(x)\mathcal{F}^{-1}\left(H(\sigma)\partial_\sigma^K\left((1-\chi(\sigma))\sigma^{N+2+\lvert s\rvert -s}(x\partial_x)^{M}\tilde{u}^+\left(\frac{x}{\sigma}\right)\right)\right)(\mathfrak{t})\\
=& \left(\frac{i}{\mathfrak{t}}\right)^{K} \int_0^{+\infty}e^{i\mathfrak{t}\sigma}\chi(x)\partial_\sigma^K\left((1-\chi(\sigma))\sigma^{N+2+\lvert s\rvert -s}(x\partial_x)^{M}\tilde{u}\left(\frac{x}{\sigma}\right)\right)\dd \sigma
\end{align*}
and we can bound the modulus by $\mathfrak{t}^{-K}x^{1-} C \int_{\inf(\text{supp}(1-\chi))}^{+\infty}\sigma^{N+1+\lvert s\rvert -s-K+}\dd \sigma$.
\end{proof}

\begin{definition}\label{defAsympt}
The previous proposition leads us to define $u_{I_+}:= \left(\frac{x\mathfrak{t} + 1}{\mathfrak{t}}\right)^{3+\lvert s\rvert -s}u^0_{I_+}(x\mathfrak{t})$.

We also define $u_{K_+} := -\chi(\mathfrak{t}^{-1})\mathcal{F}^{-1}\left(\sigma^{2+2\lvert s\rvert}\frac{\ln\lvert\sigma\rvert-i\pi H(\sigma)}{2\lvert s\rvert+1} u^{(0)}(c_f)\right)$.

\end{definition}

We recall the following Fourier transform calculation:
\begin{align*}
u_{K_+} =& -\chi(\mathfrak{t}^{-1})(-i)^{2\lvert s\rvert}\frac{(2+2\lvert s\rvert)!}{2\lvert s\rvert+1} \mathfrak{t}^{-3-2\lvert s\rvert} u^{(0)}(c_f) 
\end{align*}

\begin{lemma}
\label{remainderFTK_+}
If we define
$R:=-\chi(\mathfrak{t}^{-1})\mathcal{F}^{-1}\left((1-\chi(\sigma))\sigma^{2+2\lvert s\rvert}\frac{\ln\lvert\sigma\rvert-i\pi H(\sigma)}{2\lvert s\rvert+1} u^{(0)}(c_f)\right)$
We have the following bounds (for all $N,M,K \in \N$):
\begin{align*}
\left| \partial_{\mathfrak{t}}^{N}(x\partial_x)^{M} R\right|\leq C_{N,M,K} \left<\mathfrak{t}\right>^{-K}x^{-s-\lvert s\rvert}
\end{align*}
In particular, $R\in Z^{\infty,\infty}_{\infty,\infty,-s-\lvert s \rvert-}$.
\end{lemma}
\begin{proof}
Note that using the definition of $u^{(0)}(c_f)$, we have that for all $M \in \N$: $(x\partial_x)^{M}u^{(0)}(c_f)\leq x^{-s-\left|s\right|}$. Moreover, we have that for any $p>K$ $g(\sigma):=(1-\chi(\sigma))\sigma^{2+2\lvert s\rvert}\frac{\ln\lvert\sigma\rvert-i\pi H(\sigma)}{2\lvert s\rvert+1} \in C^{p}(\R_{\sigma}, \C)$ with support away from zero. For all $q\leq p$, we have $(\sigma\partial_{\sigma})^q g(\sigma) \leq C\sigma^{2+2\lvert s\rvert}$ which is in $\sigma^{-N+p}L^1(\R_{\sigma}, \C)$ for $p$ large enough (for example larger than $4+N+2\left|s\right|$). Therefore, we can use Lemma \ref{boundedFT} to conclude the proof.
\end{proof}

\begin{definition} \label{definitionuluh}
Let $\chi$ be a smooth compactly supported cutoff equal to $1$ on $[-1,1]$ and equal to zero on $\R\setminus [-2,2]$.
For $f\in C^{\infty}(\R_\sigma, \overline{H}_{(b)}^{\tilde{r},l})$ with $l>-\frac{3}{2}-s-\left|s\right|$, $\tilde{r}+l>-\frac{1}{2}-2s$ and $\tilde{r}>\frac{1}{2}+s$, we define the low energy part of the solution $\ul := \mathcal{F}^{-1}\chi(\sigma)R(\sigma)f$ and its high energy part $\uh := \mathcal{F}^{-1} (1-\chi(\sigma))R(\sigma)f$.
\end{definition}

\begin{prop}
\label{highEnergyPart}
Let $k\in \N$.
Let $f\in C^{\infty}(\R_\sigma, \overline{H}_{(b)}^{\tilde{r},l})$ with $l+1<-\frac{1}{2}$, $\tilde{r}+l-2k>-\frac{1}{2}-2s$ and $\tilde{r}-2k-1>\max(\frac{1}{2}+s,0)$. We assume in addition that there exists $C>0$ independent of $\sigma$ such that for all $j\leq k$ $(\partial_\sigma)^jf(\sigma)\leq C\left<\sigma\right>^{-1-p-}$.
We have that $\uh\in Z^{p-\tilde{r}+k+1, \tilde{r}-2k-1}_{k,k+l+\frac{5}{2}, l+\frac{5}{2}}$.
\end{prop}
\begin{proof}
By proposition \ref{differentiabilityRes}, for all $j\leq k$, $\partial_{\sigma}^j(1-\chi(\sigma))R(\sigma)f \in \lvert \sigma \rvert^{-p+j}L^1((-\infty, -1)\cup(1,+\infty), \overline{H}_{(b),\lvert \sigma\rvert^{-1}}^{\tilde{r}-2k-1, l+1})$. In particular, since $\left\|u\right\|_{\overline{H}_{(b)}^{\tilde{r}-2k-1, l+1}}\leq \lvert \sigma\rvert^{\tilde{r}-2k-1} \left\|u\right\|_{\overline{H}_{b,\lvert \sigma\rvert^{-1}}^{\tilde{r}-2k-1, l+1}}$ uniformly when $\lvert \sigma \rvert \geq 1$, we have for all $j\leq k$, $\partial_{\sigma}^j(1-\chi(\sigma))R(\sigma)f \in \lvert \sigma \rvert^{-p+\tilde{r}-k-1}L^1((-\infty, -1)\cup(1,+\infty), \overline{H}_{b}^{\tilde{r}-2k-1, l+1})$. We conclude using lemma \ref{boundedFT} that for all $j\leq p-\tilde{r}+k+1$:
\begin{align*}
\left\|\mathfrak{t}^{k} x^{-l-\frac{5}{2}}\partial_{\mathfrak{t}}^j \uh\right\|_{\overline{H}_{(b)}^{\tilde{r}-2k-1,-\frac{3}{2}}}\leq C_j 
\end{align*}
\end{proof}

We can now compute exactly the constant $\Im(b)$ in the expression of $R(\sigma)f$ in a neighborhood of zero by adapting an ingenious causality consideration from \cite[Remark~3.5]{hintz2022sharp}.
\begin{prop}
\label{lowEnergyPartFinal}
Let $f$ be defined as in Corollary \ref{lowEnergyRecap}.
\begin{align*}
R(\sigma)f =& \sigma^{2-s+\lvert s\rvert}\chi(x) \tilde{u}\left(\sigma\right) + \frac{\sigma^{2+2\lvert s\rvert}}{1+2\lvert s\rvert}\left(\ln\lvert\sigma\rvert-i\pi H(\sigma)\right)(\chi(x)x^{-s-\lvert s\rvert}c_f - u^{(0)}(c_f))\\
& +\lvert \sigma\rvert^{2\lvert s\rvert +2 + \epsilon}W_k(\R_\sigma, \overline{H}_{(b)}^{\tilde{r},-\frac{1}{2}-s-\lvert s\rvert-\epsilon -})+C^\infty(\R_\sigma, \overline{H}_{(b)}^{\tilde{r}-2k-1, -\frac{1}{2}-s-\lvert s\rvert-\epsilon})
\end{align*}
\end{prop}
\begin{proof}
Let $\Psi$ be a smooth cutoff localizing away from $x=0$. We consider $\Psi\mathcal{F}^{-1}(R(\sigma)f)$.
By Proposition \ref{highEnergyPart} (more precisely by the proof of the Proposition), we have:
\begin{align*}
\left\|\Psi u_h\right\|_{\overline{H}_{(b)}^{\tilde{r}-k,-\frac{3}{2}}}\leq C\left<t\right>^{-3-2\lvert s\rvert -\epsilon}
\end{align*}
We now look at the low energy part. We have:
\begin{align*}
\Psi u_l =& \mathcal{F}^{-1}\left( \chi(\sigma)\sigma^{2-s+\lvert s\rvert}\chi(x)\Psi(x) \tilde{u}(\sigma)\right.\\
& + \chi(\sigma)\sigma^{2+2\lvert s\rvert}\left(\frac{\ln\lvert\sigma\rvert}{1+2\lvert s\rvert} - 2i\Im(b)H(\sigma)\right)(\chi(x)\Psi(x)x^{-s-\lvert s\rvert}c_f - \Psi(x)u^{(0)}(c_f)) \\
& \left. +\lvert \sigma\rvert^{2\lvert s\rvert +2 + \epsilon}\chi(\sigma)W_k(\R_\sigma, \overline{H}_{(b)}^{\tilde{r},-\frac{1}{2}-s-\lvert s\rvert-\epsilon -})+C_c^\infty(\R_\sigma, \overline{H}_{(b)}^{\tilde{r}, -\frac{1}{2}-s-\lvert s\rvert-\epsilon})\right)
\end{align*}
By lemma \ref{errorLowEnergy} (or rather the analog for negative time), we have (for $\mathfrak{t} \to -\infty$):
\begin{align*}
\Psi u_l =& \mathcal{F}^{-1}\left( \chi(\sigma)\sigma^{2-s+\lvert s\rvert}\chi(x)\Psi(x) \tilde{u}(\sigma) \right.\\
&\left.+ \chi(\sigma)\sigma^{2+2\lvert s\rvert}\left(\frac{\ln\lvert\sigma\rvert}{1+2\lvert s\rvert} - 2i\Im(b)H(\sigma)\right)(\chi(x)\Psi(x)x^{-s-\lvert s\rvert}c_f - \Psi(x)u^{(0)}(c_f))  \right)\\
& + O(\lvert\mathfrak{t}\rvert^{-3-2\lvert s\rvert-\epsilon})_{\overline{H}_{(b)}^{\tilde{r},\infty}}
\end{align*}
Moreover, by \eqref{expansion_u_plus} we have 
\begin{align*}\left( \sigma^{2+2\lvert s\rvert}\left(\frac{\ln\lvert\sigma\rvert}{1+2\lvert s\rvert} - 2i\Im(b)H(\sigma)\right)\chi(x)\Psi(x)x^{-s-\lvert s\rvert}c_f \right.&\\
 \hphantom{(}\left. + \sigma^{2-s+\lvert s\rvert}\chi(x)\Psi(x) \tilde{u}(\sigma)\right)\chi(\sigma)& \in & &\sigma^{2+2\lvert s\rvert}\overline{H}_{(b)}^{-\frac{1}{2}-s-\lvert s\rvert}\\& & &+ \sigma^{2+2\lvert s\rvert}\mathcal{A}\left((-1,1)_\sigma\setminus\left\{0\right\}, \sigma^{\epsilon-}\overline{H}^{\infty, \infty}\right).
\end{align*} Therefore, by lemma \ref{boundedFT}, we have (as $\mathfrak{t}\to -\infty$):
\begin{align*}
\left\|\mathcal{F}^{-1}\left(\chi(\sigma)\chi(x)\Psi(x)\left(\sigma^{2-s+\lvert s\rvert}\tilde{u}(\sigma) + \sigma^{2+2\lvert s\rvert}\left(\frac{\ln\lvert\sigma\rvert}{1+2\lvert s\rvert} - 2i\Im(b)H(\sigma)\right)x^{-s-\lvert s\rvert}c_f\right)\right)\right\|_{\overline{H}^{\tilde{r},\infty}_{(b)}}\\
\leq C\left<t\right>^{-3-2\lvert s\rvert-\epsilon}
\end{align*}

We conclude that 
\begin{align*}
\Psi u =& \mathcal{F}^{-1}\left(-\sigma^{2+2\lvert s\rvert}\left(\frac{\ln\lvert\sigma\rvert}{1+2\lvert s\rvert} - 2i\Im(b)H(\sigma)\right)\Psi(x)u^{(0)}(c_f)\right)+O(\mathfrak{t}^{-3-2\lvert s\rvert -\epsilon})_{\overline{H}_{(b)}^{\tilde{r},\infty}} \\
\end{align*}
Near $-\infty$, we have that:
\begin{align*}
\mathcal{F}^{-1}\left(-\sigma^{2+2\lvert s\rvert}\left(\frac{\ln\lvert\sigma\rvert}{1+2\lvert s\rvert} - 2i\Im(b)H(\sigma)\right)\Psi(x)u^{(0)}(c_f)\right) = (-i)^{2\lvert s\rvert}\frac{(2+2\lvert s\rvert)!}{t^{3+2\lvert s\rvert}}\left(-\frac{\pi}{1+2\lvert s\rvert}+2\Im(b)\right)
\end{align*}

On the other hand, we know that $\Psi u$ vanish for $\mathfrak{t}$ in a neighborhood of $-\infty$. Therefore if $\Psi u^{(0)}(c_f)\neq 0$, we must have:
\begin{align*}
\Im(b) = \frac{\pi}{2(1+2\lvert s\rvert)}
\end{align*}
Note that by definition, $\Im(b)$ does not depend on $f$ and therefore, since we can always arrange $\Psi u^{(0)}(c_f) \neq 0$ for some $f$ the equality is unconditional.
If $\Psi u^{(0)}(c_f) = 0$ for all the cutoffs $\Psi$, we have $u^{(0)}(c_f) = 0$ and therefore $c_f = 0$. As a consequence, the proposition is also true.
\end{proof}

We can now compute explicitly the function $u^0_{I^+}$.
\begin{lemma}\label{explicitUZeroIPlus}
For $v\in (0,+\infty)$, we have:
\begin{align*}
u^0_{I^+}(v) = (-i)^{2+2\left|s\right|}(2\left|s\right|)!\frac{v\left((2\left|s\right|+2)v+2(\left|s\right|-s+1)\right)}{(v+2)^{2+\left|s\right|+s}(v+1)^{3+\left|s\right|-s}}c_f.
\end{align*}
\end{lemma}
\begin{proof}
Starting from the definition of $\tilde{u}^+$, we have for $X>0$:
\begin{align}\label{equ+:1}
(X^{-2-\left|s\right|+s}N_{\eff, \left|s\right|}^+X^{2+\left|s\right|+s}) X^{-2-\left|s\right|+s}\tilde{u}^+(X) = X^{-2-2\left|s\right|}
\end{align}
We define the operator
\begin{align*}
Q =& -Y^2\partial_{Y}^2 + (-2iY+2\left|s\right|+2)Y\partial_{Y} + 2i(\left|s\right|-s+1)Y -2(\left|s\right|+1)
\end{align*}
which is obtained from $X^{-2-\left|s\right|+s}N_{\eff, \left|s\right|}^+X^{2+\left|s\right|-s}$ by the change of variable $Y = X^{-1}$. Therefore \eqref{equ+:1} becomes (for all $Y>0$):
\begin{align}\label{equ+:2}
Q Y^{2+\left|s\right|-s}\tilde{u}^+(Y^{-1}) = Y^{2+2\left|s\right|}
\end{align}
Replacing $Y$ by $-Y$ in \eqref{equ+:2}, and taking the complex conjugate we obtain (for all $Y<0$):
\begin{align}\label{equ+:3}
(-1)^{-\left|s\right|-s}QY^{2+\left|s\right|-s}\overline{\tilde{u}^+}(-Y^{-1}) = Y^{2+2\left|s\right|}
\end{align}
Combining \eqref{equ+:2} and \eqref{equ+:3}, we obtain for every $Y\in \R\setminus\left\{0\right\}$:
\begin{align}\label{equ+:4}
Q\left(Y^{2+\left|s\right|-s}\tilde{u}(Y,1)\right) = Y^{2+2\left|s\right|}
\end{align}
Using Remark \ref{devUtildePlus}, we see that $Y \mapsto QY^{2+\left|s\right|-s}\tilde{u}(Y,1)$ understood in the sense of distributions is in $L^1_{loc}(\R)$. Therefore it does not contain any Dirac term at $Y=0$ and equality \eqref{equ+:4} holds globally in the sense of distributions.

Applying the inverse Fourier transform, we get:
\begin{align*}
\hat{Q}\mathcal{F}_{Y}^{-1}\left(Y^{2+\left|s\right|-s}\tilde{u}(Y,1)\right) = i^{2+2\left|s\right|}\delta_0^{(2+2\left|s\right|)}
\end{align*}
where 
\begin{align*}
\hat{Q} =& -Y(Y+2)\partial_Y^2 + 2(-Y(\left|s\right|+3)+s-\left|s\right|-3)\partial_Y -4\left|s\right|-6.
\end{align*}
Restricting to $(0,+\infty)$ we deduce:
\begin{align*}
\hat{Q}\mathcal{F}_{Y}^{-1}\left(Y^{2+\left|s\right|-s}\tilde{u}(Y,1)\right) = 0.
\end{align*}
The previous equation is hypergeometric. The regular singular point at infinity has indicial roots $2$ and $3+2\left|s\right|$. Using the development of $u^0_{I_+}$ at infinity (see Remark \ref{expansionUZeroIPlus}), we deduce that $\mathcal{F}_Y^{-1}\left(Y^{2+\left|s\right|-s}\tilde{u}(Y,1)\right)$ is equal to the unique (up to a constant complex factor) solution of order $v^{-3-2\left|s\right|}$ at infinity. This solution could be expressed by the mean of the general hypergeometric functions. We can also check directly that:
\begin{align*}
\hat{Q}\frac{(2\left|s\right|+2)v+2(\left|s\right|-s+1)}{v^{2+\left|s\right|-s}(v+2)^{2+\left|s\right|+s}} = 0
\end{align*}
and therefore, $\frac{(2\left|s\right|+2)v+2(\left|s\right|-s+1)}{v^{2+\left|s\right|-s}(v+2)^{2+\left|s\right|+s}}$ is the desired solution.
Finally, we use the development of $u^0_{I_+}$ at infinity to identify the complex factor.
\end{proof}

\begin{prop}
\label{lowEnergyPart}
Let $k\in \N$.
Let $l>\frac{1}{2}-s+\lvert s\rvert$, $\tilde{r}+l-2k-1>-\frac{1}{2}-2s$ and $\tilde{r}-2k-1>\frac{1}{2}+s$. Let $f\in C^{\infty}(\R_\sigma,\overline{H}_{(b)}^{\tilde{r}, l})$ the principal term of $\ul$ at $K_+$ is $u_{K_+}$ and the principal term of $\ul$ at $I_+$ is $u_{I_+}$. More precisely:
\begin{align*}
\chi(\mathfrak{t}^{-1})(\ul-u_{K_+}) \in Z^{E(k-3-2\lvert s\rvert -\epsilon),\tilde{r}-2k-1}_{3+2\lvert s\rvert+\epsilon, 3-s+\lvert s\rvert, -s-\lvert s \rvert-\epsilon}\\
\chi(\mathfrak{t}^{-1})(\ul-u_{I_+})\in Z^{E(k-3-2\lvert s\rvert -\epsilon),\tilde{r}-2k-1}_{(3+2\lvert s\rvert), 3-s+\lvert s\rvert +\epsilon, -s-\lvert s\rvert}
\end{align*}
\end{prop}
\begin{proof}
Using proposition \ref{lowEnergyPartFinal} and lemma \ref{errorLowEnergy}, we see that
\begin{align*}
\chi(\mathfrak{t}^{-1})\ul =& \chi(\mathfrak{t}^{-1})\mathcal{F}^{-1}\left(\sigma^{2-s+\lvert s\rvert}\chi(x) \tilde{u}\left(\sigma\right) + \sigma^{2+2\lvert s\rvert}\frac{\ln\left|\sigma\right|-i\pi H(\sigma)}{1+2\lvert s\rvert}(\chi(x)x^{-s-\lvert s\rvert}c_f - u^{(0)}(c_f))\right)\\
 &+ Z^{E(k-3-2\lvert s\rvert -\epsilon),\tilde{r}-2k-1}_{3+2\lvert s\rvert+\epsilon, 3-s+\lvert s\rvert +\epsilon, -s-\lvert s\rvert}
\end{align*}
Therefore, using lemma \ref{remainderFTI_+} and lemma \ref{remainderFTK_+}, it is enough to prove that:
\begin{align}
\label{firstClaim}
\chi(\mathfrak{t}^{-1})\mathcal{F}^{-1}\left(\chi(\sigma)\sigma^{2-s+\lvert s\rvert}\chi(x) \tilde{u}\left(\sigma\right)\right) \in Z^{E(k-3-2\lvert s\rvert -\epsilon),\tilde{r}-2k-1}_{3+2\lvert s\rvert+\epsilon, 3-s+\lvert s\rvert, -s-\lvert s \rvert-\epsilon}\\
\label{secondClaim}
\chi(\mathfrak{t}^{-1})\mathcal{F}^{-1}\left(\frac{\sigma^{2+2\lvert s\rvert}}{1+2\lvert s\rvert}\left(\ln\left|\sigma\right|-i\pi H(\sigma)\right)(\chi(x)x^{-s-\lvert s\rvert}c_f - u^{(0)}(c_f))\right) \in Z^{E(k-3-2\lvert s\rvert -\epsilon),\tilde{r}-2k-1}_{(3+2\lvert s\rvert), 3-s+\lvert s\rvert +\epsilon, (1-s-\lvert s\rvert)-}
\end{align}
The claim \eqref{firstClaim} follows from \eqref{expansion_u_plus} and Lemma \ref{boundedFT} (we even get the claim with the stronger space $Z^{\infty,\infty}_{3+2\lvert s\rvert +\epsilon, 3+\lvert s\rvert -s, -s-\lvert s\rvert -\epsilon}$).

We now show \eqref{secondClaim}
Note that by definition of $u^{(0)}(c_f)$, $\chi(x)x^{-s-\lvert s\rvert}c_f - u^{(0)}(c_f)\in \overline{H}_{(b)}^{\infty, (-\frac{1}{2}-s-\lvert s\rvert)-}$ and therefore (using the explicit Fourier transform of $\sigma^{2+2\lvert s\rvert}\frac{\ln\left|\sigma\right|-i\pi H(\sigma)}{1+2\lvert s\rvert}$):
\begin{align*}
\chi(\mathfrak{t}^{-1})\mathcal{F}^{-1}\left(\frac{\sigma^{2+2\lvert s\rvert}}{1+2\lvert s\rvert}\left(\ln\left|\sigma\right|-i\pi H(\sigma)\right)(\chi(x)x^{-s-\lvert s\rvert}c_f - u^{(0)}(c_f))\right) \in Z^{\infty, \infty}_{(3+2\lvert s\rvert), (4+s-\lvert s\rvert)-, (1-s-\lvert s\rvert)-}
\end{align*}
\end{proof}

We now compute the principal term at $\mathscr{I}_+$. We follow the argument in \cite{hintz2022sharp}. We use coordinates $v= x\mathfrak{t}$ and $\tau = \mathfrak{t}^{-1}$ which are smooth near $\mathscr{I}_+$ in the blow up space. $v=0$ is a defining function of $\mathscr{I}_+$. The normal operator of $T_s$ at $\mathscr{I}_+$ is:
\begin{align*}
A =& -2v^{-1}\left(v\partial_v - \tau\partial_\tau\right)(v\partial_v -1)
\end{align*}

We denote by $\mathcal{N} = [0,1)_v\times [0,1)_\tau\times \mathbb{S}^2$ (the $\mathbb{S}^2$ part being the boundary at $x=0$).
We denote by $\overline{H}_b^{\tilde{r}, \mu, \nu}\left(\left(\mathcal{B}_s\right)_{|_{\mathcal{N}}}\right) = v^\mu \tau^\nu \overline{H}_b^{\tilde{r},0,0}$ the usual b-Sobolev space (with b volume form) with extendible conditions at $v=1$ and $\tau = 1$ and $\dot{H}_b^{\tilde{r}, \mu, \nu}\left(\left(\mathcal{B}_s\right)_{|_{\mathcal{N}}}\right)$ the same space with supported conditions at $v=1$ and $\tau = 1$.
\begin{lemma}
\label{restrictionZ}
Let $\tilde{r}\in \N$.
The restriction to $\mathcal{N}$ is a continuous map from
$Z^{\tilde{r},\tilde{r}}_{-\infty, \alpha, \beta}$ to $\overline{H}_b^{\tilde{r}, \beta,\alpha-}$
\end{lemma}
\begin{proof}
Let $\tilde{r}\in \N$.
Let $u \in Z^{\tilde{r},\tilde{r}}_{-\infty, \alpha, \beta}$. As a distribution, for all $j,k\in \N$ with $j\leq \tilde{r}$ and $k\leq \tilde{r}$, we have $(x\partial_x)^j(\mathfrak{t}\partial_{\mathfrak{t}})^k u \in Z^{0,0}_{-\infty, \alpha, \beta}$. In particular for $j+k\leq \tilde{r}$, $v:=(-\mathfrak{t}\partial_{\mathfrak{t}}+x\partial_x)^j(x\partial_x)^k u \in Z^{0,0}_{-\infty, \alpha, \beta}$. In particular, for any $\eta>0$ we have:
\begin{align*}
I:=&\int_{\mathfrak{t}\geq 1}\int_{x\leq \mathfrak{t}^{-1}}\mathfrak{t}^{2(\alpha-\beta)-1-\eta}x^{-2\beta}\left|v(\mathfrak{t},x)\right|^2x^{-1}\dd x \dd \mathfrak{t}\\
\leq & \max(2^{2(\alpha-\beta)+N},1)\int_{\mathfrak{t}\geq 1}\mathfrak{t}^{-1-\eta}\int_{x\leq \mathfrak{t}^{-1}}(x\mathfrak{t}+1)^{-N}\left(\frac{x\mathfrak{t}+1}{\mathfrak{t}}\right)^{-2\alpha}\left(\frac{x\mathfrak{t}}{x\mathfrak{t}+1}\right)^{-2\beta}\left|v(\mathfrak{t},x)\right|^2x^{-1}\dd x \dd \mathfrak{t}\\
\leq& \max(2^{2(\alpha-\beta)+N},1)\int_{\mathfrak{t}\geq 1}\mathfrak{t}^{-1-\eta}\left\|v\right\|_{Z^{0,0}_{-N,\alpha,\beta}}\dd \mathfrak{t} \\ \leq& \frac{\max(2^{2(\alpha-\beta)+N},1)}{\eta}\left\|v\right\|_{Z^{0,0}_{-N,\alpha,\beta}}\\
\end{align*}
We perform the change of variable $(\tau, v) = (\mathfrak{t}^{-1}, x\mathfrak{t})$ in the first integral and we find:
\begin{align*}
\int_{\tau\leq 1}\int_{v\leq 1}\tau^{\eta-2\alpha} v^{-2\beta}\left|(\tau\partial_{\tau})^j(v\partial_v)^k u(\tau,v)\right|^2\frac{\dd \tau}{\tau}\frac{\dd v}{v} \leq & \frac{\max(2^{2(\alpha-\beta)+N},1)}{\eta}\left\|v\right\|_{Z^{0,0}_{-N,\alpha,\beta}}
\end{align*}
Since it is true for all $j,k$ such that $j+k\leq \tilde{r}$, we deduce:
\begin{align*}
\left\|u\right\|_{\overline{H}_b^{\tilde{r}, \beta,\alpha-\frac{\eta}{2}}\left(\left(\mathcal{B}_s\right)_{|_{\mathcal{N}}}\right)}\leq C\frac{\max(2^{2(\alpha-\beta)+N},1)}{\eta}\left\|v\right\|_{Z^{0,0}_{-N,\alpha,\beta}}
\end{align*}
\end{proof}

\begin{prop}
\label{estimateTransport}
Let $\tilde{r}\in [0,+\infty]$ and $\beta>\alpha$. 
Let $\chi:\R_v\rightarrow \R$ be a smooth cutoff equal to 1 near 0 and vanishing for $v\geq 1$.
For all $\nu \in (\alpha,\beta]$, there exists $C_\nu>0$ such that for all $u\in \overline{H}_b^{\tilde{r},-N,\beta}\left(\left(\mathcal{B}_s\right)_{_{\mathcal{N}}}\right)$.
\begin{align*}
\left\|\chi u\right\|_{\overline{H}_b^{\tilde{r},\alpha,\beta}}\leq C_{\nu}\left(\left\|u\right\|_{\overline{H}_b^{\tilde{r},-N,\beta}}+\left\|(v\partial_v-\tau\partial_\tau)u\right\|_{\overline{H}_b^{\tilde{r},\nu,\beta}}\right)
\end{align*}
\end{prop}
\begin{proof}
Let $f\in\mathcal{D}:=\left\{ \phi_{|_{(0,1)_v\times (0,1)_\tau\times \mathbb{S}^2}}, \phi\in C^{\infty}_c((0,1)_v\times (0,+\infty)_\tau \times \mathbb{S}^2\right\}$. 
We define $\tilde{H}_b^{\tilde{r},\alpha,\beta}$ the space with extendible condition at $\tau = 1$ and supported condition at $v=0$.
We consider the following map:
\[S:\begin{cases}
\mathcal{D}\rightarrow \tilde{H}_b^{\tilde{r},\alpha,\beta}\\
f\mapsto \int_{0}^{-\ln(v)} f(e^{-s},v\tau e^{s})\dd s
\end{cases}
\]
We have that $u:(v,\tau)\mapsto \int_{0}^{-\ln(v)} f(e^{-s},v\tau e^{s})\dd s$ is well defined as a smooth function on $(0,1)_v\times (0,1)_\tau \times \mathbb{S}^2$ vanishing near $v=1$. We have to prove that it is in $\dot{H}_b^{\tilde{r},\alpha,\beta}$.
By induction, we prove that for $k,j\in \N$ with $k+j\leq \tilde{r}$, there exists coefficients $a_{\mu,\nu}, b_{\mu}$ (independent of $f$) such that
\begin{align*}
(\tau\partial_{\tau})^k(v\partial_v)^{j}u = \sum_{\mu+\nu\leq \tilde{r}} a_{\mu,\nu}(v\partial_v)^\mu(\tau\partial_\tau)^\nu f + \sum_{\mu = 0}^{\tilde{r}}b_{\mu}\int_0^{-\ln(v)} ((x_2\partial_{x_2})^{\mu}f)(e^{-s}, v\tau e^s)\dd s
\end{align*}
(where we have called $x_2$ the second variable of $f$).
The $\dot{H}_b^{0,\alpha,\beta}$ norm of the first sum in the right hand side is bounded by $C\left\|f\right\|_{\dot{H}_b^{\tilde{r},\alpha,\beta}}$. Moreover we have:
\begin{align*}
I:=& \int_0^1\int_0^1 v^{-2\alpha}\tau^{-2\beta}\left|\int_0^{-\ln(v)} ((x_2\partial_{x_2})^{\mu}f)(e^{-s}, v\tau e^s)\dd s\right|^2 \frac{\dd \tau}{\tau}\frac{\dd v}{v}\\
 \leq& \int_0^1\int_0^1 v^{-2\alpha}\tau^{-2\beta}\left|\ln(v)\right|\int_0^{-\ln(v)} \left|((x_2\partial_{x_2})^{\mu}f)(e^{-s}, v\tau e^s)\right|^2\dd s \frac{\dd \tau}{\tau}\frac{\dd v}{v} \\
\end{align*}
By the change of variable $(w,z,v)=(e^{-s},v\tau e^s,v)$, we find:
\begin{align*}
I\leq& \int_{(0,1)}\int_{(w,z)\in (v,1)\times (0,1)}\left|\ln(v)\right| v^{2(\beta-\alpha)}z^{-2\beta}w^{-2\beta}\left|((x_2\partial_{x_2})^{\mu}f)(w, z)\right|^2 \frac{\dd w}{w}\frac{\dd z}{z}\frac{\dd v}{v}\\
\leq& \int_{(0,1)}\left|\ln(v)\right|v^{2\eta} \int_{(w,z)\in (0,1)^2}z^{-2\beta}w^{-2\alpha-2\eta}\left|((x_2\partial_{x_2})^{\mu}f)(w, z)\right|^2 \frac{\dd w}{w}\frac{\dd z}{z}\frac{\dd v}{v}
\end{align*}
which is true for all $\eta \in (0,\beta-\alpha]$
We deduce that for all $\nu \in (\alpha, \beta]$ there exists $C_\nu>0$ such that:
\begin{align*}
I\leq C_{\nu}\left\|f\right\|_{\tilde{H}_b^{\tilde{r},\nu,\beta}}
\end{align*}
Finally, we get $\left\|S(f)\right\|_{\tilde{H}_b^{\tilde{r},\alpha,\beta}}\leq C\left\|f\right\|_{\tilde{H}_b^{\tilde{r},\nu,\beta}}$. By density of $\mathcal{D}$ in $\dot{H}_b^{\tilde{r},\nu,\beta}$, we get that $S$ extends uniquely as a continuous linear map from $\dot{H}_b^{\tilde{r},\beta-1,\beta}$ to $\dot{H}_b^{\tilde{r},\alpha,\beta}$. Moreover, since $(v\partial_v-\tau\partial_\tau)S(f) = f$ for $f\in \mathcal{D}$, the relation stays true in the sense of distributions for $f\in \dot{H}_b^{\tilde{r},\alpha,\beta}$.
Let $u\in  H_b^{\tilde{r},-N,\beta}$ such that $(v\partial_v-\tau\partial_\tau)u \in \overline{H}_b^{\tilde{r},\nu,\beta}$. Let $\tilde{\chi}$ be defined a smooth cutoff equal to 1 on $\text{supp}\chi$ and with support in $[0,1)$. We have $(v\partial_v-\tau\partial_\tau)\chi(v) u = \chi(v)(v\partial_v-\tau\partial_\tau) u+ v\chi'(v)  u$ and since $(v\partial_v-\tau\partial_\tau)$ has no kernel in $\tilde{H}_b^{\tilde{r},-\infty,-\infty}(\mathcal{N})$, we deduce $S(\chi(v)(v\partial_v-\tau\partial_\tau) u+ v\chi'(v) u) = \chi(v) u$. In particular:
\begin{align*}
\left\|\chi(v) u\right\|_{\tilde{H}_b^{\tilde{r},\alpha,\beta}}\leq & C\left(\left\|\chi(v) (v\partial_v-\tau\partial_\tau) u\right\|_{\tilde{H}_b^{\tilde{r},\nu,\beta}} + \left\|v\chi'(v)  u\right\|_{\tilde{H}_b^{\tilde{r},\nu,\beta}}\right)
\end{align*}
Looking at the supports, it implies:
\begin{align*}
\left\|\chi(v) u\right\|_{\overline{H}_b^{\tilde{r},\alpha,\beta}}\leq & C\left(\left\|(v\partial_v-\tau\partial_\tau) u\right\|_{\overline{H}_b^{\tilde{r},\nu,\beta}} + \left\| u\right\|_{\overline{H}_b^{\tilde{r},-\infty,\beta}}\right)
\end{align*}
\end{proof}

\begin{prop} \label{radiationField}
Let $k\in \N$ such that $k\geq 8+2\left|s\right|$.
Let $\tilde{r}-2k-1>4+\max(\frac{1}{2}+s, -2s,0)$ and $l>\frac{1}{2}-s+\lvert s\rvert$. Let $f\in C^{1+p}_c(\R_{\mathfrak{t}}, H_{(b)}^{\tilde{r},l})$ with $p\geq \max(1, \tilde{r}-4-2\left|s\right|)$ and let $u$ be the solution of $T_s u = f$ vanishing for $\mathfrak{t}$ near $-\infty$. Then we have:
\begin{align*}
u-v u_{\text{rad}} \in \overline{H}_b^{\min(k-8-2\left|s\right|,\tilde{r}-2k-5), (3-s+\lvert s\rvert)-, 2-}
\end{align*}
where $u_{\text{rad}}\in \overline{H}_b^{\min(k-8-2\left|s\right|,\tilde{r}-2k-5), (3-s+\lvert s\rvert)-}(\mathscr{I}_+)$ and has a leading order term in $\tau = 0$ equal to $\tau^{-3-\left|s\right|+s}(-i)^{2+2\left|s\right|}(2\left|s\right|)!\frac{\left|s\right|-s+1}{2^{1+\left|s\right|+s}}c_f$
\end{prop}
\begin{proof}
By propositions \ref{lowEnergyPart} and \ref{highEnergyPart}, we have $u\in Z^{k-4-2\left|s\right|,\tilde{r}-2k-1}_{3+2\lvert s\rvert,3-s+\lvert s\rvert,-s-\lvert s\rvert -\epsilon}$ and by lemma \ref{restrictionZ}, we have $u\in \overline{H}_b^{\min(k-4-2\left|s\right|,\tilde{r}-2k-1),-s-\lvert s\rvert -\epsilon, (3-s+\lvert s\rvert)-}$. We also have that:
\begin{align*}
-2v^{-1}\left(v\partial_v-\tau\partial_\tau\right)(v\partial_v -1) u = f + \text{Diff}^2_b u
\end{align*}
where $\text{Diff}^2_b$ is the space of second order $b$-differential operators on $\left(\mathcal{B}_s\right)_{|_{\mathcal{N}}}$. Since $-\frac{v}{2}\left( f + \text{Diff}^2_b u\right) \in \overline{H}_b^{\min(k-4-2\left|s\right|,\tilde{r}-2k-1)-2,1-s-\lvert s\rvert -\epsilon, (3-s+\lvert s\rvert)-}$, we can use proposition \ref{estimateTransport} to get \begin{align*}(v\partial_v-1)u \in \overline{H}_b^{\min(k-4-2\left|s\right|,\tilde{r}-2k-1)-2,(1-s+\lvert s\rvert)-,(3-s+\lvert s\rvert)-}\end{align*}. By inverting $v\partial_v-1$ using the Mellin transform\footnote{See for example the proof of Proposition \ref{refinementAtScri} for a similar argument}, we get $u\in \overline{H}_b^{\min(k-4-2\left|s\right|,\tilde{r}-2k-1)-2,1-, (3-s+\lvert s\rvert)-}$. One more application of proposition \ref{estimateTransport} gives $(v\partial_v-1)u\in \overline{H}_b^{\min(k-4-2\left|s\right|,\tilde{r}-2k-1)-4,2-,(3-s+\lvert s\rvert)-}$. Using the Mellin transform and a contour deformation argument, we find that $u = vu_{\text{rad}}+\overline{H}_b^{\min(k-4-2\left|s\right|,\tilde{r}-2k-1)-4,2-,(3-s+\lvert s\rvert)-}$ with $u_{\text{rad}}\in \overline{H}_b^{\min(k-4-2\left|s\right|,\tilde{r}-2k-1)-4,(3-s+\lvert s\rvert)-}(\mathscr{I}_+)$. Moreover, combining Propositions \ref{highEnergyPart} and \ref{lowEnergyPart}, we also know that $u- u_{I_+}\in Z^{\min(k-4-2\left|s\right|,p-\tilde{r}+k),\tilde{r}-2k-1}_{(3+2\left|s\right|)-, 3-s+\left|s\right|+\epsilon, -s-\left|s\right|}$ and we conclude using the explicit expression of $u_{I_+}$ (see Lemma \ref{explicitUZeroIPlus}).
\end{proof} 

\begin{remark}
\label{BoundAtScrI}
Note that by a slight adaptation of Theorem B.2.7. in \cite{hormander2007analysis}, we get that for $\tilde{r}>\frac{1}{2}+j$ with $j\in \N$, $\dot{H}^{\tilde{r}, \mu, \nu}_b \subset Z^{j, \tilde{r}-j-\frac{1}{2}}_{\infty, \nu, \mu}$. In particular, in the setting of proposition \ref{radiationField}, we deduce that for $j<\min(k-8-2\left|s\right|,\tilde{r}-2k-5)-\frac{1}{2}$, we have $u \in Z^{j, \min(k-8-2\left|s\right|,\tilde{r}-2k-5)-\frac{1}{2}-j}_{3+2\left|s\right|, 3-s+\left|s\right|, 1-}$.
\end{remark}

\section{Contour deformation argument} \label{secContour}
We consider the forcing problem:
\begin{align*}
T_s u = f
\end{align*}
with $f \in C^{1+p}_c(\R_{\mathfrak{t}}, \overline{H}^{\tilde{r},l}_{(b)})$. By the classical hyperbolic estimate, we get that there is a unique solution $u$ to this equation which vanishes in a neighborhood of $\mathfrak{t} = -\infty$. Moreover, there exists $C>0$ such that $\left\|u\right\|_{\overline{H}^{\tilde{r},-N}_{(b)}}\leq Ce^{C\mathfrak{t}}$ (for some $N>0$). In particular, we can define the Laplace transform of $u$, $\hat{u}(\sigma):= \int e^{-i\sigma \mathfrak{t}}u(\mathfrak{t})\dd \mathfrak{t}$ on the domain $\left\{\Im(\sigma)> C\right\}$. The equation on $\hat{u}$ is:
\begin{align*}
\hat{T}_s(\sigma) \hat{u}(\sigma) = \hat{f}(\sigma)
\end{align*}
with $\hat{f}$ holomorphic from $\C$ to $\overline{H}_{(b)}^{\tilde{r},l})$. Moreover, there exists $D>0$ such that for all $k\in \N$ with there exists $C_{k}>0$ such that:
\begin{align}
\label{boundOnf3}
\left\|\partial_{\sigma_x}^{k}\hat{f}(\sigma_x + i\sigma_y)\right\|_{\overline{H}_{(b)}^{\tilde{r},l}}\leq C_{k}\left<\sigma_x\right>^{-(1+p)} e^{D\lvert \sigma_y\rvert}
\end{align}

In particular, if $l>-\frac{3}{2}-s-\left|s\right|$ $\tilde{r}>\frac{1}{2}+s$ and $\tilde{r}+l>-\frac{1}{2}-2s$, we can define $R(\sigma)\hat{f}(\sigma)$ and we have on $\Im(\sigma) = C+1$:
\begin{align*}
\hat{u}(\sigma)=& R(\sigma)\hat{f}(\sigma)\\
u(\mathfrak{t}) =& \int_{\Im(\sigma) = C+1} e^{i\sigma \mathfrak{t}}R(\sigma)\hat{f}(\sigma)\dd \sigma
\end{align*}

Using the fact that $R(\sigma)\hat{f}(\sigma)$ is holomorphic from $\Im(\sigma)>0$ to $\overline{H}_{(b)}^{\tilde{r}-1,l+1}$ (see proposition \ref{regHighFreq}) and \eqref{highFreqInStrip} in proposition \ref{defResolvent} to control the error terms, we find that for every $\epsilon>0$:
\begin{align*}
\int_{\Im(\sigma) = C+1} e^{i\sigma \mathfrak{t}}R(\sigma)\hat{f}(\sigma)\dd \sigma = \int_{\Im(\sigma) = \epsilon}e^{i\sigma \mathfrak{t}}R(\sigma)\hat{f}(\sigma)\dd \sigma
\end{align*}
Moreover using \ref{globalContinuity}, for every $\sigma \in \R$ we have $\lim\limits_{\epsilon \to 0}\left\|R(\sigma +i\epsilon)\hat{f}(\sigma+i\epsilon) - R(\sigma)f(\sigma)\right\|_{\overline{H}_{(b)}^{\tilde{r}-1,l-}} = 0$. If $\tilde{r}-1\geq 0$, it is also true in $\overline{H}_{(b)}^{0,l-}$ and we can combine \ref{globalContinuity} and \ref{boundOnf} to get that there exists $C_p>0$ independent of $\epsilon\in[0,1]$ such that: $\left\|R(\sigma +i\epsilon)\hat{f}(\sigma+i\epsilon) \right\|_{\overline{H}_{(b)}^{0, l-}}\leq C_p\left<\sigma\right>^{-(1+p)}$. In particular if $p>0$, using Lebesgue convergence theorem (see \cite{hille1996functional} for this theorem in the case of Bochner integral, here functions $\sigma \mapsto R(\sigma+i\epsilon)f(\sigma+i\epsilon)$ is continuous in $\overline{H}_{(b)}^{0,l-}$ and therefore strongly measurable in the sense of \cite{hille1996functional}), we get:
\begin{align}\label{FourierReprFormula}
u(\mathfrak{t}) = \int_{\R} e^{i\sigma\mathfrak{t}}R(\sigma)\hat{f}(\sigma) \dd \sigma
\end{align}

We can now state the two main theorems (and their corollaries concerning the Cauchy problem) which are based on this formula and the precise analysis of $R(\sigma)$.
The first theorem concerns the forcing problem with a forcing term having moderate decay at $\mathscr{I}^+$.
\begin{theorem}\label{PreciseThmLowDecay}
Let $l = -\frac{3}{2}+\alpha$ where $\alpha \in (0,1)$. Let $k\geq \alpha+s+\left|s\right|+1$. We assume that $\tilde{r}-2k-1>\max(\frac{1}{2}+s, 0)$ and $\tilde{r}-\frac{3}{2}-s-\left|s\right|>-\frac{1}{2}-2s$. Let $f\in C^{1+p}_c\left(\R_{\mathfrak{t}}, \overline{H}^{\tilde{r},l}_{(b)}\right)$ with $p> \tilde{r}-k-1$.
Let $u$ be the unique solution to the forcing problem $\hat{T}_s u = f$ which vanishes near $\mathfrak{t}=-\infty$. There exists $C>0$ such that for all $\mathfrak{t}\in \R$ and for all $j\leq \min(k-\alpha-s-\left|s\right|+1, p-\tilde{r}+k+1-)$:
\begin{align*}
\left\|(\mathfrak{t}\partial_{\mathfrak{t}})^ju(\mathfrak{t})\right\|_{\overline{H}_{(b)}^{\tilde{r}-2k-1,-\frac{3}{2}-s-\left|s\right|}}\leq C\left<\mathfrak{t}\right>^{-1-\alpha-s-\left|s\right|+}\\
\left\|(\mathfrak{t}\partial_{\mathfrak{t}})^ju(\mathfrak{t})\right\|_{\overline{H}_{(b)}^{\tilde{r}-2k-1,-\frac{1}{2}-}}\leq C\left<\mathfrak{t}\right>^{-\alpha}
\end{align*}
\end{theorem}
\begin{remark}
Using Sobolev injections, we have the uniform boundedness version:
\begin{align*}
\left\|u(\mathfrak{t},x,\omega)\right\|\leq C\min(\left<\mathfrak{t}\right>^{-1-\alpha-s-\left|s\right|+}x^{-s-\left|s\right|-},\left<\mathfrak{t}\right>^{-\alpha+}x^{1-})\\
\leq C\frac{\left<\mathfrak{t}\right>^{-\alpha+}x^{1-}}{1+(\left<\mathfrak{t}\right>x)^{1+s+\left|s\right|}}
\end{align*}
\end{remark}
\begin{proof}
We first use the Fourier representation formula \eqref{FourierReprFormula}.
Then, it is a consequence of Proposition \ref{moderateDecay} (with $l_c = -\frac{3}{2}-s-\left|s\right|-$ for the first estimate and $l_c = -\frac{1}{2}-$ for the second one) and Lemma \ref{errorLowEnergy} (to handle the low energy part) and of Proposition \ref{highEnergyPart} (to handle the high energy part).
\end{proof}

We now use Corollary \ref{FourierTrans2} and Theorem \ref{PreciseThmHighDecay} to deduce the following Corollary concerning the Cauchy problem on the hypersurface $\Sigma_0 := t_0^{-1}(\left\{0\right\})$ (see \eqref{deft_0} for the precise definition of $t_0$). We use the notations of Subsection \ref{CauchyMoreGeneral}.
\begin{coro}\label{preciseCauchyLowDecay}
Let $\alpha \in (0,1)$. Let $k\geq \alpha+s+\left|s\right|+1$. We assume that $\tilde{r}-2k-1>\max(\frac{1}{2}+s, 0)$ and $\tilde{r}-\frac{3}{2}-s-\left|s\right|>-\frac{1}{2}-2s$. Let $p> \tilde{r}-k-1$. 
Let $u_0 \in \overline{H}_b^{\tilde{r}+6+p, 1+\alpha}$ and $u_1 \in \overline{H}_b^{\tilde{r}+5+p, 1+\alpha}$.
Let $u$ be the solution of the Cauchy problem:
\begin{align*}
\begin{cases}
T_s u = 0\\
u_{|_{\Sigma_0}} = u_0\\
\rho_0\nabla^\mu t_0\partial_{\mu}u_{|_{\Sigma_0}} = u_1.
\end{cases}
\end{align*}
There exists $C>0$ such that for all $\mathfrak{t}\in \R$ and for all $j\leq \min(k-\alpha-s-\left|s\right|+1, p-\tilde{r}+k+1-)$:
\begin{align*}
\left\|(\mathfrak{t}\partial_{\mathfrak{t}})^ju(\mathfrak{t})\right\|_{\overline{H}_{(b)}^{\tilde{r}-2k-1,-\frac{3}{2}-s-\left|s\right|}}\leq C\left<\mathfrak{t}\right>^{-1-\alpha-s-\left|s\right|+}\\
\left\|(\mathfrak{t}\partial_{\mathfrak{t}})^ju(\mathfrak{t})\right\|_{\overline{H}_{(b)}^{\tilde{r}-2k-1,-\frac{1}{2}-}}\leq C\left<\mathfrak{t}\right>^{-\alpha+}
\end{align*}
\end{coro}

We go back to the forcing problem and provide a more precise result when the forcing term has higher decay:
\begin{theorem}\label{PreciseThmHighDecay}
Let $l>\frac{1}{2}-s+\left|s\right|$, $k> 8+2\left|s\right|+\frac{1}{2}$, $\tilde{r}-2k-5> \max(\frac{1}{2}+s, -2s, \frac{1}{2})$ and $p\geq \tilde{r}-4-2\left|s\right|$. Let $f\in C^{1+p}_c\left(\R_{\mathfrak{t}}, \overline{H}_{(b)}^{\tilde{r},l}\right)$. Let $u$ be the unique solution to the forcing problem $\hat{T}_s u = f$ which vanishes near $\mathfrak{t}=-\infty$. Let $\chi$ be a smooth compactly supported cutoff equal to $1$ near $0$. Then, with the notations introduced in Definition \ref{defAsympt} there exists $C>0$ and $\epsilon>0$ such that for all $j\in \N$ such that $0\leq j< \min(k-8-2\left|s\right|, \tilde{r}-2k-5)-\frac{1}{2}$:
\begin{align*}
\chi(\mathfrak{t}^{-1})\left(u- \mathfrak{p}(\mathfrak{t},x,\omega)\right)\in Z^{j, \min(k-8-2\left|s\right|, \tilde{r}-2k-5)-j-\frac{1}{2}}_{3+2\left|s\right|+\epsilon,3+\left|s\right|-s+\epsilon,1-}\\
\end{align*}
where 
\begin{align*}
\mathfrak{p}(\mathfrak{t},x,\omega):=(-i)^{2+2\left|s\right|}(2\left|s\right|)!\mathfrak{t}^{-3-2\left|s\right|}\frac{(x\mathfrak{t})^{1+\left|s\right|+s}\left((2\left|s\right|+2)x\mathfrak{t}+2(\left|s\right|-s+1)\right)}{(x\mathfrak{t}+2)^{2+\left|s\right|+s}}u^{(0)}(c_{\hat{f}})
\end{align*}
\end{theorem}

\begin{proof}
We first use the Fourier representation formula \eqref{FourierReprFormula}.
Then, it is a consequence of Proposition \ref{lowEnergyPart} (note that $\mathfrak{p}$ is asymptotic to $u_{K_+}$ near $K_+$ and to $u_{I^+}$ near $I^+$)  and Lemma \ref{errorLowEnergy} (to handle the low energy part), of Proposition \ref{highEnergyPart} (to handle the high energy part) and of remark \ref{BoundAtScrI}.
\end{proof}

Finally, we use Corollary \ref{FourierTrans1} and Theorem \ref{PreciseThmHighDecay} to deduce the following Corollary concerning the Cauchy problem on the hypersurface $\Sigma_0 := t_0^{-1}(\left\{0\right\})$ (see \eqref{deft_0} for the precise definition of $t_0$).
\begin{coro}\label{PreciseCauchyHighDecay}
Let $k> 8+2\left|s\right|+\frac{1}{2}$, $\tilde{r}-2k-5> \max(\frac{1}{2}+s, -2s, \frac{1}{2})$ and $p\geq \tilde{r}-4-2\left|s\right|$. Let $u_0 \in H^{\tilde{r}+\frac{7}{2}+p}(\Sigma_0, \mathcal{B}_s)$, $u_1\in H^{\tilde{r}+\frac{5}{2}+p}(\Sigma_0, \mathcal{B}_s)$ be compactly supported. Let $u$ be the solution of the Cauchy problem:
\[\begin{cases}
T_s u = 0\\
u_{|_{\Sigma_0}} = u_0\\
\nabla^\mu t_0 \partial_\mu u_{|_{\Sigma_0}} = u_1.
\end{cases}\]
Let $\chi$ be a smooth compactly supported cutoff equal to $1$ near $0$.
There exists $C>0$ and $\epsilon>0$ such that for all $j\in \N$ such that $0\leq j< \min(k-8-2\left|s\right|, \tilde{r}-2k-5)-\frac{1}{2}$:
\begin{align*}
\chi(\mathfrak{t}^{-1})\left(u- \mathfrak{p}(\mathfrak{t},x,\omega)\right)\in Z^{j, \min(k-8-2\left|s\right|, \tilde{r}-2k-5)-j-\frac{1}{2}}_{3+2\left|s\right|+\epsilon,3+\left|s\right|-s+\epsilon,1-}\\
\end{align*}
where 
\begin{align*}
\mathfrak{p}(\mathfrak{t},x,\omega):=(-i)^{2+2\left|s\right|}(2\left|s\right|)!\mathfrak{t}^{-3-2\left|s\right|}\frac{(x\mathfrak{t})^{1+\left|s\right|+s}\left((2\left|s\right|+2)x\mathfrak{t}+2(\left|s\right|-s+1)\right)}{(x\mathfrak{t}+2)^{2+\left|s\right|+s}}u^{(0)}(c_{\hat{f}})
\end{align*}
and the function $\hat{f}$ used to compute the constant $c_{\hat{f}}$ is the Fourier transform with respect to $\mathfrak{t}$ of a function $f$ as defined in Proposition \ref{cutoffCauchy}.
\end{coro}

\begin{remark}\label{constantWrtInitialData}
Using the freedom on the choice of the function $f$, we can express $c_{\hat{f}}$ in terms of the initial data. Indeed, by Remark \ref{ChangeOfCutoff}, we can take a sequence of function $\chi_{n} \in C^{\infty}(\R)$ converging to the Heaviside function in the sense of distributions and such that $v_n = \chi_n(t_0) u$ are as in Proposition \ref{cutoffCauchy}. We define:
\begin{align*}
f_n := T_s v_n = [T_s, \chi_n(t_0)]u.
\end{align*}
Since all the $v_n$ are equal near $\mathfrak{t} = +\infty$, the constant $c_{\hat{f}_n}$ does not depend on $n$. Moreover, $\lim\limits_{n\to +\infty} f_n$ exists in the sense of distribution and is equal to $f_{\infty}:=[T_s, H(t_0)]u$ where $H$ is the Heaviside function. Note that the function $f_{\infty}$ only depends on the initial data of the Cauchy problem. Using the definition of $c_{\hat{f}_n}$, we see that $c_{\hat{f}_{\infty}}$ is well defined and $\lim\limits_{n\to +\infty} c_{\hat{f}_n} = c_{\hat{f}_\infty}$. This provides a way to compute the constant $c_{\hat{f}}$ (and therefore $\mathfrak{p}$) only in terms of the initial data.
\end{remark}
\appendix
\section{Hyperbolic estimate with small parameter} \label{semiclassicalHyperboEstSection}
We begin by a general definition
\begin{definition}
Let $E$ be a (complex) vector bundle over $\mathfrak{M}$ and $\Theta$ be a real connection on $E$
We define the gradient operator by:
\begin{align*}
\grad_{g,\Theta}:\begin{cases}
\Gamma(E)\rightarrow \Gamma(E\otimes T\mathfrak{M})\\
u \mapsto g^{i,j} (\Theta_{\partial_j} u)\otimes \partial_i
\end{cases}
\end{align*}
where we used the Einstein summation convention and local coordinates on $\mathfrak{M}$ (but it does not depend on the choice). 

If local coordinates are fixed on an open set $U$, we define the function $J$ on $U$ by $J(x)= \sqrt{\lvert\det(g_{i,j}(x))\rvert}$.
We define the divergence operator by linearity using
\begin{align*}
\Gamma(E\otimes T\mathfrak{M})\rightarrow \Gamma(E)\\
\sum_{i=1}^{n} u^i\otimes \partial_i \mapsto \sum_{i=1}^n J^{-1} \Theta_{\partial_i} J u^i
\end{align*}
The local definition is once again independent of the choice of local coordinates.
We define the operator $\square_{g, \Theta} = \div_{g,\Theta}\grad_{g,\Theta}$.
\end{definition}
\begin{remark}
Let $\nabla_{LC}$ be the Levi-Civita connection on $\mathfrak{M}$.
Since we have the connection $\Theta\otimes\nabla_{LC}$ on $E\otimes F$ where $F$ is a tensor bundle (meaning that $F=  T\mathfrak{M}^{\otimes k} \otimes T^*\mathfrak{M}^{\otimes r}$ for some non negative integers $k$ and $r$), we can extend the definition of $\grad$ and $\div$ naturally to these bundles.
\end{remark}

We introduce the energy momentum tensor:
\begin{definition}
\label{energyMomentumTensor}
Let $u\in \Gamma^2(E)$.
We call energy momentum tensor of $u$ the tensor defined (in abstract index notation) by:
\begin{align*}
T^{\delta,\gamma}(u) = \Re (k(\Theta_\mu u, \Theta_\nu u))g^{\mu,\delta}g^{\nu, \gamma} - \frac{1}{2}g^{\delta,\gamma}g^{\mu, \nu}k(\Theta_\mu u, \Theta_\nu u) + \frac{1}{2}g^{\delta, \gamma}k(u,u)
\end{align*}
We remark that this tensor is symmetric and real.
\end{definition}

\begin{prop}
\label{divT}
We have the following expression for the divergence of $T$:
\begin{align*}
\div(T)^\delta =& g^{\mu, \delta}\Re(k(\Theta_\mu u, \square_{g,\Theta} u))+ \Re(g^{\delta, \mu}g^{\gamma,\nu}k(R^\Theta_{\gamma,\mu}u, \Theta_\nu u))\\
&+ \Re((\Theta_\gamma k)(\Theta_\mu u, \Theta_\nu u))(g^{\mu, \delta}g^{\nu, \gamma}-\frac{1}{2}g^{\delta, \gamma}g^{\mu,\nu})\\
& + \frac{1}{2}g^{\delta,\gamma}(\Theta_\gamma k)(u,u)+g^{\delta,\gamma}\Re(k(\Theta_\gamma u,u))
\end{align*}
where $R^{\Theta}_{\mu,\nu}$ is the curvature tensor associated to the connection $\Theta$.
\end{prop}
\begin{remark}
We note that in the scalar case ($E$ is a trivial bundle of rank 1), and with $\Theta$ the trivial connection and $k$ the canonical Hermitian form on $\C$ we get $\div(T)^{\delta} = g^{\mu, \delta}\Re(\partial_\mu \overline{u} \square_{g} u) + g^{\delta,\gamma}\Re(\overline{u}\partial_\delta u )$ which is the usual result in hyperbolic estimates.
\end{remark}

\begin{proof}
The proposition follows from a local computation. Let $x\in \mathfrak{M}$, we compute $\div(T)^{\delta}(x)$ in normal local coordinates centered at $x$ (we will use the fact that the coordinates are normal many times in the computation):
\begin{align*}
\div(T)^{\delta}(x) =& \partial_\gamma T^{\delta, \gamma} \\
=& g^{\mu,\delta}g^{\nu, \gamma}\Re(\partial_\gamma k(\Theta_\mu u, \Theta_\nu u)) - \frac{1}{2}g^{\delta,\gamma}g^{\mu, \nu}\Re(\partial_\gamma k(\Theta_\mu u, \Theta_\nu u))\\
&+\frac{1}{2}g^{\delta,\gamma}(\Theta_\gamma k)(u,u)+g^{\delta,\gamma}\Re(k(\Theta_\gamma u,u))
\end{align*}
Note that we used the fact that $g^{\mu,\nu}k(\Theta_\mu u, \Theta_\nu u)$ is real to introduce for free a real part in the second term.
We have:
\begin{align*}
\partial_\gamma k(\Theta_\mu u, \Theta_\nu u) = (\Theta_\gamma k)(\Theta_\mu u, \Theta_\nu u)+ k(\Theta_\gamma\Theta_\mu u, \Theta_\nu u)+k(\Theta_\mu u, \Theta_\gamma\Theta_\nu u)
\end{align*}
therefore, we get
\begin{align*}
\div(T)^{\delta}(x) =& g^{\mu, \delta}\Re(k(\Theta_{\mu} u, g^{\nu, \gamma}\Theta_\gamma \Theta_\nu u)) + g^{\mu,\delta}g^{\nu,\gamma}\Re(k(\Theta_\gamma\Theta_\mu u - \Theta_\mu \Theta_\gamma u, \Theta_\nu u)\\
&+\Re((\Theta_\gamma k)(\Theta_\mu u, \Theta_\nu u))(g^{\mu, \delta}g^{\nu, \gamma}-\frac{1}{2}g^{\delta, \gamma}g^{\mu,\nu})\\
&+\frac{1}{2}g^{\delta,\gamma}(\Theta_\gamma k)(u,u)+g^{\delta,\gamma}\Re(k(\Theta_\delta u,u))
\end{align*}
But since the coordinates are normal around $x$, we get at $x$:
\begin{align*}
\square_{g,\Theta} u (x) =& g^{\nu,\gamma}\Theta_\gamma\Theta_\nu u\\
R_{\gamma, \mu}^{\Theta} u (x) =& \Theta_\gamma \Theta_\mu u -\Theta_\mu\Theta_\gamma u
\end{align*}
and this conclude the proof.
\end{proof}

Let $\mathcal{M} = \R_t \times X$ with $X$ a smooth compact manifold with boundary of dimension $n$. Let $g$ be a smooth metric on $\mathcal{M}$ such that $g(\dd t, \dd t) = 1$. Let $E$ be a smooth complex vector bundle of rank $m$ over $X$ with connection $\nabla$ and with a smooth hermitian inner product $k$. Note that $C^{\infty}(\R_t, \Gamma(E))$ is naturally identified with $\Gamma( \pi^*_2(E))$ where $\pi_2$ is the second projection on the product $\R_t\times X$. We denote by $\Theta$, the unique connection on $\Gamma( \pi^*_2(E))$ such that for all $u\in C^{\infty}(\R_t, \Gamma(E))$:
\begin{align*}
\Theta_{\partial_t} u &= \partial_t u\\
\Theta_{V} u &= (t\mapsto \nabla_V u(t)) \text{\quad for $V\in TX$}
\end{align*} 

\begin{definition}
Let $\Sigma_t := \left\{t\right\}\times X$ with volume form $\dd \vol_t = \iota_{\grad\, t}\dd \vol_g$. We define the energy at time $t\in \R$ of $u\in C^{2}(\R_t, \Gamma^2(E))$ by
\begin{align*}
\mathcal{E}(t)[u] =& \int_{\Sigma_t} T_{\delta,\gamma} (\grad\, t)^\delta (\grad\,t)^\gamma \dd \vol_t
\end{align*}
When there is no ambiguity, we use $\mathcal{E}(t)$ for $\mathcal{E}(t)[u]$.
\end{definition}

\begin{lemma}
\label{comparisonEnergy}
Let $u\in C^{2}(\R_t, \Gamma^2(E))$. There exists $C>0$ such that
\begin{align*}
\left\|\grad_{g_{T_1},\nabla} u(t)\right\|^2_{L^2(E)} + \left\|\partial_t u(t)\right\|^2_{L^2(E)}+\left\|u(t)\right\|^2_{L^2(E)}\leq C\mathcal{E}(t)
\end{align*}
uniformly in $t\in [T_1, T_2]$. In the above expression $g_{T_1}$ is the metric induced on $X$ at time $T_1$ and the $L^2$ norm is taken with respect to the inner product of the bundle (for $u(t)$, this inner product is $k$ and for $\grad u(t)$, it is $k\otimes g_{T_1}$) and with the volume form $\dd vol_{T_1}$ (induced by the metric $g_{T_1}$).
\end{lemma}
\begin{proof}
Let $x_0\in X$, and $t_0\in [T_1,T_2]$. We take an orthonormal basis for $g$ at $(t_0,x_0)$ with the first vector equal to $\grad\, t((t_0,x_0))$ and the other vectors $(X_i)_{i=1}^n$ in $T_{x_0} X$. There exists local coordinates $(x^\mu)_{\mu=0}^n$ around $(t_0,x_0)$ such that $(x_i)_{i=1}^n$ are local coordinates on $X$ on a neighborhood $U$ of $x_0$ and such that $\partial_0(t_0, x_0) = \grad\, t$ and $\partial_i(x_0) = X_i$. To alleviate the notations, we write $\Theta_{i}$ for $\Theta_{\partial_i}$. We compute in these local coordinates:
\begin{align*}
T_{\delta,\gamma} (\grad\, t)^\delta (\grad\, t)^\gamma (x_0) =& \frac{1}{2}\Re(k(\Theta_0 u(t_0,x_0), \Theta_0 u(t_0,x_0)))+\sum_{i=0}^n \frac{1}{2}\Re(k(\Theta_i u(t_0, x_0), \Theta_i u(t_0, x_0)))\\
&+\frac{1}{2}k(u,u)\\
\dd vol_{t}(x) = \dd x^1 ... \dd x^n
\end{align*}
Note that at $(t_0, x_0)$:
\begin{align*}
\grad\, t =& \partial_t - \dd x^i(\partial_t)X_i\\
\end{align*}
Therefore,
\begin{align*}
k(\partial_t u, \partial_t u) &= k(\Theta_0 u, \Theta_0 u) - 2\dd x^i(\partial_t)\Re(k(\partial_t u, \Theta_i u))+ \dd x^i(\partial_t)\dd x^j(\partial_t)k(\Theta_i u, \Theta_j u)\\
&\leq (1+C_0)k(\Theta_0 u, \Theta_0 u) + (C_0+C_0^2n)\sum_{i=1}^n k(\Theta_i u, \Theta_i u)
\end{align*}
where $C_0 := \sup_{y\in \mathcal{M}}\left\{ y(\partial_t) \right\}$ where $\mathcal{M}:= \left\{f\in T^* [T_1,T_2]\times X : g(f,f)=1, g(f,\dd t) = 0\right\}$ is a compact set because $-g$ is positive definite on the orthogonal of $\dd t$ (and since $\partial_t$ is a smooth vector field, the function $y \in T^* [T_1,T_2]\times X \mapsto y(\partial_t)$ is smooth therefore bounded on $\mathcal{M}$).
On the other hand we have:
\begin{align*}
\left\|\grad_{g_{T_1},\nabla} u(t_0)\right\|^2_{g\otimes k}(x_0) =& \sum_{1\leq i,j\leq n}\left|g_{(T_1, x_0)}(\dd x^i, \dd x^j)\right|\Re(k(\Theta_i u, \Theta_j u))\\
\leq& \frac{C_1n}{2}\left(\sum_{i=1}^n \Re(k(\Theta_i u, \Theta_i u))\right)\\
\end{align*}
Where the constant $C_1 := (\inf_{(t,y)\in [T_1,T_2]\times S^*X_{-g_{T_1}}} \left\{-g_{t}(y,y)\right\})^{-1}$ where $SX_{-g_{T_1}}$ is the cosphere bundle associated to the riemannian metric $-g_{T_1}$ (this is well defined since $Y:=[T_1,T_2]\times SX_{-g_{T_1}}$ is compact and $-g$ is continuous on $Y$ and positive since $-g_{t}$ is Riemannian). 
Moreover:
\begin{align*}
\dd vol_{T_1} \leq C_2 \dd vol_{t}
\end{align*}
with $C' := (\inf_{(t,y)\in [T_1,T_2]\times SOX_{-g_{T_1}}}\left\{\dd vol_{t}(y)\right\})^{-1}$ where $SOX_{-g_{T_1}}$ is the bundle of oriented orthonormal bases associated to the riemannian metric $-g_{T_1}$ which is compact.
We conclude that,
\begin{align*}
\left\|\grad_{g_{T_1},\nabla} u(t)\right\|^2_{L^2(E)} + \left\|\partial_t u(t)\right\|^2_{L^2(E)}+\left\|u(t)\right\|_{L^2(E)}\leq C_3 \mathcal{E}(t)
\end{align*}
for some constant $C_3$ depending only on $n, C_0, C_1$ and $C_2$ (and therefore uniform in $[T_1, T_2]$).
\end{proof}

\begin{prop}
\label{BrutEnergyEstimate}
For all $u\in C^{2}(\R_t, \Gamma^2(E))$ and for all $T_1<T_2$, there exists $C>0$ independent of $u$ such that:
\begin{align*}
\left| \mathcal{E}(T_2)-\mathcal{E}(T_1) \right|\leq&  C\left(\int_{T_1}^{T_2} \mathcal{E}(t) \dd t + \int_{[T_1, T_2]} \left\|\square_{g,\Theta}u\right\|^2_{k}+\left|\int_{[T_1,T_2]\times \partial X} i^*(\iota_{T\grad\, t}\dd vol_g)\right|\right)
\end{align*}
where the orientation on $[T_1,T_2]\times \partial X$ is given by an outgoing vector and $i: [T_1,T_2]\times \partial X\rightarrow [T_1,T_2]\times X$ is the inclusion map.
\end{prop}
\begin{proof}
We use Stockes' theorem:
\begin{align*}
\int_{[T_1,T_2]\times X} \div_{g} T(\grad\,t) \dd vol = \mathcal{E}(t)-\mathcal{E}(0) + \int_{[T_1,T_2]\times \partial X} i^*(\iota_{T\grad\, t}\dd vol_g)
\end{align*}
Decomposing $\dd vol = \frac{\dd t}{\sqrt{g(\dd t, \dd t)}} \dd vol_t = \dd t \dd vol_t$, we have:
\begin{align*}
\int_{[T_1,T_2]\times X} \div_{g} T(\grad\,t) \dd vol = \int_{T_1}^{T_2} \int_X \div(T)^0 + T(\nabla_{LC}\grad\,t ) \dd vol_{t} \dd t
\end{align*}
Using the the compactness of $[T_1, T_2]\times X$, we get that 
\begin{align*}
\int_X \left|\div(T)^0 + T(\nabla_{LC}\grad\,t )\right| \dd vol_{t} \leq& C\left(\left\|\square_{g,\Theta}u\right\|_{L^2(E)}^2+\left\|\grad_{g_{T_1},\nabla} u(t)\right\|^2_{L^2(E)} + \left\|\partial_t u(t)\right\|^2_{L^2(E)}\right)\\
\leq& C' (\mathcal{E}(t)+ \left\|\square_{g,\Theta}u\right\|_{L^2(E)}^2)
\end{align*} where $C,C'>0$ are uniform with respect to $t\in[T_1,T_2]$.
\end{proof}

Let $X_h$ be a family of compact manifolds with boundary included in $X$ (the parameter $h$ is in some arbitrary set but we will use it with $h\in (0,1]$) and let $E_h$ be the restriction of $E$ to $X_h$. We also define the norm $\mathcal{E}_h$ as $\mathcal{E}$ but with integration on $\Sigma_t\cap \R_t\times X_h$ instead of $\Sigma_t$.
Let $u_h$ be a family of smooth sections of $E_h$.
Let $L_h:\Gamma(E)\rightarrow \Gamma(E)$ be a family of linear operator such that there exists $C>0$ such that for all $t\in [T_1,T_2]$, $\left\|L_hu_h(t)\right\|_{L^2(E)}\leq C\left\|u_h(t)\right\|_{T_1}$ (note that we only require this bound to be uniform for this particular family $u_h$).
\begin{coro}
\label{energyEstimate}
If the term $\left|\int_{[T_1,T_2]\times \partial X_h} i^*(\iota_{T\grad\, t}\dd vol_g)\right|$ vanishes for all $h$ (for example if $X_h$ has no boundary or if $\Theta u_h$ vanishes on $[T_1,T_2]\times\partial X_h$), we have for some $C'>0$ independent of $t$ and $h$ such that:
\begin{align*}
\mathcal{E}_h(t)\leq \left(\mathcal{E}_h(T_1)+ \left\|(\square_{g,\Theta}+L_h)u\right\|^2_{L^2(E_h)}\right)e^{C'(t-T_1)}\\
\end{align*}
\end{coro}
\begin{proof}
It follows from proposition \ref{BrutEnergyEstimate} applied on the manifold with boundary $X_h$. Since they are all included in the compact manifold $X$ and since $g$ is continuous on $[T_1,T_2]\times X$, the constant are uniform with respect to $h$.
\end{proof}

We will use this corollary to deduce the semiclassical energy estimate.

\subsection{Semiclassical hyperbolic estimate}
Let $\mathcal{M} = \R_t\times X$ where $X$ is a smooth compact oriented manifold of dimension $n$ or $\R^n$.
Let $E$ be a smooth complex vector bundle of rank $m$ over $X$ with connection $\nabla$ and with a smooth hermitian inner product $k$. 
Let $g$ be a smooth metric on $\mathfrak{M}:=\R_t\times X \times \R_\tau$ translation invariant with respect to $\tau$ and such that $\dd t$ is timelike. We assume without loss of generality that $g(\dd t,\dd t) = 1$. 
Note that $C^{\infty}(\R_t\times \R_\tau, \Gamma(E))$ is naturally identified with $\Gamma( \pi^*_2(E))$ where $\pi_2$ is the second projection on the product $\R_t\times X\times \R_\tau$. We define a connection $\Theta$ on $\pi^*_2(E)$ by $C^{\infty}(\R_t\times X\times \R_\tau)$-linearity using the identities (for $u\in C^{\infty}(\R_t\times\R_\tau, \Gamma(E))$):
\begin{align*}
\Theta_{\partial_t} u &= \partial_t u\\
\Theta_{\partial_\tau} u &= \partial_\tau u\\
\Theta_{V} u &= ((t,\tau)\mapsto \nabla_V u(t,\tau)) \text{\quad for $V\in TX$}
\end{align*}
We then consider an operator $P_h$ acting on $C^{\infty}(\R_t, \Gamma(E))$ by the following:
\begin{align*}
P_h u = e^{-i\frac{\tau}{h}}h^2\square_{g, \Theta} e^{i\frac{\tau}{h}} u + R_1(t)h\partial_t u+R_2(t) u
\end{align*}
 where $R_2(t)$ is a smooth family of operators in $h\Psi_h^1(X; E)$ and $R_1(t)$ is a smooth family of operators in $h\Psi_h^0(X;E)$. The first term is a priori in $C^{\infty}(\R_t\times\R_\tau, \Gamma(E))$ but is in fact independent of $\tau$ and we identify it with an element of $C^{\infty}(\R_t, \Gamma(E))$.

More explicitly, if we choose local coordinates $(x_1,...,x_n)$ in $X$, $x_0= t$, $x_{n+1} = \tau$, we have: 
\begin{align*}P_h = h^2\sum_{0\leq \mu,\nu \leq n}^n J^{-1}\Theta_{\nu}J g^{\mu,\nu}\Theta_{\mu} + h\sum_{\mu=0}^n \left(g^{\mu,n+1}\Theta_{\mu} + J^{-1}\Theta_\mu J\right)+ g^{n+1,n+1} +R_1(t)h\partial_t+R_2(t)
\end{align*} where $J := \sqrt{\left|\det g\right|}$, $R_2(t)$ is a smooth family of operators in $h\Psi_h^1(X; E)$ and $R_1(t)$ is a smooth family of operators in $h\Psi_h^0(X;E)$. 
Note that the principal part $A_h = h^2\sum_{0\leq \mu,\nu \leq n}^n J^{-1}\nabla_{\nu}J g^{\mu,\nu}\nabla_{\mu} + h\sum_{\mu=0}^n \left(g^{\mu,n+1}\nabla_{\mu} + J^{-1}\nabla_\mu J\right)+ g^{n+1,n+1}$ is invariantly defined (the expression does not depend of the chosen local coordinates and trivializations).

\begin{lemma}
\label{scHyperboEstimate}
Let $T\in (0,+\infty)$.
With the previous notations, for all $s\in \R$, there exists a constant $C_s$ depending on $g,k,s$ and on the symbol norms of $R(t)$ on $[0,T]$ such that for all $\phi\in C^{2}([0,T]; \Gamma(E))$ with $\phi(0) = \partial_t \phi(0) = 0$:
\begin{align*}
\left\|\phi\right\|_{L^2([0,T],H^s_h)}+\left\|h\partial_t \phi\right\|_{L^2([0,T],H^{s-1}_h)}\leq C_s h^{-1}\left\|P_h \phi \right\|_{L^2([0,T],H^{s-1}_h)}
\end{align*}
\end{lemma}
\begin{remark}
If we assume that $s=1$, it will be clear from the proof that the estimate holds for $\phi \in C^2([0,T]; \Gamma^2(E))$.
\end{remark}

\begin{proof}
First we reduce to the case $s=1$. We denote by $g_0$ the metric induced by $g$ on $\left\{t = \tau = 0\right\}$. Let $s>0$ and $Q_s$ be a semiclassical pseudodifferential operator on $\Gamma(E)$ with principal symbol $(1+\lvert\xi\rvert^2_{g_0})^{\frac{s}{4}} I_m$ (where $I_m$ is the identity operator in $\mathcal{L}(E,E)$). Let $G_s := Q_s^*Q_s$ where the adjoint is taken with respect to the metric $g_0$ and the fiber inner product $k$. The principal symbol of $G_s$ is $(\lvert\xi\rvert_{g_0}^{^2}+1)^{\frac{s}{2}} I_m$ and therefore, $G_s$ is semiclassically elliptic of order $s$ and since $X$ is compact, $G_s$ is Fredholm between $H_h^{r}(E)$ and $H_h^{r-s}(E)$ for every $r\in \R$. Since $G_s$ is also formally symmetric, its index is zero and since it is injective for $h$ small enough, $G_s^{-1}$ exists as a (uniformly with respect to $h$) bounded operator from $H_h^{r}(E)$ to $H_h^{r+s}(E)$ for every $r\in \R$. We show that $G_s^{-1}$ is a semiclassical pseudodifferential operator with principal symbol $(\lvert \xi \rvert_{g_0}^{s}+1)^{-1}I_m$.  Indeed, there exists a parametrix $R_s$ for $G_s$ by ellipticity. Therefore,  we have for any $N\in \N$ and $u\in \Gamma(E)$:
\begin{align*}
\left\|(R_s - G_s^{-1}) u\right\|_{H_h^N} &\leq \left\| G_s^{-1}G_s (R_s-G_s^{-1})u\right\|_{H_h^{N}}\\
 &\leq C_r\left\|(G_s R_s - I) u\right\|_{H_h^{N-s}}\\
&\leq C_r h^{N}\left\|u\right\|_{H_h^{-N}}
\end{align*}
As a consequence, $R_s-G_s^{-1}$ which was a priori only a bounded operator from $H_h^{-N}$ to $H_h^{-N+s}$ is in fact a (h uniformly) bounded operator from $h^{-N} H_h^{-N}$ to $H^{N}$ for any $N\in \N$. Therefore, it is an element of $h^{\infty}\Psi_h^{-\infty}(E)$ and $G_s^{-1} \in R_s + h^{\infty}\Psi_h^{-\infty}(E)$ is a pseudodifferential operator in $\Psi_h^{-s}(E)$ with principal symbol $(\lvert \xi \rvert_{g_0}^2+1)^{-s}I_m$. 
 
If the estimate of the lemma is true for $s=1$, we can then apply it to the operator $G_s P_h G_s^{-1}$ which is of the correct form and $G_s u$.

Let $h\in (0,1]$.
In order to get the inequality in the case $s=0$, we apply corollary \ref{energyEstimate} to the section $\tilde{u}_h = e^{\frac{i\tau}{h}}u$ on the manifold with boundary $X\times [0, 2\pi h]_\tau\subset X\times [0, 2\pi]$ between times $T_1=0$ and $T_2 = T$. We check that there exists $C>0$ such that for all $h\in (0,1]$ and all $t\in [0,T]$:
\begin{align*}
\left\|h^{-1}u(t)\right\|^2_{H^1_h(X)}+\left\|\partial_t u(t)\right\|^2_{L^2(X)}\leq C\mathcal{E}_h(t)[\tilde{u}_h]
\end{align*}
Moreover, since $\tilde{u}$ is $2\pi h$ periodic with respect to $\tau$ and since $\dd vol_g$, $g$ and $k$ are $\partial_\tau$ invariant, we have 
\begin{align*}
\int_{[0,T]\times X\times \left\{2\pi h\right\}_\tau}i^*(\iota_{\grad\, T}\dd vol_g)\dd t-\int_{[0,T]\times X\times \left\{0\right\}_\tau}i^*(\iota_{\grad\, T}\dd vol_g)\dd t = 0.
\end{align*}
We also check that with $L_h:= h^{-1}R_1(t)\partial_t  + h^{-2}R_2(t)$, we have $\left\|L_h\tilde{u}_h\right\|_{L^2}\leq C\left\|u\right\|_{H_h^1}\leq C'\left\|\tilde{u}_h\right\|_{\mathcal{H}_0}$ with the constant $C'$ independent of $u$ and $h$.
Therefore, corollary \ref{energyEstimate} gives:
\begin{align*}
\left\|h^{-1}u(t)\right\|^2_{L^2([0,T],H^1_h(E))}+\left\|\partial_t u(t)\right\|^2_{L^2([0,T], L^2(E))}\leq \mathcal{E}_h(t)[\tilde{u}_h]\leq e^{Ct}\left\|(\square_{g,\Theta} + L)\tilde{u}_h\right\|^2_{L^2}
\end{align*}
Multiplying both sides by $h^2$ and integrating on $[0,T]$:
\begin{align*}
\left\|u(t)\right\|^2_{L^2([0,T],H^1_h(E))}+\left\|h\partial_t u(t)\right\|^2_{L^2([0,T], L^2(E))}\leq C'\left(h^{-2}\left\|P_hu\right\|^2_{L^2([0,T],L^2(E))}\right)
\end{align*}
where the constant $C'$ depends on $T$ but is independent of $u$ and $h$.
\end{proof}

We now want to get estimates in the norm $\overline{H}_{h}^s$ (therefore have the same order of derivation in time and space). We need a quantitative and semiclassical version of lemma B.2.9 in \cite{hormander2007analysis} applied to the operator $P_h$. This require finer spaces with two indices of regularity. We first recall the definition of spaces $H_{(m,s)}(\R^n)$ and $\overline{H}_{(m,s)}(\R^n_+)$ following \cite[Section B.2]{hormander2007analysis}. We begin by the definition of $H_{(m,s)}(\R^n)$.
\begin{definition}
For $m,s\in \R$, we denote by $H_{(m,s)}(\R^n)$ the set of tempered distribution $u$ such that:
\begin{align*}
\left<\xi\right>^{m}\left<\xi'\right>^{s}\hat{u}\in L^2(\R^n)
\end{align*}
where $\xi = ( \xi', \xi_n) \in \R^n$ is the dual variable of $x$. When $u\in H_{(m,s)}(\R^n)$, we define the corresponding norm:
\begin{align*}
\left\|u\right\|_{H_{(m,s)}(\R^n)} = \left\|\left<\xi\right>^{m}\left<\xi'\right>^{s}\hat{u}\right|_{L^2}.
\end{align*}
We get the semiclassical version of this space $H_{h,(m,s)}(\R^n)$ by taking the semiclassical Fourier transform 
\begin{align*}\mathcal{F}_h u(\xi) := (2\pi h)^{-\frac{n}{2}}\int e^{-i\frac{x\xi}{h}} u(x)\dd x
\end{align*}
instead of the Fourier transform in the definition (the space is the same but the norm degenerates when $h\to 0$).
\end{definition}
We define $\R^n_+:= \R^{n-1}\times (0,+\infty)$ and we introduce the space $\overline{H}_{(m,s)}(\R^n_+)$:
\begin{definition}
Elements of $\overline{H}_{(m,s)}(\R^n_+)$ are distributions on $R^n_+$ which admits an extension in $H_{(m,s)}(\R^n)$. The norm is given by the infimum over all the possible extensions of the $H_{(m,s)}(\R^n)$ norm of the extension. Replacing $H_{(m,s)}(\R^n)$ by $H_{h,(m,s)}(\R^n)$ in the definition, we get the space $\overline{H}_{h,(m,s)}(R^n_+)$.
\end{definition}
We generalize these spaces in the context of sections of a rank $k$ vector bundle $[T_0,T_1]_t \times E$ where $E$ is a rank $k$ vector bundle over a smooth compact manifold (without boundaries).
\begin{definition}
We consider a covering of $X$ by a finite family of open sets $(U_i)_{i=1}^N$ such that we have local coordinates on each $U_i$ and the bundle $E$ is trivial on $U_i$. We take a partition of unity $\chi_i$ on $X$ subordinated to $(U_i)$. We take $\psi_1, \psi_2 \in C^\infty(\R, [0,1])$ such that $\psi_1+\psi_2 = 1$ and $\psi_1(t) = 0$ on $(-\infty, \frac{2T_0 + T_1}{3}]$, $\psi_2 = 0$ on $(\frac{T_0+2T_1}{3}, +\infty)$. Note that modulo a translation and/or reflection in the $t$ variable, for any distributional section $u$ of $E_{|_{(T_0,T_1)\times X}}$, we can identify $\chi_i \Psi_j u$ with a finite family $(v_{i,j}^p)_{\substack{1\leq i\leq n\\ 1\leq j\leq 2\\ 1\leq p\leq k}}$ of distributions on $R^n_+$. The space $\overline{H}_{(m,r)}(E)$ (or $\overline{H}_{(m,r)}$ for short) is the space of $u$ such that all the $v_{i,j}^p$ belongs to $\overline{H}_{(m,r)}(\R^n_+)$. The norm is given by:
\begin{align*}
\left\|u\right\|_{\overline{H}_{(m,s)}}^2 = \sum_{i,j,p} \left\|v_{i,j}^p\right\|_{\overline{H}_{(m,s)}(\R^n_+)}^2
\end{align*}
Note that a different choice for the $\chi_i$, the $\psi_j$ and the local charts and trivializations induces an equivalent norm. Replacing $\overline{H}_{(m,s)}(\R^n_+)$ by $\overline{H}_{h,(m,s)}(\R^n_+)$ in the definition, we get the space $\overline{H}_{h,(m,s)}(E)$.
\end{definition}

We use the following quantitative and semiclassical version of lemma B.2.9 in \cite{hormander2007analysis} applied to the operator $P_h$:
\begin{lemma}
\label{improveTimeReg}
Let $k\in \N$, $m,r\in \R$. If $u\in \overline{H}_{h,(m-k+1,r)}$ and $P_hu\in \overline{H}_{h,(m,r-k)}$, then $u\in \overline{H}_{h, (m+1,r-k)}$ and there exists a constant independent of $u$ such that:
\begin{align*}
\left\|u\right\|_{\overline{H}_{h, (m+1,r-k)}}\leq& C\left( \left\|u\right\|_{\overline{H}_{h, (m-k+1,r)}}+ \left\|P_hu\right\|_{\overline{H}_{h, (m,r-k)}}\right)
\end{align*}
\end{lemma}
\begin{proof}
We prove this result by induction on $k$. For $k=0$, there is nothing to prove. Assume that the result is true for some $k\in \N$. Then let $u\in \overline{H}_{h,(m-k,r)}$ such that $P_hu\in \overline{H}_{h,(m,r-k-1)}$. Using the expression of the operator $P_h$, we get:
\begin{align*}
\left\|h^2\partial_t^2 u\right\|_{\overline{H}_{h,(m-k-1,r-1)}}\leq& \left\|P_h u\right\|_{\overline{H}_{h,(m-k-1,r-1)}} + C\left\|h\partial_t u\right\|_{\overline{H}_{h,(m-k-1,r)}} + \left\|u\right\|_{\overline{H}_{h,(m-k-1,r+1)}}\\
&\leq \left\|P_h u\right\|_{\overline{H}_{h,(m,r-k-1)}} + \left\| u\right\|_{\overline{H}_{h,(m-k,r)}}
\end{align*}
Then, we use the semiclassical version of lemma B.2.3 (on a compact manifold with boundary) in \cite{hormander2007analysis} to get:
\begin{align*}
\left\|h\partial_t u\right\|_{\overline{H}_{h,(m-k,r-1)}}\leq& C\left(\left\|h^2\partial_t^2 u\right\|_{\overline{H}_{h,(m-k-1,r-1)}}+ \left\|h\partial_t u\right\|_{\overline{H}_{h,(m-k-1,r)}}\right)\\
\leq& C\left(\left\|P_h u\right\|_{\overline{H}_{h,(m,r-k-1)}} + \left\| u\right\|_{\overline{H}_{h,(m-k,r)}}\right)
\end{align*}
and the same lemma applied to $u$ gives:
\begin{align*}
\left\| u\right\|_{\overline{H}_{h,(m-k+1,r-1)}}\leq& C\left(\left\|h\partial_t u\right\|_{\overline{H}_{h,(m-k,r-1)}}+ \left\| u\right\|_{\overline{H}_{h,(m-k,r)}}\right)\\
\leq& C\left(\left\|P_h u\right\|_{\overline{H}_{h,(m,r-k-1)}} + \left\| u\right\|_{\overline{H}_{h,(m-k,r)}}\right)
\end{align*}
By the induction hypothesis, we conclude that:
\begin{align*}
\left\|u\right\|_{\overline{H}_{h,(m+1, r-k-1)}}\leq& C\left( \left\|u\right\|_{\overline{H}_{h, (m-k+1,r-1)}}+ \left\|P_hu\right\|_{\overline{H}_{h, (m,r-1-k)}}\right)\\
\leq& C\left(\left\|P_h u\right\|_{\overline{H}_{h,(m,r-k-1)}} + \left\| u\right\|_{\overline{H}_{h,(m-k,r)}}\right)
\end{align*}
\end{proof}

Using lemma \ref{improveTimeReg}, we can state the following proposition.
Let $\mathcal{H}_h^s$ be the space of distributional sections of $(-\infty, T)\times E$ which can be obtained from a restriction of a distribution in $H_h^s(\R\times E))$ and which have support in $[0,T)\times X$.
\begin{prop}
\label{finalScHyperboEstimate}
There exists $C>0$ such that for all $u\in \cup_{N\in \N}\mathcal{H}_h^{-N}$:
\begin{align*}
\left\|u\right\|_{\mathcal{H}_h^s}\leq Ch^{-1}\left\|P_h u\right\|_{\mathcal{H}_h^{s-1}}
\end{align*}
In the strong sense that if $P_h u \in \mathcal{H}_h^{s-1}$, then $u\in \mathcal{H}_h^s$ and the inequality holds.
\end{prop}
\begin{proof}
First, note that the combination of lemma \ref{scHyperboEstimate} and lemma \ref{improveTimeReg} gives the proposition if we assume $u\in \Gamma^{4}((-\infty, T]\times E)$ with support in $[0,T]\times X$. 
We will use this estimate on $P_h$ and $P_h^*$ to get the general result. We note that the formal adjoint $P_h^*$ is of the correct form to apply the estimate (and we can put zero initial data at $t=T$ instead of $t=0$). The dual space of $\mathcal{H}^s$ is $\mathcal{G}^{-s}$, the space of distributional sections of $(0, +\infty)\times E^*$ which can be obtained from a restriction of a distribution in $H_h^{-s}(\R \times E^*)$ and which have support in $(0,T]\times X$. Therefore, for $f\in \mathcal{H}^{s-1}$ and $\phi\in \Gamma([0, +\infty)\times E^*)$ with support in $[0,T]\times X$ , we have:
\begin{align*}
\left|\left< \phi, f\right>_{L^2}\right|\leq& \left\|\phi\right\|_{\mathcal{G}_h^{-s+1}}\left\|f\right\|_{\mathcal{H}_h^{s-1}}\\
\leq& Ch^{-1}\left\|P^*_h\phi\right\|_{\mathcal{G}_h^{-s}}\left\|f\right\|_{\mathcal{H}_h^{s-1}}
\end{align*}
Therefore, by Hahn-Banach, there exists $u\in \mathcal{H}_h^{s}$ such that $P_h u = f$ and $\left\|u\right\|_{\mathcal{H}^s}\leq Ch^{-1}\left\|Pu\right\|_{\mathcal{H}^{s-1}}$. By taking $f\in \Gamma((-\infty, T]\times E)$ with support in $[0,T]\times X$, the previous construction with $s$ large enough provides $u$ in $\Gamma^{N}((-\infty, T]\times E)$ with support in $[0,T]\times X$ for an arbitrary large $N\geq 4$ such that $P_h u= f$. This argument applied to $P^*_h$ gives the density of \[\left\{P^*\phi, \phi\in \Gamma^{N}([0, +\infty)\times E^*) \text{ with support in } [0,T]\times X\right\}\] in $\mathcal{G}_h^s$ for every $s$. Therefore, we have uniqueness of $u\in \mathcal{H}_h^{-M}$ such that $Pu = f \in \mathcal{H}_h^{s}$ (for any $M$ and $s$). 
Now we can prove the estimate for a general $u\in \mathcal{H}_h^{-N}$ such that $P_hu\in \mathcal{H}^s_h$. We can take a sequence $f_n$ of sections in  $\Gamma((-\infty, T]\times E)$ with support in $[0,T]\times X$ such that $\lim\limits_{n\to +\infty}f_n = P_hu$. Then, the associated $u_n$ are in $\Gamma((-\infty, T]\times E)$ with support in $[0,T]\times X$ (we can construct $u_n$ in $\Gamma^N$ for any $N$ and we have uniqueness). By the estimate, the sequence is Cauchy in $\mathcal{H}^s_h$ and by uniqueness of $u\in \mathcal{H}^{-N}_h$ such that $Pu = f$, the limit is $u$. Therefore, $u \in \mathcal{H}^s_h$ and the estimate holds. This also prove that the set of sections $u\in \Gamma((-\infty, T]\times E)$ with support in $[0,T]\times X$ is dense for the norm $\left\|u\right\|_{\mathcal{H}^s_h}+\left\|P_hu\right\|_{\mathcal{H}^{s-1}_h}$.
\end{proof}

\section{Absence of kernel for the effective normal operator}
\label{absenceOfKernel}
Let $\mathcal{N}$ be a smooth compact Riemannian manifold (without boundary).
Then $X := [0,+\infty]_x\times \mathcal{N}$ is a smooth manifold with boundary and we denote by $n$ its dimension. The boundary is given by two faces  each of them is diffeomorphic to $\mathcal{N}$ and associated with the smooth boundary defining functions $x$ or $r = x^{-1}$. Let $E$ be a smooth complex line bundle over $\mathcal{N}$ endowed with a hermitian metric $m$ and let $\tilde{E}$ be the semitrivial bundle given by $[0,+\infty]\times E$. 
\begin{definition}\label{defbSpacesTwoEnds}
Let $\tilde{r},l,\nu\in \R$. We define the space $H_b^{\tilde{r},l,\nu}(\tilde{E})$ ($H_b^{\tilde{r},l,\nu}$ for short) as the space of sections of $\tilde{E}$ which are in the usual $b$ space (with a $b$-volume form) $H^{\tilde{r},l+\frac{n}{2}}_b$ near the end $x=0$, in the space $H^{\tilde{r}}$ in any compact region of the interior of $X$ and in the space $H^{r,\nu-\frac{n}{2}}_b$ near the conic end ($x^{-1} = 0$). More concretely for $\tilde{r}\in \N$, if $\chi_1$, $\chi_2 \in C^{\infty}([0,+\infty];[0,1])$ are such that $\chi_1+\chi_2= 1$, $\chi_1 =1$ near $x=0$, $\chi_2 = 1$ near $x^{-1} = 0$, we can define the square of the norm of $H_b^{\tilde{r},0,0}(\tilde{E})$ by the following expression:
\begin{align*}
\left\|x^{-\frac{n}{2}}(x\partial_x)^{\tilde{r}}\chi_1  u\right\|^2_{L^2_b((0,+\infty);L^2(E))}+\left\|x^{-\frac{n}{2}}\chi_1 u\right\|^2_{L^2_b((0,+\infty);H^{\tilde{r}}(E))} \\
 + \left\|x^{-\frac{n}{2}}(x\partial_x)^{\tilde{r}}\chi_2  u\right\|^2_{L^2_b((0,+\infty);L^2(E))}+\left\|x^{-\frac{n}{2}}\chi_2 u\right\|^2_{L^2_b((0,+\infty);H^{\tilde{r}}(E))}
\end{align*}
The spaces for $\tilde{r} \in \R$ are then defined by interpolation and duality. Eventually, we define $H_b^{\tilde{r},l,\nu}:= \left(\frac{x}{x+1}\right)^{l}(x+1)^{-\nu}H_b^{\tilde{r},0,0}$.
 \end{definition}
 The $\frac{n}{2}$ shift corresponds to the fact that we want the index to be coherent with respect to the volume form $x^{-n-1}\dd x \dd vol_{\mathcal{N}}$ instead of a b-volume form.
Let $\beta \in i\R$, $\gamma, \beta' \in \C$, $\zeta \in \C\setminus\left\{0\right\}$ with $\Im(\zeta)\geq 0$.
We consider the following operator on $\tilde{E}$:
\begin{align*}
\tilde{P}(\zeta) =& \left(x^2D_x\right)^2+i(n-1)x^3D_x + x^2 L + \beta\left(x^3D_x + i\frac{n-2}{2}x^2\right)\\
&+\beta'x^2 - 2 \zeta x\left(xD_x + i\frac{n-1}{2} + \frac{\beta-\gamma}{2}\right)
\end{align*}
where $\Re(\beta')>\frac{\beta^2}{4}-\left(\frac{n-2}{2}\right)^2$, $L$ is a non negative elliptic formally selfadjoint (with respect to the metric $m$ and volume form $\dd vol_{\mathcal{N}}$) second order differential operator on $E$ (in particular it is selfadjoint on $L^2(E)$ with domain $H^2(E)$). We denote by $(\lambda_k)_{k\in \N}$ the sequence of its eigenvalues.

\begin{prop}
\label{trivialKernel}
We assume that for all $k\in \N$,
\begin{align*}
\tilde{a}&:= \frac{1}{2}+\frac{i\gamma}{2}+\sqrt{\left(\frac{n-2}{2}\right)^2-\frac{\beta^2}{4}+\lambda_k+\beta'}
\end{align*}
is not a non positive integer and $\Im(\zeta)\geq 0$.
If $u \in H_b^{r, l, \nu}(\tilde{E})$ is in the kernel of $\tilde{P}(\zeta)$ with $l<-\frac{1}{2}+\frac{\Im(\beta-\gamma)}{2}$ and  $r+l>-\frac{1}{2}+\frac{\Im(\beta+\gamma)}{2}$ and $\nu \in \left(1-\frac{\Im(\beta)}{2}-\Re\sqrt{-\frac{\beta^2}{4}+\left(\frac{n-2}{2}\right)^2 + \beta'}, 1-\frac{\Im(\beta)}{2}+\Re\sqrt{-\frac{\beta^2}{4}+\left(\frac{n-2}{2}\right)^2+\beta'}\right)$, then $u = 0$.
\end{prop}

\begin{prop}
\label{trivialKernel2}
We assume that for all $k\in \N$, $\frac{1}{2}-\frac{i\gamma}{2}+\sqrt{\left(\frac{n-2}{2}\right)^2-\frac{\beta^2}{4}+\lambda_k + \beta'}$ is not a negative integer and $\Im(\zeta)\leq 0$.
If $u \in H_b^{r, l, \nu}(\tilde{E})$ is in the kernel of $\tilde{P}(\zeta)$ with $l>-\frac{1}{2}+\frac{\Im(\beta-\gamma)}{2}$ and  $r+l<-\frac{1}{2}+\frac{\Im(\beta+\gamma)}{2}$ and $\nu \in \left(1-\frac{\Im(\beta)}{2}-\Re\sqrt{-\frac{\beta^2}{4}+\left(\frac{n-2}{2}\right)^2 + \beta'}, 1-\frac{\Im(\beta)}{2}+\Re\sqrt{-\frac{\beta^2}{4}+\left(\frac{n-2}{2}\right)^2+\beta'}\right)$, then $u = 0$.
\end{prop}

To prove these propositions, we first prove precise a priori results for $u$. Then we reduce the problem to a singular ordinary differential equation by diagonalizing the operator $L$. Finally, we use standard theory of the confluent hypergeometric equation to prove that $u = 0$.
\subsection{Preliminary result about \texorpdfstring{$u$}{u} under the hypothese of proposition \ref{trivialKernel}}
\begin{lemma}
\label{fredholmEstimateEff}
 Let u be defined as in proposition \ref{trivialKernel}. For all $r',l'$ satisfying the same inequalities as $r$ and $l$, we have
\[
\left\| u \right\|_{H^{r',l',\nu}_b}\leq C\left\|u\right\|_{H^{r,l,\nu}_b}
\]
In particular, $u\in H^{\infty, l, \nu}_b$.
\end{lemma}
\begin{proof}
First note that the operator $\hat{P}(\zeta)$ is $b$-elliptic near $x = +\infty$ and outside any neighborhood of the end $x=0$.
Therefore this lemma follows from the application of the estimate presented in the proof of Theorem 1.1 in \cite{vasy2020limiting}. Indeed, the operator 
\begin{align*}P(\zeta) :&= e^{i\frac{\zeta}{x}}\tilde{P}(\zeta)e^{-i\frac{\zeta}{x}}\\
&= (x^2D_x)^2 + \left(\beta+i(n-1)\right)x^3D_x + \left(L+i\beta\frac{n-2}{2} + \beta'\right)x^2 + \gamma x\zeta - \zeta^2
\end{align*} 
satisfy the three conditions presented in the proof of Proposition \ref{midFreqInfEstimate} with $a_0 := \beta$, $A := L+i\beta\frac{n-2}{2}+\beta'$, $\sigma:= \zeta$, $Q = \gamma x$.
Note that the third condition is a consequence of the following property of $L$:
There exists $(L_k)_{k=0}^N$ in $\text{Diff}^1(E)$ and $(L'_k)$ and $L''$ in $C^{\infty}(\mathcal{N})$ such that 
\begin{align*}
L = \sum_{k=0}^{N} L_k^*L_k + L'_kL_k + L''
\end{align*}
By using a finite partition of unity on the compact manifold $\mathcal{N}$, we can assume that $L$ has compact support on an open set of trivialization $U$ of $E$ ($L$ is not selfadjoint anymore after localization but it is still principally self adjoint). We write the operator in a local trivialization (in which $m$ is the canonical hermitian product on $\C$):
\begin{align*}
L = \sum_{1\leq i,j\leq n} a_{i,j}D_{y_i} D_{y_j} + V
\end{align*}
with $V\in \text{Diff}^1(U)$ (compactly supported in $U$). We can assume that $a_{i,j}=a_{j,i}$ in this decomposition. Because $L$ is principally self adjoint, $a_{i,j}$ are real valued. Let $\chi$ smooth non negative compactly supported in $U$ be such that $\chi = 1$ on the support of $L$ and $V$. We have:
\begin{align*}
&L-\sum_{1\leq i,j\leq n}\left(\frac{a_{i,j}}{2}D_{y_i} + \chi^2 D_{y_j}\right)^*\left(\frac{a_{i,j}}{2}D_{y_i} + \chi^2 D_{y_j}\right)\\
& - \sum_{i=1}^n \left(\sqrt{\sum_{j=1}^{n}\frac{a_{i,j}^2}{4}+n\chi^2}\partial_{y_i}\right)^*\left(\sqrt{\sum_{j=1}^{n}\frac{a_{i,j}^2}{4}+n\chi^2}\partial_{y_i}\right) \in \text{Diff}^1(U)
\end{align*}
Note that we use the fact that $\sqrt{\sum_{j=1}^{n}\frac{a_{i,j}^2}{4}+n\chi^2}$ is smooth on $U$. Therefore, we can take $(\frac{a_{i,j}}{2}D_{y_i} + \chi^2 D_{y_j})$ and $\sqrt{\sum_{j=1}^{n}\frac{a_{i,j}^2}{4}+n\chi^2}\partial_{y_i}$ for the $L_k$. The important fact is that they generate the space of smooth vector fields with support included in the support of $L$.
\end{proof}

\begin{lemma}
\label{asymptoticExpansion}
Let $u$ be defined as in proposition \ref{trivialKernel}. There exists a sequence $(a_k)_{k\in \N}$ of smooth sections of $E$ such that for all $N\in \N$ and all $\epsilon>0$:
\begin{align*}
u -\sum_{k=0}^{N} x^{k-\frac{1}{2}+\frac{n}{2}-i\frac{\beta-\gamma}{2}}a_k \in H^{\infty, \frac{1}{2}+N -\epsilon +\frac{\Im(\beta-\gamma)}{2}, -\infty}_b.
\end{align*}
\end{lemma}
\begin{proof}
Let $\chi$ be a smooth non negative function such that $\chi = 1$ for $x\leq 1$ and $\chi = 0$ for $x\geq 2$.
Since $\tilde{P}(\zeta)u = 0$, we have $x^{-\frac{n+2}{2}}\tilde{P}(\zeta)x^{\frac{n}{2}}(x^{-\frac{n}{2}}u) = 0$ with $v := x^{-\frac{n}{2}}u \in H^{\infty, l}_b$ (this space being considered with a usual $b$ volume form).
We denote by $N := -2\zeta\left(xD_x - \frac{i}{2} + \frac{\beta-\gamma}{2}\right)$ the normal operator of $Q:=x^{-\frac{n+2}{2}}\tilde{P}(\zeta)x^{\frac{n}{2}}$ at the boundary $x = 0$.
We have
\begin{align*}
N\chi v &= Q\chi v + (N-Q)\chi v \\
&= [Q(\zeta),\chi]v + (N-Q)\chi v.
\end{align*}
Therefore $f := N\chi v \in H^{\infty, l+1}_b$ and has support in $\left\{x\leq 2\right\}$.
We use the Mellin transform (with the convention $\mathcal{M}u(\sigma) = \int_{-\infty}^{+\infty}x^{-i\sigma-1}u(x)\dd x$) to get the following equality on $\left\{\Im(\sigma)>-l\right\}$:
\begin{align*}
-2\zeta\left(\sigma - \frac{i}{2}+\frac{\beta-\gamma}{2}\right)\mathcal{M}(\chi v)(\sigma) = \mathcal{M}f(\sigma)\\
\mathcal{M}(\chi v)(\sigma) = \frac{\mathcal{M}f(\sigma)}{-2\zeta \left(\sigma - \frac{i}{2}+\frac{\beta-\gamma}{2}\right)}.
\end{align*}
Since $\mathcal{M}f$ is holomorphic on $\left\{\Im(\sigma)>-l-1\right\}$, $\mathcal{M}u(\sigma)$ has a meromorphic extension on $\left\{\Im(\sigma)>-l-1\right\}$ with the only possible pole at $\frac{i}{2}-\frac{\beta-\gamma}{2}$. We can assume that $\frac{i}{2}-\frac{\beta-\gamma}{2} \in \left\{-l>\Im(\sigma)>-l-1\right\}$ (otherwise $\mathcal{M}u$ is holomorphic and by a contour deformation argument (the rhs has enough decay at $\Re(\sigma) = \pm \infty$), we prove that $\chi v \in H^{\infty, l+1}_b$ and $f \in H^{\infty, l+2}$ and we can iterate the process.
Therefore, the same contour deformation argument comes with a residu and we have:
\begin{align*}
\chi v(x) = 2i\pi x^{-\frac{1}{2}-i\frac{\beta-\gamma}{2}} \mathcal{M}f\left(\frac{i}{2}-\frac{\beta-\gamma}{2}\right) + R(x)
\end{align*}
with $R(x) \in H^{\infty, l+1}_b$. We use this to prove that $\mathcal{M}f$ has a meromorphic extension up to $\left\{\Im(\sigma)>-l-2\right\}$. Indeed, $\mathcal{M}([Q,\chi]v)$ is holomorphic on the whole complex plane and $(N-Q)(\chi v) =  N_2 x\chi v + S \chi v$ with $N_2$ a dilation invariant operator and $S \in x^2 \text{Diff}^2_b$ and therefore $S\chi v \in H^{\infty, l+2}_b$ and $\mathcal{M}(S\chi v)$ holomorphic on $\left\{\Im(\sigma)>-l-2\right\}$ and $\mathcal{M
}(N_2 x\chi)$ is a polynomial times a translate of $\mathcal{M}u$ by $-i$.
The iteration of the process gives a development:
\begin{align*}
\chi v = \sum_{k=0}^N h_k x^{-\frac{1}{2}+k-i\frac{\beta-\gamma}{2}} + R_N
\end{align*}
with $R_N \in H^{\infty, \frac{1}{2}+N -\epsilon +\frac{\Im(\beta-\gamma)}{2}}_b$.
\end{proof}

\subsection{Diagonalization of \texorpdfstring{$L$}{L}}
We denote by $(f_i, \lambda_i)_{i\in \N}$ an orthonormal family of eigenfunctions and eigenvalues of $L$.
We define the following operators:
\begin{align*}
\Pi_i:\begin{cases}
H^{\infty, l, \nu}_b(\tilde{E})\rightarrow C^{\infty}((0,+\infty)_x)\\
u \mapsto \left(x\mapsto \int_{\mathcal{N}}m(u(x,y),f_i(y))\dd vol(y)\right)
\end{cases}
\end{align*}
(where $\dd vol(y)$ is the volume form on $\mathcal{N}$).
\begin{lemma}
We have $\lim_{n\to +\infty} \sum_{i=0}^n(\Pi_i u)f_i =  u$ in the topology of $H^{\infty, l, \nu}$.
\end{lemma}
\begin{proof}
Since $L$ is an elliptic self-adjoint operator on sections of a bundle over a compact manifold, for all $x\in [0,+\infty)_x$, we have (in $L^2(E)$):
\begin{align}
\lim\limits_{n\to +\infty} \sum_{i=0}^n(\Pi_iu)(x)f_i = u(x,.) \notag \\
\label{dominationCondition}
\left\| \sum_{i=0}^n (\Pi_iu)(x)f_i \right\|^2_{L^2(\partial {X})}\leq \left\|u(x,.)\right\|^2_{L^2(\partial X)}.
\end{align}
We denote by $w_{l,\nu}$ a positive smooth weight on $[0,+\infty)_x$ such that $w_{l,\nu} = x^{-l}$ near $\left\{x=0\right\}$ and $w_{l,\nu} = x^{\nu}$ near $\left\{x= +\infty\right\}$. In particular for all $N\in \N$, the multiplication by $w_{l,\nu}$ is a homomorphism from $H^{N,l,\nu}_b$ to $H^{N,0,0}_b$. Note also that $w_{l,\nu}$ commutes with $\Pi_i$.
Therefore, we have pointwise convergence of $w_{l,\nu}\sum_{i=0}^n (\Pi_i u)f_i$ towards $w_{l,\nu}u$ in $L^2([0,+\infty)_x, L^2(E))$ and the bound \eqref{dominationCondition} is what we need to apply Lebesgue convergence theorem.
Therefore we have the result of the lemma for the topology of $H^{0,l,\nu}$. To get it for the topology of $H^{2N,l,\nu}_b$, we apply the previous method to $\left((xD_x)^2+\Delta_{\partial X}\right)^N w_{l,\nu} u \in H^{0,0,0}_b$ (the important point being that $\left((x\partial_x)^2+\Delta_{\partial X}\right)$ is elliptic of order $2N$ in the $b$ sense and commutes with $\Pi_i$).
\end{proof}

By orthogonality of the $f_i$, the previous lemma implies that $u = 0$ if and only if $u_i := \Pi_i u = 0$ for all $i\in \N$.
Moreover, using the lemma again and the fact that $\tilde{P}(\zeta)u = 0$, we see that $\tilde{P}_i(\zeta) u_i = 0$ for all $i\in \N$ where \[\tilde{P}_i(\zeta) = \left(x^2D_x\right)^2+i(n-1)x^3D_x + x^2 \lambda_i + \beta\left(x^3D_x + i\frac{n-2}{2}x^2\right)+\beta'x^2 - 2 \zeta x\left(xD_x + i\frac{n-1}{2} + \frac{\beta-\gamma}{2}\right).\]
 
For the rest of this section, we fix $k\in \N$ and we prove that $u_k = 0$. We use standard results in analysis of singular ODE.
First, we can change variable $r = x^{-1}$ in the equation and we get:
\begin{align*} Q_k(\zeta) =& D_r^2 - i(n-1)r^{-1}D_r + r^{-2}\lambda_k + \beta\left(-r^{-1}D_r + i\frac{n-2}{2}r^{-2}\right) + \beta'r^{-2}\\
& - 2\zeta \left(-D_r + i\frac{n-1}{2}r^{-1} + \frac{\beta-\gamma}{2}r^{-1}\right) \\
v_k(r) :=& u_k(x^{-1})\\
\tilde{P}_k u_k =& 0 \Leftrightarrow Q_k(\zeta) v_k &= 0
\end{align*}
Therefore, the equation has meromorphic coefficients (with only poles at 0), has a regular singularity at $r = 0$ and a rank one singularity at $r=\infty$.
The indicial equation at $0$ is 
\begin{align*}
-\alpha^2 + (2-n+i\beta)\alpha + \left(\lambda_k + i\beta\frac{n-2}{2} + \beta'\right)=0
\end{align*}
with two roots:
\begin{align*}
\alpha_{\pm} = 1-\frac{n}{2}-\frac{i\beta}{2}\pm \sqrt{\left(\frac{n-2}{2}\right)^2-\frac{\beta^2}{4}+\lambda_k+\beta'}
\end{align*}
(we choose the branch of the square root which extends the square root on $\R_+$ to $\C\setminus (-\infty,0]$).

We know by Frobenius theory that the space of solution of the equation $Q_k(\zeta) v_k = 0$ is of dimension 2 and consists of analytic functions with a holomorphic extension (at least) to $\C \setminus (-\infty,0]$ where 0 is a branch point. After fixing a continuous determination of $\ln(x)$ on $\C \setminus (-\infty,0]$, one of the solutions can be written as:
\begin{align*}
\sum_{n=0}^{+\infty}a_n r^{n+\alpha_+}
\end{align*}
where the radius of convergence of the serie $\sum_{n=0}^{+\infty}a_n r^n$ is infinite.
The precise expression of an other independent solution as a serie depends on whether or not $\alpha_+-\alpha_- \in \N$ but it is equivalent to $bx^{\alpha_-}$ for some $b \neq 0$ (since $\alpha_+\neq\alpha_-$).
The fact that $u\in H^{l,\nu}_b$ with \[\nu \in \left(1-\frac{\Im(\beta)}{2}-\Re\sqrt{-\frac{\beta^2}{4}+\left(\frac{n-2}{2}\right)^2 + \beta'}, 1-\frac{\Im(\beta)}{2}+\Re\sqrt{-\frac{\beta^2}{4}+\left(\frac{n-2}{2}\right)^2+\beta'}\right)\] implies that $u_k$ is $o(x^{-1+\frac{n}{2}+\frac{\Im(\beta)}{2}-\Re\sqrt{-\frac{\beta^2}{4}+\left(\frac{n-2}{2}\right)^2+\beta'}})$ when $x\rightarrow +\infty$. Therefore, since $\lambda_k \geq 0$, we have $v_k = o(r^{\alpha_-})$. We deduce that $v_k = C\sum_{n=0}^{+\infty}a_n r^{n+\alpha_+}$ for some $C\in \C$. In particular, the function $r^{-\alpha_+}v_k$ has a holomorphic extension to $\C$.

Now we show that the differential equation can be written as a confluent hypergeometric equation after the transformation $\tilde{Q}_k := r^{1-\alpha_+}Q_k r^{\alpha_+}$ followed by the change of variable $z = -2i\zeta r$.
\begin{align*}
\tilde{Q}_k &= -r\partial_r^2+(1-n-2\alpha_++i\beta -2i\zeta r)\partial_r+ \zeta (-2i\alpha_+ -\beta + \gamma - i(n-1))\\
&= z\partial_z^2 + (2\alpha_+ + (n-1)-i\beta-z)\partial_z+ \left(-\alpha_+ + \frac{i\beta}{2}-\frac{i\gamma}{2}-\frac{n-1}{2}\right) \\
&= z\partial_z^2 + (c-z)\partial_z - \tilde{a}
\end{align*}
where
\begin{align*}
c&= 2\alpha_+ + n -1 - i\beta\\
&= 1+2\sqrt{\left(\frac{n-2}{2}\right)^2-\frac{\beta^2}{4}+\lambda_k+\beta'}\\
\tilde{a}&= \alpha_+-\frac{i\beta}{2}+\frac{i\gamma}{2}+\frac{n-1}{2} \\
&= \frac{1}{2}+\frac{i\gamma}{2}+\sqrt{\left(\frac{n-2}{2}\right)^2-\frac{\beta^2}{4}+\lambda_k+\beta'}
\end{align*}
Now we assume that $\tilde{a}$ is not a non positive integer. If it is the case, $v_k$ could be non zero (at least it is true if $c$ is not a non positive integer with $c\leq a$).
Therefore, since $Q_k v_k = 0$, we have $\tilde{Q}_k \tilde{v}_k = 0$ where $\tilde{v}_k(z) = (r^{-\alpha_+}v_k)\left(\frac{iz}{2\zeta}\right)$. Moreover, $\tilde{v}_k$ is holomorphic on $\C$. Therefore, by Frobenius theory, we have $\tilde{v}_k = C\mathbf{M}(\tilde{a},c,z)$ (where $\mathbf{M}$ is the renormalized Kummer's function, defined for example by (9.04) in \cite{olver1997asymptotics}). Note that $z\mapsto\mathbf{M}(\tilde{a},c,z)$ is a non zero entire function except when $a,c\in -\N$ and $c\leq a$ where it is $0$. Using (10.09) and (10.10) in \cite{olver1997asymptotics}, we see that we have the following asymptotic equivalent for $\mathbf{M}(\tilde{a},c,z)$ when $\lvert z\rvert \rightarrow \infty$:
\begin{align*}
\mathbf{M}(\tilde{a},c,z) &\sim \frac{e^{\tilde{a}i\pi}z^{-\tilde{a}}}{\Gamma(c-\tilde{a})}+\frac{e^zz^{c-\tilde{a}}}{\Gamma(\tilde{a})} &\text{when } -\frac{\pi}{2}-\delta<\text{arg}(z)<-\frac{\pi}{2}+\delta \\
\mathbf{M}(\tilde{a},c,z) &\sim \frac{e^zz^{c-\tilde{a}}}{\Gamma(\tilde{a})} &\text{  when } -\frac{\pi}{2}+\delta<\text{arg}(z)<\frac{\pi}{2}-\delta \\
\mathbf{M}(\tilde{a},c,z) &\sim \frac{e^{-\tilde{a}i\pi}z^{-\tilde{a}}}{\Gamma(c-\tilde{a})}+\frac{e^zz^{c-\tilde{a}}}{\Gamma(\tilde{a})} &\text{when } \frac{\pi}{2}-\delta<\text{arg}(z)<\frac{\pi}{2}+\delta \\
\mathbf{M}(\tilde{a},c,z) &\sim \frac{e^{-\tilde{a}i\pi}z^{-\tilde{a}}}{\Gamma(c-\tilde{a})} &\text{  when } \frac{\pi}{2}+\delta < \text{arg}(z) <\frac{3}{2}\pi -\delta
\end{align*}
where $z^{-a}$ and $z^{c-a}$ are fixed using $\ln(z) = \ln\lvert z\rvert + i \text{arg(z)}$ (pay attention to the fact that in the last asymptotics, arg is not the principal argument).

We now end the proof of proposition \ref{trivialKernel}. When $r\in (0,+\infty)$ and $\zeta$ in the upper half complex plane, we have $z = -2i\zeta r$ has principal argument $\text{arg}(z) \in [-\frac{\pi}{2}, \frac{\pi}{2}]$. Therefore, if $\tilde{a}$ is not a non positive integer $\mathbf{M}(\tilde{a},c,-2i\zeta r)$ has no asymptotic expansion when $r\rightarrow +\infty$ (because of the exponential term in the previous asymptotic expansion). Moreover under the hypotheses of proposition \ref{trivialKernel}, we know that $v_k$ has an asymptotic expansion (comes from the fact that $u$ has such an expansion by lemma \ref{asymptoticExpansion}). Therefore, $v_k= 0$ for all $k$ and $u = 0$.

We now end the proof of proposition \ref{trivialKernel2}. If $r\in (0,+\infty)$ and $\zeta$ in in the lower half complex plane, we have $z = -2i\zeta r$ has argument in $[\frac{\pi}{2}, \frac{3\pi}{2}]$ and if $c-\tilde{a}$ is not a non positive integer, $\mathbf{M}(\tilde{a},c,z)$ has a non zero $z^{-\tilde{a}}$ term in its expansion and since $\lvert r^{\alpha_+}z^{\tilde{a}}r^{-\frac{n-2}{2}+\Im\left(\frac{\beta-\gamma}{2}\right)}\rvert^2r^{n-1}$ is not integrable at infinity, and since $v_k \in H^{r, l, \nu}_b(\tilde{E})$ with $l>-\frac{n-2}{2}+\Im\left(\frac{\beta-\gamma}{2}\right)$, we must have $v_k = 0$.

\begin{coro}
\label{invertibilityNormalEff}
Let $\zeta \in \C\setminus\left\{0\right\}$ and $\Im(\zeta)\geq 0$. We assume that for all $k\in \N$,
\begin{align*}
\frac{1}{2}+i\frac{\gamma}{2}+\sqrt{\left(\frac{n-2}{2}\right)^2 - \frac{\beta^2}{4} + \lambda_k + \beta'}\notin \N \\
\frac{1}{2}-i\frac{\overline{\gamma}}{2}+\sqrt{\left(\frac{n-2}{2}\right)^2 - \frac{\beta^2}{4} + \lambda_k + \overline{\beta'}}\notin \N
\end{align*}
Then $\tilde{P}(\zeta)$ is invertible from $\left\{u\in H^{r,l,\nu}_b: \tilde{P}(\zeta)u\in H^{r, l+1, \nu-2}\right\}$ to $H^{r,l+1, \nu-2}_b$ where $l<-\frac{1}{2}+\frac{\Im(\beta-\gamma)}{2}$ and  $r+l>-\frac{1}{2}+\frac{\Im(\beta+\gamma)}{2}$ and $\nu \in \left(1-\frac{\Im(\beta)}{2}-\Re\sqrt{-\frac{\beta^2}{4}+\left(\frac{n-2}{2}\right)^2 + \beta'}, 1-\frac{\Im(\beta)}{2}+\Re\sqrt{-\frac{\beta^2}{4}+\left(\frac{n-2}{2}\right)^2+\beta'}\right)$.
\end{coro}
\begin{proof}
We already know that we have Fredholm estimates for $\tilde{P}(\zeta)$ (see the proof of lemma \ref{fredholmEstimateEff}). The triviality of the kernel follows from proposition \ref{trivialKernel}. The adjoint operator $\tilde{P}(\zeta)^*$ has the same form than $\tilde{P}$ but with $\zeta$ replaced by $\overline{\zeta}$, $\beta$ replaced by $\overline{\beta} = -\beta$, $\gamma$ replaced by $\overline{\gamma}$, $\beta'$ replaced by $\overline{\beta'}$. The second non integer coincidence condition therefore correspond to the condition in proposition \ref{trivialKernel2}. Moreover, $\tilde{P}(\zeta)^*$ goes from $H^{-r, -l-1, 2-\nu}_b$ to $H^{r,l,\nu}_b$. Therefore the hypotheses of proposition \ref{trivialKernel2} are satisfied since:
\begin{align*}
-l-1 >& -\frac{1}{2} - \frac{\Im(\beta-\gamma)}{2} = -\frac{1}{2}+ \frac{\Im(\overline{\beta}-\overline{\gamma})}{2}\\
-r-l-1 <&-\frac{1}{2} + \frac{\Im(\overline{\beta}+\overline{\gamma})}{2}\\
2-\nu \in& \left(1-\frac{\Im(\overline{\beta})}{2}-\Re\sqrt{-\frac{\overline{\beta}^2}{4}+\left(\frac{n-2}{2}\right)^2 + \overline{\beta'}}, 1-\frac{\Im(\overline{\beta})}{2}+\Re\sqrt{-\frac{\overline{\beta}^2}{4}+\left(\frac{n-2}{2}\right)^2+\overline{\beta'}}\right)
\end{align*}

\end{proof}

\section{Proofs of Propositions \ref{basicConormalEstimate}, \ref{refinementAtScri}}\label{ProofCauchyConormal}

The goal of this section is to provide detailed proofs of Propositions \ref{basicConormalEstimate}, \ref{refinementAtScri} and Corollary \ref{FourierTrans2}. The proof of proposition \ref{basicConormalEstimate} is based on energy estimates on the model of what is done in \cite[Section~4.1]{hintz2020stability}. The estimate on $U_3$ (defined below) follows less closely \cite{hintz2020stability} but remains in the same spirit. The proof of Proposition \ref{refinementAtScri} is then based on an idea presented in \cite[Section~5.1]{hintz2020stability}.

For the proof it will be convenient to have a concrete defintion for $t_0$. Moreover, we also need a time coordinate $\tt$ whose level sets are transverse to $\mathscr{I}^+$ and $\mathfrak{H}$ and such that $\dd \tt$ remains uniformly timelike with respect to $\tilde{G}$ up to $\mathscr{I}^+$ and up to $r_+-2\epsilon$. Therefore we introduce:
\begin{align}
t_0 := t_*+h(r) \label{deft_0}\\
\tt := t_*+F(r) \notag\\
\end{align}
where $h(r)$ and $F(r)$ are defined as follows. Let $\chi_0\in C^{\infty}(\R, [0,1])$ be a smooth cutoff such that $\chi_0 = 1$ on $(-\infty, 3M]$ and $\chi_0 = 0$ near $+\infty$. Let $\chi_1 \in C^{\infty}(\R, [0,1])$ such that $\chi_1 = 0$ on $(-\infty, 5M)$ and $\chi_1 = 1$ near $+\infty$. Moreover, we can assume\footnote{This assumption is used to ensure $\mathfrak{t}+1\leq t_0$ and be able to figure out the relative position of some particular level sets of $\mathfrak{t}$ and $\Sigma_0$ later in the proof.} that $\sup \supp(\chi_0)<\inf \supp(\chi_1)$ and $\chi_0\leq \psi_0$, $\chi_1\leq \psi_1$ (where $\psi_0,\psi_1$ were used to define $L(r)$ in \eqref{defLr}). We define:
\begin{align*}
h(r) :=& \begin{cases}-T(r)-\int_r^{+\infty} \chi_0(r)\left(\frac{a^2+r^2}{\Delta_r}-\frac{M^2}{a^2+r^2}\right) \dd r \text{ if $r>r_+$}\\
-\int_{3M}^{+\infty}\chi_0(r)\frac{a^2+r^2}{\Delta_r}\dd r+\int_{r}^{+\infty}\chi_0(r)\frac{M^2}{a^2+r^2}\dd r
\end{cases}\\
F(r) =& h(r) + \int_{-\infty}^r \chi_1(r) \left(\frac{M^2}{\Delta_r}-\frac{(a^2+r^2)}{\Delta_r}\right) \dd r
\end{align*}
Note that $h(r)$ and $F(r)$ are smooth across the event horizon. Moreover, $\tt$ is smooth up to $\mathscr{I}^+$. 

By construction, we have:
\begin{align*}
 \mathfrak{t}+1\leq \tt \leq t_0\leq t 
\end{align*}

An explicit computation provides:
\begin{align*}
\tilde{G}(\dd t_0, \dd t_0) =& -a^2\sin^2\theta + \frac{(a^2+r^2)^2}{\Delta_r}-\chi_0(r)^2 \frac{(a^2+r^2)^2}{\Delta_r}\left(1-\frac{M^2\Delta_r}{(a^2+r^2)^2}\right)^2 \\
\geq& -a^2\sin^2\theta + M^2\left(2-\frac{M^2\Delta_r}{(a^2+r^2)^2}\right) 
\end{align*}
Since on $r>r_+-2\epsilon$ (for $\epsilon$ sufficiently small), we have:
\begin{align*}
M^2 \leq r^2\leq a^2+r^2 \\
\Delta_r\leq a^2+r^2
\end{align*}
we deduce $\tilde{G}(\dd t_0, \dd t_0)\geq -a^2\sin^2\theta + M^2$ and $\dd t_0$ is timelike on $\mathcal{M}_{2\epsilon}$.
Similarly, since $\dd \tt = \dd t_0$ for $r< \inf(\supp \chi_1)$ we have that $\dd \tt$ is timelike when $r<\inf(\supp \chi_1)$. For $r\geq \inf(\supp \chi_1)$, we have:
\begin{align*}
\tilde{G}(\dd \tt, \dd \tt) =& -a^2\sin^2\theta + \frac{(a^2+r^2)^2}{\Delta_r} - \chi_1(r)^2\frac{(a^2+r^2-M^2)^2}{\Delta_r}\\
\geq& -a^2\sin^2\theta + \frac{(a^2+r^2)M^2}{\Delta_r}\left(2-\frac{M^2}{a^2+r^2}\right)\\
\geq& M^2-a^2
\end{align*}
This proves that the $\dd \tt$ is uniformly timelike with respect to $\tilde{G}$.

Let $\alpha>0$, $\beta>0$ and $T>\frac{1}{2}$ be constants to be chosen later. Note that since $\mathfrak{t}+1\leq \tt$, we have $\left\{\tt \leq \frac{1}{2}\right\}\subset \mathcal{U}$.
Our goal to prove Proposition \ref{basicConormalEstimate} is to get estimates on $\left\{ r\geq r_+-\frac{\epsilon}{2}, t_0\geq 0, \tt\leq T\right\}$ for all $T>\frac{1}{2}$. To achieve that, we glue four estimates together:
\begin{itemize}
\item An estimate on $U_1 := \left\{ t_0\geq 0, \rho_I \geq \alpha, \rho_0\leq \beta\rho_I \right\}$
\item An estimate on $U_2 := \left\{ \tt \leq \frac{1}{2}, \rho_I\leq 2\alpha\right\}$
\item An estimate on $U_3 := \left\{ r_+-\frac{\epsilon}{2}\leq r,\frac{1}{4} \leq \tt \leq T\right\}$
\item An estimate on $K := \left\{ r_+-\frac{\epsilon}{2}\leq r, \tt \leq \frac{1}{2}, t_0\geq 0, \rho_I\geq \alpha, \rho_0\geq \frac{\beta\alpha}{2} \right\}$
\end{itemize}
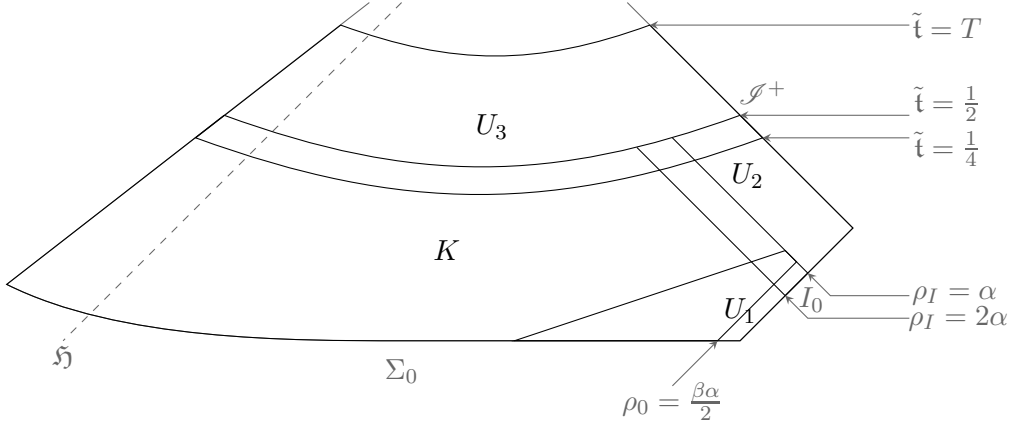
\begin{figure}
\begin{center}
\begin{tikzpicture}[scale = 1.5]
\draw [black!60] (-2.286, 3)--(-5.5, 0.5).. controls +(1, -0.5) and +(-1,0) .. (-2,0) node[yshift = -0.4cm] {$\Sigma_0$}--(1,0)--(2,1) node [midway, xshift = 0.2cm, yshift = -0.2cm] {$I_0$}--(0,3) node [midway, xshift = 0.3cm, yshift = 0.3cm] {$\mathscr{I}^+$};
\draw [dashed, black!60](-2,3)--(-5,0) node[yshift = -0.2cm] {$\mathfrak{H}$};
\draw (1.6,0.6)--(1.4, 0.8)--(-1,0)--(1,0) node [yshift = 0.4cm] {$U_1$}--cycle;
\draw (-3.573, 2)to [in= -160, out =-20](1, 2);
\draw (1, 2) node[yshift = -0.8cm, xshift = 0.1cm] {$U_2$}--(2,1)--(1.4, 0.4)--(0.085, 1.715);
\draw (0.2, 2.8) to [out = -160, in = -20](-2.54,2.8)--(-3.825, 1.8)to [in= -160, out =-20](1.2, 1.8) --cycle;
\draw (-1.2, 1.9)node {$U_3$};
\draw (-3.573, 2)--(-5.5, 0.5).. controls +(1, -0.5) and +(-1,0) .. (-2,0)--(0.8, 0)--(1.5, 0.7)--(0.4, 1.8);
\draw (-1.6, 0.8) node {$K$};
\draw [stealth-, black!60] (1.2, 1.8)--(2.5, 1.8) node[xshift = 0.5cm, yshift = -0.15cm] {$\tt = \frac{1}{4}$};
\draw [stealth-, black!60] (1,2)--(2.5,2) node[xshift = 0.5cm, yshift = 0.15cm] {$\tt = \frac{1}{2}$};
\draw [stealth-, black!60] (1.6, 0.6)--(1.8, 0.4)--(2.5, 0.4) node[xshift = 0.6cm, yshift = 0.05] {$\rho_I = \alpha$};
\draw [stealth-,black!60] (1.4, 0.4)--(1.6, 0.2)--(2.5, 0.2) node[xshift = 0.65cm, yshift = -0.05] {$\rho_I = 2\alpha$};
\draw [stealth-,black!60] (0.8, 0)--(0.4, -0.4) node[yshift = -0.2cm] {$\rho_0 = \frac{\beta \alpha}{2}$};
\draw [stealth-,black!60] (0.2, 2.8)--(2.5, 2.8) node[xshift = 0.5cm] {$\tt = T$};
\end{tikzpicture}
\end{center}
\caption{Representations of the sets $U_1, U_2, U_3$ and $K$}
\end{figure}
Note that $\rho_I$ is bounded\footnote{This follows from, $-x\mathfrak{t} = \frac{h(r)+L(r)-t_0}{r}$} on $t_0\geq 0$, therefore, for every $\eta>0$, by taking $\beta$ small enough, we have $\rho_0\leq \eta$ on $U_1$. We denote by $\beta\mapsto\eta(\beta)$ a positive function with $\lim\limits_{\beta \to 0}\eta(\beta) = 0$ and such that $\rho_0\leq \eta(\beta)$ on $U_1$. Moreover since $\rho_0 = \frac{1}{L(r)-t_0+h(r)}$ on $U$, if $t_0\geq 0$, $\frac{1}{L(r)+h(r)}\leq \rho_0$ and therefore for all $C>0$, there exists $\eta>0$ such that $r>C$ if $\rho_0\leq \eta$. Injecting this in the definition of $\rho_I$, we deduce that for $\beta$ small enough, $\rho_I\leq \frac{3}{2}$ on $U_1$ (in the sequel we assume that $\beta$ is small enough to ensure this property).
Note that $K$ is compact included in $\mathcal{M}_{2\epsilon}$ (and can be included into a relatively compact hyperbolic region $\mathcal{R}$).
Since $t_0$ is a global time function on $\mathcal{M}_{2\epsilon}\cap \left\{t_0>-1\right\}$, by classical hyperbolic theory, we get:
\begin{lemma}
For all $s \in \R$, there exists $C>0$ such that:
\begin{align*}
\left\| u\right\|_{\overline{H}^{s+1}(K)}\leq C\left\|u_0\right\|_{\overline{H}^{s+1}(\mathcal{R}\cap \Sigma_0)} + C\left\|u_1\right\|_{\overline{H}^{s}(\mathcal{R}\cap \Sigma_0)}
\end{align*}
\end{lemma}

\subsection{Estimate on \texorpdfstring{$U_1$}{U\_ 1}}
On $U_1$, we have the following expression for the operator $T_s$:
\begin{align*}
T_s = \rho_I(2-\rho_I)\partial_{\rho_I}^2 -2\rho_0\partial_{\rho_0}\partial_{\rho_I}- \slashed{G} + \rho_0\text{Diff}_b^2 + \text{Diff}^1_b
\end{align*}
Therefore, the inverse metric writes:
\begin{align}
\label{metricU1}
\tilde{G} =\rho_I(2-\rho_I)\partial_{\rho_I}^2 -2\rho_0\partial_{\rho_0}\partial_{\rho_I}- \slashed{G} + \rho_0 \text{Diff}_b^2
\end{align}
and is therefore a $b$-metric on $U_1$ (note that the only part of the boundary intersecting $U_1$ is $I_0$).

\begin{lemma}
\label{derivativeBLC}
If $g$ is a b-metric (Lorentzian or Riemannian) on some manifold with boundary $\mathcal{M}$ (and we call $\rho$ the defining function of the boundary), and if $X \in \Gamma\left({}^bT^*\mathcal{M}^{\otimes k}\otimes {}^b T\mathcal{M}^{\otimes k'}\right)$, then $\nabla X \in {}^bT^*\mathcal{M}^{\otimes (k+1)}\otimes {}^b T\mathcal{M}^{\otimes k'}$ where $\nabla$ is the Levi-Civita connection associated to $g$.
\end{lemma}
\begin{proof}
By Leibniz rule and since for $Y \in \Gamma({}^bT\mathcal{M})$ and $f\in C^{\infty}(\mathcal{M})$ we have $Yf \in C^{\infty}(\mathcal{M})$, it is enough to prove that if $\rho$ is a local defining function of the boundary and $(y_i)_{i=1}^n$ are local coordinates on the boundary:
\begin{align}
\nabla_{\rho\partial_\rho} \partial_{y_i} =& a_0(\rho,y)\rho \partial_\rho + \sum_{i=1}^n a_i(\rho,y) \partial_{y_i} \label{firstbDerivative}\\
\nabla_{\rho\partial_\rho} (\rho \partial_\rho) =& a_0(\rho,y)\rho \partial_\rho + \sum_{i=1}^n a_i(\rho,y) \partial_{y_i}\\
\nabla_{\partial_{y_i}} \partial_{y_i} =& a_0(\rho,y)\rho \partial_\rho + \sum_{i=1}^n a_i(\rho,y) \partial_{y_i}\\
\nabla_{\partial_{y_i}} (\rho\partial_{\rho}) =& a_0(\rho,y)\rho \partial_\rho + \sum_{i=1}^n a_i(\rho,y) \partial_{y_i} \label{lastbDerivative}
\end{align}
where coefficients $(a_\mu)_{\mu=0}^{n}$ are smooth up to the boundary (and can be different in each equality). From these equalities, we can deduce the corresponding ones on  $\frac{\dd \rho}{\rho}, \dd y_i$ since for $Z \in \Gamma {}^bT\mathcal{M}$, $\omega \in \Gamma({}^bT^*\mathcal{M})$ and $X\in \Gamma({}^bT\mathcal{M})$, $\nabla_Z \omega(X) = Z(\omega(X)) - \omega(\nabla_Z X)$ (which is smooth up to the boundary since by the previous computation $\nabla_Z X \in \Gamma({}^bT\mathcal{M})$). Then using the definition of the tensor product connection, we get the full result.
To prove \eqref{firstbDerivative}-\eqref{lastbDerivative}, we only need to prove that $\Gamma^\beta_{\alpha,\gamma} = \rho^{\delta^\beta_0 - \delta^\alpha_0 - \delta^\gamma_0}a(\rho,y)$ with $a$ smooth up to the boundary. Since we have
\begin{align*}
\Gamma^\beta_{\alpha,\gamma} = \frac{1}{2}G^{\beta,\theta}\left(\partial_{\alpha}g_{\theta,\gamma} + \partial_{\gamma}g_{\theta,\alpha}-\partial_{\theta}g_{\alpha,\gamma}\right)
\end{align*}
we conclude the argument using the fact that $g\in {}^bT^*\mathcal{M}^{\otimes 2}$ and $G\in {}^bT\mathcal{M}^{\otimes 2}$.
\end{proof}

Since we do not know the behavior of $u$ near the boundary, we first prove an estimate on $U_1^c := U_1\cap \left\{t+r\leq c\right\}$ for $c>c_0$ with constants uniform with respect to $c$. We then take the limit $c\to +\infty$ and Fatou's lemma to get the estimate on $U_1$.
Note that since $U_1^c$ is relatively compact in $\mathcal{M}_{\epsilon}$, by classical hyperbolic theory we can approximate $u$ by smooth functions $u_n$ such that on $U_1^c$, $\lim\limits_{n\to +\infty}u_n = u$ in $\overline{H}^{s+1}(U_1^c)$ and $\lim\limits_{n\to +\infty}(u_n)_{|_{\Sigma_0\cap U_1^c}} = u_0$ (in $H^{s+1}(\Sigma_0\cap U_1^c)$) and $\lim\limits_{n\to +\infty}(\partial_t u_n)_{|_{\Sigma_0\cap U_1^c}} = u_1$ (in $H^{s+1}(\Sigma_0\cap U_1^c)$). Therefore, we can assume for the estimates on $U_1^c$ that $u$ is smooth.
We define the energy-momentum tensor:
\begin{align*}
T^{\delta,\gamma}(u) = \Re (\mathbb{m}(\Theta_\mu u, \Theta_\nu u))\tilde{G}^{\mu,\delta}\tilde{G}^{\nu, \gamma} - \frac{1}{2}\tilde{G}^{\delta,\gamma}\tilde{G}^{\mu, \nu}\mathbb{m}(\Theta_\mu u, \Theta_\nu u) + \frac{1}{2}\tilde{G}^{\delta, \gamma}\mathbb{m}(u,u).
\end{align*}
This expression is an adaptation of the classical energy-momentum tensor for the scalar wave equation with the introduction of the metric $\mathbb{m}$ (due to the fact that $u$ is valued in $\mathcal{B}_s$) and of the last term which is necessary to control the $L^2$ norm of $u$.
To limit the amount of notations, we sometimes use the name of coordinates as index on tensors or vectors. For example $T^{\rho_0,\rho_0}$ is used to denote the component $T(\dd \rho_0, \dd\rho_0)$. Moreover, we use the same notation for a bilinear form and for the associated quadratic form. For example $\tilde{G}(\dd\rho_I)$ stands for $\tilde{G}(\dd \rho_I, \dd \rho_I)$.
We define the vector fields:
\begin{align*}
V= \nabla^\mu \rho_I \\
W = \rho_0^{-2a_0}e^{A\rho_I}V
\end{align*}
Note that $\tilde{G}(\dd\rho_I) \geq \frac{\alpha}{4}$ on $U_1$ if $\beta$ is small enough since $\alpha\leq \rho_I\leq \frac{3}{2}$ (and $\rho_0\leq \eta(\beta)$)  on $U_1$.
We compute:
\begin{align*}
\div(T^{\mu,\nu}W_\mu) = e^{A\rho_I}\rho^{-2a_0}\left(\div(T)^{\rho_I} - T^{\mu,\nu}\nabla_{\mu}V_{\nu} - \frac{2a_0}{\rho_0} T^{\rho_0,\rho_I} + AT^{\rho_I,\rho_I}\right)
\end{align*}

By lemma \ref{derivativeBLC}, $\nabla \dd \rho_I \in \overline{\Gamma}({}^bT^*U_1^{\otimes 2})$ and therefore we have:
\begin{align*}
\tilde{G}^{\mu,\alpha}\tilde{G}^{\nu,\beta}\nabla_\alpha V_{\beta} \in \overline{\Gamma}({}^bTU_1^{\otimes 2})\\
\tilde{G}^{\mu,\nu} \nabla_{\mu}V_{\nu} \in C^{\infty}(\overline{U}_1)
\end{align*}
Since $U_1$ is relatively compact, we deduce the bounds:
\begin{align*}
\left|T^{\mu,\nu}\nabla_{\mu}V_{\nu}\right|\leq C\left( \mathbb{m}(\partial_{\rho_I}u) + \mathbb{m}(\partial_{\rho_0}u) + \slashed{G}^{\omega_i,\omega_j}\mathbb{m}(\Theta_{\omega_i}u, \Theta_{\omega_j}u)+\mathbb{m}(u)\right)
\end{align*}
where $(\omega_0, \omega_1)$ are unspecified local coordinates on the sphere (note that the expression $\slashed{G}^{\omega_i,\omega_j}\mathbb{m}(\Theta_{\omega_i}u, \Theta_{\omega_j}u)$ does not depend of the choice of these local coordinates).
Similarly since $\frac{\dd \rho_0}{\rho_0}\otimes \dd \rho_I \in \overline{\Gamma}({}^bT^*U_1^{\otimes 2})$, we get:
\begin{align*}
\left|\frac{2a_0}{\rho_0}T^{\rho_0,\rho_I}\right|\leq C\left( \mathbb{m}(\partial_{\rho_I}u) + \mathbb{m}(\partial_{\rho_0}u) + \slashed{G}^{\omega_i,\omega_j}\mathbb{m}(\Theta_{\omega_i}u, \Theta_{\omega_j}u)+\mathbb{m}(u)\right)
\end{align*}

Finally, using the explicit expression of $\tilde{G}$ \eqref{metricU1}, we get:
\begin{align*}
T^{\rho_I,\rho_I} =& \rho_I^2(2-\rho_I)^2\mathbb{m}(\partial_{\rho_I}u)-2\rho_I(2-\rho_I)\Re \mathbb{m}(\partial_{\rho_I}u,\rho_0\partial_{\rho_0}u)+\mathbb{m}(\rho_0\partial_{\rho_0}u)\\
&-\frac{1}{2}\rho_I(2-\rho_I)\left(\rho_I(2-\rho_I)\mathbb{m}(\partial_{\rho_I}u) -2\Re \mathbb{m}(\rho_0\partial_{\rho_0}u,\partial_{\rho_I}u)-\slashed{G}^{\omega_i,\omega_j}\mathbb{m}(\Theta_{\omega_i}u,\Theta_{\omega_j}u)-\mathbb{m}(u)\right)+R \\
=& \frac{1}{2}\rho_I^2(2-\rho_I)^2\mathbb{m}(\partial_{\rho_I}u) - \rho_I(2-\rho_I)\Re \mathbb{m}(\partial_{\rho_I}u,\rho_0\partial_{\rho_0}u)+\mathbb{m}(\rho_0\partial_{\rho_0}u)\\
& + \frac{1}{2}\rho_I(2-\rho_I)\slashed{G}^{\omega_i,\omega_j}\mathbb{m}(\Theta_{\omega_i}u,\Theta_{\omega_j}u)+\frac{1}{2}\rho_I(2-\rho_I)\mathbb{m}(u)\\
\left| R\right|\leq& C\rho_0\left( \mathbb{m}(\partial_{\rho_I}u) + \mathbb{m}(\partial_{\rho_0}u) + \slashed{G}^{\omega_i,\omega_j}\mathbb{m}(\Theta_{\omega_i}u, \Theta_{\omega_j}u)+\mathbb{m}(u)\right)
\end{align*}
Using that:
\begin{align*}
\left|\rho_I(2-\rho_I)\Re \mathbb{m}(\partial_{\rho_I}u,\rho_0\partial_{\rho_0}u)\right|\leq \frac{3}{4}\mathbb{m}(\rho_0\partial_{\rho_0}u) + \frac{\rho_I^2(2-\rho_I)^2}{3}\mathbb{m}(\partial_{\rho_I}u)
\end{align*}
We get that for $\beta$ small enough, 
\begin{align*}
T^{\rho_I,\rho_I} \geq C\left( \mathbb{m}(\partial_{\rho_I}u) + \mathbb{m}(\partial_{\rho_0}u) + \slashed{G}^{\omega_i,\omega_j}\mathbb{m}(\Theta_{\omega_i}u, \Theta_{\omega_j}u)+\mathbb{m}(u)\right)
\end{align*}
on $U_1$ for some constant $C>0$ independent of $u$.
Using the fact that the curvature $R^{\Theta}_{\mu,\nu}$ and $\Theta \mathbb{m}$ are bounded on $U_1$, we get:
\begin{align*}
\left|\div(T)^{\rho_I} \right|\leq C\left(\mathbb{m}(\partial_{\rho_I}u)+\mathbb{m}(\rho_0\partial_{\rho_0}u)+ \slashed{G}^{\omega_i,\omega_j}\mathbb{m}(\Theta_{\omega_i}u,\Theta_{\omega_j}u)+\mathbb{m}(\square_{g,\Theta} u)+\mathbb{m}(u)\right)
\end{align*}
We check that $T_s - \square_{g,\Theta} \in \text{Diff}^1_b$ on $U_1$ and therefore:
\begin{align*}
\left|\div(T)^{\rho_I} \right|\leq C\left(\mathbb{m}(\partial_{\rho_I}u)+\mathbb{m}(\rho_0\partial_{\rho_0}u)+ \slashed{G}^{\omega_i,\omega_j}\mathbb{m}(\Theta_{\omega_i}u,\Theta_{\omega_j}u)+\mathbb{m}(u)+\mathbb{m}(T_s u)\right)
\end{align*}
Finally, we see that if we choose $A$ large enough, we get on $U_1$:
\begin{align*}
\div(T^{\mu,\nu}W_\mu) \geq e^{A\rho_I}\rho^{-2a_0}\left(AC\left( \mathbb{m}(\partial_{\rho_I}u) + \mathbb{m}(\partial_{\rho_0}u) + \slashed{G}^{\omega_i,\omega_j}\mathbb{m}(\Theta_{\omega_i}u, \Theta_{\omega_j}u)+\mathbb{m}(u)\right)-C\mathbb{m}(T_s u)\right)
\end{align*}
where $C>0$ is a constant independent of $u$.
Finally, we can apply the Stockes theorem on $U_1^c$ which has boundaries included in:
\begin{itemize}
\item $\left\{t_0=0\right\}$, an outward normal is $-\dd t_0$ which is timelike past oriented.
\item $\left\{\rho_I=\alpha\right\}$, an outward normal is $-\dd \rho_I$ which is timelike future oriented.
\item $\left\{\rho_0 =\beta\rho_I \right\}$, an outward normal is $\dd \rho_0 - \beta\dd\rho_I$ which is timelike future oriented if $\beta$ is small enough.
\item $\left\{t+r= c\right\}$, an outward normal is $\dd t+\dd r$. We have $\tilde{G}(\dd t+\dd r) = 4Mr + O(1)$ when $r\to +\infty$ and therefore the normal is timelike future oriented on the boundary of $U_1^c$ (provided $r$ is large enough which is the case if $\beta$ is small enough).
\end{itemize}
Since $\dd \rho_I$ is timelike past oriented, we deduce that (uniformly with respect to $c$ in a neighborhood of $+\infty$):
\begin{align*}
\int_{U_1^c}e^{A\rho_I}\rho^{-2a_0}AC\left( \mathbb{m}(\partial_{\rho_I}u) + \mathbb{m}(\partial_{\rho_0}u) + \slashed{G}^{\omega_i,\omega_j}\mathbb{m}(\Theta_{\omega_i}u, \Theta_{\omega_j}u)+\mathbb{m}(u)\right)\dd vol_g\\
 - C\int_{U_1^c}\rho_0^{-2a_0}e^{A\rho_I}\mathbb{m}(T_s u)\dd vol_g\leq \int_{U_1^c}\div(T^{\mu,\nu}W_{\mu})\dd vol_g \\
\leq \int_{\partial U_1^c \cap \left\{t_0=0\right\}} -\rho_0^{-2a_0}e^{A\rho_I}T(\dd \rho_I, \frac{\dd t_0}{\sqrt{\tilde{G}(\dd t_0)}}) \dd vol_{\Sigma_0}
\end{align*}
Since $\dd \rho_I$ and $\frac{\dd t_0}{\sqrt{\tilde{G}(\dd t_0)}}$ are uniformly timelike $b$-one forms when restricted to $U_1 \cap \left\{t_0=0\right\}$ (we have $\dd t = -\left(1+\frac{a^2+r^2}{\Delta_r \rho_I}\right)\frac{\dd \rho_0}{\rho_0^2}-\frac{a^2+r^2}{\Delta_r}\frac{\dd \rho_I}{\rho_0\rho_I^2}$  and since $g$ is a $b$-metric, we have (letting $c\to +\infty$):
\begin{align}
A\int_{U_1^c}e^{A\rho_I}\rho^{-2a_0}\left( \mathbb{m}(\partial_{\rho_I}u) + \mathbb{m}(\partial_{\rho_0}u) + \slashed{G}^{\omega_i,\omega_j}\mathbb{m}(\Theta_{\omega_i}u, \Theta_{\omega_j}u)+\mathbb{m}(u)\right)\dd vol_g \notag \\
 \leq C\int_{U_1^c}\rho_0^{-2a_0}e^{A\rho_I}\mathbb{m}(T_s u)\dd vol_g + Ce^{\frac{3A}{2}}\left(\left\|u_1\right\|_{\rho_0^{a_0}H^1_b(U_1\cap \Sigma_0)}+\left\|u_1\right\|_{\rho_0^{a_0}L^2_b(U_1\cap \Sigma_0)}\right) \label{estimateI0} 
\end{align}
with constants independent of $A$ and $u$.
Now we assume that $T_s u = 0$, in particular, in this case for all $k\geq 1$ $(\rho_0\nabla t_0)^k u$ has a trace on $\Sigma_0\cap U_1$ and we can express 
\begin{align*}((\rho_0 \nabla t_0)^k u)_{|_{\Sigma_0\cap U_1}}= \text{Diff}^{k-1}_b(\Sigma_0\cap U_1)((\rho_0\nabla t_0) u)_{|_{\Sigma_0\cap U_1}} + \text{Diff}^{k}_b(\Sigma_0\cap U_1) u_{|_{\Sigma_0\cap U_1}}.
\end{align*}

We want to get the estimate in $H^k_b(U_1)$ with $k\geq 1$, we use the fact that for $Z\in \text{Diff}^k_b$ defined on $U_1$, $[T_s, Z]\in \text{Diff}^{k+1}_b(U_1)$ as follows:
Let $(Z_i)_{i=1}^N$ be operators in $\text{Diff}^k_b$ such that $\sum_{i=1}^N \left\|Z_i u\right\|_{L^2_b(U_1)}$ is equivalent to $\left\|u\right\|_{H^k_b(U_1)}$.
The estimate \eqref{estimateI0} applied to $Z_l u$ gives:
\begin{align*}
A\int_{U_1^c}e^{A\rho_I}\rho^{-2a_0}\mathbb{m}(\partial_{\rho_I}Z_l u)+\mathbb{m}(\rho_0\partial_{\rho_0}Z_l u)+\slashed{G}^{\omega_i,\omega_j}\mathbb{m}(\Theta_{\omega_i}Z_l u, \Theta_{\omega_j}Z_l u)+\mathbb{m}(Z_l u)\dd vol_g\\
\leq  C\int_{U_1^c}\rho_0^{-2a_0}e^{A\rho_I}[T_s,Z_l] u\dd vol_g +Ce^{\frac{3A}{2}}\left(\left\|u_0\right\|_{\rho_0^{a_0}H^{k+1}_b(U_1\cap \Sigma_0)}+\left\|u_1\right\|_{\rho_0^{a_0}H^{k}_b(U_1\cap \Sigma_0)}\right)
\end{align*}
We sum the estimates from $l=1$ to $N$ and if we take $A$ large enough so that the left-hand-side absorb the term $C\int_{U_1^c}\rho_0^{-2a_0}e^{A\rho_I}[T_s,Z_l] u\dd vol_g$ we get:
\begin{lemma}\label{estU1}
There exists a constant $C>0$ such that for all $u$ as in proposition \ref{basicConormalEstimate}:
\begin{align*}
\left\| u\right\|_{\rho_0^{a_0}H^{k+1}_b(U_1)}\leq C\left(\left\|u_0\right\|_{\rho_0^{a_0}H^{k+1}_b(U_1\cap \Sigma_0)}+\left\|u_1\right\|_{\rho_0^{a_0}H^{k}_b(U_1\cap \Sigma_0)}\right)
\end{align*}
in the strong sense that the left-hand side is finite whenever the right-hand side is finite and the inequality holds.
\end{lemma}

\subsection{Estimate on \texorpdfstring{$U_2$}{U\_ 2}}
We denote by ${}^{\mathscr{I}}TU_2$ the bundle with smooth local sections $\rho_0\partial_{\rho_0}$, $\rho_I\partial_{\rho_I}$, $\rho_I^{\frac{1}{2}} \slashed{\partial}$ (where $\slashed{\partial}$ denotes a smooth vector field on $\mathbb{S}^2$) and by ${}^{\mathscr{I}}T^*U_2$ the dual bundle. We can then define the Sobolev space $H_{\mathscr{I},b}^1$ with norm:
\begin{align*}
\left\|u\right\|_{H_{\mathscr{I},b}^1}^2 = \int_{U_2}\mathbb{m}(\rho_0\partial_{\rho_0} u)+\mathbb{m}(\rho_I\partial_{\rho_I}u) + \rho_I \slashed{G}^{\omega_i,\omega_j}\mathbb{m}(\Theta_{\omega_i} u, \Theta_{\omega_j}u)\frac{\dd \rho_0}{\rho_0}\frac{\dd \rho_I}{\rho_I}\dd \vol_{\mathbb{S}^2}
\end{align*}
Let $Z_1, Z_2, Z_3 \in \Gamma(T\mathbb{S}^2)$ be vector fields spanning $T\mathbb{S}^2$ (as a $C^{\infty}(\mathbb{S}^2)$-module).
For $k\in \N$, we define by induction the norm:
\begin{align*}
\left\|u\right\|_{H_{\mathscr{I},b}^{1,k+1}}^2 := \left\|u\right\|_{H_{\mathscr{I},b}^{1,k}}^2+\left\|\rho_0\partial_{\rho_0}u\right\|_{H_{\mathscr{I},b}^{1,k}}^2 + \left\|\rho_I\partial_{\rho_I}u\right\|_{H_{\mathscr{I},b}^{1,k}}^2+ \left\|\rho_0\partial_{\rho_0}u\right\|_{H_{\mathscr{I},b}^{1,k}}^2 + \sum_{i=1}^3 \left\|Z_i u\right\|_{H_{\mathscr{I},b}^{1,k}}^2
\end{align*}

Let $u \in H^{\tilde{r}+1}_{loc}(\mathcal{M}_{2\epsilon})$ with $\rho_0^{-a_0}\rho_I^{-a_I}T_s\rho_I u \in \overline{H}^{\tilde{r}}_b$. At first, we assume that $u=0$ on $\left\{\rho_I\geq \frac{3}{2}\alpha\right\}$. As before, for energy estimates on relatively compact sets, we can assume that $u$ is smooth. We do an energy estimate with the vector field $W = \rho_0^{-2a_0}\rho_I^{-2a_I}(-(1+c_V)\rho_I\partial_{\rho_I}+\rho_0\partial_{\rho_0})$ where $c_V$ is a constant to be chosen later.
Note that on $U_2$, we have:
\begin{align*}
\tilde{G} =& \frac{1}{\rho_I}\left(2\rho_I^2\partial_{\rho_I}^2 - 2\rho_0\rho_I\partial_{\rho_I}\partial_{\rho_0} - \rho_I\slashed{G} + \rho_I\text{Diff}^2_{\mathscr{I},b}\right)\\
\tilde{g} =& \rho_I\left(-\frac{2\dd \rho_0^2}{\rho_0^2} - 2\frac{\dd \rho_0\dd \rho_I}{\rho_0\rho_I} - \frac{1}{\rho_I}\slashed{g} + \rho_I\Gamma\left({}^{\mathscr{I}}T^*U_2^{\otimes 2}\right)\right)
\end{align*}
We consider:
\begin{align*}
T_{\mu,\nu} =& \mathbb{m}(\Theta_\mu u, \Theta_\nu u) -\frac{1}{2}\tilde{g}_{\mu,\nu}\tilde{G}^{\gamma, \delta}\mathbb{m}(\Theta_\gamma u, \Theta_\delta u) +\frac{1}{2\rho_I}\tilde{g}_{\mu,\nu}\mathbb{m}(u)\\
\div(T_{\mu,\nu}W^\mu) =& \div(T)_\mu W^\mu + T_{\mu,\nu} \pi^{\mu,\nu}\\
\pi^{\mu,\nu}:=& \frac{1}{2}\left(\nabla^\mu W^\nu + \nabla^{\nu}W^\mu\right)\\
=&\frac{1}{2}\left(\tilde{G}^{\mu, \alpha} \partial_\alpha W^\nu + \tilde{G}^{\alpha,\nu}\partial_\alpha W^\mu - W^{\alpha}\partial_\alpha \tilde{G}^{\mu,\nu}\right)
\end{align*}

A tedious computation provides:
\begin{align*}
T_{\mu,\nu}\pi^{\mu,\nu} =& \rho_0^{-2a_0}\rho_I^{-2a_I-1}\left((2(a_I-a_0)c_V+2a_I)\mathbb{m}(\rho_I\partial_{\rho_I}u)-4a_I\Re \mathbb{m}(\rho_0\partial_{\rho_0}u,\rho_I\partial_{\rho_I}u)+2a_I \mathbb{m}(\rho_0\partial_{\rho_0}u)\right.\\
&\left.+\frac{2(a_I-a_0)+2a_Ic_V-c_V-1}{2}\rho_I\slashed{G}^{\omega_i,\omega_j}\mathbb{m}(\Theta_{\omega_i}u,\Theta_{\omega_j}u+\mathbb{m}(u))\right) + R\\
\left|R\right|\leq& C\rho_I^{-2a}\rho_0^{-2a_0}\left(\mathbb{m}(\rho_I\partial_{\rho_I}u)+\mathbb{m}(\rho_0\partial_{\rho_0}u)+\rho_I\slashed{G}^{\omega_i,\omega_j}\mathbb{m}(\Theta_{\omega_i}u,\Theta_{\omega_j}u)+\mathbb{m}(u)\right)
\end{align*}

If we assume that $c_V>0$, $a_I<a_0$ and $a_I<0$, we get:
\begin{align*}
&(2(a_I-a_0)c_V+2a_I)\mathbb{m}(\rho_I\partial_{\rho_I}u)\\ &-4a_I\Re \mathbb{m}(\rho_0\partial_{\rho_0}u,\rho_I\partial_{\rho_I}u) +2a_I \mathbb{m}(\rho_0\partial_{\rho_0}u)\leq -c_1 \mathbb{m}(\rho_I\partial_{\rho_I}u) - c_2 \mathbb{m}(\rho_0\partial_{\rho_0}u)\\
&c_1 := (a_0-a_I)c_V >0\\
&c_2 := -2a_I-\frac{16a_I^2}{-8a_I+4(a_0-a_I)}>0
\end{align*}

We also check that (see Proposition \ref{divT} for a similar computation) that there exists $C_0>0$ and $C>0$ such that on $U_2$ we have:
\begin{align}
\left|\div(T)_\mu W^\mu \right|\leq& C_0\rho_I^{-2a_I}\rho_0^{-2a_0}\left(\mathbb{m}(\rho_I\partial_{\rho_I}u)+\mathbb{m}(\rho_0\partial_{\rho_0}u)+\rho_I\slashed{G}^{\omega_i,\omega_j}\mathbb{m}(\Theta_{\omega_i}u,\Theta_{\omega_j}u)\right. \notag\\
&\left.+\mathbb{m}(u)\right)+C\rho_I^{-2a_I-1}\rho_0^{-2a_0}\mathbb{m}(\rho_I\square_{\Theta,\tilde{g}}u) \label{estimateDivU2}
\end{align}

Moreover, $\rho_I^{-1}T_s\rho_I-\square_{\Theta,\tilde{g}} \in \text{Diff}^1_{\mathscr{I}}$ and therefore, if $\alpha$ is small enough we can replace $\rho_I\square_{\Theta,g}$ by $T_s\rho_I$ in \eqref{estimateDivU2} at the cost of increasing $C_0$.
As a consequence, if we choose $\alpha$ small enough, we have the following inequality on $U_2$ for some constants $c_3>0$ and $C>0$:
\begin{align}
\div(T_{\mu,\nu}\pi^{\mu,\nu})\leq& -c_3\rho_0^{-2a_0}\rho_I^{-2a_I-1}\left(\mathbb{m}(\rho_I\partial_{\rho_I}u)+\mathbb{m}(\rho_0\partial_{\rho_0}u)+\rho_I\slashed{G}^{\omega_i,\omega_j}\mathbb{m}(\Theta_{\omega_i}u,\Theta_{\omega_j}u)\right. \notag\\
&\left.+\mathbb{m}(u)\right)+C\rho_I^{-2a_I-1}\rho_0^{-2a_0}\mathbb{m}(T_s\rho_I u) \label{firstIneqU2}
\end{align}
where constant $c_3$ does not depend on $\alpha$ (as long as it is smaller than some small threshold $\alpha_0>0$).

For $0<c<\alpha$ and $c'$ large enough, we apply Stockes formula on $U_2^{c,c'} := U_2 \cap \left\{\rho_I\geq c, t+r\leq c'\right\}$ which is relatively compact in $\mathcal{M}_{\epsilon}$ (and we will take $c\to 0$ and $c'\to \infty$).
We have boundaries included in :
\begin{itemize}
\item $\left\{ \rho_I = c\right\}$ with exterior normal $-\dd \rho_I$ which is timelike and future oriented.
\item $\left\{ \tt = \frac{1}{2} \right\}$ with exterior normal $\dd \tt$ which is timelike and future oriented.
\item $\left\{ t+r = c'\right\}$ with exterior normal $\dd t + \dd r$ which is timelike and future oriented when $r$ is large (and on the boundary of $U_2$ we have $\tt \leq \frac{1}{2}$ and since for some constant $C_0>0$, $t-r-C_0\leq \tt \leq \frac{1}{2}$, we deduce that on the boundary $t+r=c'$, $r\geq \frac{c'-C_0-\frac{1}{2}}{2}$ is as large as we want if we chose $c'$ large enough).
\item $\left\{ \rho_I = 2\alpha\right\}$ with exterior normal $\dd \rho_I$ which is timelike and past oriented.
\end{itemize}

Since $u$ has support in $\left\{\rho_I\leq \frac{3}{2}\alpha\right\}$, the boundary term corresponding to $\left\{\rho_I = 2\alpha\right\}$ vanishes. The other boundary terms are non negative and we deduce:
\begin{align*}
-\int_{U_2^{c,c'}} \div(T_{\mu,\nu}W^{\mu}) \dd vol_{g'} \leq 0
\end{align*}
Combining the previous inequality with \eqref{firstIneqU2}, we get:
\begin{align}
\int_{U_2^{c,c'}}\rho_0^{-2a_0}\rho_I^{-2a_I-1}\left(\mathbb{m}(\rho_I\partial_{\rho_I}u)+\mathbb{m}(\rho_0\partial_{\rho_0}u)+\rho_I\slashed{G}^{\omega_i,\omega_j}\mathbb{m}(\Theta_{\omega_i}u,\Theta_{\omega_j}u)+\mathbb{m}(u)\right)+\rho_I^{-2a_I-1}\rho_0^{-2a_0}\mathbb{m}(T_s \rho_I u) \notag\\
\leq \frac{1}{c_3}\int_{U_2^{c,c'}}\rho_I^{-2a_I-1}\rho_0^{-2a_0}\mathbb{m}(T_s\rho_I u)\dd \vol_g \label{secondEstimateU2}
\end{align}
Since the constant $c_3$ does not depend on $c,c'$, we can use Fatou's lemma and prove that \eqref{secondEstimateU2} holds with $U_2^{c,c'}$ replaced by $U_2$.

Now we introduce $Z_0 = \rho_I\partial_{\rho_I}$, $Z_1= \rho_0\partial_{\rho_0}$, and $Z_2,Z_3,Z_4 \in \text{Diff}^1(\mathcal{B}(s))$ generating the $C^\infty(\mathbb{S}^2)$ module $\text{Diff}^1(\mathcal{B}_s)$. We then have for $1\leq i\leq 4$, $[Z_i,\rho_I^{-1}T_s\rho_I] \in \text{Diff}^2_b$ and $[Z_0, \rho_I^{-1}T_s\rho_I]=\rho_I^{-1}T_s\rho_I + \text{Diff}^2_b$. We use this property together with \eqref{secondEstimateU2} (applied to $Z_iu$) to prove by induction that for all $k\in \N$, there exists $\alpha_k>0$ such that for $\alpha = \alpha_k$ in the definition of $U_2$ and for all $u$ with support in $\left\{\rho_I\leq \frac{3}{2}\alpha_k\right\}$:
\begin{align*}
\left\|u\right\|_{\rho_I^{a_I}\rho_0^{a_0}\overline{H}_{\mathscr{I},b}^{1,k}(U_2)}\leq \frac{1}{c_3}\left\|\mathbb{m}(T_s\rho_I u)\right\|_{\rho_I^{a_I}\rho_0^{a_0}\overline{H}^{k}_b(U_2)}
\end{align*}
Note that a priori $\alpha_k$ has to depend on $k$ since at each step we have to take $\alpha$ small enough so that the term arising from the commutator can be absorbed into the left-hand side.

Let $k\in \N$. Let $\chi\in C^{\infty}(\R,[0,1])$ be equal to $1$ on $(-\infty, \alpha_k)$ and zero on $(\frac{3}{2}\alpha_k, +\infty)$.
We do not assume anymore that $u$ has support in $\left\{\rho_I\leq \frac{3}{2}\alpha_k\right\}$, we use the previous argument on $\chi u$ combined with estimate on $U_1$ and $K$ to control the commutator term $\left\|[T_s\rho_I, \chi]u\right\|_{\rho_0^{a_0}\rho_I^{a_I}\overline{H}^{k}_b(U_2)}\leq C\left(\left\|u\right\|_{\overline{H}^{k+1}(K)}+\left\|u\right\|_{\overline{H}^{k+1}_b(U_1)}\right)$ (since these estimates do not impose any restriction on $\alpha$, we can apply them with $\alpha = \alpha_k$) and we conclude:
\begin{lemma}\label{estU2}
For all $k\in \N$, there exists constants $C_k>0$ and $\alpha_k>0$ such that for all $u$ as in proposition \ref{basicConormalEstimate} with $\tilde{r}\geq k$:
\begin{align*}
\left\|u\right\|_{\rho_I^{a_I+1}\rho_0^{a_0}\overline{H}_{\mathscr{I},b}^{1,k}(U_2)}\leq C_k\left(\left\|u_0\right\|_{\rho_0^{a_0}\overline{H}^{k+1}_b(\Sigma_0)}+\left\|u_1\right\|_{\rho_0^{a_0}\overline{H}^{k}_b(\Sigma_0)}\right)
\end{align*}
where we have taken $\alpha = \alpha_k$ in the definition of $U_2$.
\end{lemma}

\subsection{Estimate on \texorpdfstring{$U_3$}{U\_ 3}}
We define ${}^{\mathscr{I}}T U_3$ and $H_{\mathscr{I},b}^{1,k}$ which are the natural extension to the corresponding spaces defined on $U_2$ but since $U_3$ is not included in $U$, we use vector fields $x\partial_x$, $\partial_{\mathfrak{t}}$ instead of $\rho_I\partial_{\rho_I}$ and $\rho_0\partial_{\rho_0}$. We also introduce the bundle ${}^{\beta}TU_3$ which is generated by the sections $x^{\frac{1}{2}}\partial_{\tt}$, $x\partial_x$ and $x^{\frac{1}{2}}\slashed{\partial}$ (where $\slashed{\partial}$ denotes a smooth vector field on $\mathbb{S}^2$).
We compute:
\begin{align*}
\tilde{G} =& x^{-1}\left(-2\partial_{\tt}x\partial_x -x \slashed{G} + x\Gamma\left({}^{\mathscr{I}}TU_3^{\otimes 2}\right)\right)\\
\tilde{G} =& x^{-1}\left((2M^2-a^2\sin^2\theta)x\partial_{\tt}^2-2\partial_{\tt}x\partial_x - x\slashed{G} + x^{\frac{1}{2}}\Gamma\left({}^{\beta}TU_3^{\otimes 2}\right)\right)
\end{align*}
As before, we obtain the estimate on $U_3$ as limit of estimate on the relatively compact set $U_{3,\eta}:=U_3\cap \left\{\frac{x}{1+\tilde{t}}\geq \eta\right\}$.
We define the norm $\mathfrak{H}^1(U_3)$ by 
\begin{align*}
\left\|u\right\|^2_{\mathfrak{H}^1(U_3)} = \int_{U_3} x\mathbb{m}(\partial_{\tt}u) + x^2\mathbb{m}(\partial_{x}u) + x\slashed{G}^{i,j}\mathbb{m}(\Theta_{i}u, \Theta_{j}u) + \mathbb{m}(u)\frac{\dd x}{x}\dd \tt \dd^2 \omega
\end{align*}
the norm $\mathfrak{H}^1(U_{3,\eta})$ with the same expression but integration on $U_{3,\eta}$.
For $(A_i)_{i=1}^3$ a set of vector fields spanning the $C^{\infty}(\mathbb{S}^2)$ module of vector fields on $\mathbb{S}^2$, we define:
\begin{align*}
\left\|u\right\|^2_{\mathfrak{H}^{1,k}(U_3)} = \left\|u\right\|^2_{H_b^k(U_3)} + \sum_{i=1}^3\left\|x^{\frac{1}{2}}\Theta_{A_i}u\right\|^2_{H_b^k(U_3)} + \left\|x^{\frac{1}{2}}\partial_\tt u\right\|^2_{H_b^k(U_3)} + \left\|x\partial_x u\right\|^2_{H_b^k(U_3)}
\end{align*}

Let $u \in H^{\tilde{r}+1}_{loc}(\mathcal{M}_{2\epsilon})$ with $x^{-a_I}T_s x u \in \overline{H}^{\tilde{r}}_b$. At first, we assume that $u=0$ on $\left\{\tilde{t}\leq \frac{3}{8}\right\}$. As before, for energy estimates on relatively compact sets, we can assume that $u$ is smooth. We prove that for $A$ large enough:
\begin{align*}
\left\|x^{-a_I}e^{-\frac{A\tt}{2}}\partial_\tt u\right\|^2_{L^2_b(U_{3,\eta})} + A\left\|x^{-a_I}e^{-\frac{A\tt}{2}} u\right\|_{\mathfrak{H}^1(U_{3})}^2\leq& C\left\|e^{-\frac{A\tt}{2}}x^{-a_I}T_s xu\right\|^2_{L^2_b(U_3)}.
\end{align*}
Then, we prove by induction on $k$ that for $u\in H^{1,k}_{(loc)}(\mathcal{M}_+)$ with $x^{-a_I}T_s xu \in H^k_b(U_3)$:
\begin{align*}
\left\|x^{-a_I}e^{-\frac{A\tt}{2}}\mathbb{m}(\Theta_\tt u)\right\|^2_{H^k_b(U_{3,\eta})} + A\left\|x^{-a_I}e^{-\frac{A\tt}{2}} u\right\|^2_{\mathfrak{H}^{1,k}(U_{T})}\leq& C\left\|e^{-\frac{A\tt}{2}}x^{-a_I}T_s xu\right\|^2_{H^k_b(U_3)}.
\end{align*}

We define the energy momentum tensor as follows:
\begin{align*}
T^{\delta,\gamma}(u) = \Re (\mathbb{m}(\Theta_\mu u, \Theta_\nu u))\tilde{G}^{\mu,\delta}\tilde{G}^{\nu, \gamma} - \frac{1}{2}\tilde{G}^{\delta,\gamma}\tilde{G}^{\mu, \nu}\mathbb{m}(\Theta_\mu u, \Theta_\nu u) + \frac{1}{2}\tilde{G}^{\delta, \gamma}x^{-1}\mathbb{m}(u)
\end{align*}
We define the vector fields:
\begin{align*}
V :=& \left(-\partial_{\tt}+x\partial_x\right) \in \Gamma({}^{\mathscr{I}}TU_3)\\
W :=& e^{-A\tt}x^{-2a_I}V
\end{align*}
Then we define the energy one form:
\begin{align*}
J^W_\mu := T_{\nu,\mu}W^{\mu}
\end{align*}
The divergence is given by:
\begin{align*}
K^W(u)= \div(T)_\mu W^{\mu} + T_{\nu,\mu}\nabla^{\nu} W^{\mu}
\end{align*}
Moreover, we have (see Proposition \ref{divT}):
\begin{align*}
\div(T)_\mu W^{\mu} :=& e^{-A\tt}x^{-2a_I}\left(\vphantom{\frac{1}{2}}\Re(\mathbb{m}(\Theta_{V} u, \square_{\tilde{g},\Theta} u))+ V^{\mu}\Re(\tilde{G}^{\gamma,\nu}\mathbb{m}(R^\Theta_{\gamma,\mu}u, \Theta_\nu u))+ \Re((\Theta_\gamma \mathbb{m})(\Theta_{V} u, \Theta_\nu u))\tilde{G}^{\nu, \gamma}\right.\\
& \left.-\frac{1}{2}\Re((\Theta_{V} \mathbb{m})(\Theta_\mu u, \Theta_\nu u))\tilde{G}^{\mu,\nu} + \frac{1}{2}(\Theta_V x^{-1}\mathbb{m})(u,u)+\Re(x^{-1}\mathbb{m}(\Theta_{V} u,u))\right)
\end{align*}
Using that $\mathbb{m}$ and $\Theta$ are smooth up to $x=0$ and independent of $\tt$ and $x$ (in particular $\Theta_V \mathbb{m} = 0$ and $\Theta_V x^{-1}\mathbb{m} = -x^{-1}\mathbb{m}$ is of order $x^{-1}$) and Cauchy-Schwarz estimates, we get for all $\epsilon>0$, there exists $C_\epsilon>0$ (independent of $A$ and $u$) such that:
\begin{align*}
\left|\div(T)_\mu W^{\mu}\right|\leq& e^{-A\tt}x^{-2a_I-1}\left(\epsilon \mathbb{m}(\partial_{\tt}u)+C_\epsilon x^2\mathbb{m}(\Theta_{\partial_x}u)+C_\epsilon x\slashed{G}^{i,j}\mathbb{m}(\Theta_{i}u, \Theta_{j}u)+C_\epsilon\mathbb{m}(u)\right. \\
&\hphantom{e^{-A\tt}x^{-2a_I-1}}\left. \vphantom{\slashed{G}^{i,j}} + C_\epsilon \mathbb{m}(x\square_{\tilde{g},\Theta} u)\right) + R(\tt, x, \omega)\\
\left|R(\tt, x, \omega)\right| \leq& Ce^{-A\tt}x^{-2a_I}\left( \mathbb{m}(\partial_{\tt}u)+x^2\mathbb{m}(\Theta_{\partial_x}u)+x\slashed{G}^{i,j}\mathbb{m}(\Theta_{i}u, \Theta_{j}u)+\mathbb{m}(u) \right)
\end{align*}
where $C$ depends only on the metric. In the sequel, we use this estimate with $\epsilon = -a_I>0$.

In order to replace $\square_{\tilde{g},\Theta}$ by $x^{-1}T_s x$ in the estimate (allowing the constant $C$ to change), we check that $\square_{\tilde{g},\Theta} - x^{-1}T_s x \in \text{Diff}^1_{\mathscr{I}^+}$. We compute (using trivialisation $\mathcal{T}_m$):
\begin{align*}
J =& \frac{\sin\theta}{1+x^2a^2\cos^2\theta} = \sin(\theta)(1+O(x^2))\\
\square_{\tilde{g},\Theta} - x^{-1}T_s x =& J^{-1}\Theta_\mu J\tilde{G}^{\mu,\nu} \Theta_\nu -x^{-1}T_s x\\
=& \tilde{G}^{\mu,\nu}\partial_\mu \partial_\nu + \partial_\mu(J\tilde{G}^{\mu,\nu})\partial_\nu + is\frac{\cos\theta}{\sin^2\theta}\partial_\phi - s^2\cotan^2\theta - x^{-1}T_s x\\
 = &\left(s+xC^{\infty}\left(\left[0,\frac{1}{r_+-\epsilon}\right]_x\times\mathbb{S}^2\right)\right)\\
&+\left(2s(2M-ia\cos\theta)+xC^{\infty}\left(\left[0,\frac{1}{r_+-\epsilon}\right]_x\times\mathbb{S}^2\right)\right)\partial_{\tt}\\
& \left(2s+xC^{\infty}\left(\left[0,\frac{1}{r_+-\epsilon}\right]_x\times\mathbb{S}^2\right)\right)x\partial_{x}\\
 & +\left(2as+xC^{\infty}\left(\left[0,\frac{1}{r_+-\epsilon}\right]_x\times\mathbb{S}^2\right)\right)x\partial_{\phi} - \frac{2 a^{2} x^{2} \sin{\left(\theta \right)} \cos{\left(\theta \right)}}{a^{2} x^{2} \cos^{2}{\left(\theta \right)} + 1}\partial_{\theta}
\end{align*}
Moreover, we have:
\begin{align*}
T_{\nu,\mu}\nabla^{\nu} W^{\mu} =& \frac{1}{2}\left(\tilde{G}^{\mu,\alpha}\nabla_\alpha W^\nu +\tilde{G}^{\nu,\alpha}\nabla_\alpha W^\mu\right)\Re \mathbb{m}(\Theta_\mu u, \Theta_\nu u)\\ & - \frac{1}{4}\tilde{g}_{\mu,\nu} \left(\tilde{G}^{\mu,\alpha}\nabla_\alpha W^\nu + \tilde{G}^{\nu,\alpha}\nabla_\alpha W^\mu \right)\left(\tilde{G}^{\kappa,\delta}\Re\mathbb{m}(\Theta_\kappa u, \Theta_\delta u)-x^{-1}\mathbb{m}(u)\right)
\end{align*}
Using the defintion of the Levi-Civita connection, we get:
\begin{align*}
\pi^{\mu,\nu}:=\tilde{G}^{\mu,\alpha}\nabla_\alpha W^\nu + \tilde{G}^{\alpha,\nu}\nabla_{\alpha}W^\mu = \tilde{G}^{\mu, \alpha} \partial_\alpha W^\nu + \tilde{G}^{\alpha,\nu}\partial_\alpha W^\mu - W^{\alpha}\partial_\alpha \tilde{G}^{\mu,\nu}
\end{align*}

We get the following intermediary computation results:
\begin{align*}
\tilde{G}^{x,\alpha}\partial_\alpha W^x =& Ax^{-2a_I+1}e^{-A\tt}+ Ae^{-A\tt}O(x^{-2a_I+2})\\
\tilde{G}^{\tt, \alpha}\partial_\alpha W^x =& -(-2a_I+1)x^{-2a_I}e^{-A\tt}+Ae^{-A\tt}O(x^{-2a_I+1})\\
\tilde{G}^{x,\alpha}\partial_\alpha W^{\tt} =& -Ax^{-2a_I}e^{-A\tt}+Ae^{-A\tt}O(x^{-2a_I+1})\\
\tilde{G}^{\tt, \alpha}\partial_\alpha W^{\tt} =& (-2a_I+(2M^2-a^2\sin^2\theta)Ax) x^{-2a_I-1}e^{-A\tt}+e^{-A\tt}O(x^{-2a_I}+Ax^{-2a_I+\frac{1}{2}})
\end{align*}
All the other terms of this form only contribute to the error part (i.e have a decay as a component of a section of $Ae^{-A\tt}x^{-2a_I}{}^{\mathscr{I}}TU_3\otimes {}^{\mathscr{I}}TU_3$).
Note that $W^{\alpha}\partial_\alpha \tilde{G}^{\mu,\nu}$ decays like the components of a section of $e^{-A\tt}x^{-2a_I}{}^{\mathscr{I}}TU_3\otimes {}^{\mathscr{I}}TU_3$ and are therefore in the error term.
We compute:
\begin{align*}
\tilde{g}_{\mu,\nu}\pi^{\mu,\nu} =& (A-2a_I+1)x^{-2a_I}e^{-A\tt}+Ae^{-A\tt}O(x^{-2a_I +1})
\end{align*}

We conclude:
\begin{align*}
T_{\nu,\mu}\nabla^{\nu} W^{\mu} =& x^{-2a_I-1}e^{-A\tt}\left(\left(-2a_I + Ax\frac{2M^2-a^2\sin^2\theta}{2}\right) \mathbb{m}(\Theta_\tt u) + Ax^2 \mathbb{m}(\Theta_x u)\right.\\&\left. + \frac{1}{2}(A-2a_I+1)\slashed{G}^{\omega_i,\omega_j}x \mathbb{m}\left(\Theta_{\omega_i}u,\Theta_{\omega_j}u\right)+\frac{1}{2}(A-2a_I+1)\mathbb{m}(u)\right)+R(\tt, x, \omega)\\
\left|R(\tt, x, \omega)\right|\leq& (A+C(Ax^{\frac{1}{2}}+1))x^{-2a_I}e^{-A\tt}\frac{2M^2-a^2\sin^2\theta}{4}\mathbb{m}(\Theta_t u)\\&+Cx^{-2a_I}Ae^{-A\tt}\left(x^2 \mathbb{m}(\Theta_x u) + \slashed{G}^{\omega_i,\omega_j}x \mathbb{m}\left(\Theta_{\omega_i}u,\Theta_{\omega_j}u\right)+\mathbb{m}(u)\right)
\end{align*}
where $C$ depends on $a_I$ and on the bounds of the metric coefficients (and their derivatives) but (crucially) not on $A$ (nor on $T$ since $\tilde{G}$ is stationary).

We can find $\epsilon>0$ and $A_0>0$ (both independent of $u$) such that there exists a constant $C>0$ such that for all $A>A_0$ the following inequality holds on $\left\{ x\leq \epsilon\right\}$:
\begin{align}
K_W \geq& \frac{x^{-2a_I-1}e^{-A\tt}}{2}\left(\left(-a_I + Ax\frac{2M^2-a^2\sin^2\theta}{4}\right) \mathbb{m}(\Theta_\tt u) + \frac{A}{2}x^2 \mathbb{m}(\Theta_x u)\right.\notag \\ &\left.+ \frac{1}{4}(A-2a_I+1)\slashed{G}^{\omega_i,\omega_j}x \mathbb{m}\left(\Theta_{\omega_i}u,\Theta_{\omega_j}u\right)+\frac{1}{4}(A-2a_I+1)\mathbb{m}(u)-C\mathbb{m}(T_sx u)\right)
\label{energyEstimateNearNullInf}
\end{align}
and, even restricting $\epsilon$, we can require that $V$ is uniformly timelike on $\left\{ \frac{\epsilon}{2}\leq x\leq \epsilon\right\}$

This is the part of the estimate near $x=0$. 
We can modify $V$ on $\left\{x>\epsilon\right\}$ so that it is uniformly timelike on $\left\{ x>\epsilon\right\}$, smooth and stationnary. From now on, $V$ denotes such a modified vector field. Note that on $\left\{x>\epsilon\right\}$, we have:
\begin{align*}
\left|\div(T)_\mu W^{\mu}\right| \leq& Ce^{-A\tt}\left(\mathbb{m}(\square_{\tilde{g},\Theta} u))+ \mathbb{m}(\Theta_\tt u) + \mathbb{m}(\Theta_x u) + \slashed{G}^{\omega_i,\omega_j}\mathbb{m}(\Theta_{\omega_i}u,\Theta_{\omega_j}u) + \mathbb{m}(u)\right)
\end{align*}
for some constant $C$ (independent of $A$ and $u$).
\begin{align*}
T_{\nu,\mu}\nabla^{\nu} W^{\mu} =& T_{\nu,\mu}V^{\mu}\nabla^{\nu}\left(x^{-2a_I}e^{-A\tt}\right) + x^{-2a_I}e^{-A\tt}T_{\mu,\nu}\nabla^{\nu}V^{\mu}\\
\left\|T_{\mu,\nu}\nabla^{\nu}V^{\mu}\right\| \leq& C\left(\mathbb{m}(\Theta_t u) + \mathbb{m}(\Theta_x u) + \slashed{G}^{\omega_i,\omega_j} \mathbb{m}(\Theta_{\omega_i}u,\Theta_{\omega_j}u) + \mathbb{m}(u)\right)
\end{align*}
Since $V$ and $\frac{-2a_I\nabla x}{Ax}-\nabla \tt$ are uniformly timelike (for $A$ large enough) on $\left\{x>\epsilon\right\}$, we get as usual in hyperbolic estimates:
\begin{align*}
 T_{\nu,\mu}V^{\mu}\nabla^{\nu}\left(x^{-2a_I}e^{-A\tt}\right)\geq C_\epsilon Ae^{-At}\left(\mathbb{m}(\Theta_t u)+ \mathbb{m}(\Theta_x u)+ \slashed{G}^{\omega_i,\omega_j}\mathbb{m}(\Theta_i u, \Theta_j u) + \mathbb{m}(u)\right)
\end{align*}
Therefore, if we choose $A$ large enough, we have the following inequality on $\left\{x>\epsilon\right\}$:
\begin{align*}
K_W \geq Ce^{-A\tt}\left((A\mathbb{m}(\Theta_{\tt} u) + A\mathbb{m}(\Theta_x u) + A\slashed{G}^{\omega_i,\omega_j}\mathbb{m}(\Theta_{\omega_i}u,\Theta_{\omega_j}u)+A\mathbb{m}(u)-\mathbb{m}(T_s xu)\right)
\end{align*}
with $C>0$ independent of $u$, $\tt$ and $A$. Combining this estimate with \eqref{energyEstimateNearNullInf}, we get that there exists a constant $C>0$ independent of $A$ (as long as it is large enough) and $u$ such that for all $(\tt, x,\omega)$:
\begin{align*}
K_W \geq Cx^{-2a_I-1}e^{-A\tt}\left((1+Ax)\mathbb{m}(\Theta_\tt u) + Ax^2 \mathbb{m}(\Theta_x u) + Ax\slashed{G}^{\omega_i,\omega_j}\mathbb{m}(\Theta_{\omega_i}u,\Theta_{\omega_j}u)+ A\mathbb{m}(u)-\mathbb{m}(T_s xu)\right)
\end{align*}

We can then apply the Stockes formula on the domain $U_{3,\eta}$ for $\eta$ so small that $\tilde{G}(\dd \frac{x}{1+\tt}) = \frac{x}{(1+\tt)^3} + O\left(\frac{x^2}{(1+\tt)^2}\right)>0$ on the boundary face $\left\{ \frac{x}{1+\tt} = \eta, \frac{1}{4}\leq \tt\leq T\right\}$. Therefore, the outward normal is causal future oriented on each boundary face of $U_{3,\eta}$ except at the face $\left\{\tt =\frac{1}{4}, \frac{x}{1+\tt}\geq \eta\right\}$ but $u=0$ near this face. Therefore, we get (for some constant $C>0$ independent of $u$ and $T$):
\begin{align*}
\left\|x^{-a_I}e^{-\frac{A\tt}{2}}\partial_\tt u\right\|^2_{L^2_b(U_{3,\eta})} + A\left\|x^{-a_I}e^{-\frac{A\tt}{2}} u\right\|_{\mathfrak{H}^1(U_{3,\eta})}\leq C\left(\left\|e^{-\frac{A\tt}{2}}x^{-a_I}T_s xu\right\|_{L^2_b(U_3)}\right)
\end{align*}

By the monotone convergence theorem we can take $\eta\to 0$ and this concludes the first part of the estimate.

We now prove the higher regularity statement for $k>0$ by induction. We assume that $u \in H^{k+1}_{loc}(\mathcal{M}_{2\epsilon})$ with $x^{-a_I}T_s x u \in \overline{H}^{k}_b$ and that the high regularity statement is true at order $k-1$. As before, we also assume that $\tt\geq \frac{3}{8}$ on $\supp(u)$.
We introduce the set of commutators $\mathcal{C}:= \left\{Z_0:=\partial_{\tt}, Z_1:=\tt \partial_{\tt}-x\partial_x, Z_2, Z_3, Z_4\right\}$ where $Z_2, Z_3, Z_4 \in \text{Diff}^1(\mathcal{B}_s)$ are generators of the $C^\infty(\mathbb{S}^2)$ module $\text{Diff}^1(\mathcal{B}_s)$ modulo $\text{Diff}^0(\mathcal{B}_s)$.
We consider $v_i = Z_i u$ which satisfy the induction hypothesis. Therefore, we have:
\begin{align*}
\left\|x^{-a_I}e^{-\frac{A\tt}{2}}\partial_\tt v_i\right\|^2_{H^{k-1}_b(U_{T})} + A\left\|x^{-a_I}e^{-\frac{A\tt}{2}} v_i\right\|^2_{\mathfrak{H}^{1,k-1}(U_{T})}\leq& C\left(\left\|e^{-\frac{A\tt}{2}}x^{-a_I}T_s xv_i\right\|^2_{H^{k-1}_b(U_3)}\right).
\end{align*}
We also have the estimate of order $k-1$ on $u$ instead of $v_i$:
\begin{align*}
\left\|x^{-a_I}e^{-\frac{A\tt}{2}}\partial_\tt u\right\|^2_{H^{k-1}_b(U_{T})} + A\left\|x^{-a_I}e^{-\frac{A\tt}{2}} u\right\|^2_{\mathfrak{H}^{1,k-1}(U_{T})}\leq& C\left(\left\|e^{-\frac{A\tt}{2}}x^{-a_I}T_s xu\right\|^2_{H^{k-1}_b(U_3)}\right).
\end{align*}
Considering all these estimates, we deduce that $x^{-a_I}u \in H^{1,k}_{\mathscr{I},b}(U_3)$ and $\partial_\tt u \in H^{k}_b(U_3)$. In particular, considering the definition of $H^{1,k}_{\mathscr{I},b}(U_3)$, we deduce that $x^{-a_I+1}u\in H^{k+1}_b(U_3)$ and $\left\|x^{-a_I+1}u\right\|_{H^{k+1}_b(U_3)}\leq C\left\|x^{-a_I}u\right\|_{H^{1,k}_{\mathscr{I},b}(U_3)}$ with $C$ independent of $u$.
Note that for all $0\leq i\leq 4$, $[x^{-1}T_s x, Z_i]\in \text{Diff}^2_b$ and therefore we get:
\begin{align*}
\left\|x^{-a_I}e^{-\frac{A\tt}{2}}\partial_\tt v_i\right\|^2_{H^{k-1}_b(U_{T})} + A\left\|x^{-a_I}e^{-\frac{A\tt}{2}} v_i\right\|^2_{\mathfrak{H}^{1,k-1}(U_{T})}\leq& C\left(\left\|e^{-\frac{A\tt}{2}}Z_i x^{-a_I}T_s x u\right\|^2_{H^{k-1}_b(U_3)}\right.\\ & \hphantom{C(}\left.+ \left\|e^{-\frac{A\tilde{t}}{2}}x^{-a_I+1}u\right\|^2_{H^{k+1}_b(U_3)}\right)
\end{align*}
Summing all the estimates (and the one on $u$) and taking $A$ large enough, we find:
\begin{align*}
\left\|x^{-a_I}e^{-\frac{A\tt}{2}}\partial_\tt u\right\|^2_{H^{k}_b(U_{T})} + A\left\|x^{-a_I}e^{-\frac{A\tt}{2}} u\right\|^2_{\mathfrak{H}^{1,k}(U_{T})}\leq& C\left(\left\|e^{-\frac{A\tt}{2}}x^{-a_I}T_s xu\right\|^2_{H^{k}_b(U_3)}\right)
\end{align*}

In particular, using the fact that $\frac{1}{2}\leq \tt \leq T$ on $U_3$, we get:
\begin{align*}
\left\|u\right\|^2_{x^{a_I}H^{1,k}_{\mathscr{I},b}(U_3)} \leq Ce^{AT}\left\|T_s xu\right\|^2_{H^k_b(U_3)}
\end{align*}
where constants $A$ and $C$ does not depend on $T$ (but could depend on $k$).

We apply the previous bound to $\chi(\tt)u$ where $\chi$ is a cutoff equal to $1$ on $\left\{\tt\geq \frac{1}{2} \right\}$ and equal to zero on $\left\{ \tt \leq \frac{3}{8}\right\}$. We use estimates on $K$ and $U_2$ to control the term $\left\|[T_s x, \chi]u\right\|_{H^k_b(U_3)}\leq C\left(\left\|u\right\|_{\overline{H}^{k+1}_b(K)}+\left\|u\right\|_{\overline{H}^{1,k}_b(U_2)}\right)$.
Therefore we get:
\begin{lemma}\label{estU3}
For all $k\in \N$, there exists a constant $C>0$ such that for all $u$ as in proposition \ref{basicConormalEstimate} with $\tilde{r}\geq k$ and all $T>\frac{1}{2}$:
\begin{align*}
\left\|u\right\|^2_{x^{a_I+1}H^{1,k}_{\mathscr{I},b}(U_3)} \leq Ce^{CT}&\left(\left\|T_s u\right\|^2_{\rho_I^{a_I}\rho_0^{a_0}\overline{H}^{k}_b(U_2\cup K \cup U_1)}+\left\|T_s u\right\|^2_{x^{a_I}\overline{H}_b^{k}(U_3)}+ \left\|u_0\right\|^2_{\rho_0^{a_0}\overline{H}^{k+1}_b(\Sigma_0)}\right.\\
&\hphantom{(}\left.+\left\|u_1\right\|^2_{\rho_0^{a_0}\overline{H}^{k}_b(\Sigma_0)}\right)
\end{align*}
\end{lemma}

\subsection{Proof of Proposition \ref{basicConormalEstimate}}

\begin{lemma}\label{traceBSpace}
There exists $C>0$ such that for all $v\in H^{1,\tilde{r}}_{\mathscr{I},b}$ and for all $T>0$, $v \in C^{0}([0,T]_{\tt}, \overline{H}_b^{\tilde{r}})$ and:
\begin{align*}
\sup_{\tt \in [0,T]} \left\|v(\tt)\right\|_{\overline{H}^{\tilde{r}}_b}\leq C\left\|v\right\|_{H^{1,\tilde{r}}_{\mathscr{I},b}(\left\{ 0 \leq \tt \leq T\right\})}
\end{align*}
\end{lemma}
\begin{proof}
First, we define two open sets of the form $(-1, +\infty)_\tt \times \left(0, \frac{1}{r_+-\frac{3\epsilon}{2}}\right)_x\times U_i$ where $U_i$ is diffeomorphic to a disc for $i\in \left\{1,2\right\}$ which cover $(-1, +\infty)_\tt \times \left(0, \frac{1}{r_+-\epsilon}\right)_x\times\mathbb{S}^2$. Using a subordinated partition of unity, local trivializations (associated to local coordinates $(\tt, x, y_1, y_2)$) and after the change of variable $x= e^{-y_0}$, it is enough to prove the statement: if $v\in H^{\tilde{r}}((-1,+\infty)_\tt\times \R_y^3)$ and $\partial_\tt v \in H^{\tilde{r}}((-1,+\infty)_\tt\times \R_y^3)$ (we call $\mathcal{Y}$ the space of such functions), then $v \in C^{0}(\R_{\tt}, H^{\tilde{r}}(\R_y^3))$ and:
\begin{align*}
\sup_{\tt \in \R} \left\|v(\tt)\right\|_{H^{\tilde{r}}(\R_y^3)}\leq C\left(\left\|v\right\|_{H^{\tilde{r}}((-1,+\infty)_\tt\times\R_y^3)} + \left\|\partial_{\tt}v\right\|_{H^{\tilde{r}}((-1,+\infty)_\tt\times\R_y^3)}\right)
\end{align*}
for a constant $C>0$ independent of $v$. We first assume that $v$ is smooth and compactly supported. We denote by $\hat{v}$ the spatial Fourier transform of $v$ and by $\mathcal{F}(v)$ the total Fourier transform of $v$. We have:
\begin{align*}
\left\|v(0)\right\|^2_{H^{\tilde{r}}(\R_y^3)} =& \int \left<\xi\right>^{2\tilde{r}}\left|\hat{v}(0, \xi)\right|^2 \dd \xi\\
\hat{v}(0,\xi) =& \frac{1}{2\pi}\int \int e^{-i\tt\sigma}\hat{v}(\tt,\xi)\dd \tt \dd \sigma.
\end{align*}
Thus,
\begin{align*}
\left|\hat{v}(0,\xi)\right|^2\leq & \frac{1}{4\pi^2}\int (1+\sigma^2)^{-1}\dd \sigma \int (1+\sigma^2) \left|\int e^{-i\tt\sigma}\hat{v}(\tt,\xi)\dd \tt\right|^2 \dd \sigma\\
\leq& \frac{1}{4\pi}\left(\int \left|\mathcal{F}(v)(\sigma,\xi)\right|^2\dd \sigma+ \int \left|\mathcal{F}(\partial_t v)(\sigma,\xi)\right|^2\dd \sigma\right)\\
\left\|v(0)\right\|^2_{H^{\tilde{r}}(\R_y^3)}\leq& \frac{1}{4\pi}\left(\left\|v\right\|_{H^{\tilde{r}}(\R_y^3)}^2 + \left\|\partial_t v\right\|_{H^{\tilde{r}}(\R_y^3)}\right).
\end{align*}
The previous computation proves that the trace on $\left\{\tt = 0\right\}$ extends continuously to a linear map from $\mathcal{Y}$ to $H^{\tilde{r}}(\R_y^3)$. We call this map $\gamma_0$. We define $\gamma_{\tt} v = \gamma_{0} (\tau_{\tt}v)$ where $\tau_{\theta}v$ is the pullback of $v$ by the diffeomorphism $(t,y)\mapsto (t+\theta, y)$ ($\gamma_{\tt}$ is the continuous extension of the trace on $\left\{\tt\right\}\times R_y^3$). Note that since $\tt \mapsto \tau_{\tt}v$ is continuous from $\R$ to $\mathcal{Y}$, we get that $\tt \mapsto \gamma_\tt v$ is continuous from $\R$ to $H^{\tilde{r}}(\R^3_y)$.
\end{proof}

\begin{proof}[Proof of Proposition \ref{basicConormalEstimate}]
In order to have weights defined on $\textbf{M}_\epsilon$, we introduce \begin{align*}
\tilde{\rho}_I := \chi_0(\tt)\rho_I + (1-\chi_0(\tt))x\\
\tilde{\rho}_0:= \chi_1(\rho_0)\rho_0 + (1-\chi_1(\rho_0))
\end{align*} where $\chi_0(\tt)=1$ when $\tt \leq -2$, $\chi_0(\tt) = 0$ when $\tt\geq -1$, $\chi_1(\rho_0) = 1$ when $\rho_0\leq 1$ and $\chi(\rho_0)=0$ on $\textbf{M}_{\epsilon}\setminus \left\{\rho_0\leq 2\right\}$. 
We define $U_4 := \left\{t_0\geq 0, \tt \leq T, x\leq \frac{1}{r_+-\epsilon}\right\}$.
Combining Lemmas \ref{estU1}, \ref{estU2} and \ref{estU3}, we get:
\begin{align*}
\left\|u\right\|^2_{\tilde{\rho}_I^{a_I+1}\tilde{\rho}_0^{a_0}H^{1,k}_{\mathscr{I},b}(U_4)}\leq Ce^{CT}\left(\left\|u_0\right\|^2_{\rho_0^{a_0}\overline{H}^{k+1}_b(\Sigma_0)}+\left\|u_1\right\|^2_{\rho_0^{a_0}\overline{H}^{k}_b(\Sigma_0)}\right)
\end{align*}
By Lemma \ref{traceBSpace}, we deduce that for all $T>0$, $u\in C^0\left([0,T]_{\tt}, x^{a_I+1}\overline{H}^{\tilde{r}}_b\right)$ and there exists $C>0$ such that:
\begin{align*}
\sup_{t\in [0,T]} \left\|x^{-a_I-1}u(t)\right\|_{\overline{H}_b^{\tilde{r}}} \leq Ce^{CT}\left(\left\|u_0\right\|^2_{\rho_0^{a_0}\overline{H}^{k+1}_b(\Sigma_0)}+\left\|u_1\right\|^2_{\rho_0^{a_0}\overline{H}^{k}_b(\Sigma_0)}\right)
\end{align*}
If $\tilde{r}\geq k$ for $k\in \N \setminus\left\{0\right\}$ , we can apply Lemma \ref{traceBSpace} to $\partial_{\mathfrak{t}}^k u$ which is (locally with respect to $\mathfrak{t}$) in $H_{\mathscr{I},b}^{1,\tilde{r}-k}$. We deduce $u\in C^k(\R_{\mathfrak{t}}, \overline{H}^{\tilde{r}-k}_b)$
\end{proof}

\subsection{Proof of Proposition \ref{refinementAtScri}}
For $\tilde{r}\in \N$, $a_0,a_I \in \R$, we define the space $\mathcal{H}_{b}^{\tilde{r}, 0, 0}$ as the space of sections of $\mathcal{B}(s,s)$ defined on $U$ with support in $\left\{\rho_I\leq 1\right\}$ and with index of $b$ regularity (with respect to the boundary $\left\{\rho_I=0\right\}\cup\left\{\rho_0=0\right\}$) equal to $\tilde{r}$. We then define $\mathcal{H}_{b}^{\tilde{r},a_I,a_0} := \rho_0^{a_0}\rho_I^{a_I}\mathcal{H}_{b}^{\tilde{r},0,0}$ and $\overline{\mathcal{H}}_{b}^{\tilde{r},a_I,a_0}$ the space of restriction to $\left\{\rho_0\leq 1\right\}$ of elements of $\mathcal{H}_{b}^{\tilde{r},a_0,a_I}$.
\begin{prop}[Elementary transport estimate]
\label{estimateTransport1}
Let $\tilde{r}\in \N$ and $\beta>\alpha$. 
Let $\chi:\R_{\rho_I}\rightarrow \R$ be a smooth cutoff equal to 1 near 0 and vanishing for $\rho_I\geq 1$.
For all $\nu \in (\alpha,\beta]$ and all $N\in \N$, there exists $C_{\nu,N}>0$ such that for all $u\in \mathcal{H}_b^{\tilde{r},-N,\beta}$:
\begin{align*}
\left\|\chi u\right\|_{\overline{\mathcal{H}}_{b}^{\tilde{r},\alpha,\beta}}\leq C_{\nu,N}\left(\left\|u\right\|_{\overline{\mathcal{H}}_{b}^{\tilde{r},-N,\beta}}+\left\|(\rho_I\partial_{\rho_I}-\rho_0\partial_{\rho_0})u\right\|_{\overline{\mathcal{H}}_{b}^{\tilde{r},\nu,\beta}}\right)
\end{align*}
in the strong sense that if the right-hand side is finite, $\chi u \in \overline{\mathcal{H}}_{b}^{\tilde{r},\alpha,\beta}$ and the inequality holds.
\end{prop}
\begin{proof}
Let $f\in\mathcal{D}:=\left\{ \phi_{|_{(0,1)_{\rho_I}\times (0,1)_{\rho_0}\times \mathbb{S}^2}}, \phi\in C^{\infty}_c((0,1)_{\rho_I}\times (0,+\infty)_{\rho_0} \times \mathbb{S}^2)\right\}$. 
We consider the following map:
\[S:\begin{cases}
\mathcal{D}\rightarrow \overline{\mathcal{H}}_{b}^{\tilde{r},\alpha,\beta}\\
f\mapsto \int_{0}^{-\ln(\rho_I)} f(e^{-s},\rho_I\rho_0 e^{s})\dd s
\end{cases}
\]
Note that we formally have $(\rho_I\partial_{\rho_I}-\rho_0\partial_{\rho_0})S(f) = f$.
We have that $u:(v,\rho_0)\mapsto \int_{0}^{-\ln(v)} f(e^{-s},v\rho_0 e^{s})\dd s$ is well defined as a smooth section of $(0,1)_{\rho_I}\times (0,1)_{\rho_0} \times \mathcal{B}_s$ vanishing near $v=1$. We have to prove that it belongs to $\overline{\mathcal{H}}_{b}^{\tilde{r},\alpha,\beta}$.
By induction, we prove that for $k,j\in \N$ with $k+j\leq \tilde{r}$, there exists coefficients $a_{\mu,\nu}, b_{\mu}$ (independent of $f$) such that
\begin{align*}
(\rho_0\partial_{\rho_0})^k(\rho_I\partial_{\rho_I})^{j}u = \sum_{\mu+\nu\leq \tilde{r}} a_{\mu,\nu}(\rho_I\partial_{\rho_I})^\mu(\rho_0\partial_{\rho_0})^\nu f + \sum_{\mu = 0}^{\tilde{r}}b_{\mu}\int_0^{-\ln(\rho_I)} ((x_2\partial_{x_2})^{\mu}f)(e^{-s}, \rho_I\rho_0 e^s)\dd s
\end{align*}
(where we have called $x_2$ the second variable of $f$).
The $\overline{\mathcal{H}}_{b}^{0,\alpha,\beta}$-norm of the first sum in the right hand side is bounded by $C\left\|f\right\|_{\overline{\mathcal{H}}_{b}^{\tilde{r},\alpha,\beta}}$. Moreover we have:
\begin{align*}
I:=& \int_0^1\int_0^1 \rho_I^{-2\alpha}\rho_0^{-2\beta}\left|\int_0^{-\ln(\rho_I)} ((x_2\partial_{x_2})^{\mu}f)(e^{-s}, \rho_I\rho_0 e^s)\dd s\right|^2 \frac{\dd \rho_0}{\rho_0}\frac{\dd \rho_I}{\rho_I}\\
 \leq& \int_0^1\int_0^1 \rho_I^{-2\alpha}\rho_0^{-2\beta}\left|\ln(\rho_I)\right|\int_0^{-\ln(\rho_I)} \left|((x_2\partial_{x_2})^{\mu}f)(e^{-s}, \rho_I\rho_0 e^s)\right|^2\dd s \frac{\dd \rho_0}{\rho_0}\frac{\dd \rho_I}{\rho_I} \\
\end{align*}
By the change of variable $(w,z,\rho_I)=(e^{-s},\rho_I\rho_0 e^s,\rho_I)$, we find:
\begin{align*}
I\leq& \int_{(0,1)}\int_{(w,z)\in (\rho_I,1)\times (0,1)}\left|\ln(\rho_I)\right| \rho_I^{2(\beta-\alpha)}z^{-2\beta}w^{-2\beta}\left|((x_2\partial_{x_2})^{\mu}f)(w, z)\right|^2 \frac{\dd w}{w}\frac{\dd z}{z}\frac{\dd \rho_I}{\rho_I}\\
\leq& \int_{(0,1)}\left|\ln(\rho_I)\right|\rho_I^{2\eta} \int_{(w,z)\in (0,1)^2}z^{-2\beta}w^{-2\alpha-2\eta}\left|((x_2\partial_{x_2})^{\mu}f)(w, z)\right|^2 \frac{\dd w}{w}\frac{\dd z}{z}\frac{\dd \rho_I}{\rho_I}
\end{align*}
which is true for all $\eta \in (0,\beta-\alpha]$
We deduce that for all $\nu \in (\alpha, \beta]$ there exists $C_\nu>0$ such that:
\begin{align*}
I\leq C_{\nu}\left\|f\right\|_{\overline{\mathcal{H}}_{b}^{\tilde{r},\nu,\beta}}
\end{align*}
Finally, we get $\left\|S(f)\right\|_{\overline{\mathcal{H}}_{b}^{\tilde{r},\alpha,\beta}}\leq C\left\|f\right\|_{\overline{\mathcal{H}}_{b}^{\tilde{r},\nu,\beta}}$. By density of $\mathcal{D}$ in $\overline{\mathcal{H}}_{b}^{\tilde{r},\nu,\beta}$, we get that $S$ extends uniquely as a continuous linear map from $\overline{\mathcal{H}}_{b}^{\tilde{r},\beta-1,\beta}$ to $\overline{\mathcal{H}}_{b}^{\tilde{r},\alpha,\beta}$. Moreover, since $(\rho_I\partial_{\rho_I}-\rho_0\partial_{\rho_0})S(f) = f$ for $f\in \mathcal{D}$, the relation stays true in the sense of distributions for $f\in \overline{\mathcal{H}}_{b}^{\tilde{r},\alpha,\beta}$.
Let $u\in \rho_I^{-N}\rho_0^\beta H_b^{\tilde{r}}(U)$ such that $(\rho_I\partial_{\rho_I}-\rho_0\partial_{\rho_0})u \in \rho_I^\nu\rho_0^\beta H_b^{\tilde{r}}(U)$. We have $(\rho_I\partial_{\rho_I}-\rho_0\partial_{\rho_0})\chi(\rho_I) u = \chi(\rho_I)(\rho_I\partial_{\rho_I}-\rho_0\partial_{\rho_0}) u+ \rho_I\chi'(\rho_I)  u$ and since $(\rho_I\partial_{\rho_I}-\rho_0\partial_{\rho_0})$ has no kernel in $\overline{\mathcal{H}}_{b}^{\tilde{r},-\infty,-\infty}$ (since all the integral curves of $\rho_I\partial_{\rho_I}-\rho_0\partial_{\rho_0}$ starting in $\left\{\rho_I\leq 1, \rho_0\leq 1\right\}$ intersect transversally the boundary $\rho_I= 1$), we deduce $S(\chi(\rho_I)(\rho_I\partial_{\rho_I}-\rho_0\partial_{\rho_0}) u+ \rho_I\chi'(\rho_I) u) = \chi(\rho_I) u$ on $\left\{\rho_I \leq 1, \rho_0\leq 1\right\}$. In particular:
\begin{align*}
\left\|\chi(\rho_I) u\right\|_{\overline{\mathcal{H}}_{b}^{\tilde{r},\alpha,\beta}}\leq & C\left(\left\|\chi(\rho_I) (\rho_I\partial_{\rho_I}-\rho_0\partial_{\rho_0}) u\right\|_{\overline{\mathcal{H}}_{b}^{\tilde{r},\nu,\beta}} + \left\|\rho_I\chi'(\rho_I)  u\right\|_{\overline{\mathcal{H}}_{b}^{\tilde{r},\nu,\beta}}\right)
\end{align*}
Looking at the supports, it implies:
\begin{align*}
\left\|\chi(\rho_I) u\right\|_{\overline{\mathcal{H}}_{b}^{\tilde{r},\alpha,\beta}}\leq & C\left(\left\|(\rho_I\partial_{\rho_I}-\rho_0\partial_{\rho_0}) u\right\|_{\overline{\mathcal{H}}_{b}^{\tilde{r},\nu,\beta}} + \left\| u\right\|_{\overline{\mathcal{H}}_{b}^{\tilde{r},-\infty,\beta}}\right)
\end{align*}
\end{proof}

\begin{proof}[Proof of Proposition \ref{refinementAtScri}]
We first use Proposition \ref{basicConormalEstimate} with $a_0-2<a_I<0$ (by hypothesis, $a_0=1+\alpha<2$).
We exploit the fact that $T_s - \frac{2}{\rho_I}(\rho_I\partial_{\rho_I}-\rho_0\partial_{\rho_0})(\rho_I\partial_{\rho_I}-1)\in Diff_b^2$. Since $T_su= 0$, we get:
\begin{align*}
(\rho_I\partial_{\rho_I}-\rho_0\partial_{\rho_0})(\rho_I\partial_{\rho_I}-1) u \in \rho_I^{a_I+2}\rho_0^{a_0}H^{\tilde{r}-2}_b(U)
\end{align*}
We apply Lemma \ref{estimateTransport1} and get $w:=(\rho_I\partial_{\rho_I}-1)\chi(\rho_I) u \in \overline{\mathcal{H}}_{b}^{\tilde{r}-2, a_0-, a_0}$ where $\chi\in C^{\infty}([0,+\infty),[0,1])$ is a smooth cutoff equal to 1 near zero and with support in $[0,1)$. We can use the Mellin transform with respect to $\rho_I$ which send $\overline{\mathcal{H}}_{b}^{\tilde{r}-2, a_0-, a_0}$ to a space of holomorphic functions from $\left\{\Im(\lambda)>-a_0+\right\}$ to $\overline{H}_b^{\tilde{r}-2}([0,1]_{\rho_0}\times \mathcal{B}_s)$ (where in this last space the sections are extendible at $\rho_0 = 1$ and have b regularity $\tilde{r}-2$ at $\rho_0 = 0$). Therefore:
\begin{align*}
(i\lambda-1)\mathcal{M}(\chi u)(\lambda) = g(\lambda)
\end{align*}
where $g$ is holomorphic on $\left\{\Im(\lambda)>-1-\alpha+\right\}$.Therefore, $\mathcal{M}(\chi u)(\lambda)$ is meromorphic on the same set with only one possible simple pole at $\lambda = -i$. Moreover, by Plancherel, we have the following bound for $\lambda \in \R+i(\epsilon-a_0)$ (with $\epsilon>0$):
\begin{align*}
\sum_{j=0}^{\tilde{r}-2}\int_{\R+i(\epsilon-a_0)}\Re(\lambda)^{2j}\left\| g(\lambda)\right\|^2_{\overline{H}_b^{\tilde{r}-2-j}([0,1]_{\rho_0}\times \mathcal{B}_s)}\dd \lambda\leq& \sum_{j=0}^{\tilde{r}-2}\int_{0}^{1}\left\|\rho_I^{-a_0+\epsilon}(\rho_I\partial_{\rho_I})^jw(\rho_I)\right\|^2_{\overline{H}_b^{\tilde{r}-2-j}([0,1]_{\rho_0}\times \mathcal{B}_s)}\frac{\dd\rho_I}{\rho_I}\\
\leq& C\left\|w\right\|_{\overline{\mathcal{H}}_{b}^{\tilde{r}-2,a_0-\epsilon, a_0}}.
\end{align*}
We can also get the bound (for $\lambda \in \R+i(\epsilon-a_0)$):
\begin{align*}
\left\| g(\lambda)\right\|^2_{\overline{H}_b^{\tilde{r}-2}([0,1]_{\rho_0}\times \mathcal{B}_s)}\leq& \int_{0}^{1}\rho_I^{\frac{\epsilon}{2}}\left\|\rho_I^{-a_0+\frac{\epsilon}{2}}w(\rho_I)\right\|_{\overline{H}_b^{\tilde{r}-2-j}([0,1]_{\rho_0}\times \mathcal{B}_s)}\frac{\dd\rho_I}{\rho_I}\\
\leq& \frac{1}{\epsilon}\left\|w\right\|_{\overline{\mathcal{H}}_{b}^{\tilde{r}-2,a_0-\frac{\epsilon}{2}, a_0}}
\end{align*}
We deduce that we can use a contour deformation argument to get (the limit are understood in the sense of distributions):
\begin{align*}
\chi u =& \lim\limits_{M\to +\infty}\frac{1}{2\pi}\int_{[-M,M]+i(\epsilon-a_I)}\rho_I^{i\lambda}\frac{g(\lambda)}{i\lambda-1}\dd \lambda\\
=& \lim\limits_{M\to +\infty}\frac{1}{2\pi}\int_{[-M,M]+i(\epsilon-1-\alpha)}\rho_I^{i\lambda}\frac{g(\lambda)}{i\lambda-1}\dd \lambda + g(-i)\rho_I
\end{align*}
Therefore, $\chi u = \chi(\rho_I) g(-i)\rho_I+ \overline{\mathcal{H}}_{b}^{\tilde{r}-2,1+\alpha-\epsilon,a_0}$.
We now investigate the dependency of $g(-i)$ with respect to $\rho_0$. Using Lemma \ref{basicConormalEstimate}, for $k\in \N$ with $k\leq \tilde{r}-2$ and the fact that $\rho_0=-\frac{1}{\tt}$, we deduce that $\rho_0\mapsto w(\rho_0)$ belongs to $\rho_0^{a_0}C^{k}([0,1]_{\rho_0}, \overline{H}^{\tilde{r}-2-k, 1+\alpha-}_b)$.
\begin{align*}
g(-i) = \int_0^{+\infty}w(\rho_0,\rho_I)\dd\rho_I 
\end{align*}
we deduce $g(-i) \in C^{k}\left([0,1]_{\rho_0}, H^{\tilde{r}-2-k}(\mathcal{B}_s)\right)$ and therefore $\chi u - \chi(\rho_I) g(-i)\rho_I \in C^{k}\left([0,1]_{\rho_0}, \overline{H}^{\tilde{r}-2-k, 1+\alpha-}_b\right)$. We define $\chi_1 \in C^{\infty}(\R_{\mathfrak{t}}, [0,1])$ equal to 1 near $+\infty$ and such that $\supp(\chi_1')\subset [\mathfrak{t}_0,\mathfrak{t}_1]\subset (-\infty, -1)$. Then if we define $v := \chi_1(\mathfrak{t})u$ we get $v\in C^k(\R_{\mathfrak{t}}, \overline{H}^{\tilde{r}-2-k,1-}_{b})$. Moreover, 
\begin{align*}
T_s v =& [T_s, \chi_1]u\\
=& \left(2a_{\mathfrak{t},\mathfrak{t}}\chi_1'\partial_{\mathfrak{t}} + a_{\mathfrak{t},\varphi}\chi_1'\partial_{\varphi}+a_{\mathfrak{t},r}\chi_1'\partial_r+a_{\mathfrak{t}}\chi_1'+a_{\mathfrak{t},\mathfrak{t}}\chi_1''\right)u
\end{align*}
Since $2a_{\mathfrak{t},\mathfrak{t}}\chi_1'\partial_{\mathfrak{t}} + a_{\mathfrak{t},\varphi}\chi_1'\partial_{\varphi}+a_{\mathfrak{t},r}\chi_1'\partial_r+a_{\mathfrak{t}}\chi_1'+a_{\mathfrak{t},\mathfrak{t}}\chi_1'' = 2x^{-1}\chi_1'\left(-x\partial_x+1\right)+\text{Diff}^1_b(\mathfrak{t}^{-1}((\mathfrak{t}_0,\mathfrak{t}_1))$, and $2x^{-1}\chi_1'\left(-x\partial_x+1\right)g(-i)\rho_I = 0$ we get $T_s v \in C^{k-1}(\R_{\mathfrak{t}}, \overline{H}^{\tilde{r}-3-k, \alpha}_b)$ which conclude the proof of proposition \ref{refinementAtScri}.
\end{proof}

\bibliography{biblio.bib}
\end{document}